\documentclass[11pt]{imsart}

\usepackage{xr}

\usepackage{tikz}

\usepackage{amsthm,amsmath,amssymb}
\usepackage{parskip}
\usepackage[numbers]{natbib}
\usepackage{hypernat}
\usepackage{marvosym}
\usepackage[hang,small,bf]{caption}
\usepackage{mathtools}

\usepackage[colorlinks]{hyperref}
\usepackage{hyperref}
\hypersetup{
    colorlinks,%
    citecolor=blue,%
    filecolor=black,%
    linkcolor=blue,%
    urlcolor=black
}
\usepackage{hypernat}
\usepackage{enumerate}

\usepackage{color}
\usepackage{bbm}
\usepackage{amsthm}
 \usepackage{amsfonts}
\usepackage{amsmath}
\usepackage{amssymb}
\usepackage{thmtools} 
\usepackage{xcolor}
\usepackage{subcaption}
\usepackage{textcomp}

\usepackage{lineno}

\usepackage[numbers]{natbib}
\usepackage{tikz}
\usepackage{mathtools}
\mathtoolsset{showonlyrefs}

\newtheoremstyle{exampstyle}
{5pt} 
{5pt} 
{\it} 
{} 
{\bfseries} 
{.} 
{.5em} 
{} 

\theoremstyle{exampstyle}

\newtheorem{theorem}{Theorem}[section]
\newtheorem{proposition}[theorem]{Proposition}
\newtheorem{lemma}[theorem]{Lemma}
\newtheorem{corollary}[theorem]{Corollary}
\newtheorem{remark}{Remark}[section]
\declaretheorem[name={Definition}, sibling=remark]{definition}

\startlocaldefs

\newcommand{\eat}[1]{}

\DeclareMathOperator*{\argmin}{\arg\!\min}
\DeclareMathOperator*{\argmax}{\arg\!\max}

\renewcommand{\bar}[1]{\overline{#1}}
\renewcommand{\hat}[1]{\widehat{#1}}
\renewcommand{\tilde}[1]{\widetilde{#1}}
\newcommand{\E}{\mathbb{E}}
\renewcommand{\P}{\mathbb{P}}

\newcommand{\R}{\mathbb{R}}
\newcommand{\A}{\mathcal{A}}

\newcommand{\B}{\mathcal{B}}
\newcommand{\s}{\mathcal{S}}

\renewcommand{\b}[1]{\boldsymbol{\mathbf{#1}}}
\renewcommand{\vec}[1]{\b{#1}}

\usepackage{hyperref}

\newcommand{\conv}{\mathrm{Conv}}
\newcommand{\Vor}{\mathrm{Vor}}

\newcommand{\X}{\mathcal{X}}
\newcommand{\Y}{\mathcal{Y}}

\newcommand{\dom}{\operatorname{dom}}
\newcommand{\N}{\mathbb{N}}

\newtheorem{claim}{Claim}[section]
\theoremstyle{plain}

\def\beq{\begin{equation}}
\def\eeq{\end{equation}}
\def\ba{\begin{enumerate}[(a)]}
\def\bei{\begin{enumerate}[(i)]}
\def\be{\begin{enumerate}[(1)]}
\def\ee{\end{enumerate}}
\def\bi{\begin{itemize}}
\def\ei{\end{itemize}}
\def\beg{\begin{eg}}
\def\eeg{\end{eg}}
\def\bd{\begin{defn}}
\def\ed{\end{defn}}
\def\bt{\begin{thm}}
\def\et{\end{thm}}
\def\bl{\begin{lemma}}
\def\el{\end{lemma}}
\def\bfac{\begin{fact}}
\def\efac{\end{fact}}

\def\bc{\begin{corollary}}
\def\ec{\end{corollary}}
\def\bp{\begin{prop}}
\def\ep{\end{prop}}
\def\bo{\begin{observe}}
\def\eo{\end{observe}}
\def\bas{\begin{assumption}}
\def\eas{\end{assumption}}
\def\RR{\mathbb{R}}

\def\NN{\mathbb{N}}

\endlocaldefs

\begin{document}

\begin{frontmatter}
\title{Multivariate Ranks and Quantiles using Optimal Transport: Consistency, Rates, and Nonparametric Testing}

\runtitle{Multivariate Ranks and Quantiles using Optimal Transport}

\begin{aug}
\author[A]{\fnms{Promit} \snm{Ghosal}\ead[label=e1]{promit@mit.edu}}
\and
\author[B]{\fnms{Bodhisattva} \snm{Sen}\thanksref{t2}$^,$\ead[label=e2]{bodhi@stat.columbia.edu}}
\thankstext{t2}{Supported by NSF grant DMS-2015376.}

\runauthor{Ghosal and Sen}

\address[A]{Department of Mathematics, Massachusetts Institute of Technology,
\printead{e1}}

\address[B]{Department of Statistics, 
Columbia University, 
\printead{e2}}

\end{aug}

\begin{abstract}
In this paper we study multivariate ranks and quantiles, defined using the theory of optimal transport, and build on the work of Chernozhukov et al.~\cite{Cher17} and Hallin et al.~\cite{dCHM}. We study the characterization, computation and properties of the multivariate rank and quantile functions and their empirical counterparts. We derive the uniform consistency of these empirical estimates to their population versions, under certain assumptions. In fact, we prove a Glivenko-Cantelli type theorem that shows the asymptotic stability of the empirical rank map in any direction. Under mild structural assumptions, we provide global and local rates of convergence of the empirical quantile and rank  maps. We also provide a sub-Gaussian tail bound for the global $L_2$-loss of the empirical quantile function. Further, we propose tuning parameter-free multivariate nonparametric tests --- a two-sample test and a test for mutual independence --- based on our notion of multivariate quantiles/ranks. Asymptotic consistency of these tests are shown and the rates of convergence of the associated test statistics are derived, both under the null and alternative hypotheses.     
\end{abstract}



\end{frontmatter}

\section{Introduction}\label{sec:Intro}
Suppose that $X$ is a random vector in $\R^d$, for $d \ge 1$, with distribution $\nu$. When $d=1$, the rank and quantile functions of $X$ are defined as $F$ and $F^{-1}$ (the inverse\footnote{$F^{-1}(p) := \inf \left\{x\in {\mathbb  {R}}:p\leq F(x)\right\}$.} of $F$), respectively, where $F$ is the cumulative distribution function of $X$. Moreover, when $d=1$, quantile and rank functions and their empirical counterparts are ubiquitous in statistics and form the backbone of what is now known as classical nonparametrics (see e.g.,~\cite{Lehmann75} and the references therein) and are important tools for inference (see e.g.,~\cite{HW03} and the references therein). In this paper we study many properties of {\it multivariate} (empirical) ranks and quantiles defined using the theory of optimal transport (OT), as introduced in~\cite{Cher17}. 

Unlike the real line, the $d$-dimensional Euclidean space $\R^d$, for $d \ge 2$, has no natural ordering. This has been a major impediment in defining analogues of quantiles and ranks in $\R^d$, for $d \ge 2$. Several notions of multivariate quantiles have been proposed in the statistical literature --- some based on data depth ideas (see e.g.,~\cite{Oja83, Liu92, Zou03}) and some based on geometric ideas (see e.g.,~\cite{Chaud96, Kol97, HPS10}); see~\cite{Serfling10} and~\cite{dCHM} for recent surveys on this topic. However, most of these notions do not enjoy the numerous appealing properties that make univariate ranks and quantiles so useful. For example, most of these notions can lead to multivariate quantiles that may take values outside the support of the distribution $\nu$.

To motivate the notions of ranks and quantiles based on the theory of OT (the subject of our study) let us first consider the case when $d=1$. Suppose that $X \sim \nu$ has a continuous distribution function $F$. An important property of the one-dimensional rank function $F$ is that $F(X) \sim \mu$ where $\mu \equiv $ Uniform$([0,1])$, i.e., $F$ {\it transports} (see \eqref{eq:PushMeasure-1} for the formal definition) the distribution $\nu$ to $\mu$. Similarly, the quantile function $F^{-1}$ (which is the inverse of the rank map) transports $\mu$ to $\nu$, i.e., $F^{-1}(U) \sim X$ where $U \sim \mu$. In fact, it can be easily shown that the quantile function $F^{-1}$ (or $F$) is the unique monotone nondecreasing map that transports $\mu$ to $\nu$ (or $\nu$ to $\mu$). Moreover, if $\nu$ has finite second moment, it can be shown that $F^{-1}$ is the almost everywhere (a.e.) unique map (on $[0,1]$) that transports $\mu$ to $\nu$ and minimizes the expected squared-error cost, i.e., 
\begin{equation}\label{eq:Q-F}
F^{-1} = \arg \min_{T: T(U) \sim \nu} \E [(U - T(U))^2], \qquad \mbox{where} \;\; U \sim \mu
\end{equation}
and the minimization is over all functions $T$ that transport $\mu$ to $\nu$ (and thus the connection to OT); see Section~\ref{sec:Q-R} for the details. The rank function $F$ also minimizes the expected squared-error cost where now one considers maps that transport $\nu$ to $\mu$.

The multivariate quantile and rank functions using OT essentially extend the above properties of univariate rank and quantile functions. Now let $\mu$ be an absolutely continuous probability measure with respect to (w.r.t.) Lebesgue measure on $\RR^d$ ($d \ge 1$) and supported on a compact convex set $\s$; e.g., we can take $\mu$ to be Uniform$([0,1]^d)$ or uniform on the ball of radius one around $0 \in \R^d$. We often refer to $\mu$ as the \emph{reference distribution} and will define quantiles relative to this reference measure (when $d=1$ we usually take $\mu$ to be Uniform$([0,1])$). Let $\nu$ be another probability measure in $\RR^d$ which we term as the \emph{target distribution}; we think of $\nu$ as the population distribution of the observed data. We define the {\it multivariate quantile function} $Q: \s \to \R^d$ of $\nu$ w.r.t.~$\mu$ as the solution to the following optimization problem:
\begin{equation}\label{eq:Multi-Q-F}
Q := \arg \min_{T: T(U) \sim \nu} \E [\|U - T(U)\|^2], \qquad \mbox{where} \;\; U \sim \mu,
\end{equation}
and the minimization is over all functions $T:\s \to \R^d$ that transport $\mu$ to $\nu$; cf.~\eqref{eq:Q-F} and see Section~\ref{sec:Q-R} for the details. Here $\|\cdot \|$ denotes the usual Euclidean norm in $\R^d$. Moreover, if $\nu$ does not have a finite second moment, the above optimization problem might not be meaningful but the notion of multivariate quantiles (using OT) can still be defined as follows. By Brenier-McCann's theorem (see Theorem~\ref{thm:Brenier}), there exists an a.e.~unique map $Q: \s \to \R^d$ --- which we define as the quantile function of $\nu$ (w.r.t.~the reference measure $\mu$) --- that is the gradient of a convex function and transports $\mu$ to $\nu$; i.e., $Q(U) \sim \nu$ where $U \sim \mu$. Further, it is known that when~\eqref{eq:Multi-Q-F} is meaningful, the above two notions yield the same function $Q$. Note that when $d=1$, the gradient of a convex function is a monotone nondecreasing function and thus the above two characterizations of the quantile function $Q$ are the exact analogues of the one-dimensional case described in the previous paragraph.


Although the rank function can be intuitively thought of as the inverse of the quantile function, such an inverse might not always exist --- especially when $\nu$ is a discrete probability measure (which arises naturally when defining the empirical rank map). In Section~\ref{sec:Q-R} we tackle this issue and use the notion of the Legendre-Fenchel transform (see Section~\ref{sec:prelim}) to formally define the rank function. Indeed, if the reference and the target distributions are absolutely continuous, this notion of rank function is the inverse of the quantile function almost everywhere (a.e.).  Furthermore, it can be shown that (see Proposition~\ref{thm:Q_prop}; also see~\cite[Theorem 1]{CF19}), under mild regularity conditions, the quantile and rank functions are continuous bijections (i.e., homeomorphisms) between the (interiors of the) supports of the reference and target distributions and they are inverses of each other. It is worth noting that when $d=1$, a continuous bijective rank map corresponds to the distribution function being continuous and strictly increasing. 

In Section~\ref{sec:Q-R} we describe some important properties of the defined multivariate quantile and rank functions; also see Section~\ref{sec:Q-R-Prop}. 
 For example, in Lemma~\ref{cor:RankProp} we show that, under appropriate conditions, the rank map approaches a limit, along every ray, that depends only on the geometry of $\s$ (and not on $\nu$); this plays a crucial role in proving the uniform convergence result for the empirical rank map in Theorem~\ref{thm:GCProp}. Some useful properties of the multivariate quantile and rank functions, including: (i) equivariance under orthogonal transformations when the reference distribution is spherically symmetric, and (ii) decomposition/splitting into marginal quantile/rank functions when $X \sim \nu$ has mutually independent sub-vectors and $\mu =$ Uniform($[0,1]^d$), are given in Section~\ref{sec:Euiv-Inv}. 
 Note that the choice of the reference distribution $\mu$ affects the properties of the multivariate ranks/quantiles; see Remark~\ref{rem:Choice-mu} for a discussion on this. 

Given $n$ i.i.d.~random vectors $X_1,\ldots, X_n \sim \nu$ in $\R^d$, in Sections~\ref{sec:Emp_Q_R}, we discuss the characterization and properties of the empirical quantile and rank maps --- which are defined via~\eqref{eq:Multi-Q-F} but with $\nu$ replaced by the empirical distribution of the data. Thus, the computation of the empirical quantile map reduces to a semi-discrete OT problem; see Section~\ref{sec:Comp} for the details where we show that the empirical quantile map can be computed by solving a convex optimization problem with $n$ variables. An attractive property of the empirical ranks, when $d=1$, that makes ranks useful for statistical inference, is that they are distribution-free. Lemma~\ref{lem:Uniform} shows that a distribution-free version of empirical multivariate ranks can be obtained by external randomization (also see Lemma~\ref{lem:Distfree}). Although our approach of defining multivariate quantiles/ranks via the theory of OT has many similarities with those of~\cite{Cher17, dCHM} and~\cite{BSS18} there are subtle and important differences; in Section~\ref{sec:Comparison} we discuss these connections. 

The main statistical contributions of this paper are divided in the three sections ---  Sections~\ref{sec:UnifConv},~\ref{sec:Gl-Lo-Rate}, and~\ref{sec:Goodness-Fit-Test}. In the following we highlight some of the main results in these sections and their novelties. 

\textbf{(I) Uniform convergence of empirical quantile/rank maps}: In Section~\ref{sec:UnifConv} we state our first main theoretical result on the almost sure (a.s.) uniform convergence of the empirical quantile and rank maps to their population counterparts. An informal statement of this result (Theorem~\ref{thm:GCProp}) is given below.
Suppose that $\mu$ is supported on a compact convex set $\s\subset \RR^d$ with non-empty interior. Let $\Y$ be the support of $\nu$ and let $\{\hat{\nu}_n\}_{n\ge 1}$ be a sequence of random probability distributions converging weakly to $\nu$ a.s. Suppose that the quantile map $Q$ of $\nu$ (w.r.t.~$\mu$) is a continuous bijection from $\mathrm{Int}(\s)$ (the interior of $\s$)  to $\mathrm{Int}(\Y)$. Then, with probability (w.p.) 1, the empirical quantile and rank maps corresponding to $\hat{\nu}_n$ (w.r.t.~$\mu$) --- $\hat{Q}_n$ and $\hat{R}_n$ --- converge  uniformly to $Q$ and $R \equiv Q^{-1}$, respectively, over compacts inside $\mathrm{Int}(\s)$ and $\mathrm{Int}(\Y)$. Moreover, if $\s \subset \R^d$ is a strictly convex set (see Definition~\ref{defn:St-Cvx}) then $\hat{R}_n$ converges uniformly to $R = Q^{-1}$ over the whole of $\RR^d$ a.s.; furthermore, w.p.~1, the tail limit of $\hat{R}_n$ stabilizes along any direction. 
We mention below two main novelties of the above result. 

\noindent \textbf{(a)} One of the main consequences of Theorem~\ref{thm:GCProp} is the a.s.~convergence of the empirical rank function $\hat{R}_n$ on the whole of $\RR^d$, under the strong convexity condition on the support $\s$ of $\mu$. This can indeed be thought of as a  generalization of the famous Glivenko-Cantelli theorem for rank functions when $d>1$. Moreover, our result does not need any boundedness assumption on the support of $\nu$ and even applies when the second moment of $\nu$ is not finite. This is a major improvement over the corresponding results in~\cite[Theorem 3.1]{Cher17} and~\cite[Theorem 2.3]{BSS18}. 
Furthermore, unlike in~\cite{dCHM}, $\mu$ can be any absolutely continuous distribution supported on a compact convex domain with minor restrictions on its boundary. Note that for Theorem~\ref{thm:GCProp} to hold we need to assume that $Q$ is a homeomorphism; in particular, if $\nu$ has a convex support with a bounded density then the above holds; see e.g.,~Proposition~\ref{thm:Q_prop}. 
 

\noindent \textbf{(b)} Our result (see \eqref{eq:Asymptot} of Theorem~\ref{thm:GCProp}) implies that when the population rank map is a homeomorphism, the tail limits of the estimated rank maps $\hat{R}_n$ depend neither on $\nu$ nor on $\mu$; rather they depend on the geometry of $\s$ --- the support of the reference distribution $\mu$. This is reminiscent of the case when $d=1$ where the limits of the distribution (rank) function towards $-\infty$ and $+\infty$ are always $0$ and $1$, respectively (irrespective of $\nu$).   





\textbf{(II) Rate of convergence of empirical quantile/rank maps}: Theorem~\ref{thm:GCProp} naturally leads to the question: ``What are the rates of convergence of the empirical quantile/rank maps --- $\hat{Q}_n$ and $\hat{R}_n$?". We study this question in detail in Section~\ref{sec:Gl-Lo-Rate}. We first introduce the following notation:
\begin{equation}\label{eq:rate} 
r_{d,n} := \begin{cases} n^{-1/2} & d =1, 2,3, \\ n^{-1/2} \log n & d = 4, \\
n^{-2/d} & d > 4. \end{cases}
\end{equation} 
\noindent \textbf{(a)} In Theorem~\ref{thm:Q-Rate} we provide upper bounds on the $L_2$-global risk of the empirical quantile map $\hat{Q}_n$. In particular, we show that, for all $n \ge 1$, $$\E \left[\int  \| \hat{Q}_n  - Q \|^2 d \mu  \right] \le C \,  r_{d,n},$$ where $C>0$ is a constant that depends only on $\mu$ and $\nu$. This result is proved using Lemma~\ref{ppn:RateProp2}, which is of independent interest, and gives a quantitative stability estimate for OT maps in the semi-discrete setting. Note that the rates obtained in Theorem~\ref{thm:Q-Rate} are strictly better than those obtained for OT maps in~\cite[Theorem 1.1]{Berman2018} and~\cite[Section 4]{LN2020}. Furthermore, in Theorem~\ref{thm:Q-Rate} we also give a sub-Gaussian tail bound for $\int \| \hat{Q}_n  - Q \|^2 d \mu$. We believe that Theorem~\ref{thm:Q-Rate} gives the exact rate of convergence for the empirical quantile map $\hat{Q}_n$ when $d > 4$; see~\cite{HR19} where the conjecture of this optimality of the rate $n^{-2/d}$ (for $d > 4$) is made. Furthermore, Theorem~\ref{thm:Q-Rate} holds under minimal structural assumptions on $\nu$ --- we only assume strong convexity of the underlying potential function (see~\eqref{eq:Quantile} below).

\noindent \textbf{(b)} In Theorem~\ref{thm:RateProp1}, under appropriate assumptions, we give an upper bound on the risk of the sample rank map, i.e., we show that, for all $n \ge 1$, 
$$\E \left[ \frac{1}{n} \sum_{i=1}^n \big\| \hat R_n(X_i) - R(X_i)\big\|^2 \right] \le K \,  r_{d,n},$$ where $K>0$ is a constant that depends only on $\mu$ and $\nu$. Deriving such a rate result for the multivariate sample rank map $\hat R_n$ is a bit more tricky as $\hat R_n$ is not an OT map per se, but is defined via the Legendre-Fenchel transform (see Section~\ref{sec:Q-R} for the details). 

\noindent \textbf{(c)} 
We address the local uniform rate of convergence of the empirical quantile and rank maps in Theorem~\ref{thm:RateTheo}. The pointwise rate of convergence of empirical rank/quantile maps, defined via the theory of OT, is indeed a hard problem when $d>1$ and not much is known in the literature. Under similar assumptions as in Theorems~\ref{thm:Q-Rate} and~\ref{thm:RateProp1}, we show that  $\hat Q_n$ and $\hat R_n$ converge locally uniformly to $Q$ and $R$, respectively, at the rate $r^{1/(d+2)}_{d,n}$. We consider Theorem~\ref{thm:RateTheo} as a first step towards understanding the local behavior of transport maps. The proof of this result uses Theorem~\ref{thm:Q-Rate} and a correspondence result between the local uniform and local $L_2$ rates of convergence of the empirical rank and quantile functions (see Proposition~\ref{ppn:RateProp1}) 
that could be of independent interest. 


\textbf{(III) Applications to nonparametric testing}: In Section~\ref{sec:Goodness-Fit-Test}, we investigate two statistical applications of the multivariate rank and quantile functions studied in this paper --- we propose methodology for multivariate two-sample goodness-of-fit testing (in Section~\ref{sec:2S-Goodness-Fit-Test}) and testing for mutual independence (in Section~\ref{sec:IndepTest}). 
Both the proposed tests are tuning parameter-free. Applying the uniform convergence results of Theorem~\ref{thm:GCProp}, we prove the consistency of these proposed tests, i.e., the power of these tests converges to 1 under fairly general assumptions on the underlying distributions (see Propositions~\ref{lem:Power1} and~\ref{lem:Power2}). Moreover, using the results in Section~\ref{sec:Gl-Lo-Rate} we provide rates of convergence of the test statistics (for both the testing problems), under both the null and alternative hypotheses; see Propositions~\ref{thm:Two-S-Rate},~\ref{thm:Two-S-Alt},~\ref{lem:Rate2} and~\ref{lem:Rate2-Alt}. This leads to omnibus consistent nonparametric tests that are computationally feasible, and being rank based, do not depend on moment assumptions on the underlying distribution(s). 

Although we state most of our results in terms of multivariate quantile and rank functions, many of the results have immediate implications in estimation of OT maps. Indeed, in recent years there has been a deluge of work at the intersection of statistics and the theory of OT; see e.g.,~\cite{RW18-MLE,WB19, RW19, PZ19, Klatt20, Ramdas17, delB-JMA-19, delB-AoP-19,Peyre2019} and the references therein.

The paper is organized as follows. We introduce notation and some basic notions from convex analysis and the theory of OT in Section~\ref{sec:prelim}. Section~\ref{sec:Q-R} defines the multivariate quantile and rank maps and their empirical counterparts and investigates some of their properties, including computation. The asymptotic results on the uniform a.s.~convergence of the empirical quantile and rank maps are given in Section~\ref{sec:UnifConv}. Global and local rates of convergence of the empirical quantile/rank maps are given in Section~\ref{sec:Gl-Lo-Rate}. The two statistical applications in nonparametric testing are given in Section~\ref{sec:Goodness-Fit-Test}. All the proofs of the main results, additional (technical) results, further remarks and discussions are relegated to the Appendices~\ref{sec:Proofs}--\ref{Appendix-B}. 

\section{Preliminaries}\label{sec:prelim}
We start with some notation and recall some important concepts from convex analysis that will be relevant for the rest of the paper. For $u, v \in \R^d$, we use $\langle u , v \rangle$ to denote the dot product of $u$ and $v$ and $\|\cdot\|$ denotes the usual Euclidean norm in $\R^d$. For $y_1, \ldots, y_k \in \R^d$ we write $\conv(y_1,\ldots, y_k)$ to denote the convex hull of $\{y_1, \ldots, y_k\} \subset \R^d$. A {\it convex polyhedron} is the intersection of finitely many closed half-spaces. A {\it convex polytope} is the convex hull of a finite set of points. The interior, closure and boundary of a set $\X \subset \R^d$ will be denoted by Int($\X$), Cl($\X$), and Bd($\X$), respectively. The Dirac delta measure at $x$ is  denoted by $\delta_x$. For $\delta >0$ and $x \in \R^d$, $B_{\delta}(x) := \{y \in \R^d: \|y - x\| < \delta\}$ denotes the open ball of radius $\delta$ around $x$. The set of natural numbers will be denoted by $\N$.

The {\it domain} of a function $f : \R^d \to \R \cup \{+\infty\}$, denoted by $\dom(f)$, is the set $\{x \in \R^d: f(x) < + \infty\}$. A function $f$ is called {\it proper} if $\dom(f) \ne \emptyset$. A function $f \in L^{\infty}(\mathcal{S})$, where $\s \subset \R^d$, if $\sup_{x \in \s} |f(x)| < \infty$. We say that $f$ is lower semi-continuous (l.s.c.) at $x_0 \in \R^d$ if $\liminf _{x\to x_{0}}f(x)\geq f(x_{0})$. For a proper function $f : \R^d \to \R \cup \{+\infty\}$, the {\it Legendre-Fenchel dual} (or convex conjugate or simply the dual) of $f$ is the proper function $f^* : \R^d \to \R \cup \{+\infty\}$ defined by
\begin{equation}\label{eq:LF-dual}
f^*(y) := \sup_{x \in \R^d} \left\{ \langle x,y\rangle - f(x) \right\}, \qquad \mbox{for all } y \in \R^d.
\end{equation}
It is well known that $f^*$ is a proper, l.s.c.~convex function. The Legendre-Fenchel duality theorem says that for a proper l.s.c.~convex function $f$, $(f^*)^* = f$. Throughout the paper, we will assume that all the convex functions that we will be dealing with are l.s.c.

Given a convex function $f : \R^d \to \R \cup \{+\infty\}$ we define the \emph{subdifferential} set of $f$ at $x\in \dom(f)$ by $$\partial f(x) := \big\{\xi \in \R^d: f(x) + \langle y-x, \xi\rangle \le f(y), \quad \mbox{for all } y \in \R^d \big\}.$$ Any element in $\partial f(x)$ is called a \emph{subgradient} of $f$ at $x$. The subdifferential $\partial f(x)$ is empty if $f(x) = + \infty$ and nonempty if $x \in $ Int$(\dom (f))$.  If $f$ is differentiable at $x$ then $\partial f(x) = \{\nabla f(x)\}$. 
A convex function is a.e.~differentiable (w.r.t.~Lebesgue measure) on Int($\dom(f)$). As a consequence, a convex function is continuous in the interior of its domain. For a convex function $f : \R^d \to \R \cup \{+\infty\}$ we sometimes just write $\nabla f(x)$ to denote the (sub)-differential of $f$ at $x$ with the understanding that when $f$ is not differentiable at $x$ we can take $\nabla f(x)$ to be any point in the set $\partial f(x)$. This avoids the need to deal with the set-valued function $\partial f$. However, sometimes we will need to view $\partial f$ as a multi-valued mapping, i.e., a mapping from $\R^d$ into the power set of $\R^d$, and we will use the notation $\partial f$ in that case. We will find the following results useful (see e.g.,~\cite[Proposition 2.4]{V03}).

\begin{lemma}[Characterization of subdifferential]\label{lem:SubD}
Let $f:\R^d \to \R \cup \{+ \infty\}$ be a proper l.s.c.~convex function. Then for all $x,y \in \R^d$, 
\begin{equation}\label{eq:Charac-Sub}
\langle x,y\rangle = f(x) + f^*(y) \Longleftrightarrow y \in \partial f(x) \Longleftrightarrow x \in \partial f^*(y). 
\end{equation}
\end{lemma}
Lemma~\ref{lem:SubD} shows a one-to-one relation between the subdifferential set of a convex function and its Legendre-Fenchel dual. 

\begin{definition}[Strongly convex function]
A function $f:\R^d \to \R \cup \{+\infty\}$ is {\it strongly convex} with parameter $\lambda >0$ if for all $x,y \in \dom(f)$,
\begin{equation}\label{eq:St-Cvx}
f(y) \ge f(x)+ \nabla f(x)^\top (y-x) + \frac{\lambda}{2} \|y-x\|^2.
\end{equation}
\end{definition}

\begin{definition}[Set convergence]\label{bd:SetConv} 
Let $K_1\subset K_2 \subset \ldots$ be an increasing sequence of sets in $\RR^d$. We say that $K_n$ \emph{increases} to $K\subset \RR^d$, and write $K_n \uparrow K$, if for any compact set $A \subset \mathrm{Int}(K)$ there exists $n_0=n_0(A) \in \N$ such that $A \subseteq K_n$ for all $n\geq n_0$.
\end{definition}
The above notion is slightly stronger than just assuming $K_1 \subset K_2 \subset \ldots$ and $\lim \inf_{n \to \infty} K_n = K$.

A supporting hyperplane of a closed convex set $S \subset \R^{d}$ is a hyperplane that has both of the following two properties: (i) $S$ is entirely contained in one of the two closed half-spaces bounded by the hyperplane, and (ii) $S$ intersected with the hyperplane is nonempty. 
\begin{definition}[Strictly convex set]\label{defn:St-Cvx}
A convex set $S \subset \R^d$ is said to be strictly convex if any supporting hyperplane to $\mathrm{Cl}(S)$ touches $\mathrm{Cl}(S)$ at only one point.
\end{definition}

Let $\mu$ and $\nu$ be two Borel probability measures supported on $\s \subset \R^d$ and $\Y\subset \R^d$ respectively. The goal of OT (Monge's problem), under the squared Euclidean loss, is to find a measurable transport map $T \equiv T_{\mu;\nu} : \s \to \Y$ solving the (constrained) minimization problem
\begin{equation}\label{eq:Meas_Trans-1}
	\inf_T \int  \|u - T(u)\|^2 d\mu(u)  \qquad \quad \mbox{subject to }\quad T\#\mu  = \nu
\end{equation}
where the minimization is over $T$ (a {\it transport map}), a measurable map from $\s$ to $\Y$, and $T\#\mu$ is the {\it push forward} of $\mu$ by $T$, i.e., 
\begin{align}\label{eq:PushMeasure-1}
T\#\mu(B) = \mu(T^{-1}(B)), \qquad \mbox{for all } B \subset \Y\;\; \mbox{Borel}. 
\end{align}
A map $T_{\mu;\nu}$ that attains the infimum in~\eqref{eq:Meas_Trans-1} is called an OT map  from $\mu$ to $\nu$. We state an important result in this theory, namely Brenier-McCann's theorem (\cite{B91, McCann95}). This result will be very useful to us; see Section~\ref{sec:OT} 
for a  brief introduction to the field of OT.

\begin{theorem}[Brenier-McCann theorem]\label{thm:Brenier}
Let $\mu$ and $\nu$ be two Borel probability measures on $\R^d$. Suppose further that $\mu$ has a Lebesgue density. Then there exists a convex function $\psi: \R^d \to \R \cup \{+\infty\}$ whose gradient $G =  \nabla \psi : \R^d \to \R^d$ pushes $\mu$ forward to $\nu$. In fact, there exists only one such $G$ that arises as the gradient of a convex function, i.e., $G$ is unique $\mu$-a.e. Moreover, if $\mu$ and $\nu$ have finite second moments, $G$ uniquely minimizes Monge's problem~\eqref{eq:Meas_Trans-1}.
\end{theorem}

\section{Quantile and rank maps in $\R^d$ when $d \ge 1$}\label{sec:Q-R}
Suppose that $X \sim \nu$ is supported on $\Y \subset \R^d$. Let $\mu$ be a known absolutely continuous distribution on $\R^d$ (i.e., $\mu$ has a density w.r.t.~Lebesgue measure on $\R^d$) with support $\s$ --- a compact convex subset of $\R^d$ with nonempty interior; e.g., we can take $\mu$ to be Uniform$([0,1]^d)$. Other natural choices of $\mu$ are the uniform distribution on the unit ball $B_1(0)$ in $\R^d$~\cite{Cher17}, and the spherical uniform distribution ($V$ has the spherical uniform distribution if $V = L \varphi$ where $\varphi$ is  uniformly distributed on the unit sphere around $0 \in \R^d$ and $L \sim $ Uniform($[0,1]$), and $L$ and $\varphi$ are mutually independent); see~\cite{dCHM, Figalli18}. 

In the following we define the multivariate {\it quantile} and {\it rank} maps for $\nu$ w.r.t.~the distribution $\mu$ using the theory of OT. We first define the quantile function for $\nu$ and then use it to define the rank map. 
Our approach is essentially the same as outlined in~\cite{Cher17} although there are some important and subtle differences; see Section~\ref{sec:Comparison} for a discussion. 

\begin{definition}[Quantile function]\label{def:Quantile} 
The quantile function of the probability measure $\nu$ (w.r.t.~$\mu$) is defined as the $\mu$-a.e.~unique map $Q :\s \to \R^d$ which pushes $\mu$ to $\nu$ and has the form
\begin{equation}\label{eq:Quantile}
Q := \nabla \psi
\end{equation} 
where $\psi:\R^d \to \R\cup \{+ \infty\}$ is convex. We call $\psi$ a potential function.
\end{definition}
\begin{remark}[Uniqueness of $Q$]\label{rem:QUnique} 
As the convex function $\psi$ in Definition~\ref{def:Quantile} need not be differentiable everywhere, there is a slight ambiguity in the definition of $Q$. When $\psi$ is not differentiable, say at $u \in \s$, we can define $Q(u)$ to be any element of the subdifferential set $\partial \psi(u)$ (see Section~\ref{sec:prelim} for its formal definition). As a convex function is differentiable a.e.~(on its domain) this convention does not affect the $\mu$-a.e.~uniqueness of $Q$. Further, this convention bypasses the need to define quantiles as a multi-valued map. 
\end{remark}
The existence and $\mu$-a.e.~uniqueness of the quantile map $Q(\cdot)$, for any probability measure $\nu$ on $\R^d$, is guaranteed by the Brenier-McCann theorem (see Theorem~\ref{thm:Brenier}). Further, by Theorem~\ref{thm:Brenier}, if $\nu$ has finite second moment, then $Q(\cdot)$ can be expressed as in~\eqref{eq:Multi-Q-F}.
As discussed in the Introduction, the above notion of quantiles obviously extend our usual definition of quantiles when $d=1$; see Section~\ref{sec:d=1} 
for a more detailed discussion. 

\begin{remark}[Non-uniqueness of $\psi$]
Although $Q$ is $\mu$-a.e.~unique it is easy to see that $\psi$ (as in Definition~\ref{def:Quantile}) is not unique; in fact, $\psi(\cdot) + c$ where $c\in \R$ is a constant would also suffice (as $\partial (\psi+c) = \partial \psi$). Further, we can change $\psi(\cdot)$ outside the set $\s$ and this does not change $Q$ (as $Q$ has domain $\s$). For this reason, we will consider
\begin{equation}\label{eq:Conv-Cvx}
\psi(u) = +\infty, \qquad \mbox{for } u \in \R^d\setminus \s.
\end{equation} 
The above convention will be useful in the subsequent discussion.  
\end{remark}

\begin{definition}[Rank map]\label{def:Rank} 
Recall the convex potential function $\psi:\R^d \to \R \cup \{+\infty\}$ whose gradient yields the quantile map (see~\eqref{eq:Quantile}; also see~\eqref{eq:Conv-Cvx}). We define the rank function $R:\R^d \to \s$ of $\nu$ (w.r.t.~$\mu$) as  
\begin{equation}\label{eq:Rank}
R := \nabla \psi^*
\end{equation} 
where $\psi^*:\R^d \to \R \cup \{+\infty\}$ is the Legendre-Fenchel dual of the convex function $\psi$, i.e., 
$\psi^*(x) := \sup_{u \in \R^d} \{\langle x, u\rangle - \psi(u)\}$, for $x \in \R^d$. Note that $\psi^*$ is also referred to as the dual potential of $Q \equiv \nabla \psi$.
\end{definition}
 A few remarks are in order now.
\begin{remark}[The domain of $R$]  
The rank map $R(x)$ is finite for all $x \in \R^d$; cf.~the quantile map $Q(\cdot)$ which is $\mu$-a.e.~uniquely defined. This follows from the fact that $\psi^*(x) <\infty$ for every $x \in \R^d$; see Lemma~\ref{lem:psi-finite}. 
If $\psi^*$ is not differentiable at $x$ (say) we can define $R(x)$ to be any element in the subdifferential set $\partial \psi^*(x)$.
\end{remark}
 

\begin{remark}[The range of the rank map]
Using Lemma~\ref{lem:SubD} one can argue that $R(x) \in \s$ for every $x \in \R^d$. This follows from the fact that $R(x) \in \partial \psi^*(x)$ exists for every $x \in \R^d$, and by~\eqref{eq:Charac-Sub}, $u \in \partial \psi^*(x) \Leftrightarrow  x \in \partial \psi(u).$ Note that as $\partial \psi(u)$ exists we must have $\psi(u) <+\infty$, which in turn implies that $u \in \s$ (as $\psi(u) = +\infty$, for $u \in \R^d\setminus \s$ by~\eqref{eq:Conv-Cvx}).
\end{remark}

\begin{remark}[When $\psi^*$ is not differentiable]
As $\psi^*(x) <\infty$ for every $x \in \R^d$, $\psi^*$ has a gradient a.e. Thus, $R(x)$ is uniquely defined for a.e.~$x$. For $x\in \R^d$ where $\psi^*(x)$ is not differentiable, $R(x)$ is not uniquely defined. Although for such an $x$ we can define $R(x)$ to be any element in the subdifferential set $\partial \psi^*(x)$ (as was done in~\cite{Cher17}), in Section~\ref{sec:Ranks} we give a randomized choice of $R(x)$ that leads to the map $R$ having appealing statistical properties.
\end{remark}

Absolute continuity of $\nu$ is a sufficient condition for the rank map $R$ and the quantile map $Q$ to be the essential inverses of one another, i.e.,
 \begin{equation}\label{eq:Ess-Inv}
\quad R \circ Q (u)= u,\quad\mbox{for}\;\mu\mbox{-a.e.}~u, \qquad \mbox{and} \qquad Q \circ R (x)= x,\quad\mbox{for}\;\nu\mbox{-a.e.}~x,
\end{equation} 
and $R\# \nu= \mu$ (see e.g.,~\cite[Theorem 2.12 and Corollary 2.3]{V03}). This justifies the definition of $R$ via~\eqref{eq:Rank}. Observe that the rank map, as in Definition~\ref{def:Rank}, clearly extends the notion of the distribution function beyond $d=1$. There is an intimate connection between the quantile map and the celebrated Monge-Amp\`{e}re differential equation; see e.g.,~\cite[Lemma 4.6]{V03} (also see \cite{Ca1,PF14,CF19}). In Section~\ref{sec:Q-R-Prop} 
we discuss a few other interesting properties of quantile/rank and potential functions. 


Although we know that $R = Q^{-1} \;\mu$-a.e.~when $\nu$ is absolutely continuous, one may ask if the equality holds everywhere (as opposed to a.e.). Several results have been obtained in this direction that provide sufficient conditions for such an equality. Caffarelli (see \cite{Ca1,Ca2,Ca3}) showed that when $\mathcal{S}$ and $\mathcal{Y}$ are two bounded convex sets in $\RR^d$ and $\mu$ and $\nu$ are absolutely continuous with positive densities (on their supports), then, the corresponding OT maps $T:\mathcal{S}\to \mathcal{Y}$ (such that $T\# \mu = \nu$) and $T^{*}:\mathcal{Y}\to \mathcal{S}$ (such that $T^{*}\# \nu= \mu$) are continuous homeomorphisms and $T^*= T^{-1}$ everywhere in $\mathcal{Y}$; see \cite[Pages 317--323]{V09} for other sufficient conditions. 
In Proposition~\ref{thm:Q_prop} below we give another such sufficient condition that is particularly useful in statistical applications. As pointed out by an anonymous referee the main result in the recent paper~\cite{CF19} implies Proposition~\ref{thm:Q_prop}-(a); Proposition~\ref{thm:Q_prop}-(b) can then be derived as a consequence. 
 However, for the convenience of the reader we provide a self-contained proof of Proposition~\ref{thm:Q_prop} in Section~\ref{sec:Q_prop}.
 
 \begin{proposition}\label{thm:Q_prop} 
Let $\mathcal{S} \subset \R^d$ be a compact convex set and $\mathcal{Y} \subset \R^d$ be a convex set. Let $\mu$ be a probability measure supported on $\mathcal{S}$ such that the density of $\mu$ (w.r.t.~Lebesgue measure) is bounded away from zero and bounded above by a constant (on $\s$). Let $\nu$ be a distribution supported on $\mathcal{Y}$ with density $p_{\mathcal{Y}}$ satisfying the following: there exists a sequence of convex compact sets $\{K_n\}_{n\geq 1}$ with $K_n\uparrow \mathcal{Y}$ and constants $\{\lambda_n, \Lambda_n\}_{n\ge1} \subset \R$ such that 
\begin{align}\label{eq:sandwitch}
0<\lambda_n\leq p_{\mathcal{Y}}(x) \leq \Lambda_{n}, \qquad \quad \mbox{for all }\; x\in K_n. 
\end{align}
Let $\psi:\R^d\to \R\cup \{+\infty\}$ be a convex function such that $\psi(x) = +\infty$ for $x\notin \mathcal{S}$, $\partial \psi(\mathrm{Int}(\mathcal{S})) = \mathrm{Int}(\mathcal{Y})$ and $\partial \psi \# \mu = \nu$. Let $\psi^{*}:\R^d \to \R\cup \{+\infty\}$ be the Legendre-Fenchel dual of $\psi$. Then:
\begin{enumerate}
\item[(a)] $\nabla \psi^{*}$, restricted to $\mathrm{Int}(\mathcal{Y})$, is a homeomorphism from $\mathrm{Int}(\mathcal{Y})$ to $\mathrm{Int}(\mathcal{S})$. 

\item[(b)] $\nabla \psi$ is a homeomorphism from $\mathrm{Int}(\mathcal{S})$ to $\mathrm{Int}(\mathcal{Y})$. Furthermore, we have $\nabla \psi = (\nabla \psi^{*})^{-1}$ in $\mathrm{Int}(\mathcal{S})$.  
\end{enumerate}      
\end{proposition}

%

\begin{remark}[Convexity of $\mathcal{S}$ and $\mathcal{Y}$]\label{rem:Cvx-S-Y}
Convexity of the domains, $\mathcal{S}$ and $\mathcal{Y}$, is one of the important conditions for the existence of continuous transport maps. Cafarrelli constructed a counterexample (see e.g.,~\cite[pp.~283--285]{V09}) where he showed that the transport map may fail to be continuous when the two measures are absolutely continuous with bounded densities on two smooth and simply connected non-convex domains. 
\end{remark}

\begin{remark}\label{rem:Cond-11} 
The condition \eqref{eq:sandwitch} in Proposition~\ref{thm:Q_prop} does not necessarily require $\mathcal{Y}$ to be compact. For example, any unimodal density supported on a convex domain $\mathcal{Y}\subset \RR^d$ satisfies \eqref{eq:sandwitch};  in particular, this includes the family of all absolutely continuous multivariate normal distributions. 
\end{remark}
Similar to the univariate distribution function, the one-dimensional projection of the rank map $R$ along any direction, is nondecreasing (see Lemma~\ref{lem:Monotonicity} 
for a formal statement of this result). A univariate distribution function is not only nondecreasing but takes the value 0 or 1 as one approaches $-\infty$ or $+\infty$, irrespective of $\nu$. Under mild assumptions on $\s$ and $R$, we show in the following lemma (proved in Section~\ref{sec:RankProp}) 
that $R(\cdot)$ is continuous on the whole of $\RR^d$ and it approaches a limit along every ray that depends only on the geometry of $\s$ and not on the measure $\nu$. 
\bl\label{cor:RankProp} 
Let $\mathcal{S} \subset \RR^d$ be a strictly convex compact set (as in Definition~\ref{defn:St-Cvx}). Let $\mu$ and $\nu$ be two probability measures on $\s$ and $\Y \subset \RR^d$, respectively, where $\Y$ has nonempty interior. Let $R$ be the rank map of $\nu$ w.r.t. $\mu$. Suppose that $R$ is a homeomorphism from $\mathrm{Int}(\Y)$ to $\mathrm{Int}(\mathcal{S})$. Then, for any $x\in \RR^d$, $\lim_{\lambda \to  +\infty} R(\lambda x) = \argmax_{v\in \mathcal{S}} \langle x, v\rangle.$
\el
Note that the above required condition on $\s$ is certainly satisfied, for example, when $\s$ is the unit ball in $\R^d$, i.e., $\s = B_1(0)$; unfortunately when $\s = [0,1]^d$, the condition is not satisfied. Lemma~\ref{cor:RankProp} has a simple interpretation when $\s = B_1(0)$ --- in this case  $\lim_{\lambda \to  +\infty} R(\lambda x) = \argmax_{v \in \s} \langle x, v \rangle = \frac{x}{\|x\|}$ if $x \ne 0$; cf.~\cite[Corollary 3.1]{delBarrio2020}. This generalizes the fact that for a distribution function $F$ on $\R$, $F(x) \to 1$ as $x \to +\infty$ and $F(x) \to 0$ as $x \to -\infty$.

\subsection{The sample quantile and rank maps}\label{sec:Emp_Q_R}
As before, we fix an absolutely continuous distribution $\mu$ with compact convex support $\s \subset \R^d$. Given a random sample $X_1,\ldots, X_n$ from a distribution $\nu$ (on $\R^d$), we now consider estimating the population quantile and rank maps $Q$ and $R$, respectively (w.r.t.~$\mu$). We simply define the sample versions of the quantile and rank maps as those obtained by replacing the unknown distribution $\nu$ with its empirical counterpart $\hat \nu_n$ --- the empirical distribution of the data, i.e., 
$$
\hat \nu_n(A) = \frac{1}{n} \sum_{i=1}^n \delta_{X_i}(A), \qquad \mbox{ for any Borel set $A \subset \R^d$}.
$$
\subsubsection{Empirical quantile function}\label{sec:Q_n}
By Theorem~\ref{thm:Brenier} there exists an $\mu$-a.e.~unique map $\hat Q_n$ which pushes $\mu$ to $\hat \nu_n$ and can be expressed as
\begin{equation}\label{eq:SampQ_n}
\hat Q_n = \nabla \hat \psi_n,
\end{equation} 
where $\hat \psi_n:\R^d \to \R \cup \{+\infty\}$ is a convex function. Further, by Theorem~\ref{thm:Brenier}, $\hat Q_n$ can be computed via:
\begin{equation}\label{eq:Emp_Q}
\hat Q_n = \argmin_T \int  \|u - T(u)\|^2 d\mu(u)  \qquad \quad \mbox{subject to }\quad T\#\mu  = \hat \nu_n.
\end{equation}
Note that $\hat Q_n \equiv \nabla \hat \psi_n$ is $\mu$-a.e.~unique; when $\hat \psi_n$ is not differentiable, say at $u$, we can define $\hat Q_n(u)$ to be any point in $\partial \hat \psi_n(u)$.  As $\hat Q_n \equiv \nabla \hat \psi_n$ pushes $\mu$ to $\hat \nu_n$, $\hat \psi_n$ is a convex function whose gradient takes $\mu$-a.e.~finitely many values --- in the set $\X := \{X_1,\ldots, X_n\}$. Thus $\hat \psi_n$ is piecewise linear (affine), and hence, there exists $\hat h = (\hat h_1,\ldots, \hat h_n) \in \R^n$ (unique up to adding a scalar multiple of $(1,\ldots, 1) \in \R^n$) such that $\hat{\psi}_n:\R^d \to \R \cup \{+\infty\}$ can be represented as
\begin{equation}\label{eq:Emp_G}
\hat{\psi}_n(u) := \begin{cases} \max \limits_{i=1,\ldots, n} \{u^\top X_i + \hat h_i\} & \mbox{for } u \in \s \\ +\infty & \mbox{for } u \notin \s. \end{cases}
\end{equation}
The vector $\hat h$ can be computed by solving a convex optimization problem;  see Section~\ref{sec:Comp} for the details. 
\begin{remark}[Form of the subdifferential set $\partial \hat{\psi}_n(u)$]\label{rem:Q-d} As $\hat{\psi}_n$ is  piecewise linear (affine) and convex (and thus a finite pointwise maximum of affine functions), we can explicitly write its subdifferential, i.e., for any $u \in \s$, $$\partial \hat{\psi}_n(u) =  \conv(\{X_i: \langle u,X_i \rangle + \hat h_i = \hat{\psi}_n(u)\}).$$
\end{remark}
\label{sec:Emp-Q-R}
\begin{figure}
\centering
\includegraphics[scale=0.9]{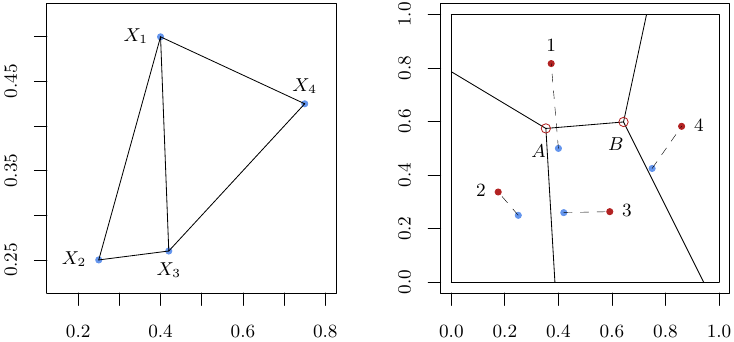}
\caption{The left plot shows a data set with four 2-dimensional points $X_1,X_2,X_3$ and $X_4$.  The right plot shows the four cells (each with area 1/4) marked 1, 2, 3, 4, and the four data points in blue (appropriately scaled to lie in $[0,1]^2$) along with dashed lines connecting each data point to the centroid (in red) of the corresponding cell. The two points $A$ and $B$ in the right plot correspond to the intersection of three cells --- $1, 2, 3$ and $1, 3, 4$.}
\label{fig:Q-map}
\end{figure}
As $\hat Q_n(u) \in \partial \hat{\psi}_n(u)$, for any $u \in \mathrm{Int}(\s)$, $\hat Q_n(u)$ belongs to the convex hull of the data. The function $\hat Q_n = \nabla \hat{\psi}_n$ induces a cell decomposition of $\s$: Each cell is a convex set and is defined as
\begin{equation}\label{eq:W_i}
W_i(\hat h) := \{u \in \s: \nabla \hat{\psi}_n(u) = X_i\}.
\end{equation}
In defining $W_i(\hat h)$  we only consider points $u \in \s$ where $\hat{\psi}_n$ is differentiable; see~\cite{Gu12} for more details. Note that, for a.e.~sequence $X_1,\ldots, X_n$, each cell $W_i(\hat h)$ has $\mu$ measure $1/n$ and $\cup_{i=1}^n W_i(\hat h) \subset \s$. Figure~\ref{fig:Q-map} illustrates this with four points $X_1,X_2,X_3$ and $X_4$ and $\mu = $ Uniform$([0,1]^2)$. Each point in the four cells (labelled 1, 2, 3  and 4) is mapped to the corresponding data point ($X_1,X_2,X_3$ and $X_4$) by the sample quantile function $\hat Q_n \equiv \nabla \hat \psi_n$. The convex function $\hat \psi_n$ is not differentiable at the boundary of the 4 cells (marked by the black lines in the right panel of Figure~\ref{fig:Q-map}). Remark~\ref{rem:Q-d=1} 
illustrates the above ideas when $d=1$ and $\mu = $ Uniform$([0,1])$.

\subsubsection{Empirical rank map}\label{sec:R_n}
Let us define $\hat{\psi}_n^*:\R^d \to \R \cup \{+ \infty\}$ to be the Legendre-Fenchel dual of $\hat{\psi}_n$, i.e., 
 \begin{equation}\label{eq:Sup}
\quad  \hat{\psi}^{*}_n(x) :=\sup_{y\in \RR^d}\big\{\langle x,y\rangle - \hat{\psi}_n(y)\big\} = \sup_{u\in \s}\big\{\langle x,u\rangle - \hat{\psi}_n(u)\big\}, \quad \mbox{for } x \in \R^d.
\end{equation} 
We define the {\it multivariate sample rank function} $\hat{R}_n:\R^d \to \s$ as 
\begin{equation}\label{eq:SampR_n}
\hat{R}_n  := \nabla  \hat{\psi}^{*}_n.
\end{equation} 
Lemma~\ref{lem:RAltDef} 
gives an alternate expression for $\hat{R}_n$ which was used in~\cite[Definition 3.1]{Cher17}. Further, Remark~\ref{rem:R-d=1} 
shows that when $d=1$, $\hat{R}_n$ is not defined uniquely at the data points. Note that the non-uniqueness of the rank function when $d=1$ was finessed by enforcing right-continuity, which is hard to do as we go beyond $d=1$. Indeed, for any $d \ge 1$, $\hat R_n(X_{i})$ could be defined as any element in the (closure of the) cell $W_i(\hat h)$; this follows from Lemma~\ref{lem:SubD}. Figure~\ref{fig:Q-map} illustrates this when $\mu = $ Uniform$([0,1]^2)$. We see in Figure~\ref{fig:Q-map} that any point in the interior of the triangle formed by $X_1, X_2$ and $X_3$ (or $X_1, X_3$ and $X_4$) is mapped to the point $A$ (or $B$) by the sample rank map $\hat R_n$. The following result (Lemma~\ref{lem:Char3},  proved in Section~\ref{sec:Char3}, 
formalizes this observation when $d \ge 2$ and provides a way of finding the empirical rank map at any given point.
\begin{lemma}\label{lem:Char3} 
Fix $x\in \RR^d$. Suppose that for $i_1,\ldots, i_{d+1} \subset \{1,\ldots, n\}$: (i) $x\in \mathrm{Int}\big(\mathrm{Conv}\big(X_{i_1}, \ldots , X_{i_{d+1}}\big)\big)$, and (ii) there exists a unique $u\in \s$ such that $u = \mathrm{Cl}(W_{i_1}(\hat h)) \cap \ldots \cap \mathrm{Cl}(W_{i_{d+1}}(\hat h))$ (see~\eqref{eq:W_i}). Then, $u$ is the unique point in $\s$ such that $x\in \partial \hat{\psi}_n (u)$. Furthermore, $\partial\hat{\psi}^{*}_n(x) = u = \hat{R}_n(x)$.   
\end{lemma}

\subsubsection{Empirical ranks}\label{sec:Ranks}
By the ``ranks" of the data points we mean the rank function evaluated at the data points. When $d=1$ and the underlying distribution is continuous, the usual ranks, i.e., $\{\mathbb{F}_n(X_i)\}_{i=1}^n$ (here $\mathbb{F}_n$ is the empirical distribution function), are identically distributed on the discrete set $\{1/n,2/n,\ldots, n/n\}$ with probability $1/n$ each. As a consequence, the usual ranks are {\it distribution-free} (in $d=1$), i.e., the distribution of $\mathbb{F}_n(X_i)$ does not depend on the distribution of $X_i$. We may ask: ``Does a similar property hold for the multivariate ranks $\hat R_n(X_i)$?".

From the discussion in Section~\ref{sec:R_n} it is clear that the multivariate ranks $\hat R_n(X_{i})$ are non-unique. In fact, we can choose $\hat R_n(X_i)$ to be any point in the set $W_i(\hat h)$ (see~\eqref{eq:W_i}). In the sequel we will use a special choice of $\hat R_n(X_i)$ which will lead to a distribution-free notion. We define $\hat R_n(X_i)$ as a random point drawn from the distribution $\hat \mu_i$, i.e., for $i \in \{1,\ldots, n\}$,
\begin{equation}\label{eq:Rank-X_i}
\hat R_n(X_i)|X_1,\ldots, X_n \sim \hat \mu_i
\end{equation}
where $$\hat \mu_i: B \mapsto n \mu(W_i(\hat h) \cap B) \quad \mbox{ for any Borel } \; B \subset \R^d.$$ Note that $\hat \mu_i$ is a Borel probability measure supported on the cell $W_i(\hat h)$ as $\mu(W_i(\hat h)) = n^{-1}$. 
When $\mu$ is the uniform distribution on $\s$, an equivalent representation of~\eqref{eq:Rank-X_i} is
$\hat R_n(X_i)|X_1,\ldots, X_n \sim \mathrm{Uniform}(W_i(\hat h)).$ Thus, our choice of the empirical ranks $\{\hat R_n(X_i)\}_{i=1}^n$ is random. However, this external randomization leads to the following interesting consequence --- the multivariate ranks are marginally  distribution-free. This is formalized in the following lemma (proved in Section~\ref{pf:Uniform}). 
\begin{lemma}\label{lem:Uniform}
Suppose that $X_1,\ldots, X_n$ are i.i.d.~$\nu$, an absolutely continuous distribution on $\R^d$. Then, for any $i =1,\ldots, n$, $\hat R_n(X_i) \sim \mu.$
\end{lemma}
Compare Lemma~\ref{lem:Uniform} with the result $R(X) \sim \mu$ where $R$ is the (population) rank map of $X \sim \nu$ (as $R$ pushes forward $\nu$ to $\mu$). If we do not want a randomized choice of ranks, then we can define  $\hat R_n(X_i) := \max_{u \in \mathrm{Cl}(W_i(\hat h))}\|u\|;$ the above choice is convenient for computational purposes.

\begin{remark}[Choice of the reference distribution $\mu$]\label{rem:Choice-mu}
As may have been clear from the above discussion, the concept of multivariate (empirical) ranks and quantiles, based on OT, depends on the choice of the reference distribution $\mu$. For example, the choice of a spherically symmetric $\mu$ (e.g., Uniform($B_1(0)$)) leads to quantile maps that are equivariant under orthogonal transformations (which can be an useful property when studying multivariate depth, outlyingness, etc.), whereas the choice of $\mu = $ Uniform$([0,1]^d)$ guarantees factorization into lower dimensional marginals under independence (that may be more appropriate for measuring association/independence between the marginals of $\nu$); see Section~\ref{sec:Q-R-Prop} 
for formal results in this regard. We would like to point out here that $\s$, the support of $\mu$, can play an important role in determining the behavior of the rank map $\hat{R}_n(\cdot)$; we will see in Theorem~\ref{thm:GCProp} that the choice $\s = [0,1]^d$ could lead to inconsistent estimation of $\hat{R}_n(x)$ for $x$ near the boundary of $\mathcal{Y}$ (see Remark~\ref{rem:Suff-Cond} for further discussion). Thus, when $d >1$, the choice of $\mu$ should be dictated by the application at hand. 

The two plots in Figure~\ref{fig:4Plots} show the cell decompositions corresponding to the uniform measures on $[0,1]^2$ and $B_1(0) \subset \R^2$ respectively. As in $d=1$, we believe that the use of appropriate score functions can mitigate the dependence of multivariate rank-based procedures on the reference distribution $\mu$; see e.g.,~\cite[Chapters 13 and 15]{vdV98} for the usefulness of a score-based approach when $d=1$. Indeed, the recent papers~\cite{Shi2020, deb2021efficiency} illustrate the flexibility and power of incorporating score functions in the definition of multivariate rank-based tests; we expect this to be a fertile area of research in the future.

\end{remark} 

\subsection{Computation of the sample quantile and rank functions}\label{sec:Comp}
The computation of the empirical quantile function $\hat Q_n$ (in~\eqref{eq:Emp_Q}) reduces to a semi-discrete OT problem. There are several methods proposed in the literature to solve the semi-discrete OT problem. Oliker and Pr\"{u}ssner~\cite{OP88} proposed one of the earliest algorithms in this regard relying on coordinate-wise increments; also see~\cite{Caffarelli99,  Merigot2020}. Although this algorithm has convergence guarantees (see~\cite{Kitagawa14}) it is quite slow in practice. Recently, fast algorithms for solving~\eqref{eq:Emp_Q} have been proposed that typically rely on the formulation of the semi-discrete OT problem as an unconstrained convex optimization problem which is then solved using a (damped) Newton or quasi-Newton method; see e.g.,~\cite{Auren98, Merigot11, Levy15, Kita-Meri-19}. See~\cite{Merigot2020} for a detailed account of many of the algorithms cited above. 

In the following we outline our approach to computing $\hat Q_n$ (see the \texttt{R} package \url{https://github.com/Francis-Hsu/testOTM}~\cite{OTM}). We use Newton-type algorithms proposed in the papers~\cite{Kita-Meri-19, Merigot11, Levy15} and implemented in the Geogram\footnote{\url{http://alice.loria.fr/software/geogram/doc/html/}} package. These algorithms are experimentally efficient and converge globally with linear rate; see~\cite{Kita-Meri-19}. Our approach is  different from the ``gradient algorithm'' used in~\cite[Section 4]{Cher17} to solve the semi-discrete problem.

The computation of $\hat Q_n$ leads to a ``partition'' of $\s$ into $n$ convex sets, each with volume $1/n$ (i.e., the $W_i(\hat h)$'s in~\eqref{eq:W_i}), and involves what is usually called the power diagram~\cite{Auren87} --- a type of weighted Voronoi diagram. Recall that $\X := \{X_1,\ldots, X_n\}$. Let $w = (w_1,\ldots, w_n) \in \R^n$ be a given (weight) vector. The {\it power diagram} of $(\X,w)$ is the decomposition of the set $\s$ into a finite number of cells, one for each element in $\X$, defined by (for $i=1,\ldots, n$) $$\Vor^w_\X(i) := \left\{u \in \s: \|u-X_i\|^2 - w_i \le \|u - X_j\|^2 - w_j, \; \mbox{for all } \; j \ne i \right\}.$$ Note that if the weights are all zero, this coincides with the usual Voronoi diagram. The computation of the power diagram is a classical problem in computational geometry, for which there exists very efficient software, such as CGAL\footnote{\url{http://www.cgal.org}} or Geogram. Two-dimensional power diagrams can be constructed by an algorithm that runs in time $O(n \log n)$. More generally, $d$-dimensional power diagrams (for $d > 2$) may be constructed by an algorithm with worst case complexity $O(n^{\lceil d/2\rceil })$; see e.g.,~\cite{Auren87},~\cite[Chapter 6]{Auren13}.


Given the power diagram of $(\X,w)$ we can define the {\it power map} $T_\X^w: \s \to \X$ such that $T_\X^w(u) := X_i$ if $u \in \Vor^w_\X(i)$. This map is well-defined $\mu$-a.e. (except on the boundary of the power cells). A weight vector $w \in \R^n$ is called {\it adapted} to $(\mu, \hat \nu_n)$ if for every $i = 1,\ldots, n$, one has $\mu(\Vor^w_\X(i)) = \int_{\Vor^w_\X(i)} d\mu(u) = {n}^{-1}$.~\cite[Theorem 2]{Merigot11} shows that finding a weight vector adapted to $(\mu, \hat \nu_n)$ amounts to finding the global minimum of the convex function
\begin{equation}\label{eq:L}
L(w) := -\sum_{i=1}^n \left[\frac{w_i}{n} + \int_{\Vor^w_\X(i)} \big(\|u - X_i\|^2 - w_i\big) d \mu(u) \right], \quad \mbox{for }\; w \in \R^n;
\end{equation}
also see~\cite{Auren98}. Moreover,~\cite[Theorem 2]{Merigot11} shows that the power map $T_\X^{\hat w}$, where $\hat w$ is the global minimizer of $L(\cdot)$, is the OT map between $\mu$ and $\hat \nu_n$, i.e., $$\hat Q_n = T_\X^{\hat w}.$$ Thus, we have to minimize $L(\cdot)$ in~\eqref{eq:L} to obtain $\hat w$, which will yield $\hat Q_n$. Note that the gradient and Hessian of $L(\cdot)$ can be easily computed; see e.g.,~\cite{Levy15}. This makes Newton-type algorithms especially attractive in computing ${\hat w}$. 




The potential function $\hat \psi_n$, as defined in~\eqref{eq:Emp_G}, can also be recovered from the above optimization problem. Let $\hat h_i := \frac{1}{2}(\hat w_i - \|X_i\|^2)$, for $i=1,\ldots, n$. Now, we can easily see that the convex function thus defined by~\eqref{eq:Emp_G} has gradients that coincides with the quantile map $\hat Q_n$; see e.g.,~\cite[Section 3.4]{Merigot11}. The computation of the dual potential $\hat \psi_n^*$ (as defined in~\eqref{eq:Sup}) and the empirical rank map $\hat R_n$ now follows easily; see Remark~\ref{rem:hat-R_n} (in Section~\ref{sec:Q-R-Prop}) 
for the details.
\begin{figure}
\includegraphics[scale=0.9]{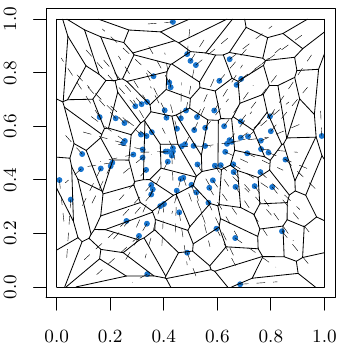}\hspace{0.3in}
\includegraphics[scale=0.9]{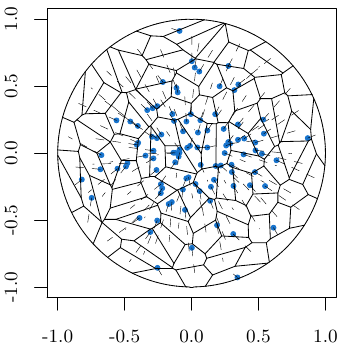}
\caption{Left plot: Shows the cell decomposition of $\s = [0,1]^2$ (each with area $1/n$ where $n=100$) induced by the estimated quantile function $\hat Q_n$   (w.r.t.~$\mu = $ Uniform$([0,1]^2)$) where the data points are drawn i.i.d. from $N_2((0,0),I_2)$ (and appropriately scaled to lie in $[0,1]^2$) along with dashed lines indicating which cell corresponds to which data point. Right plot: Shows the corresponding cell decomposition of $\s = B_1(0)$ --- the ball of radius 1 around $(0,0) \in \R^2$ --- induced by $\hat Q_n$ (w.r.t.~$\mu = $ Uniform$(B_1(0))$).}
\label{fig:4Plots}
\end{figure} 

\begin{remark}[Computation of the sample ranks]\label{rem:Rank_n}
For computing the sample rank $\hat R_n(X_i)$, for $i=1,\ldots, n$, we advocate the use of a randomized choice where we define $\hat R_n(X_i)$ as any point in the set $W_i(\hat h)$, chosen according to the probability measure in~\eqref{eq:Rank-X_i}. When $\mu$ is the uniform distribution on a convex polytope $\s$ (e.g., $[0,1]^d$), this computation is especially simple as then $W_i(\hat h)$ is also a convex polytope whose vertices are already provided by our algorithm, and thus, uniform sampling can be carried out easily (e.g., via rejection sampling on the smallest hyper-rectangle containing $W_i(\hat h)$). In fact, the above approach is much more broadly applicable, as in practice the computer always approximates $\s$ by a convex polytope (see e.g.,~Figure~\ref{fig:4Plots}).
\end{remark}

\begin{figure}
		\includegraphics[scale=0.78]{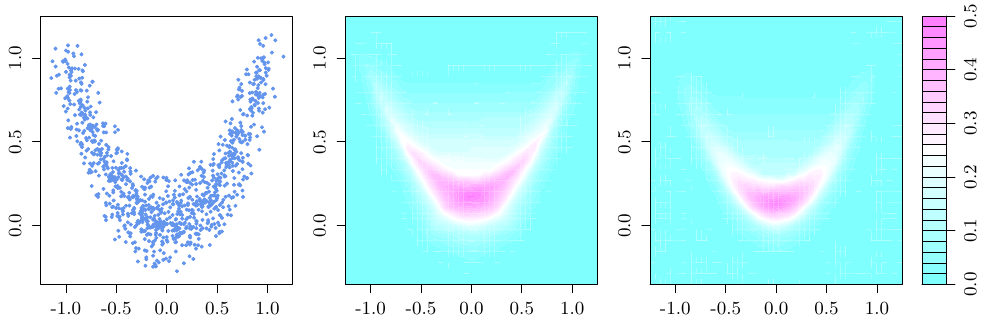}
		\caption{Left panel: A random sample of size $n = 1000$ drawn from the banana-shaped distribution. Middle panel: The estimated depth function --- defined as $\hat D_n(x) := 1/2 - \|\hat R_n(x) - (1/2)\mathbf{1}\|_\infty$ for $x \in \R^d$ (see~\cite{Cher17}), where $\mathbf{1} = (1,1,\ldots, 1) \in \R^d$ --- using $\mu =$ Uniform$([0,1]^2)$. Right panel: The estimated depth function --- defined as $\hat D_n(x) := \pi^{-1}(\theta - \cos \theta \sin \theta)$ where $\theta = \arccos (\|\hat R_n(x)\|)$ --- w.r.t.~$\mu =$ Uniform$(B_1(0))$; see~\cite[Section 5.6]{RR1999}.}
	\label{fig:Banana}
\end{figure}

The two plots in Figure~\ref{fig:4Plots} show the cell decompositions of $[0,1]^2$ and $B_1(0) \subset \R^2$, obtained from the semi-discrete OT problem; see Section~\ref{sec:Add-Plots} 
for more plots of this kind. In particular, we can directly visualize the empirical quantile maps for the two settings. As the (empirical) rank function (taking values in $\R^d$) is a bit difficult to visualize, in Figure~\ref{fig:Banana} we plot the estimated depth functions for the banana-shaped distribution when $n=1000$; cf.~\cite[Figure 2]{Cher17} where the authors motivate the use of multivariate ranks/quantiles based on OT using this data. The banana-like geometry of the data cloud is correctly picked up by the non-convex contours in Figure~\ref{fig:Banana}. We also provide  depth function plots for other distributions in Section~\ref{sec:Add-Plots}. 

\subsection{Comparison with Chernozhukov et al.~\cite{Cher17}, Hallin et al.~\cite{dCHM} and Boeckel et al.~\cite{BSS18}}\label{sec:Comparison}
In the papers~\cite{Cher17},~\cite{dCHM} and~\cite{BSS18} the authors use ideas from the theory of OT to define multivariate quantiles and ranks. Although our approach is similar in spirit to that of~\cite{Cher17} there are subtle and important differences. As opposed to~\cite{Cher17} and~\cite{BSS18}, the quantile map here is defined based on McCann's result (see Theorem~\ref{thm:Brenier}) which extended Brenier's theorem to general probabilities, without the need for a second moment. 
Whereas~\cite{Cher17} studied multivariate quantiles and ranks to obtain notions of statistical depth we study quantiles and ranks to aid us to construct nonparametric goodness-of-fit and mutual independence tests.



The approach to defining multivariate ranks and quantiles proposed and studied in~\cite{dCHM} and~\cite{BSS18} are quite different from ours. In the papers~\cite{dCHM} and~\cite{BSS18}, the authors solve a discrete-discrete OT problem, compared to our semi-discrete approach (further, in~\cite{BSS18} the authors only consider target distributions supported on a compact subset of $\R^d$). Thus, to define the empirical rank map this approach involves the choice of $n$ representative points inside the set $\s$ (that approximates the measure $\mu$) to solve the discrete-discrete OT problem (between the sample data points and the $n$ chosen points). Thus the ``ranks'' of the data points are forced to be the points in the chosen grid. This approach immediately leads to many attractive features for the empirical ranks, e.g., the distribution-freeness of the ranks. However this approach does not automatically give rise to a quantile function (or quantile contours) and special smoothing interpolation is required. In comparison, our approach has the drawback of leading to non-unique ranks at the data points. In a sense, our approach yields an elegant and useful notion of quantiles while the approach of~\cite{dCHM} (and~\cite{BSS18}) yields a notion of ranks with  attractive properties.

\section{Uniform convergence of empirical quantile and rank maps}\label{sec:UnifConv}
The rank and quantile functions in one dimension enjoy many interesting asymptotic properties. For example, if $X_1, \ldots , X_n\sim \nu$, where $\nu$ is a distribution on $\RR$, then by the Glivenko-Cantelli theorem, the empirical rank function (which is the empirical distribution function when $d=1$) converges uniformly to the population rank function a.s. Similarly, for $d=1$, the empirical quantile function converges uniformly (on compacts $[a,b] \subset (0,1)$) to the population quantile function, when the underlying distribution function is continuous and strictly increasing. One may wonder if such results also hold for the multivariate empirical quantile/rank maps studied in this paper. In Theorem~\ref{thm:GCProp} below we show that this is indeed the case. 

Suppose that $\nu$ is absolutely continuous with support $\Y \subset \RR^d$; here $\nu$ is the target distribution. Let $\mu$ be an absolutely continuous distribution supported on a compact convex set $\s \subset \R^d$. Let $Q$ and $R$ be the quantile and rank maps of $\nu$ (w.r.t.~$\mu$); as in~\eqref{eq:Quantile} and~\eqref{eq:Rank} respectively. 
Let $X_1, X_2, \ldots, X_n \stackrel{i.i.d.}{\sim} \nu$. Let $\{\hat{\nu}_n\}_{n\ge1}$ be a sequence of random probability distributions (computed from $X_1,\ldots, X_n$) such that $\hat{\nu}_n$ converges weakly to $\nu$ a.s., i.e.,
\begin{equation}\label{eq:Conv-Emp-Meas}
\hat{\nu}_n \stackrel{d}{\to} \nu \quad \mbox{a.s.}
\end{equation}
We can take $\hat{\nu}_n$ to be the empirical distribution obtained from the first $n$ data points, i.e., $\hat{\nu}_n = \frac{1}{n} \sum_{i=1}^n \delta_{X_i}$; in this case we know that~\eqref{eq:Conv-Emp-Meas} holds (see e.g.,~\cite[Theorem~11.4.1]{Dudley}). Denote the multivariate quantile/rank functions for $\hat{\nu}_n$ by $\hat{Q}_n$ and $\hat{R}_n$. In particular, when the underlying potential functions (see Definition~\ref{def:Quantile}) are not differentiable, we define $\hat{Q}_n$ and $\hat{R}_n$ to be any point in the corresponding subdifferential set. The following is a main result of this paper (see Section~\ref{proof:thm-GC} 
for its proof). 

\begin{theorem}\label{thm:GCProp}
Consider the notation introduced above and suppose that~\eqref{eq:Conv-Emp-Meas} holds. Suppose that $Q:\mathrm{Int}(\s) \to \mathrm{Int}(\Y)$ is a homeomorphism\footnote{See Proposition~\ref{thm:Q_prop} for sufficient conditions.}. 
Let $K_1\subset \mathrm{Int}(\s)$ and $K_2\subset \mathrm{Int}(\Y)$ be any two compact sets. 
\begin{enumerate}
\item[(a)]  Then,  
\begin{equation}\label{eq:GC}
\sup_{u\in K_1} \|\hat{Q}_n(u) - Q(u)\| \stackrel{a.s.}{\rightarrow} 0. 
\end{equation}  
\item[(b)] Further,
\begin{align}\label{eq:GC-Rank}
\sup_{x\in K_2} \|\hat{R}_n(x) - R(x)\| \stackrel{a.s.}{\rightarrow} 0.
\end{align}

\item[(c)] Suppose that $\mathcal{S}$ is a strictly convex compact set (as in Definition~\ref{defn:St-Cvx}). Let $\{\lambda_n\}_{n\ge 1}\subset \RR$ be a sequence such that $\lambda_n\to \infty$ as $n\to \infty$. Then, 
\begin{align}\label{eq:GC-RankFiner}
\sup_{x\in \RR^d} \|\hat{R}_n(x) - R(x)\| \stackrel{a.s.}{\rightarrow} 0, \qquad \mbox{and } \hspace{.75in} \\
\qquad \lim_{\lambda_{n}\to \infty}\hat{R}_n(\lambda_n x) \stackrel{a.s.}{=} \argmax_{v \in \s} \langle x, v\rangle, \quad \mbox{for all}\;\; x\in \RR^d.  \label{eq:Asymptot}
\end{align}


\end{enumerate}
\end{theorem}

 Theorem~\ref{thm:GCProp}-$(a)$ (i.e.,~\eqref{eq:GC}) extends the uniform convergence of the empirical quantile function (on compacts in the interior of $[0,1]$) beyond $d=1$. Theorem~\ref{thm:GCProp}-$(b)$ (i.e.,~\eqref{eq:GC-Rank}) shows the uniform convergence of the estimated rank map on any compact set inside $\mathrm{Int}(\mathcal{Y})$. One may notice that Theorem~\ref{thm:GCProp}-$(a)$ and~$(b)$ improve upon the result of \cite[Theorem 3.1]{Cher17} where the authors prove a similar convergence result for the estimated quantile/rank maps under the assumption of compactness of $\Y$. 
 In~\cite[Theorem 2.2.1]{dCHM} a result similar to~\eqref{eq:GC-RankFiner} is given for the empirical rank map arising from a discrete-discrete OT problem, when the reference measure is the spherical uniform distribution; also see~\cite[Theorem 2.3]{BSS18} for a similar result where the authors only consider a compactly supported $\nu$. In~\cite[Proposition 6]{Zemel2019} the authors prove a local uniform convergence result for the empirical quantile map, under additional finite second moment assumptions on $\nu$. Theorem~\ref{thm:GCProp}-$(c)$ (see \eqref{eq:GC-RankFiner}) can be thought of as the proper generalization of  the Glivenko-Cantelli theorem beyond $d=1$ where we show the a.s.~convergence of the estimated rank map uniformly over the whole of $\R^d$. 

To prove Theorem~\ref{thm:GCProp}, one needs to develop tools that deal with convergence of (sub)-gradients of a sequence of convex functions and their Legendre-Fenchel duals. These tools are summarized in three deterministic lemmas in Section~\ref{sec:UnifConv-App} 
--- Lemmas~\ref{KeyLemma},~\ref{KeyLemma2} and~\ref{lem:UnifConv} --- and could be of independent interest.


 \begin{remark}[On the sufficient condition for~\eqref{eq:GC-RankFiner}]\label{rem:Suff-Cond} In~\eqref{eq:GC-RankFiner} we show that the empirical rank map converges to the population rank function uniformly on $\R^d$, under the strict convexity assumption on $\s$. This sufficient condition is certainly satisfied, for example, when $\s$ is the unit ball in $\R^d$, i.e., $\s = B_1(0)$. Unfortunately when $\s = [0,1]^d$, this condition is not satisfied. 
\end{remark}

\begin{remark}[Necessity of $Q$ to be a homeomorphism]
One of the main assumptions in Theorem~\ref{thm:GCProp} is that the population quantile $Q$ is a homeomorphism; for $d=1$ this corresponds to assuming that the distribution function is continuous and strictly increasing. 
It is actually a necessary condition for showing the uniform convergence of $\hat{Q}_n$ (the sample quantile function) to $Q$; in fact, more generally, for a sequence of (sub)-gradients of convex functions. To see this, consider the example of a sequence of convex functions $\phi_n:\RR \to \RR$ defined as $\phi_n(x) := (x^2+n^{-1})^{1/2}$. 
As $n\to \infty$, $\phi_n(x)$ converges pointwise to $\phi(x) := |x|$. However, the subdifferential set of the function $\phi(x)$ at $x=0$ is equal to $[-1,1]$ whereas $\phi^{\prime}_n(0)=0$ for all $n\geq 1$. Hence, $\phi^{\prime}_n(\cdot)$ does not converge uniformly on any compact set containing $0$. 
\end{remark}

\begin{remark}[When is $Q$ a homeomorphism?]\label{rem:Q-Homeo} In Proposition~\ref{thm:Q_prop} we provide a sufficient condition on the density of $\nu$, supported on a convex set, which ensures that the quantile map $Q$ will be a homeomorphism; also see Remarks~\ref{rem:Cvx-S-Y} and~\ref{rem:Cond-11}. Recently, in~\cite[Proposition~4.5 and Corollary~4.6]{KM17} some results are provided that show that $Q$ can be a homeomorphism even when the support of $\nu$ is a union of convex domains. 
\end{remark}

\begin{remark}[Connection to~\cite{VS18}]
The recent paper~\cite{VS18} implies ``graphical convergence" of the estimated quantile maps (see~\cite[Theorem 4.2 and Corollary 4.4]{VS18}). Their result does not need absolute continuity of $\nu$ and no restrictions are placed on the supports of the measures $\mu$ and $\nu$. However, graphical convergence, which implies a form of local uniform convergence, is weaker than uniform convergence on compacta stated in Theorem~\ref{thm:GCProp}. Moreover, since the sample rank map $\hat{R}_n$ is not strictly a transport map, it is not clear if~\cite{VS18} implies any notion of convergence for $\hat{R}_n$. 
\end{remark}

\section{Rate of convergence of empirical quantile/rank maps}\label{sec:Gl-Lo-Rate}
In this section we study the global and local rates of convergence of the empirical quantile/rank maps. Section~\ref{sec:Global-Rate} provides upper bounds on the global $L_2$-risk of the empirical quantile map whereas Section~\ref{sec:Rank-Rate} provides analogous results for the empirical rank map. In Section~\ref{sec:RateSec} we provide a result that 
gives a local uniform rate of convergence for the empirical quantile/rank maps.
\subsection{Global rate of convergence for the empirical quantile map $\hat Q_n$} \label{sec:Global-Rate}

We first state a lemma (see Lemma~\ref{ppn:RateProp2} below; proved in Section~\ref{pf:RateProp2}) 
that upper  bounds the $L_2$-distance between two OT maps using the difference of the corresponding 2-Wasserstein distances (and a remainder term). Note that the {\it 2-Wasserstein distance} between $\mu$ and $\nu$ is defined as $$W_2(\mu, \nu) := \inf_{\pi \in \Pi(\mu,\nu)} \Big(\int \|u-x\|^{2} d \pi(u,x) \Big)^{{1}/{2}},$$ where $\Pi (\mu, \nu)$ denotes the collection of all joint distributions (couplings) $\pi$ with marginal distributions $\mu$ and $\nu$; see Section~\ref{sec:OT} 
for  more details. 

\begin{lemma}\label{ppn:RateProp2}
Let $\mu$, $\nu$ and $\tilde \nu$ be three probability measures on $\RR^d$ such that $\int \|x\|^2 d\mu(x) <+\infty$, $\int \|x\|^2 d\nu(x) <+\infty$ and $\int \|x\|^2 d {\tilde \nu}(x) <+\infty$. Also, let $\psi$ and $\tilde \psi$ be two convex functions such that $\nabla \psi \# \mu = \nu$ and $\nabla {\tilde \psi} \# \mu = \tilde \nu$ respectively.  Suppose $\psi^*:\R^d \to \R \cup \{+\infty\}$, the Legendre-Fenchel dual of $\psi$, is strongly convex with parameter $\lambda >0$. Then, letting $g(x) :=  \frac{\|x\|^2}{2} - \psi^*(x)$, 
\begin{equation}\label{eq:Global-Diff}
\int  \|\nabla {\tilde \psi} -  \nabla \psi  \|^2 d \mu  \le \frac{1}{\lambda} \left[ \left\{W_2^2(\mu,\tilde \nu) - W_2^2(\mu, \nu)\right\} + 2 \int g \, d (\nu - \tilde \nu)\right].
\end{equation}
\end{lemma}
The above lemma, which is of independent interest, gives a quantitative stability estimate for OT maps in the semi-discrete setting. Although the stability of OT maps has recently been studied by many authors (see e.g.,~\cite{Gigli2011, Berman2018, LN2020, HR19}) we could not find such an explicit upper bound, under such minimal assumptions on the underlying distributions. Moreover, as we illustrate in Theorem~\ref{thm:Q-Rate} below (proved in Section~\ref{pf:Q-Rate}), 
Lemma~\ref{ppn:RateProp2} can be used to obtain rates for OT maps that are strictly better than those obtained in~\cite[Theorem 1.1]{Berman2018} and~\cite[Section 4]{LN2020}. It is worth pointing out that the starting point of the proof of Lemma~\ref{ppn:RateProp2} is based on an observation in~\cite[Proposition 3.3]{Gigli2011}. 

\begin{theorem}\label{thm:Q-Rate}
Let $\mu$ be an absolutely continuous probability measure supported on a compact convex set $\mathcal{S} \subset \R^d$. 
Let $X_1,\ldots, X_n \stackrel{i.i.d.}{\sim} \nu$, where $\nu$ is an absolutely continuous distribution on $\R^d$ with population quantile map $Q \equiv \nabla \psi$ (see~\eqref{eq:Quantile}); here $\psi$ is a convex function. Suppose that $\psi^*:\R^d \to \R \cup \{+\infty\}$ (the Legendre-Fenchel dual of $\psi$) is strongly convex. Then, for all $n \ge 1$, with $\hat Q_n$ being the empirical quantile map (see~\eqref{eq:SampQ_n}),
\begin{equation}\label{eq:Q-Rate} 
\E \left[\int  \| \hat{Q}_n  - Q \|^2 d \mu  \right] \le C \,  r_{d,n}
\end{equation} 
where $C \equiv C(\mu,\nu) >0$ is a constant that depends on $\mu$ and $\nu$, and $r_{d,n}$ is defined in~\eqref{eq:rate}.  Furthermore, there exists $c>0$ such that, for all $s \ge 0$,
\begin{align}\label{eq:Prob_bound}
\mathbb{P}\Big(\int \|\hat{Q}_n - Q\|^2 d\mu \geq C r_{d,n} +n^{-1/2}s\Big)\leq \exp(-cs^2). 
\end{align}
\end{theorem}

We believe that the above result gives the exact rate of convergence for the $\hat{Q}_n$ when $d > 4$; see~\cite{HR19} where the authors mention ``... In this case, one formally recovers
the rate $n^{-2/d}$ and we conjecture that this is the minimax rate of estimation in the context where the transport map $T_0$ is only assumed to be the gradient of a strongly convex function with Lipschitz gradient...". Note that in Theorem~\ref{thm:Q-Rate} we just assume strong convexity of the dual potential associated with the quantile map $Q$. We would also like to point out here that, even when $d=1$, without some assumptions on $\nu$ it is impossible to derive rates of convergence for $\hat Q_n$ as in~\eqref{eq:Q-Rate}; see e.g.,~\cite{BGU05, BL19}. 

 The left hand side in~\eqref{eq:Q-Rate}  is obviously an upper bound for $\E[W^2_2(\hat \nu_n,\nu)]$ and, as a consequence, $$\E[W^2_2(\hat \nu_n,\nu)] \le \E \left[ \int \|\hat Q_n - Q\|^2 d \mu \right] \le C r_{d,n}.$$ Compare this with~\cite[Theorem 1]{FG15} which yields $\E[W^2_2(\hat \nu_n,\nu)] \le C r_{d,n}$, when $\nu$ has a finite moment of order $q > 4$. Thus, Theorem~\ref{thm:Q-Rate} is an improvement of the result in~\cite{FG15}, under the strong convexity assumption on the potential function $\psi^*$.

The proof of Theorem~\ref{thm:Q-Rate} utilizes the stability result of the empirical quantile map $\hat{Q}_n$ obtained in Lemma~\ref{ppn:RateProp2}. As one may note, Lemma~\ref{ppn:RateProp2} bounds the $L_2$-loss of $\hat{Q}_n$ by the difference (up to a smaller order term) between two $2$-Wasserstein distances, under minimal structural assumptions on the dual potential of $Q$. The rate of convergence in Theorem~\ref{thm:Q-Rate} is then obtained by analyzing the expected value of the difference of the Wasserstein distances using empirical process theory; see e.g., the proof of~\cite[Theorem 2]{Chizat2020}. 

%

\subsection{Global rate of convergence for the empirical rank map $\hat R_n$}\label{sec:Rank-Rate} 

Deriving a rate of convergence for the multivariate sample rank $\hat R_n$ map is a bit more tricky. Note that $\hat R_n$ is not an OT map per se, but is defined via the Legendre-Fenchel dual of the potential function $\hat \psi_n$ (see~\eqref{eq:SampR_n}). Also, the sample ranks $\hat R_n(X_i)$'s are not uniquely defined (see Section~\ref{sec:Ranks}). In this subsection we consider the randomized choice of the empirical ranks (as in~\eqref{eq:Rank-X_i}). In the following result we give an upper bound on the risk of the sample rank map (see Section~\ref{pf:RateProp1} 
for its proof). 

\begin{theorem}\label{thm:RateProp1}
Let $\mu$ be an absolutely continuous probability measure supported on a compact convex set $\mathcal{S} \subset \R^d$. 
Let $X_1,\ldots, X_n$ be i.i.d.~from an absolutely continuous distribution $\nu$ on $\R^d$ with compact support and rank map $R \equiv \nabla \psi^*$, where $\psi:\s \to \R$ is assumed to be strongly convex. For $i =1,\ldots, n$, let $\hat R_n(X_i)$ be defined as in~\eqref{eq:Rank-X_i}. Then, for all $n \ge 1$,
\begin{equation}\label{eq:Rate-Rank}
\E \left[ \frac{1}{n} \sum_{i=1}^n \big\| \hat R_n(X_i) - R(X_i)\big\|^2 \right] \le K \,  r_{d,n}
\end{equation} 
where $K \equiv K(\mu,\nu) >0$ depends on $\mu$ and $\nu$,  and $r_{d,n}$ is defined in~\eqref{eq:rate}.
\end{theorem}
The expectation on the left side of~\eqref{eq:Rate-Rank} averages over the external randomization in the definition of the empirical ranks. Note that for Theorem~\ref{thm:RateProp1} to hold we need to assume that $\nu$ has compact support, in addition to the strong convexity of $\psi$. Although a formal result on the optimality of the upper bound in Theorem~\ref{thm:RateProp1} is beyond the scope of this paper, we believe that the obtained bounds are optimal when $d > 4$. 

\subsection{Local uniform rate of convergence}\label{sec:RateSec}

Theorem~\ref{thm:RateTheo} below (proved in Section~\ref{pf:RateTheo}
), provides a local uniform rate of convergence of the empirical quantile/rank maps. In the following result, when the underlying potential functions are not differentiable, we define $\hat{Q}_n$ and $\hat{R}_n$ to be any point in the corresponding subdifferential sets.


\begin{theorem}\label{thm:RateTheo}
Let $\mu$ be an absolutely continuous distribution with a bounded non-vanishing density supported on a compact convex set $\mathcal{S} \subset \R^d$. Let $X_1, \ldots , X_n$ be i.i.d.~$ \nu$ absolutely continuous and supported on a convex set $\mathcal{Y}\subset \RR^d$ with population quantile map $Q \equiv \nabla \psi$ (see~\eqref{eq:Quantile}), where $\psi$ is a convex function.   Suppose that $Q$ is a homeomorphism from $\mathrm{Int}(\mathcal{S})$ to $\mathrm{Int}(\mathcal{Y})$. Assume that $\psi^{*}$ and $\psi$ are strongly convex inside $\mathcal{Y}$ and $\mathcal{S}$ respectively. Fix $u_0\in \mathrm{Int}(\s)$ and $\delta_0\equiv\delta_0(u_0)>0$ such that $B_{\delta_0}(u_0) \subset \s$ and $B_{\delta_0}(\nabla \psi(u_0))\subset \mathcal{Y}$. Then, there exists a constant $C \equiv C(\mu,\nu, u_0)>0$, depending only on $\mu,\nu$ and $u_0$, such that, for all $n\geq 1$, 
\begin{align}\label{eq:RateAc}
\qquad \mathbb{E}\left[\sup_{u\in B_{\delta_0/3}(u_0)} \|\hat{Q}_n(u)- Q(u)\|\right]\leq C \,r^{\frac{1}{d+2}}_{d,n},
\end{align} 
and  
\begin{align}\label{eq:RateAc2}
\qquad \mathbb{E}\left[\sup_{x\in B_{\delta_0/6}(\nabla \psi(u_0))} \|\hat{R}_n(x)- R(x)\|\right]\leq C \, r^{\frac{1}{d+2}}_{d,n}.
\end{align}
\end{theorem}
The proof of Theorem~\ref{thm:RateTheo} is built on Proposition~\ref{ppn:RateProp1} 
which connects the local uniform rate of convergence of $\hat{Q}_n$ and $\hat{R}_n$ with the local $L_2$-rate of convergence of $\hat{Q}_n$. Theorem~\ref{thm:Q-Rate} is then used to upper bound this local $L_2$-rate of convergence of $\hat{Q}_n$.

To the best of our knowledge, the above result is the first attempt to study the local uniform behavior of transport maps. However, it is not clear to us whether the above bounds are tight  when $d\ge 2$. We believe that it may be possible to improve our rate of convergence result under further assumptions on $\nu$. We hope to address this in future work.  

\section{Applications to nonparametric testing}\label{sec:Goodness-Fit-Test}
\subsection{Two-sample goodness-of-fit testing in $\R^d$}\label{sec:2S-Goodness-Fit-Test}
Suppose that we observe $X_1,\ldots, X_m$ i.i.d.~$\nu_X$ and $Y_1,\ldots, Y_n$ i.i.d.~$\nu_Y$, where $m, n \ge 1$, and $\nu_X$ and $\nu_Y$ are unknown absolutely continuous distributions on $\R^d$. We also assume that both the samples are drawn mutually independently. In this section we consider the two-sample equality of distribution hypothesis testing problem: 
\begin{equation}\label{eq:2-Sample-Test}
H_0: \nu_X = \nu_Y\quad \mbox{ versus }\quad H_1: \nu_X \ne \nu_Y.
\end{equation}
The two-sample problem for multivariate data has been extensively studied, beginning with the works of~\cite{Weiss60, Bickel68}. Several graph based methods have been proposed in the literature for this problem; see e.g.,~\cite{FR79, Sc86, Rosenbaum05, B15} and the references therein. Also see~\cite{BF04, SR13, EquivRKHS13} for distance and kernel based methods for the two-sample problem. Recently, the theory of OT and Wasserstein distances have been used to construct  goodness-of-fit tests for~\eqref{eq:2-Sample-Test}; see e.g., \cite{BSS18, Ramdas17, delB-AoP-19, Hallin2020, DebSen2019}. In the following we propose a tuning-free method that uses the (estimated) multivariate quantile/rank maps defined in Section~\ref{sec:Q-R}.

Let $\mu$ be an absolutely continuous distribution supported on a compact convex set $\s \subset \R^d$ having a density (w.r.t.~Lebesgue measure), e.g., $\mu = $ Uniform$([0,1]^d)$ or $\mu = $ Uniform$(B_1(0))$.  Let $\hat Q_X$ and $\hat Q_Y$ be the sample quantile maps estimated from the $X_i$'s and $Y_j$'s respectively (w.r.t.~$\mu$). Let $\hat R_{X,Y}$ be the empirical rank map of the pooled sample  $X_1,\ldots, X_m, Y_1,\ldots, Y_n$ (w.r.t.~$\mu$). As in Section~\ref{sec:Ranks}, we define the rank at any data point as a randomized value (as in~\eqref{eq:Rank-X_i}). 
We use the following test statistic for testing~\eqref{eq:2-Sample-Test}:
\begin{align}\label{eq:2-S-TS} 
T_{X,Y} & := \int_\s \|\hat R_{X,Y}(\hat Q_X(u)) - \hat R_{X,Y}(\hat Q_Y(u))\|^2 d \mu(u). 
\end{align} 
Exactly computing $T_{X,Y}$ is possible as the above integral reduces to a finite sum; see Section~\ref{sec:Simul-2S} 
for the details. One can also easily approximate $T_{X,Y}$ using Monte Carlo.

We reject $H_0$ when $T_{X,Y}$ is large. To motivate the form of the above test statistic consider the one-sample Cram\'{e}r-von Mises statistic when $d=1$. Let $\mathbb{F}_n$ be the empirical distribution of the data (when $d=1$) and $F$ be the true distribution function (assumed to be absolutely continuous). Then the Cram\'{e}r-von Mises statistic can be written as $$ \int_{\R} \{\mathbb{F}_n(x) - F(x)\}^2 dF(x) = \int_0^1 \{\mathbb{F}_n(F^{-1}(u)) - u\}^2 du.$$ Indeed,~\eqref{eq:2-S-TS} is similar to the right side of the above display; however as we are now in the two-sample case, $F^{-1}$ is unknown and is replaced by the sample quantile function.   

The connection to the Cram\'{e}r-von Mises statistic above immediately raises the following question: Is $T_{X,Y}$ distribution-free under $H_0$ (as the Cram\'{e}r-von Mises statistic when $d=1$)? Unfortunately, we do not know the exact answer to this question. 
In the following lemma (proved in Section~\ref{pf:Distfree}) 
we show that $\hat R_{X,Y}(\hat Q_{X}(U))$ and $\hat R_{X,Y}(\hat Q_{Y}(U))$ (as in~\eqref{eq:2-S-TS}) are both marginally distribution-free and distributed as $\mu$ under $H_0$.

\bl\label{lem:Distfree} 
Suppose that $\nu_X = \nu_Y$. Then $\hat R_{X,Y}(\hat Q_{X}(U)) \sim \mu$ and $\hat R_{X,Y}(\hat Q_{Y}(U)) \sim \mu$, and hence their distributions do not depend on $\nu_X \equiv \nu_Y$. 
\el

\begin{remark}[Finding the critical value of $T_{X,Y}$]\label{rem:Crit-Value} 
Although we have shown (in Lemma~\ref{lem:Distfree}) that $\hat R_{X,Y}(\hat Q_{X}(U)) \sim \mu$ and $R_{X,Y}(\hat Q_{Y}(U)) \sim \mu$ (and thus both quantities are distribution-free) it is not immediately clear if the test statistic $T_{X,Y}$ in~\eqref{eq:2-S-TS} is distribution-free, under $H_0$. In Section~\ref{sec:Simul-2S} 
we provide simulation evidence that suggests that a properly normalized version of $T_{X,Y}$ may be asymptotically distribution-free, at least when $d=2$. In any case, the critical value of the test can always be computed by conditioning on the observed samples and using the following permutation principle: Under $H_0$, $X_1,\ldots, X_m, Y_1, \ldots, Y_n$ are i.i.d.~and thus we can consider any partition of the $m+n$ data points into two sets of sizes $m$ and $n$ and recompute our test statistic to simulate its null distribution. This is indeed the most common approach in these nonparametric testing problems as it avoids the need to use asymptotic distributions and leads to exact tests; see e.g.,~\cite{Hoeffding1952, Kim2020, TSH-2005}.
\end{remark}

The following result (proved in Section~\ref{pf:Power1}) 
shows that our proposed test has asymptotic power 1 when $\nu_X \ne \nu_Y$.

\begin{proposition}[Consistency]\label{lem:Power1} 
Suppose that $H_0: \nu_X = \nu_Y \equiv \nu$ holds. Assume that $\nu$ is supported on a domain $\mathcal{Y} \subset \R^d$ such that the quantile map $Q:\mathrm{Int}(\s)\to \mathrm{Int}(\Y)$ is a homeomorphism. Also, assume that $m,n\to \infty$ such that $\frac{m}{m+n}\to \theta\in (0,1)$. Then, under $H_0$, as $m,n\to \infty$,
\begin{align}\label{eq:Consis1}
T_{X,Y} \stackrel{a.s.}{\longrightarrow} 0.
\end{align}
Now, suppose that $X_1, \ldots , X_m\stackrel{i.i.d}{\sim} \nu_X$ and $Y_1, \ldots , Y_n\stackrel{i.i.d}{\sim} \nu_Y$ where $\nu_X \ne \nu_Y$ are two distinct probability measures supported on domains $\mathcal{Y}_X$ and $\mathcal{Y}_Y$ respectively. Denote the quantile maps of the distributions $\nu_X$, $\nu_Y$ and $\theta \nu_X+(1-\theta)\nu_Y$ by $Q_X$, $Q_Y$ and $Q_{X,Y}$ respectively. Assume that $Q_X,Q_Y$ and $Q_{X,Y}:\mathrm{Int}(\s)\to \mathrm{Int}(\Y_X \cup \Y_Y)$ are homeomorphisms. Then, as $m,n\to \infty$,
\begin{equation}\label{eq:Consis2}
T_{X,Y} \stackrel{a.s.}{\longrightarrow} c := \int_{\s}\|R_{X,Y}(Q_{X}(u))-R_{X,Y}(Q_{Y}(u))\|^2 d\mu(u), 
\end{equation}  
where $c>0$ and $R_{X,Y}$ is the rank function for the measure $\theta \nu_X+(1-\theta)\nu_Y$. 
\end{proposition}

The following two results (proved in Sections~\ref{pf:Two-S-Rate} and~\ref{pf:Two-S-Alt}) 
provide rates of convergence of $T_{X,Y}$ under the null and alternative hypotheses. The proofs of these results are built on Theorems~\ref{thm:Q-Rate} and~\ref{thm:RateProp1}.

\begin{proposition}\label{thm:Two-S-Rate}
Suppose that $H_0:\nu_{X} = \nu_{Y} \equiv \nu$ holds. Assume that $\nu$ is absolutely continuous and supported on a compact domain $\mathcal{Y}\subset \mathbb{R}^{d}$. Further, we assume that the convex potential $\psi$ of the quantile map $Q$ of $\nu$ (w.r.t.~$\mu$) is strongly convex. Under $H_0$, if $\min\{m,n\}/(m+n) \ge \theta \in (0,1)$, then 
\begin{equation}\label{eq:TBound}
\mathbb{E}[T_{X,Y}]\leq C \, r_{d,m+n}
\end{equation}
where $C \equiv C(\mu,\nu,\theta) >0$ depends on $\mu, \nu$ and $\theta$, and $r_{d,n}$ is defined in~\eqref{eq:rate}. 
\end{proposition}
 
\begin{proposition}\label{thm:Two-S-Alt}
Suppose that $\nu_{X} \neq \nu_{Y}$, where $\nu_{X}$ and $\nu_{Y}$ are compactly supported. Recall the notation from Proposition~\ref{lem:Power1}. For convenience, we will assume that the pooled sample size $N$ is fixed and that $m|N \sim $ Binomial$(N,\theta)$, where $\theta \in (0,1)$. Further, we assume that the convex potential functions $\psi_X, \psi_Y, \psi_{X,Y}$ of the quantile maps $Q_X, Q_Y, Q_{X,Y}$ are strongly convex. Then, we have 
\begin{align}\label{eq:TBound}
\mathbb{E}[|T_{X,Y}-c|]\leq \, C \, r_{d,N}^{1/2}
\end{align}
where $C \equiv C(\mu,\nu_X,\nu_Y,\theta) >0$, and $r_{d,n}$ is defined in~\eqref{eq:rate}. 
\end{proposition}
A detailed study of the finite sample performance and the asymptotic weak limit of the above test in beyond the scope of the present paper. We plan to pursue this in a future paper. As mentioned before, $T_{X,Y}$ is inspired by the form of the Cram\'{e}r-von Mises (one-sample) goodness-of-fit statistic. One can, of course, use other test statistics based on the empirical quantile/rank maps for testing~\eqref{eq:2-Sample-Test}. A key observation for constructing such tests is to realize that, under $H_0$, $\hat R_{X,Y}(X_1), \ldots, \hat R_{X,Y}(X_m), \hat R_{X,Y}(Y_1), \ldots, \hat R_{X,Y}(Y_n)$ are {\it exchangeable} and are all marginally distributed as $\mu$.

\subsection{Mutual independence testing}\label{sec:IndepTest}
Let $X = (X^{(1)}, \ldots, X^{(k)}) \sim \nu$ be a random vector in $\R^d$ where $k \ge 2$ and $X^{(j)}\sim \nu_j$ is a random vector in $\R^{d_j}$, for $j = 1, \ldots, k$, with $\sum_{j=1}^k d_j = d$. In this subsection we consider the problem of testing the mutual independence of $X^{(1)},\ldots, X^{(k)}$. Specifically, we consider testing whether $\nu$ is equal to the product measure $\nu_1 \otimes \ldots \otimes \nu_k$, for some $\nu_1,\ldots, \nu_k$, i.e., 
\begin{equation}\label{eq:Ind-Test}
H_0: \nu = \nu_1 \otimes \ldots \otimes \nu_k  \qquad \mbox{ versus } \qquad H_1: \nu \ne \nu_1 \otimes \ldots \otimes \nu_k,
\end{equation} 
when we observe i.i.d.~data from $\nu$. This is again a fundamental problem in statistics and there has been many approaches investigated in the literature; see e.g.,~\cite{Blomqvist50, BlumEtAl61},~\cite[Chapter 8]{HW99}, ~\cite{DebSen2019} and the references therein. The use of kernel (see e.g.,~\cite{GrettonKernelMeasInd05, EquivRKHS13, Lyons13, Pfister18}) and distance covariance (see e.g.,~\cite{SzekelyBDCov09, SzekelyCorrDist07, 4-Axioms-MS19, SR13, Shi2019,  DebSen2019}) based methods have become very popular for this problem. Also see~\cite{Berrett19,Weihs18} and the references therein for some recent other approaches to testing~\eqref{eq:Ind-Test}. We use our multivariate quantile and rank functions to construct a tuning parameter-free consistent test for~\eqref{eq:Ind-Test}.

For simplicity of notation, let us assume that $k=2$. As we will see, the extension to $k > 2$ is straightforward. Let $\{Z_i \equiv (X_i,Y_i): 1 \le i \le n\}$ be i.i.d.~{$\nu$},  assumed to be absolutely continuous on $\R^{d_X} \times \R^{d_Y}$; here $d_X, d_Y \ge 1$ and $d_X + d_Y =d$. Further, we assume that $X \sim \nu_X$ and $Y \sim \nu_Y$. We want to test the hypothesis of mutual independence between $X$ and $Y$, i.e., 
\begin{equation}\label{eq:Test-2}
\qquad {H_0: \nu = \nu_X \otimes \nu_Y}\qquad \mbox{ versus }\qquad {H_1: \nu \ne \nu_X \otimes \nu_Y}.
\end{equation}
Let ${\mu_X }= $ Uniform($[0,1]^{d_X}$),  $ {\mu_Y} = $  Uniform($[0,1]^{d_Y}$) and let {$\mu := \mu_X \otimes \mu_Y$} =  Uniform($[0,1]^d$). We define $\hat R: \R^d \to [0,1]^d$ and $\hat Q:[0,1]^d \to \R^d$ to be the empirical rank and quantile maps of the joint sample $(X_1,Y_1),\ldots, (X_n, Y_n)$. Let ${\hat R_X}: \R^{d_X} \to \R^{d_X}$ be the empirical rank map of $X_1,\ldots, X_n$; similarly let {$\hat R_Y: \R^{d_Y} \to \R^{d_Y}$} be the sample {rank map} obtained from $Y_1,\ldots, Y_n$. Define ${\tilde R := (\hat R_X, \hat R_Y)}:\R^d \to [0,1]^d$. We consider the following test statistic: 
\begin{equation}\label{eq:Ind-TS}
T_{n} := {\int_{[0,1]^d} \|\hat R(\hat Q(u)) - \tilde R(\hat Q(u))\|^2 d u } = \frac{1}{n} \sum_{i=1}^n \|\hat R(Z_i) - \tilde R(Z_i)\|^2.
\end{equation}
Note that the above integral reduces to a finite average as $\hat Q(\cdot)$ can only take $n$ distinct values a.s. We reject the null hypothesis in~\eqref{eq:Test-2} when $T_{n}$ is large. As in Section~\ref{sec:2S-Goodness-Fit-Test}, the critical value of the test can be computed using the permutation principle: We take a random permutation of $\sigma$ of $\{1,\ldots, n\}$ and consider the permuted data set $\{(X_i, Y_{\sigma(i)})\}_{i=1}^n$. The (conditional) null distribution of $T_n$ can be computed by considering the permutation distribution of $T_n$ (i.e., computed from the data $\{(X_i, Y_{\sigma(i)})\}_{i=1}^n$, as $\sigma$ varies). 

The following result, proved in Section~\ref{sec:Pf:lem:Power2}, 
describes the asymptotic behavior of the proposed test statistic under the null and alternative hypotheses; in particular, it shows that the power of the test converges to 1, as the sample size $n$ increases. 
 
\begin{proposition}[Consistency]\label{lem:Power2} 
We have $\hat R(\hat Q(U)) \sim \mu$, where $U \sim \mu = $ Uniform$([0,1]^d)$. Suppose $H_0$ holds in~\eqref{eq:Test-2}, i.e., $\nu = \nu_X \otimes \nu_Y$. Then, $\tilde  R(\hat Q(U)) \sim \mu$. Assume further that $\nu_{X}$ and $\nu_{Y}$ are two probability measures supported on the domains $\mathcal{Y}_X\subset \RR^{d_{X}}$ and $\mathcal{Y}_{Y}\subset \RR^{d_{Y}}$ respectively. Denote the quantile maps of the measures $\nu_{X}$, $\nu_{Y}$ and $\nu$ w.r.t. the measures $\mathrm{Uniform}([0,1]^{d_X})$, $\mathrm{Uniform}([0,1]^{d_Y})$ and $\mathrm{Uniform}([0,1]^{d})$ by $Q_X$, $Q_Y$ and $Q$  respectively, where $d=d_X+d_Y$. Assume that $Q_X:(0,1)^{d_X}\to\mathrm{Int}(\mathcal{Y}_X)$,  $Q_X:(0,1)^{d_Y}\to\mathrm{Int}(\mathcal{Y}_Y)$ and  $Q:(0,1)^{d}\to\mathrm{Int}(\mathcal{Y}_X \times \mathcal{Y}_Y)$ are homeomorphisms. Then, under $H_0$, as $n\to\infty$,
\begin{align}\label{eq:Consis21}
T_{n} \stackrel{a.s.}{\longrightarrow} 0.
\end{align}
Now suppose that $\nu \ne \nu_X \otimes \nu_Y$. Let $\bar{R}=(R_X, R_Y)$ where $R_X$ and $R_Y$ are the rank maps of $\nu_X$ and $\nu_Y$ respectively. Then, 
$$T_{n} \stackrel{a.s.}{\longrightarrow} c:= \int_{[0,1]^d} \|u - \bar{R}(Q(u))\|^2 \, d u, \quad \mbox{ as  } n \to \infty.$$
\end{proposition}


The following two results (proved in Sections~\ref{pf:Rate2} and~\ref{pf:Rate2-Alt}) 
provide rates of convergence of $T_{n}$ under the null and alternative hypotheses.

\begin{proposition}\label{lem:Rate2} 
Suppose $H_0: \nu = \nu_X \otimes \nu_Y$ holds, where $\nu_X$ and $\nu_Y$ are compactly supported absolutely continuous distributions on $\R^{d_X}$ and $\R^{d_Y}$ with quantile maps $Q_X$ and $Q_Y$. Further, assume that the convex potentials $\psi_X$ and $\psi_Y$ of $Q_X$ and $Q_Y$ are strongly convex. Then, for $d = d_X + d_Y$, 
$$\E[T_n] \le C \,r_{d,n} $$ where $C \equiv C(\mu,\nu) >0$ depends on $\mu$ and $\nu$, and $r_{d,n}$ is defined in~\eqref{eq:rate}. 
\end{proposition}

\begin{proposition}\label{lem:Rate2-Alt} 
Suppose that $\nu \ne \nu_X \otimes \nu_Y$, where $\nu_{X}$ and $\nu_{Y}$ are compactly supported. Recall the notation from Proposition~\ref{lem:Power2}. Further, we assume that the convex potential functions $\psi_X, \psi_Y, \psi$ of the quantile maps $Q_X, Q_Y, Q$ are strongly convex functions. Then, we have 
\begin{align}\label{eq:TBound2}
\mathbb{E}[|T_{n}-c|]\leq \, Cr^{1/2}_{d,n}.
\end{align}
where $C \equiv C(\mu,\nu) >0$ and $r_{d,n}$ is defined in~\eqref{eq:rate}. 
\end{proposition}

\section*{Acknowledgements}
The authors are extremely grateful to Peng Xu for creating the R-package \url{https://github.com/Francis-Hsu/testOTM} (see~\cite{OTM}) for the computation of all the estimators studied in this paper. In particular, all the plots in the paper are obtained from his R-package. The second author would like to thank Nabarun Deb, Adityanand Guntuboyina, Marc Hallin, Johan Segers for helpful discussions. The authors also acknowledge the numerous insightful comments by the two anonymous referees that helped improve the paper.

%

\section*{Appendices~\ref{sec:Proofs}--\ref{Appendix-B}}
\appendix
This Appendices~\ref{sec:Proofs}--\ref{Appendix-B} contain proofs of all results in the main paper, other auxiliary results (with their proofs) alluded to in the main paper 
and further discussions. 

Appendix~\ref{sec:OT} gives a brief introduction to the theory of optimal transportation. 
In Appendix~\ref{sec:Q-R-Prop} we describe some important properties of the defined multivariate quantile and rank functions $Q(\cdot)$ and $R(\cdot)$ respectively. In Appendix~\ref{sec:Euiv-Inv} we study equivariance/invariance properties of the quantile and rank functions. A few remarks about the empirical quantile and rank functions when $d=1$ is given in Appendix~\ref{sec:d=1}. Some additional plots that supplement the discussion in the main paper are provided in Appendix~\ref{sec:Add-Plots}. Appendix~\ref{sec:AuxRes} states some auxiliary results about convex functions and their subgradients that will be useful in our main proofs. The proofs of all the results in Appendices~\ref{sec:OT}-\ref{sec:AuxRes} is given in Appendix~\ref{pf:Proofs}.

The main results of Section~\ref{sec:Q-R} 
are proved in Appendix~\ref{pf:Q-R}. Proof of results in Section~\ref{sec:UnifConv}  is given in Appendix~\ref{sec:UnifConv-App}. In Appendix~\ref{pf:Global-Rate} we provide the proofs of the results in Sections~\ref{sec:Global-Rate} and~\ref{sec:Rank-Rate}, whereas Appendix~\ref{pf:Loc-RateSec} gives the proof and other theoretical results needed to derive the main result in Section~\ref{sec:RateSec}. The proofs of the results in Section~\ref{sec:Goodness-Fit-Test} on nonparametric testing is given in Appendix~\ref{pf:Goodness-Fit-Test}. In Section~\ref{sec:Simul} we provide simulation studies that illustrate the finite sample and asymptotic behaviors of the test statistics proposed in Section~\ref{sec:Goodness-Fit-Test}. 
Additional technical results stated in Appendix~\ref{proof:thm-GC} are relegated to Appendix~\ref{Appendix-B}.


\section{Some Discussion and Auxiliary results}\label{sec:Proofs}

\subsection{A brief introduction to optimal transport}\label{sec:OT}
\subsubsection{Monge's problem}\label{sec:Monge}
Let $\mu$ and $\nu$ be two Borel probability measures on measurable spaces $(\X,\B_\X)$ and $(\Y,\B_\Y)$, respectively. Let $c: \X \times \Y \to [0,\infty]$ be a measurable loss function: $c(x,y)$ represents the cost of transporting $x$ to $y$. For example, when $\X= \Y = \R^d$, we can take $c : \R^{2d} \to [0,\infty]$ to be the quadratic (or $L_2$) loss function $$c(x, y) = \|x - y\|^2.$$ The goal of optimal transport (Monge's problem) is to find a measurable transport map $T \equiv T_{\mu;\nu} : \X \to \Y$ solving the (constrained) minimization problem
\begin{equation}\label{eq:Meas_Trans}
	\inf_T \int_{\X} c(x,T(x)) d\mu(x)  \qquad \quad \mbox{subject to }\quad T\#\mu  = \nu
\end{equation}
where the minimization is over $T$ (a {\it transport map}), a measurable map from $\X$ to $\Y$, and $T\#\mu$ is the {\it push forward} of $\mu$ by $T$, i.e., 
\begin{align}\label{eq:PushMeasure}
T\#\mu(B) = \mu(T^{-1}(B)), \qquad \mbox{for all } B \in \B_\Y. 
\end{align}
Another equivalent formulation of the constraint in~\eqref{eq:Meas_Trans} is: $\int f d T\#\mu = \int f\circ T d\nu$ for every measurable function $f : \Y \to \R$; see~\cite[Chapter 1]{V03}.

A map $T_{\mu;\nu}$ that attains the infimum in~\eqref{eq:Meas_Trans} is called an {\it optimal transport} map, in short, an optimal transport, of $\mu$ to $\nu$.  Note that the above optimization problem is highly non-linear and can be ill-posed (as no admissible $T$ may exist; for instance if $\mu$ is a Dirac delta measure and $\nu$ is not). Moreover, the infimum in~\eqref{eq:Meas_Trans} may not be attained, i.e., a limit of transport maps may fail to be a transport map. This problem was formalized by the French mathematician Gaspard Monge in 1781 (\cite{Monge1781}) and not much progress was made for about 160 yrs! 

\subsubsection{Kantorovich relaxation: Primal problem}
Let $\Pi (\mu, \nu)$ be the collection of joint distributions (couplings) $\pi$ of random variables $(X, Y) \in \X \times \Y$ such that $X\sim \mu$ \mbox{and} $Y\sim \nu$. Thus any $\pi \in \Pi$ admits $\mu$ and $\nu$ as marginals on $\X$ and $\Y$ respectively (i.e., $\pi(A \times \Y) = \mu(A)$ and $\pi(\X \times B) = \nu(B)$ for all measurable sets $A \subset \X$ and $B \subset \Y$). {\it Kantorovich relaxation} of Monge's problem solves the following optimization problem:
\begin{equation}\label{eq:Kanto}
	\min_{\pi \in \Pi(\mu,\nu)}  \int_\X \int_\Y c(x,y) d\pi(x,y).
\end{equation}
Note that the above is indeed an infinite-dimensional linear program (as the objective is linear and the constraints are linear equalities). It is well-known that Kantorovich's relaxation (i.e.,~\eqref{eq:Kanto}) yields an optimal solution  under the assumption that the cost function $c(\cdot,\cdot) \ge 0$ is l.s.c.; see e.g.,~\cite[Proposition 2.1]{GM96}. 

Further,~\eqref{eq:Kanto} is a relaxation of~\eqref{eq:Meas_Trans} as every transport map $T$ yields a coupling --- take $\pi = (id, T)\#\mu$ which yields $$\int_\X \int_\Y c(x,y) d\pi(x,y) = \int_{\X} c(x,T(x)) d\mu(x).$$ Thus, $$ \min_{\pi \in \Pi(\mu,\nu)}  \int_\X \int_\Y c(x,y) d\pi(x,y) \le \min_{T: T\#\mu = \nu} \int_{\X} c(x,T(x)) d\mu(x).$$ Moreover, under generals assumptions on the measures $\mu$ and $\nu$ (e.g., $\mu$ is absolutely continuous) one can show that Monge's problem and Kantorovich's relaxation have the same minimum value and Kantorovich's relaxation has a solution of the form $\pi = (id, T)\#\mu$ which yields a solution to Monge's problem; see~\cite{Gangbo99}.

The {\it Wasserstein distance} between the two probability measures $\mu$ and $\nu$ (defined on two subsets on $\R^d$ --- $\X$ and $\Y$) w.r.t.~$L_p$-cost function is: 
\begin{align}\label{eq:W_p}
W_p(\mu, \nu) = \inf_{\pi \in \Pi(\mu,\nu)} \Big(\int_{\mathcal{X} \times \Y} \|x-y\|^{p} d \pi(x,y) \Big)^{\frac{1}{p}}.
\end{align}
Clearly, the optimization problem in the definition of  $W_p$ is the same as the Kantorovich relaxation for the $L_{p}$-cost function. It can be shown that, if $\X = \Y$ is compact, then $W_p$, for any $p\geq 1$, metrizes the space of probability measures on $\mathcal{X}$ (see \cite[Chapter~8]{V03}). 

In the recent 20-30 years there has been a lot of interest and progress in this topic of optimal transportation; we refer the interested reader to the books~\cite{V03, V09, Ambro08} for a comprehensive introduction to this fascinating field of mathematics. 

\subsection{Properties of the quantile and rank maps}\label{sec:Q-R-Prop}
In this section we describe some important properties of the defined quantile and rank functions $Q(\cdot)$ and $R(\cdot)$ respectively. 

We start with some notation. For $A \subset \R^d$, we use the following notation: 
\begin{eqnarray*}
\partial f(A) & := &\{x \in \R^d: x \in \partial f(u) \;\; \mbox{for some}\;\; u \in A\}, \qquad \mbox{and} \\ (\partial f)^{-1}(A) &:= & \{u \in \R^d: x \in \partial f(u)\;\; \mbox{for some}\;\; x \in A\}.
\end{eqnarray*}

\begin{remark}\label{rem:dualrem}
For a proper l.s.c.~convex function $f:\R^d \to \R \cup \{+ \infty\}$, using \eqref{eq:Charac-Sub}, one can see that
\begin{align}
\big\{y\in \RR^d: y\in \partial f(x) \text{ for some }x\in \RR^d\big\} &= \big\{y\in \RR^d: x\in \partial f^{*}(y)\text{ for some }x\in \RR^d \big\}\nonumber\\&= \big\{y\in \RR^d: f^{*}(y)< +\infty\big\}
\end{align}
where the last equality follows since $\partial f^{*}(y) = \emptyset$ if and only if $f^{*}(y) = +\infty$. 
\end{remark}

As $Q \in \partial \psi$, and $R \in \partial \psi^*$ are (sub)-gradients of two convex functions that are convex conjugates of each other, we have the following as a direct consequence of Lemma~\ref{lem:SubD}; see Section~\ref{sec:SubD-1} for a proof, and~\cite[Chapter 2]{V03} for related results.
\begin{lemma}\label{lem:SubD-1}
We have 
 \begin{align}\label{eq:EssR-Inv}
x\in \partial \psi(\partial \psi^*(x)), \;\mbox{for} \; x \in \R^d \qquad \mbox{and} \qquad u\in \partial \psi^*(\partial \psi(u)),\;\mbox{for} \;u \in \s.
 \end{align}
Moreover, for every Borel set $B \subset \R^d$, we have $(\partial \psi)^{-1}(B)= \partial \psi^*(B)$. 
\end{lemma}
There is an intimate connection between the quantile map and the celebrated Monge-Amp\`{e}re differential equation; see e.g.,~\cite[Lemma 4.6]{V03} (also see \cite{Ca1,PF14,CF19}). One can find similar results in \cite{Ca1,PF14,CF19}. For completeness, in the following, we present the main connection as a lemma (Lemma~\ref{lem:QRrelation}) and prove it in Section~\ref{sec:QRrelation}; one can find a similar proof in \cite[Theorem~3.6]{PF14}.

\bl\label{lem:QRrelation} 
Suppose that $\mu$ and $\nu$ are absolutely continuous (w.r.t.~Lebesgue measure) with densities $f$ and $g$ respectively. Let $Q$ and $R$ be the quantile and rank maps of $\nu$ (w.r.t.~$\mu$), as defined in~\eqref{eq:Quantile} and \eqref{eq:Rank}. For any Borel set $A\subset \RR^d$, define 
$\rho_{Q}(A) := \int_{Q(A)} dx$. Then, $\rho_{Q}$ satisfies the Monge-Amp\`{e}re differential equation (see~\cite{PF14}), i.e., 
$
\rho_{Q}(A) = \int_{A}\frac{f(x)}{g(Q(x))} dx. 
$
\el

In many statistical applications it is often useful to know the regularity of the quantile and rank maps. Due to Lemma~\ref{lem:QRrelation}, the regularity of our multivariate notions of quantiles and ranks follow from the regularity theory of the solution of the Monge-Amp\'ere equation which has been extensively studied by many authors in the past (see e.g.,~\cite{Cafa90,PF13,PF15,F10,FRV11,GO17}).    


Recall the dual potential function $\psi^*$, as defined in~\eqref{eq:Rank}. The following lemma (proved in Section~\ref{sec:psi-finite}) shows that $\psi^*(x) < \infty$, for all $x \in \R^d$.
\begin{lemma}\label{lem:psi-finite}
$\psi^*(x) < \infty$, for all $x \in \R^d$.
\end{lemma}

Next we illustrate that the rank map (as defined in~\eqref{eq:Rank}) possesses many properties similar to the univariate distribution function. Lemma~\ref{lem:Monotonicity} (proved in Section~\ref{sec:Monotonicity}), states that the one-dimensional projection of the rank map, along any direction, is nondecreasing.

\bl[Monotonicity of the rank map]\label{lem:Monotonicity}
Let $R$ be the rank map of $\nu$ w.r.t. $\mu$.  For $x,y\in \R^d$, define $R_{x,y}:\RR\to \RR$ as 
\begin{align*}
R_{x,y}(t):=(x-y)^{\top}R(tx+(1-t)y), \qquad t \in \R.
\end{align*}
Then, $R_{x,y}(\cdot)$ is a nondecreasing function. 
\el
Recall the definition of empirical quantiles and ranks as in~\eqref{eq:SampQ_n} and~\eqref{eq:SampR_n} respectively. The following result, proved in Section~\ref{Pf:lem:ExtremePt}, expresses the value of $ \hat{\psi}^{*}_n$ at the data points in terms of $\hat {h}_i$ (see~\eqref{eq:Emp_G}); cf.~\cite[Corollary 2.1]{Gu12}.
\begin{lemma}\label{lem:ExtremePt}
Fix $i \in \{1,\ldots, n\}$. Consider $\hat{\psi}^{*}_n$ as defined in~\eqref{eq:Emp_G}. Then, $\hat{\psi}^{*}_n(X_i) = -\hat {h}_i$.
\end{lemma}

It is a simple fact that the Legendre-Fenchel dual of a piecewise affine convex function is also convex piecewise affine. In fact, for $x \in \conv(X_1, \ldots, X_n)$,
\begin{equation}\label{eq:Dual-2}
\hat{\psi}^{*}_n(x) = \min \left\{\sum_{i=1}^n t_i \hat{\psi}^{*}_n(X_i): t_i \ge 0, \sum_{i=1}^n t_i =1, \sum_{i=1}^n t_i X_i = x\right\};
\end{equation}
see e.g.,~\cite[Theorem 2.2.7]{H94}. Thus, $\hat R_n(\cdot)$ a.e.~takes finitely many distinct values (as it is a gradient of a piecewise  affine convex function). 

The following result (proved in Section~\ref{sec:RAltDef}) expresses $\hat{R}_n$, the sample rank map, in an alternate form which was used in~\cite[Definition 3.1]{Cher17}.

\bl\label{lem:RAltDef} 
Consider the multivariate sample rank function $\hat{R}_n$ defined in \eqref{eq:SampR_n}. Then, for all $y\in \RR^d$, 
\begin{equation}\label{eq:R_n-Comp}
\hat{R}_n(y) \in \argmax_{u\in \s}\{\langle u, y\rangle - \hat{\psi}_n(u)\}.
\end{equation}
\el

\begin{remark}[Computation of $\hat{R}_n$]\label{rem:hat-R_n}
To compute $\hat{R}_n(y)$ for arbitrary $y \in \R^d$, we use~\eqref{eq:R_n-Comp}. Recall $\hat{\psi}_n^*$, as  defined via~\eqref{eq:Sup}; i.e., $$\hat{\psi}^*(y)= \sup_{u \in \s} \{\langle u, y\rangle - \hat{\psi}_n(u)\} = \sup_{u \in \s} \left[\langle y, u \rangle -\max_{i=1,\ldots,n}\{u^\top X_i + \hat{h}_i\}\right],$$ where we have used the expression of $\hat{\psi}_n$, given in~\eqref{eq:Emp_G} which involves the $\hat h_i$'s obtained above. First let us discuss the computation of $\hat{\psi}^*(y)$. As the power cells $\{W_i(\hat h)\}_{i=1}^n$ form a partition of $\s$, finding the supremum over $\s$ is equivalent to finding first the suprema over each $W_i(\hat h)$, and then finding the overall supremum. This allows us to write: $$\hat{\psi}_n^*(y)=\max_{i=1,\dotsc,n} \left[\sup_{u\in W_i(\hat h)} \big\{u^\top(y-X_i) \big\}-\hat{h}_i \right].$$
As $W_i(\hat h)$ is encoded as a convex polytope, $\sup_{u \in W_i(\hat h)}\left\{u^\top(y-X_i)\right\}\) is a linear program. By the  vertex principle, its maximizer must be one of the vertices of $W_i(\hat h)$. Thus the supremum in the above display can be computed by finding the maximum value over the finitely many vertices of $W_i(\hat h)$, for each $i$. In the process of computing $\hat{\psi}_n^*(y)$, the value of $\hat{R}_n(y)$ is immediately obtained via~\eqref{eq:R_n-Comp} (as we vary over the finitely many vertices of $W_i(\hat h)$).
\end{remark}

\subsection{Equivariance/invariance properties of quantile and rank functions}\label{sec:Euiv-Inv}
Next we study equivariance/invariance properties of quantile and rank functions under different transformations on $X$. The first result in this direction, proved in Section~\ref{Pf:lem:LinTrans}, shows that the quantile/rank function of $Y := c X + b$, where $X \sim \nu$ is a random vector in $\R^d$ ($b \in \R^d$ and $c >0$), can be easily obtained from the quantile/rank function of $X$. 
\begin{lemma}\label{lem:LinTrans}
Suppose that $X \sim \nu$ where $\nu$ is a distribution on $\R^d$. Let $\mu$ be an absolutely continuous distribution on $\R^d$ with support $\s$. Suppose that $c >0$ is a scalar and $b \in \R^d$. Let $Y := c X + b$. Let $Q_X:\s \to \R^d$ and $R_X: \R^d \to\R^d$ be the quantile and rank maps of $X$ (w.r.t.~$\mu$). Let $Q_Y:\s \to \R^d$ and $R_Y: \R^d \to \R^d$ be the quantile and rank maps of $Y$ (w.r.t.~$\mu$). Then, for  $\mu$-a.e.~$u$ and for a.e.~$y \in \R^d$, $$Q_Y(u) = c\,Q_X(u) + b, \qquad  \mbox{and} \qquad R_Y(y) = R_X((y -b)/c).$$
\end{lemma}

The above lemma holds for any probability measure $\nu$ (discrete or continuous); it justifies the fact that we can rescale the data (by adding a constant vector and multiplying by a positive scalar) in Figures~\ref{fig:Q-map} and~\ref{fig:4Plots} and the cell decomposition (of $\s$) does not change (see~\eqref{eq:W_i}) as the transformed (piecewise linear and convex) potential function is obtained by adding a constant to a positive multiple of the previous (piecewise linear) potential function.

One may ask if it is possible to relate the quantile/rank functions of $Y := AX$, where $A_{d \times d}$ is a matrix, to those of $X$. The following result (proved in Section~\ref{Pf:lem:LinTrans-2}) shows that the quantile map is equivariant under orthogonal transformations if $\mu$ is spherically symmetric; also see~\cite[Corollary 2.12]{Cuesta-2013}.
 \begin{lemma}[{\cite[Corollary 2.12]{Cuesta-2013}}]\label{lem:LinTrans-2}
 Suppose that $A$ is an orthogonal matrix, i.e., $A A^\top  = A^\top A = I_d$. Let $\mu$ be a spherically symmetric absolutely continuous distribution on $\R^d$ (e.g., the uniform distribution on the unit ball around $0 \in \R^d$). Let us denote by $\psi_X$ the potential function linked to the random variable $X\sim \nu$, i.e., $\nabla \psi_X \# \mu = \nu$ and $\psi_X$ is convex. Then a potential function of $Y:=AX$ is given by $\psi_Y(u) = \psi_X(A^\top u)$, for $u \in \R^d$. As a consequence, for  $\mu$-a.e.~$u$ and for a.e.~$y \in \R^d$, $$Q_Y(u) = A Q_X(A^\top u), \qquad \mbox{and} \qquad R_Y(y) = A R_X(A^\top y).$$
\end{lemma}

For the next result (proved in Section~\ref{pf:prop:Indep}) we take $\mu =$ Uniform$([0,1]^d)$; note that the choice of $\mu =$ Uniform$([0,1]^d)$ has been studied before (see e.g.,~\cite{D14, Galichon16}). This has implications in testing for mutual independence between random vectors; see Section~\ref{sec:IndepTest} of the main paper for more details. 

\begin{proposition}[{\cite[Theorem 2.9]{Cuesta-2013}}]\label{prop:Indep}
Suppose that $\nu$ is a distribution on $\R^d$ and let $\mu= $ Uniform$([0,1]^d)$. Suppose that $X =(X_1,X_2, \ldots, X_k) \sim \nu$ where $k \ge 2$, $X_i \sim \nu_i$, for $i=1,\ldots, k$, are random vectors in $\R^{d_i}$ (here $d_1 + \ldots + d_k = d$). Let $Q$ and $Q_i$ be the quantile maps of $X$ and $X_i$, for $i = 1,\ldots, k$, respectively (w.r.t.~$\mu$ and $\mu_i = $ Uniform$([0,1]^{d_i})$). Let $R$ and $R_i$, for $i = 1,\ldots,k$, be the corresponding rank maps. If $X_1, \ldots, X_k$ are mutually independent then 
\begin{equation}\label{eq:Indep-Q}
Q(u_1,\ldots, u_k) = (Q_1(u_1), \ldots, Q_k(u_k)), \;\;\; \mbox{for $\mu$-a.e. } (u_1,\ldots,u_k) \in \R^{d},
\end{equation}
and
\begin{equation}\label{eq:Indep}
R(x_1,\ldots, x_k) = (R_1(x_1), \ldots, R_k(x_k)), \;\;\; \mbox{for a.e.~}(x_1,\ldots,x_k) \in \R^{d}.
\end{equation}
 Conversely, suppose that~\eqref{eq:Indep-Q} or~\eqref{eq:Indep} holds. Then $X_1,\ldots,X_k$ are mutually independent.
\end{proposition}

\subsection{The empirical quantile and rank functions when $d=1$}\label{sec:d=1}
\begin{remark}[Quantile function]\label{rem:Q-d=1}
Suppose that $d=1$ and that $\mu = $ Uniform$([0,1])$ and $\nu$ is a continuous distribution on $\R$. Let $X_1,\ldots, X_n$ be an i.i.d.~sample from $\nu$ and let $X_{(1)} < \ldots < X_{(n)}$ be the order statistics. Then the sample quantile function $\hat Q_n$, defined via~\eqref{eq:SampQ_n}, is given by
\begin{equation}\label{eq:Emp_Q-1}
	\hat Q_n(u) = X_{(i)}, \qquad \mbox{if } u \in \left({(i-1)}/{n},{i}/{n}\right), \mbox{ for } i = 1,\ldots, n,
\end{equation}
and $\hat Q_n(0) = X_{(1)}$ and $\hat Q_n(1) = X_{(n)}$. This follows from the fact that $\hat Q_n = \nabla \hat \psi_n$ where $\hat \psi_n$ can be expressed as
\begin{equation}\label{eq:Emp_psi-1}
	\hat \psi_n(u) = \max_{i=1,\ldots, n} \{X_{(i)} u + \hat h_{(i)}\}, \qquad \mbox{for } u \in [0,1],
\end{equation}
and $\hat h_{(i+1)} - \hat h_{(i)} = i(X_{(i)}-X_{(i+1)})/n$ (after simple algebra), for $i=1,\ldots, n-1$. Without loss of generality one can take $\hat h_{(1)}=0$. Observe that, we are free to define $\hat Q_n(i/n)$ as any point in the interval $[ X_{(i)},X_{(i+1)}]$, for $i=1,\ldots, n-1$. Thus, the sample quantile map obtained via~\eqref{eq:Emp_Q} is essentially the same as the usual quantile function, except at the points $i/n$, for $i=1,\ldots, n-1$. 
\end{remark}

\begin{remark}[Rank function]\label{rem:R-d=1} Let us first illustrate the above notion of the rank function $\hat{R}_n $ when $d=1$. From the form of $\hat{\psi}^{*}_n(\cdot)$ and the related discussion in~\cite{Gu12}, it follows that the graph $\{(x,\hat \psi_n^*(x)): x \in \conv(X_{(1)},\ldots, X_{(n)})\}$ of $\hat \psi_n^*$ is the lower boundary of the convex hull $\conv((X_{(1)}, - \hat h_{(1)}), \ldots, (X_{{(n)}}, - \hat h_{(n)}))$. The $i$-th divided difference of the set of points in the convex hull is $\frac{\hat h_{(i)} -\hat h_{(i+1)}}{X_{(i+1)} - X_{(i)}} = \frac{i}{n};$ cf.~\eqref{eq:Emp_psi-1}. As the divided differences are strictly increasing, the piecewise linear interpolation of the points $(X_{(1)}, - \hat h_{(1)}), \ldots, (X_{{(n)}}, - \hat h_{(n)})$ gives a convex function, which is indeed the graph of $\{(x,\hat \psi_n^*(x)): x \in \conv(X_{(1)},\ldots, X_{(n)})\}$. Moreover, this shows that $\hat R_n = \nabla \hat \psi_n^*$ is given by 
\begin{equation}\label{eq:Emp_F-1}
	\hat R_n(x) = \frac{i}{n}, \qquad \mbox{if } x \in \left(X_{(i)},X_{(i+1)}\right), \quad \mbox{for } i = 0,1,\ldots, n;
\end{equation}
here, by convention, $X_{(0)} = -\infty$ and $X_{(n+1)} = +\infty$. We are free to define $\hat R_n(X_{(i)})$ as any point in the interval $[(i-1)/n,i/n]$, for $i=1,\ldots, n$; cf.~$\mathbb{F}_n(X_{(i)}) = i/n$ where $\mathbb{F}_n$ is the empirical distribution function of the $X_i$'s. Thus,  we essentially get back the usual notion of the empirical distribution function when $d=1$ except that now the ranks of the data points $X_i$'s are not uniquely defined.
\end{remark}

\subsection{Some additional plots}\label{sec:Add-Plots}
\begin{figure}
\includegraphics{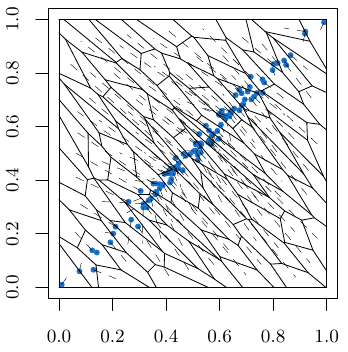}\hspace{0.2in}
\includegraphics{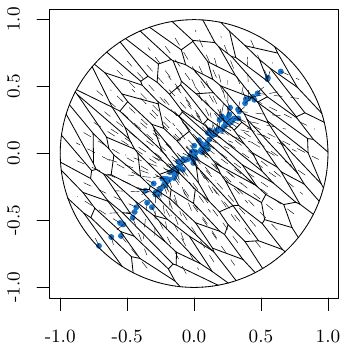}

\caption{Left plot: Shows the cell decomposition of $\s = [0,1]^2$ (each with area $1/n$ where $n=100$) induced by the estimated quantile function $\hat Q_n$ (w.r.t.~$\mu = $ Uniform$([0,1]^2)$) where the data points are drawn i.i.d. from $N_2((0,0),\Sigma)$ where $\Sigma_{1,1} = \Sigma_{2,2} = 1$ and $\Sigma_{1,2} = \Sigma_{2,1} = 0.99$ (and appropriately scaled to lie in $[0,1]^2$) along with dashed lines indicating which cell corresponds to which data point. Right plot: Shows the corresponding cell decomposition for $\s = B_1(0)$ --- the ball of radius 1 around $(0,0) \in \R^2$ --- induced by the estimated quantile function $\hat Q_n$ (w.r.t.~$\mu = $ Uniform$(B_1(0))$); here the data points are scaled to lie in $B_1(0)$.}
\label{fig:4Plots-Corr-Nor}
\end{figure}

\begin{figure}
		\includegraphics[scale=0.80]{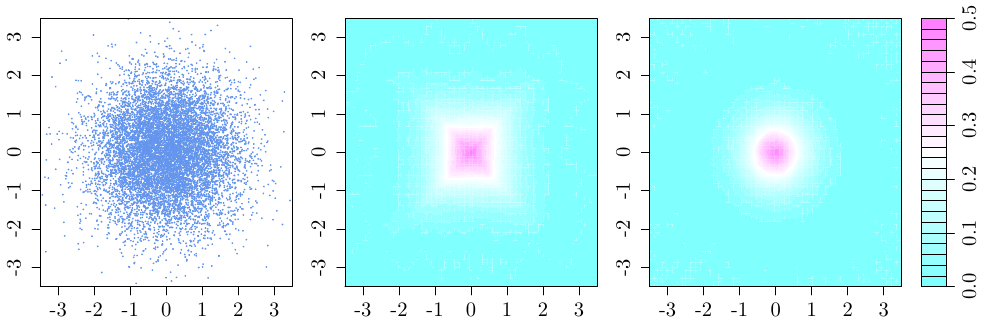}
		\caption{Left panel: A random sample of size $n = 10^4$ drawn from   $N_2((0,0),I_2)$. Middle panel: The estimated depth function --- defined as $\hat D_n(x) := 1/2 - \|\hat R_n(x) - (1/2)\mathbf{1}\|_\infty$ (for $x \in \R^d$; see~\cite{Cher17}) where $\hat R_n(\cdot)$ is the estimated rank function and $\mathbf{1} = (1,1,\ldots, 1) \in \R^d$ --- using $\mu =$ Uniform$([0,1]^2)$ as the reference distribution. Right panel: The estimated depth function --- defined as $\hat D_n(x) = \pi^{-1}(\theta - \cos \theta \sin \theta)$ where $\theta = \arccos (\|\hat R_n(x)\|)$ --- w.r.t.~$\mu =$ Uniform$(B_1(0))$; see~\cite[Section 5.6]{RR1999}.}
	\label{fig:Normal}
	\end{figure}

\begin{figure}
		\includegraphics[scale=0.80]{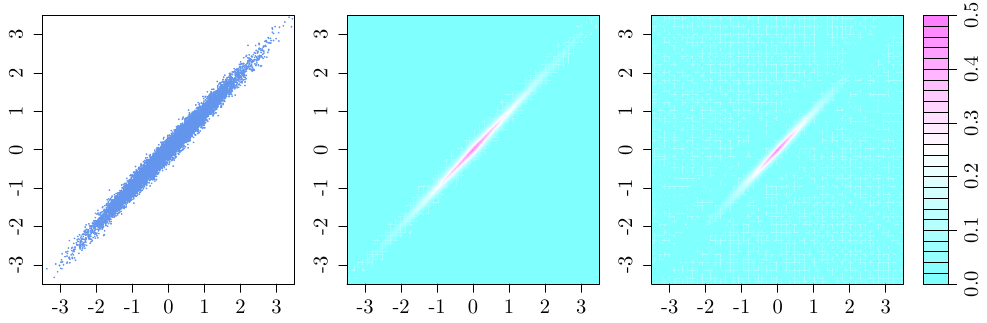}
		\caption{Left panel: A random sample of size $n = 10^4$ drawn from $N_2((0,0),\Sigma)$ where $\Sigma_{1,1} = \Sigma_{2,2} = 1$ and $\Sigma_{1,2} = \Sigma_{2,1} = 0.99$ . Middle panel: The estimated depth function --- defined as $\hat D_n(x) := 1/2 - \|\hat R_n(x) - (1/2)\mathbf{1}\|_\infty$ (for $x \in \R^d$; see~\cite{Cher17}) where $\hat R_n(\cdot)$ is the estimated rank function and $\mathbf{1} = (1,1,\ldots, 1) \in \R^d$ --- using $\mu =$ Uniform$([0,1]^2)$ as the reference distribution. Right panel: The estimated depth function --- defined as $\hat D_n(x) = \pi^{-1}(\theta - \cos \theta \sin \theta)$ where $\theta = \arccos (\|\hat R_n(x)\|)$ --- w.r.t.~$\mu =$ Uniform$(B_1(0))$; see~\cite[Section 5.6]{RR1999}.}
	\label{fig:Normal-Sig}
	\end{figure}

In this subsection we provide some additional plots that supplement the discussion in the main paper. In particular, we provide plots that: (i) illustrate the cell decomposition induced by the estimated quantile function $\hat Q_n$ when the reference distribution is uniform on $\s = [0,1]^2$ and $\s = B_1(0)$ (see Figure~\ref{fig:4Plots-Corr-Nor}); (ii) show depth functions for different distributions when the reference distribution is uniform on $\s = [0,1]^2$ and $\s = B_1(0)$ (see Figures~\ref{fig:Normal} and~\ref{fig:Normal-Sig}).

\subsection{Some auxiliary results}\label{sec:AuxRes}
In this subsection we state some simple results about convex functions and their subgradients that will be useful in our proofs. 
Many of the results below follow from~\cite[Chapter 2]{V03}; however for the convenience of the reader we give their complete proofs. Our first result, proved in Section~\ref{sec:Injectivity}, gives a sufficient condition for the injectivity\footnote{A function is called injective if and only if it is one-to-one.} of the gradient map of a strictly convex function.  

\bl\label{lem:Injectivity} 
Let $\mathcal{S}$ be a convex set in $\R^d$ and let $\psi:\mathcal{S}\to \RR $ be a strictly convex function. Assume that the gradient map $\nabla\psi:\mathcal{S}\to \RR^d$ is continuous. Then, $\nabla \psi$ is an injective function.  
\el
The following lemma, proved in Section~\ref{Pf:lem:Bd}, shows how the growth of a convex function depends on the diameter of its subdifferential set.   

\begin{lemma}\label{lem:Bd}
Let $f:\RR^d\to \RR\cup \{+\infty\}$ be a convex function. Fix $A \subseteq \mathrm{dom}(f)$. Suppose that $\sup\{\|\partial f(y)\|: y\in A\}<K$ for some absolute constant $K$. Then, $|f(x)-f(x_0)|\leq K \|x-x_0\|$ for all $x, x_0\in A$.    
\end{lemma}

Our next result, proved in Section~\ref{sec:SubGIneq}, provides a sufficient condition for a vector to be a subgradient of a convex function at a point.  
\bl\label{lem:SubGIneq} 
Let $\psi:\RR^d\to \RR\cup\{+\infty\}$ be a convex function. Fix $x_0\in \RR^d$. Suppose that there exists $\xi\in \RR^d$ such that 
\begin{align}\label{eq:SubGIneq}
\psi(x)\geq \psi(x_0)+ \langle \xi,x-x_0\rangle, \qquad \mbox{for all $x\in B_{\delta}(x_0)$},  
\end{align}
 where $B_{\delta}(x_0)$ is the open ball of radius $\delta>0$ around $x_0 \in \R^d$. Then, $\xi\in \partial \psi(x_0)$. 
\el

The following result follows directly from the definition of the Legendre-Fenchel dual of a convex function.
\bl\label{lem:Affinenv} 
Let $\psi:\RR^d\to \RR \cup \{+\infty\}$ be convex function and $c\in \RR$. Define $\phi:\RR^d \to \RR\cup \{+\infty\}$ by $\phi(y) := \psi(y)+c$ for all $y\in \RR^d$. Let $\psi^{*}$ and $\phi^{*}$ be the Legendre-Fenchel duals of $\psi$ and $\phi$ respectively. Then, $\phi^{*}(y) = \psi^{*}(y)-c$ for all $y \in \RR^d$.
\el

We end this section with the following lemma which is proved in Section~\ref{sec:CvxEq}.

\bl\label{lem:CvxEq} 
Let $\mathcal{S}$ and $\mathcal{Y}$ be two convex subsets of $\RR^d$ with non-empty interiors. Let $\phi:\s \to \R \cup \{+\infty\}$ and $\psi:\s \to \R \cup \{+\infty\}$ be two convex functions such that $\nabla \phi:\mathrm{Int}(\mathcal{S})\to \mathrm{Int}(\mathcal{Y})$ is continuous and $\nabla\phi=\nabla\psi$ a.e. in $\mathrm{Int}(\mathcal{S})$ (w.r.t. Lebesgue measure). Then, $\psi$ is differentiable everywhere in $\mathrm{Int}(\mathcal{S})$ and $\nabla\phi(y)=\nabla\psi(y)$ for all $y\in\mathrm{Int}(\mathcal{S})$. 
\el

\section{Proofs of results in Appendix~A}\label{pf:Proofs}

\subsection{Proof of Lemma~\ref{lem:SubD-1}}\label{sec:SubD-1}
The first assertion follows directly from Lemma~\ref{lem:SubD}. To see the second assertion, note that if $x\in \partial \psi^{*}(y)$ for some $y\in B$, then, $y\in \partial \psi(x)$ which implies that $x\in (\partial \psi)^{-1}(B)$ and hence, $\partial \psi^{*}(B)\subset (\partial \psi)^{-1}(B)$. Using a similar argument, one can show that $(\partial \psi)^{-1}(B)\subset \partial \psi^{*}(B)$. \qed

\subsection{Proof of Lemma~\ref{lem:QRrelation}}\label{sec:QRrelation}
Let $\Theta_{f}$ and $\Theta_{g}$ be the set of points where $f$ and $g$ are positive. Denote the sets of points where $Q$ and $R$ are  differentiable by $D_{Q}$ and $D_R$. Note that 
\begin{align}
\rho_{Q}(A) &=\int_{Q(A)\cap \Theta_{g}} \frac{1}{g(x)} g(x)dx = \int_{Q(A)\cap \Theta_{g}\cap D_{R}} \frac{1}{g(x)} d\nu(x) \\&= \int \limits_{Q^{-1}(Q(A)\cap \Theta_{g}\cap D_{R})} \frac{1}{g(Q(x))} d\mu(x) = \int \limits_{Q^{-1}(Q(A)\cap \Theta_{g}\cap D_{R})\cap D_{Q}} \frac{f(x)}{g(Q(x))}d x.\label{eq:MnAm}
\end{align}
The first two equalities follow since the complements of $\Theta_{g}$ and $D_{R}$ have zero Lebesgue measure. The third equality follows by noting that $Q\# \mu =\nu$. The last equality holds because the complement of $D_{Q}$ has zero Lebesgue measure and $\mu$ is absolutely continuous with density $f$. Now, we claim and prove that 
\begin{align}\label{eq:QR}
\qquad A \supseteq Q^{-1}(Q(A)\cap \Theta_{g}\cap D_{R})\cap D_{Q} \supseteq A \cap Q^{-1}(\Theta_{g}\cap D_{R})\cap D_{Q}. 
\end{align}
We start with proving that $Q^{-1}(Q(A)\cap \Theta_{g}\cap D_{R})\cap D_{Q} $ is contained in $A$. Fix $y\in Q^{-1}(Q(A)\cap \Theta_{g}\cap D_{R})\cap D_{Q} $. Since $y\in D_{Q}$, $\partial \psi(y)$ is singleton, i.e., there exists $x\in Q(A)\cap \Theta_{g}\cap D_{R}$ such that $Q(y)=x$. Note that $x\in D_{R}$. This in combination with the fact $Q(y)=x$ and Lemma~\ref{lem:SubD} implies that $\partial \psi^*(x)$ is also a singleton set and hence, $R(x)=y$. However, if there exists $y^{\prime}\neq y \in A$ such that $x\in Q(y^{\prime})$, then, by Lemma~\ref{lem:SubD}, we must have $y^{\prime}\in R(x)$ which contradicts $R(x)=y$. Hence, $y$ must belong to the set $A$. 

Now, we turn to show that $A \cap Q^{-1}(\Theta_{g}\cap D_{R})\cap D_{Q}$ is contained in $ Q^{-1}(Q(A)\cap \Theta_{g}\cap D_{R})\cap D_{Q}$. Fix $y \in A \cap Q^{-1}(\Theta_{g}\cap D_{R})\cap D_{Q}$. It suffices to show that $Q(y)\in Q(A)\cap \Theta_{g}\cap D_{R}$. Since $y\in A$, trivially, $Q(y)\in Q(A)$. Furthermore $\partial \psi(y)$ is a singleton set because $y\in D_{Q}$. Combining this with $y\in Q^{-1}(\Theta_{g}\cap D_{R})$ yields $Q(y)\in \Theta_{g}\cap D_{R}$. Therefore, we have $Q(y)\in Q(A)\cap \Theta_{g}\cap D_{R}$. This completes the proof of \eqref{eq:QR}. 

Applying \eqref{eq:QR} in \eqref{eq:MnAm} yields 
\begin{align}
\int_{A \cap Q^{-1}(\Theta_{g}\cap D_{R})\cap D_{Q}} \frac{f(x)}{g(Q(x))} dx\leq \rho_{Q}(A)\leq \int_{A}\frac{f(x)}{g(Q(x))} dx.
\end{align}
Observe that the complement of $Q^{-1}(\Theta_{g}\cap D_{R})\cap D_{Q}$ has zero measure w.r.t. $\mu$ because $\mu(Q^{-1}(\Theta_{g}\cap D_{R}))= \nu(\Theta_{g}\cap D_{R})=1$ and $\mu(D_{Q})=1$. As a consequence, we get  
\begin{align}
 \rho_{Q}(A)=\int_{A}\frac{f(x)}{g(Q(x))} dx.
\end{align} 
This completes the proof. \qed

\subsection{Proof of Lemma~\ref{lem:psi-finite}}\label{sec:psi-finite}
Fix $x \in \R^d$. We want to show $\psi^{*}(x)<\infty$. First, observe that $\psi^*(x) = \sup_{y \in \R^d} \{\langle x, y\rangle - \psi(y)\} = \sup_{y \in \s} \{\langle x, y\rangle - \psi(y)\}$. The last equality is a consequence of the fact that for $y \notin \s$, $\psi(y) = + \infty$ (by our convention~\eqref{eq:Conv-Cvx}).

Note that $\psi^{*}(x)\le \sup_{y\in \s} \langle y,x \rangle  - \inf_{y\in \s}\psi(y)$. We know that as $\s$ is a compact set, $\sup_{y\in \s} \langle y,x \rangle <\infty$. Thus, it suffices to show that $\inf_{y\in \s}\psi(y)>-\infty$. To show this, fix a point $y_0 \in \s$ such that: (a) $\psi(y_0)>-\infty$, and (b) there exists $z\in \R^{d}$ satisfying $z\in \partial \psi(y_0)$. Note that the above two properties can be easily satisfied as $\psi$ is a proper convex function. Then, by the convexity of $\psi$, $\psi(y)\ge \psi(y_0)+ \langle y-y_0,z \rangle $. Note that as $\inf_{y\in \s} [{\psi(y_0)+ \langle y-y_0,z \rangle}] >-\infty$, we have $\inf_{y \in \s}\psi(y)>-\infty$. This completes the proof. \qed

\subsection{Proof of Lemma~\ref{lem:Monotonicity}}\label{sec:Monotonicity}
Fix $x,y\in \R^d$. Let $\psi$ be a convex function on $\s$ such that $\nabla  \psi \# \mu = \nu$ and let $R = \nabla  \psi^*$ where $\psi^{*}$ is the Legendre-Fenchel dual of $\psi$. Consider the function $\psi^{*}_{x,y}:\RR\to \RR$ defined as $\psi^{*}_{x,y}(t) := \psi^{*}(tx+(1-t)y)$ . Since $\psi^{*}$ is a convex function in $\RR^d$ and $tx+(1-t)y$ is a linear function w.r.t. $t$, $\psi^{*}_{x,y}(t)$ is a convex function w.r.t. $t$. By the calculus of subgradients, for $t\in \RR$, $\partial \psi^{*}_{x,y}(t) = (x-y)^{\top}\partial \psi^{*}(tx+(1-t)y)$. Since the (sub)-gradient of a one-dimensional convex function is monotonic, $(x-y)^{\top}\nabla \psi^{*}(tx+(1-t)y)$ is a monotone function  of $t$, which implies the desired result.
\qed

\subsection{Proof of Lemma~\ref{lem:ExtremePt}}\label{Pf:lem:ExtremePt}
Note that for any $x \in \R^d$, from the definition of $\hat{\psi}_n^*$ (see~\eqref{eq:Sup}), we have $\hat{\psi}_n^*(x) \ge \langle x,y\rangle - \hat{\psi}_n(y)$, for all $y \in \R^d$. Fix $i \in \{1,\ldots, n\}$. As (see~\eqref{eq:Emp_G})
\begin{align}
\hat{\psi}_n(u)= \begin{cases}
\max\big\{\langle u, X_1\rangle +\hat{h}_1, \ldots , \langle u, X_n\rangle +\hat{h}_n\big\}\geq \langle u, X_i\rangle + \hat{h}_i & \text{if }u\in \s,\\
+\infty &  \text{if }u\in \s^c.
\end{cases} \nonumber
\end{align}
we have $
-\hat{h}_i \geq \langle y, X_i\rangle - \hat{\psi}_n(y)$, for all $y \in \RR^d$. Thus, $-\hat{h}_i \ge \sup_{y \in \R^d} \{\langle y, X_i\rangle - \hat{\psi}_n(y)\} =:\hat{\psi}_n^*(X_i)$. Note that, $W_i(\hat{h})$ is convex set of dimension $d$ and for all $u \in W_i(\hat h)$ (see~\eqref{eq:W_i}), we have $\hat{\psi}_n(u) = \langle u, X_i\rangle + \hat h_i$ which implies that $- \hat h_i = \langle y, X_i\rangle - \hat{\psi}_n(y) \le \hat{\psi}_n^*(X_i)$ (by definition). This completes the proof. \qed

\subsection{Proof of Lemma~\ref{lem:LinTrans}}\label{Pf:lem:LinTrans}
Let $\psi_X:\R^d \to \R \cup \{+\infty\}$ be a convex function such that $\nabla \psi_X\# \mu = \nu$. Further, we can assume that $\psi_X(u) = +\infty$ for $u \notin \s$. Thus, from the definition of $Q_X$ and $R_X$, we have $Q_X(u) = \nabla \psi_X(u)$, for $\mu$-a.e.~$u$, and $R_X(x) = \nabla \psi_X^*(x)$ for a.e.~$x$. Let $\psi_Y:\R^d \to \R\cup \{+\infty\}$ be as defined  
\begin{equation}\label{eq:psi_Y}
\psi_Y(u) := c\, \psi_X(u) + \langle b, u\rangle,  \qquad \mbox{for all } u \in \R^d.
\end{equation}
Then $\psi_Y(\cdot)$ is a convex function, and for a.e.~$y \in \R^d$, $$\nabla \psi_Y(u) = c \nabla \psi_X(u) + b.$$ Let $U \sim \mu$. Further, using the facts that $Y := c X + b$ and $\nabla \psi_X(U) \sim \nu$, we get $$\nabla \psi_Y(U) =  c \nabla \psi_X(U) + b$$ which has the same distribution as $Y$. Thus, $\nabla \psi_Y(\cdot)$ is the gradient of a convex function that pushes forward $\mu$ to the distribution of $Y$.  Therefore, for $\mu$-a.e.~$u$, $$Q_Y(u) = \nabla \psi_Y(u) = c Q_X(u) + b,$$ which yields the first result.

Next, using the form of $\psi_Y$ for a.e.~$y \in \R^d$, we have $$\psi_Y^*(y) := \sup_{u \in \R^d} \{\langle y, u\rangle - \psi_Y(u) \} = \sup_{u \in \R^d} \{\langle y, u\rangle - c\, \psi_X(u) - \langle b, u\rangle \}  = c \psi_X^*((y-b)/c).$$ Thus, for a.e.~$y \in \R^d$,  $$R_Y(y) := \nabla \psi_Y^*(y)  = \nabla \psi_X^*((y-b)/c)  =: R_X((y-b)/c)$$ which yields the second result. \qed

\subsection{Proof of Lemma~\ref{lem:LinTrans-2}}\label{Pf:lem:LinTrans-2}
Let $\psi_Y(u) = \psi_X(A^\top u)$, for $u \in \R^d$. Then $\nabla \psi_Y(u) = A \nabla \psi_X(A^\top u)$, for $u \in \R^d$. By definition, $\nabla \psi_Y$ is the gradient of a convex function (as $\psi_X$ is a convex function). Moreover, for $U \sim \mu$ and $X \sim \nu$, $$\nabla \psi_Y(U) = A  \nabla \psi_X(A^\top U) \stackrel{d}{=} A\nabla \psi_X(U) \stackrel{d}{=} A X,$$ where we have used the fact $A^\top U \stackrel{d}{=} U$, as $U$ is spherically symmetric. Thus, $\nabla \psi_Y$ is the gradient of a convex function that transports $U$ to $AX$, and this, by Lemma~\ref{thm:Brenier}, completes the proof of the first part. 

Now, for $y \in \R^d$, $$\psi_Y^*(y) =  \sup_{u \in \s} \{u^\top y - \psi_Y(u)\} = \sup_{u \in \s} \{(A^\top u)^\top A^\top y - \psi_X(A^\top u)\} = \psi^*_X(A^\top y).$$ Thus, $\nabla \psi_Y^*(y) = A \nabla \psi_X^*(A^\top y)$ which completes the proof. \qed

\subsection{Proof of Proposition~\ref{prop:Indep}}\label{pf:prop:Indep}

We will prove the result when $k=2$. The proof for $k\ge 3$ is exactly similar. By Proposition~\ref{thm:Q_prop} we can find convex functions $\psi_i:\R^{d_i} \to \R \cup \{+\infty\}$, $i=1,2$, such that $\psi_i(u_i) = +\infty$ for $u_i \in \R^{d_i} \setminus [0,1]^{d_i}$ and $$Q_i(u_i) = \nabla \psi_i(u_i),\qquad \mbox{for a.e.}~u_i \in [0,1]^{d_i}.$$ Let us define the function $\psi:\R^d \to \R \cup \{+\infty\}$ as 
\begin{equation}\label{eq:psi_1_2}
\psi(u_1,u_2) := \psi_1(u_1)  + \psi_2(u_2), \qquad \mbox{for all }\;\;\; (u_1,u_2) \in \R^{d_1} \times \R^{d_2}.
\end{equation} 
Observe that, as defined above, $\psi(\cdot)$ is a convex function and, for a.e.~$(u_1,u_2) \in [0,1]^{d_1} \times [0,1]^{d_2},$ $$\nabla \psi(u_1,u_2) = (\nabla \psi_1(u_1), \nabla \psi_2(u_2)) = (Q_1(u_1), Q_2(u_2)).$$ Obviously $\nabla \psi: [0,1]^d \to \R^d$ is the gradient of a convex function. Let $U_i \sim $ Uniform$([0,1]^{d_i})$, for $i=1,2$ be independent. As $Q_1$ and $Q_2$ are the quantile maps of $X_1 \sim \nu_1$ and $X_2 \sim \nu_2$, and $X_1$ and $X_2$ are independent, we have $$\nabla \psi(U_1,U_2) = (Q_1(U_1), Q_2(U_2)) \sim \nu_1  \times \nu_2 = \nu.$$ 
Thus, $\nabla \psi$ pushes forward $\mu$ to $\nu$. As both $Q$ and $\nabla \psi$ are: (i) gradients of convex functions, and (ii) transport $\mu$ to $\nu$,~\eqref{eq:Indep-Q} now follows from the a.s.~uniqueness of such a transport map (see Theorem~\ref{thm:Brenier}).

Recall the definitions of $\psi_1$, $\psi_2$, and $\psi$ from above. Note that, for $x_i \in \R^{d_i}$, for $i=1,2$, $$\psi_i^*(x_i) = \sup_{y_i \in \R^{d_i}} \{\langle x_i, y_i \rangle - \psi_i(y_i)\}. $$ Further, using~\eqref{eq:psi_1_2}, for $(x_1,x_2) \in \R^{d_1} \times \R^{d_2}$, $$\psi^*(x_1,x_2) = \sup_{(y_1,y_2)\in \R^{d_1} \times \R^{d_2}} \{\langle x_1, y_1 \rangle + \langle x_2, y_2 \rangle - \psi(y_1,y_2)\} = \psi_1^*(x_1) + \psi_2^*(x_2). $$ Therefore, for $(x_1,x_2) \in \R^{d_1} \times \R^{d_2} $, $$\nabla \psi^*(x_1,x_2) = (\nabla \psi_1^*(x_1), \nabla \psi_2^*(x_2))^\top $$ which yields the desired result as $R(x) = \nabla \psi^*(x)$ for a.e.~$x$ and $R_i(x) = \nabla \psi^*_i(x_i)$ for $\lambda_{d_i}$-a.e.~$x_i$, for $i=1,2$.

 Now suppose that~\eqref{eq:Indep-Q} holds. Let $Q = \nabla \psi$, where $\psi:\R^d \to \R\cap \{+\infty\}$ is a convex function. Let $U = (U_1, U_2) \sim $  Uniform$([0,1]^{d})$ where $U_i \sim $ Uniform$([0,1]^{d_i})$, for $i=1,2$, and are independent. Note that as $Q \#\mu = \nu$ we have $Q(U) \sim \nu$. But, by~\eqref{eq:Indep-Q}, $Q(U) \stackrel{d}{=} (Q_1(U_1), Q_2(U_2)) \sim \nu$. As $U_1$ and $U_2$ are independent, $Q_1(U_1)$ and $Q_2(U_2)$ are independent which in turn implies that $(X_1,X_2) \sim \nu$ factors as a product measure. Thus $X_1$ and $X_2$ are independent.

Finally, suppose that~\eqref{eq:Indep} holds. We will show that $X_1$ and $X_2$ are independent. We prove this by contradiction. Recall $Q\# \mu = \nu$. Let $\tilde{Q}:\RR^d\to \RR^d$ be such that $\tilde{Q}(u_1,u_2)= (Q_1(u_1), Q_2(u_2))$ $\mu$-a.e. Recall that \eqref{eq:Indep-Q} implies the independence of $X_1$ and $X_2$. Thus, $\tilde{Q}\#\mu \neq \nu$ when $X_1$ and $X_2$ are not independent. Let $\psi$ and $\tilde{\psi}$ be two convex functions such that $Q=\nabla \psi$ and $\tilde{Q}= \nabla \tilde{\psi}$. Suppose that there exists a Borel set $B \subset \RR^d$ such that   
 \begin{align}\label{eq:ContraIn}
 \mu\Big((\partial \psi)^{-1}(B)\Big) \neq \mu((\partial \tilde{\psi})^{-1}(B)). 
\end{align} 
Owing to Lemma~\ref{lem:SubD-1}, we have $Q^{-1}(B)= R(B)$ and $\widetilde{Q}^{-1}(B) = \widetilde{R}(B)$ where $\widetilde{R}:\RR^d\to \RR^d$ is defined by $\widetilde{R}(x_1,x_2):=(R_1(x_1), R_2(x_2))$. Due to \eqref{eq:Indep}, $R(B)= \widetilde{R}(B)$. Hence, the inequality of \eqref{eq:ContraIn} should be equality which contradicts $\tilde{Q}\#\mu \neq \nu$. Hence, \eqref{eq:Indep} implies that $X_1$ and $X_2$ are independent.

\subsection{Proof of Lemma~\ref{lem:Injectivity}}\label{sec:Injectivity}
Since $\psi$ is strictly convex in $\s$, for any $x \ne y\in \mathcal{S}$, we have $$ \psi(x)>\psi(y) +\langle \nabla \psi(y), x-y\rangle, \quad \mbox{and} \quad -\psi(x)> -\psi(y)- \langle \nabla \psi(x), x-y\rangle. $$ Adding both sides of two inequalities yields 
\begin{align}
\langle \nabla\psi(x) -\nabla\psi(y), x-y\rangle >0 , \quad \forall x \ne y\in \mathcal{S}. \label{eq:FundIneq}  
\end{align} 
Now, \eqref{eq:FundIneq} shows that if $x\neq y$, then, $\nabla \psi(x)\neq \nabla \psi(y)$. This proves that $\nabla \psi$ is an injective function in $\mathcal{S}$.
\qed

\subsection{Proof of Lemma~\ref{lem:Bd}}\label{Pf:lem:Bd}
From the definition of the subdifferential, 
\begin{align}\label{eq:UpLowBd-1}
\langle z_1, x-x_0\rangle\leq f(x) - f(x_0)\leq \langle z_2, x-x_0\rangle
\end{align}
where $z_1\in \partial f(x_0)$ and $z_2\in \partial f(x)$. Owing to \eqref{eq:UpLowBd-1}, 
\begin{align}
|f(x)-f(x_0)|\leq \max\{|\langle z_1, x-x_0\rangle|, |\langle z_2, x-x_0\rangle|\}\leq K \|x-x_0\|
\end{align}
where the last inequality follows from  Cauchy-Schwarz inequality and recalling that $\max\{\|z_1\|, \|z_2\|\}\leq K$. 
\qed

\subsection{Proof of Lemma~\ref{lem:SubGIneq}}\label{sec:SubGIneq}
Assume that there exists $z\in \RR^d$ such that 
 \begin{equation}\label{eq:RevSub}
 \psi(z)< \psi(x_0)+\langle u, z-x_0\rangle.
 \end{equation} 
Choose $\alpha\in (0,1)$ such that $z^{\prime} := \alpha z+ (1-\alpha)x_0 \in B_{\delta}(x_0)$. Then,
 \begin{align}\label{eq:Contra2}
 \psi(z^{\prime})\leq \alpha \psi(z) +(1-\alpha)\psi(x_0)< \psi(x_0) +\langle u, z^{\prime} - x_0\rangle
 \end{align}
 where the first inequality follows from the convexity of $\psi$ and the second inequality is obtained by using~\eqref{eq:RevSub} and the definition of $z^{\prime}$. Since \eqref{eq:Contra2} contradicts \eqref{eq:SubGIneq} (as $z^{\prime} \in B_{\delta}(x_0)$), there cannot exist $z$ satisfying \eqref{eq:RevSub}. \qed 
 
 \subsection{Proof of Lemma~\ref{lem:CvxEq}}\label{sec:CvxEq}

Let us define $\mathfrak{B} := \{y\in \mathrm{Int}(\s): \nabla\phi(y)=\nabla\psi(y)\}$. We need to show that $\mathfrak{B}=\mathrm{Int}(\s)$. We prove this by contradiction. Suppose that there exists $x\in \mathrm{Int}(\s)\backslash \mathfrak{B}$. Let us first assume that $\partial \psi(x)$ is a singleton set. Since the Lebesgue measure of $\mathrm{Int}(\s)\backslash \mathfrak{B}$ is $0$, $\mathfrak{B}$ is dense in $\mathrm{Int}(\s)$. Hence, there exists a sequence $\{x_n\}\subset \mathfrak{B}$ such that $x_n\to x$ as $n\to \infty$. Note that $\nabla\phi(x_n)=\nabla\psi(x_n)$ for all $n \ge 1$. Since $\nabla\phi$ is continuous in $\mathrm{Int}(\mathcal{S})$, therefore, $\nabla\phi(x_n)$ converges to $\nabla\phi(x)$. Since any limiting point of $\nabla\psi(x_n)$ belongs to $\partial \psi(x)$, therefore, $\nabla\phi(x)\in \partial \psi(x)$. As $\partial \psi(x)$ is a singleton set, so, $\nabla\psi(x)=\nabla\phi(x)$. 

In order to complete the proof, it suffices now to assume that $\partial \psi(x)$ has more than one point. Note that $\partial \psi(x)$ is a closed convex set. Let us fix $z_1,z_2\in \partial \psi(x)$ such that $\|z_1-z_2\|\geq \min\{\mathrm{diam}(\partial \psi(x)),1\}$. For any $\delta>0$, we claim that there exist $\epsilon>0$ and $w_1,w_2\in B_{\epsilon}(x)\cap \mathfrak{B}$ such that $\|\nabla\psi(w_i)-z_i\|\leq \delta$ for $i=1,2$. The proof of this claim is very similar to the proof of Claim~\ref{cl:EqCvx}. So, we omit the details here. Due to this claim there exists two sequences $\{x^{(1)}_n\}_{n\geq 1},\{x^{(2)}_n\}_{n\geq 1}\subset \mathfrak{B}$ such that $x^{(i)}_n\to x$ as $n\to \infty$ for $i=1,2$, but, $\nabla\phi(x^{(1)}_n)\to z_1$ and $\nabla\phi(x^{(2)}_n)\to z_2$. But, this contradicts the continuity of $\nabla \phi$ at $x$. Hence, the result follows.     \qed

\section{Proof of the results in Section~3}\label{pf:Q-R}
\subsection{Proof of Proposition~\ref{thm:Q_prop}}\label{sec:Q_prop}
{\textbf{Proof of} (a)}: Our proof of this proposition will be similar in spirit to the proof of Theorem~1.1 of \cite{Figalli18}. 
We denote Lebesgue measure of any Borel set $E\subset \RR^d$ by $|E|$. Define the probability measure $\rho_{\psi^{*}}$ as 
\begin{align}
\rho_{\psi^{*}}(E) := \big|\partial \psi^{*}(E)\big|,
\end{align}
for any Borel set $E\subset \mathcal{Y}$. By Lemma~\ref{lem:SubD}, we have $\partial \psi^{*}(E) =(\partial \psi)^{-1}(E)$. This implies that
$\rho_{\psi^{*}}(E) = |(\partial \psi)^{-1}(E)|.$
Fix any open bounded set $\Omega \subset \mathcal{Y}$. Now, we claim and prove that, there exist $0<\lambda\leq \Lambda $ such that  
\begin{align}\label{eq:UpLowBd-2}
\lambda |B|\leq \rho_{\psi^{*}}(B)\leq \Lambda |B| 
\end{align}
for any Borel set $B\subset \Omega$. To see \eqref{eq:UpLowBd-2}, fix a Borel set $B\subset \Omega$. Let $n\in \NN$ be such that $\Omega\subset K_n$. Using the definition of the measure $\rho_{\psi^{*}}$ and observing that $(\partial \psi)^{-1}(B)\subset \mathcal{S}$, we may write  
\begin{align}
\rho_{\psi^{*}}(B) = \int_{ (\partial \psi)^{-1}(B)} dx = \int_{(\partial \psi)^{-1}(B)}p_{\mathcal{S}}(x)\frac{1}{p_{\mathcal{S}}(x)} dx \label{eq:myRight}
\end{align}
where $p_{\mathcal{S}}(\cdot)$ is the density of the probability measure $\mu$ supported on $\mathcal{S}$. According to our assumption, $p_{\mathcal{S}}$ is upper and lower bounded by positive constants everywhere in $\mathcal{S}$. 
Applying the upper and lower bound on $p_{\mathcal{S}}(\cdot)$ into the right side of \eqref{eq:myRight} and invoking the relation $$\int_{(\partial \psi)^{-1}(B)}p_{\mathcal{S}}(x) dx =\mu((\partial \psi)^{-1}(B))=\nu(B)$$ yields 
\begin{align}
c\nu(B)\leq \rho_{\psi^{*}}(B)\leq C \nu(B)\label{eq:RhoUpLow}
\end{align}
for some constants $0<c\leq C$. Since $B\subset K_n$, we get $\lambda_n |B|\leq \nu(B)\leq \Lambda_n |B|$ via \eqref{eq:sandwitch}. Substituting these inequalities into \eqref{eq:RhoUpLow} yields \eqref{eq:UpLowBd-2}.  

If $\psi^{*}$ is a strictly convex function on $\Omega$, then, by \cite[Corollary~4.21]{Figalli} and \eqref{eq:UpLowBd-2}, $\psi^{*}$ belongs to the class\footnote{We denote the class of functions whose $k$th derivative is $\alpha$-H\"older continuous by $\mathcal{C}^{k,\alpha}$.} $\mathcal{C}^{1,\alpha}$, for some $\alpha >0$ in $\Omega$; thus $\nabla \psi^{*}$ is continuous everywhere in $\Omega$. By repeating the same argument we can show that $\nabla\psi^{*}$ will be continuous in any open bounded set $\Omega \subset \mathcal{Y}$ if $\psi^{*}$ is strictly convex on $\mathcal{Y}$. This shows the continuity of $\nabla \psi^{*}$ in $\mathrm{Int}(\mathcal{Y})$. Once $\nabla \psi^{*}$ is continuous in $\mathrm{Int}(\mathcal{Y})$, we may continuously extend it to the boundary of $\mathcal{Y}$. Hence, under the condition of strict convexity of $\psi^{*}$, the map $\nabla\psi^{*}$ will be continuous everywhere in $\mathcal{Y}$ and then, by Lemma~\ref{lem:Injectivity}, $\nabla\psi^{*}$ will be injective. Combining continuity and injectivity of the map $\nabla \psi^{*}$ with the fact\footnote{This follows by combining our assumption $\partial \psi(\mathrm{Int}(\mathcal{S})) = \mathrm{Int}(\mathcal{Y})$ with Lemma~\ref{lem:SubD}.} that $\nabla \psi^{*}(\mathrm{Int}(\mathcal{Y})) = \mathrm{Int}(\mathcal{S})$ implies that $\nabla \psi^{*}$ is a homeomorphism from $\mathrm{Int}(\mathcal{Y})$ to $\mathrm{Int}(\mathcal{S})$. So, it suffices to show that $\psi^{*}$ is strictly convex everywhere in $\mathrm{Int}(\mathcal{Y})$.  We will show this last statement in the rest of the proof. 

Suppose first that $\mathcal{Y}$ is a bounded set. 
Note that \eqref{eq:UpLowBd-2} shows that the required condition for \cite[Corollary~4.11]{Figalli} is satisfied for the convex function $\psi^{*}$ in the interior of $\mathcal{Y}$, thereby implying that $\psi^{*}$ will be strictly convex in $\mathrm{Int}(\mathcal{Y})$. 

 In the rest of the proof, we assume $\mathcal{Y}$ is unbounded. Let $y\in \mathrm{Int}(\mathcal{Y})$ and fix $\theta \in \partial \psi^{*}(y)$. Define $$\ell(z) := \psi^{*}(y)+ \langle \theta, z-y\rangle , \forall z\in \RR^d, \quad\text{and}\quad  \Sigma := \{z\in \RR^d:\psi^{*}(z)=\ell(z)\}.$$ If $\Sigma$ is singleton, then, $\psi^{*}$ is strictly convex at $y$. Suppose $\Sigma$ is not a singleton set. As \eqref{eq:UpLowBd-2} is satisfied for $\psi^{*}$ for any convex set $K_n$ (with $n\in \NN$), therefore, \cite[Theorem~4.10]{Figalli} shows that there is no \emph{exposed point}\footnote{An exposed point of convex set $\mathcal{A}$ is a point $x\in \mathcal{A}$ where some linear functional attains its strict maximum over $\mathcal{A}$.} of $\Sigma$ in the compact set $K_n$. If $\Sigma$ does not have any exposed point, then, there are only two possible ways (see \cite[Theorem~A.10]{Figalli}) in which $\Sigma$ may not be a singleton set. Those are written as follows: 
\begin{enumerate}
\item $\Sigma$ contains a full-line in $\mathcal{Y}$. 
\item $\Sigma$ contains a half-line starting from $y$.  
\end{enumerate}
In the following, we show that neither of these two conditions can hold under assumption~\eqref{eq:sandwitch}.

Let us first suppose that $\Sigma$ contains a full-line in $\mathcal{Y}$. Fix $z\in \mathrm{Int}(\mathcal{Y})$ and $\theta_1 (\neq \theta)\in \partial\psi^{*}(z)$. Since $\Sigma$ contains a full line, therefore, there exists a vector $e_0$ such that $y+te_0\in \Sigma$ for all $t\in \RR$. As $\theta \in \partial \psi^{*}(y+te_0)$, hence, because of convexity, one may write 
\begin{align}
\langle \theta - \theta_1, y+te_0 - z\rangle\geq 0.
\end{align}
Letting $t$ to $-\infty$ and $+\infty$, we see that 
\begin{align}
\langle \theta -\theta_1, e_0\rangle \leq 0, \quad \langle \theta -\theta_1, e_0\rangle \geq 0
\end{align}
which implies that $\langle \theta -\theta_1, e_0\rangle =0$.
Since the choice $z$ was arbitrary, $\partial \psi^{*}(\mathrm{Int}(\mathcal{Y}))$ is contained in the hyperplane $\{u:\langle \theta -u, e_0\rangle = 0\}$. This contradicts the fact that $\partial\psi^*\# \nu =\mu$ (follows from Lemma~\ref{lem:QRrelation}) because $\mu(\{u:\langle \theta -u, e_0\rangle = 0\}) =0$. Thus, $\Sigma$ cannot contain a full line.


Now, it remains to show that $\Sigma$ does not contain a half-line starting from $y$. We prove this by contradiction. Let $\Sigma$ contain a half-line starting from $y$. For any $y_1\in \mathrm{Int}(\mathcal{Y})$ and $\theta_1\in \partial \psi^{*}(y_1)$, 
  \begin{align}\label{eq:SubGrad}
  \langle \theta_1 - \theta, y_1-y\rangle\geq 0.  
\end{align}   
As $\Sigma$ contains a half-line starting from $y$, therefore, $y+te_1 \in \Sigma$ for all $t\geq 0$ where $e_1$ is some vector in $\RR^d$. As $\theta \in \partial \psi^{*}(y+te_1)$ (for $t\ge 0$) plugging $y+t e_1$ in place of $y$ into \eqref{eq:SubGrad} yields 
\begin{align}
\langle \theta_1-\theta, y_1 - y -te_1\rangle\geq 0 
\end{align}
which, by letting $t$ to $+\infty$, implies  
\begin{align}\label{eq:HalfLineSep}
\partial \psi^{*}(\mathrm{Int}(\mathcal{Y}))\subset H:=\Big\{ z:\langle z - \theta, e_1\rangle \leq 0\Big\}.
\end{align}
 Since $\mathrm{Int}(\mathcal{S})\subset\partial\psi^{*}(\mathrm{Int}(\mathcal{Y})) $ (which follows from the assumption $\partial \psi(\mathrm{Int}(\s)) = \mathrm{Int}(\Y)$ and Lemma~\ref{lem:SubD}), \eqref{eq:HalfLineSep} implies that $\mathrm{Int}(\mathcal{S}) \subset H$ and $\theta\in \mathrm{Bd}(\mathcal{S}\cap H)$. The last fact follows from the observations: (i) $\theta \in \s$, and (ii) $\theta \in \mathrm{Bd}(H)$ (which can be seen by taking $z=\theta$ in~\eqref{eq:HalfLineSep}). 
 Without loss of generality, we may assume that $\|e_1\|= 1$. Let $e_1,e_2, e_3,\ldots , e_d$ be a orthonormal set of vectors in $\RR^d$. For any vector $z\in \RR^d$, define 
 \begin{align}
 \|z\|_{\hat{e}_1} := \sum_{i=2}^{d} \langle z, e_i\rangle^2. 
\end{align}  
 For any $\tau>0$ small, let us define  
 \begin{align}
 U_{\tau}&:=\mathrm{Int}(\mathcal{Y})\cap \Big\{x\in \RR^d: \langle x-y, e_1\rangle \geq 0, \|x-y\|_{\hat{e}_1}\leq \tau \langle x-y, e_1\rangle\Big\}, \\
  V_{\tau} &:=\mathrm{Int}(\mathcal{S})\cap \Big\{z\in \RR^d: -\tau  \|z-\theta\|_{\hat{e}_1}\leq \langle z-\theta, e_1\rangle \leq 0\Big\}.
 \end{align}
Now, we claim and prove that $\partial \psi^{*}(U_{\tau})\subset V_{\tau}$. Let $\hat{z}\in U_{\tau}$ and fix $z\in \partial \psi^{*}(U_{\tau})$. Using convexity of $\psi^{*}$, we get 
\begin{align}
\langle  z - \theta, \hat{z} - y -te_1\rangle\geq 0, \qquad \forall \; t\geq 0. 
\end{align} 
Setting $t=0$ we get $\langle z -\theta, \hat{z} - y\rangle \geq 0$ and letting $t\to +\infty$ we see $\langle z-\theta, e_1\rangle \leq 0$. We may now write 
\begin{align}
\qquad 0\leq \langle z -\theta, \hat{z} - y\rangle  = \sum_{i=1}^{d}\langle z -\theta, e_i\rangle \cdot\langle e_i, \hat{z} - y\rangle. \label{eq:SplitInProd}
\end{align}
Note that as $\hat{z}\in U_{\tau}$, $\langle z -\theta, e_1\rangle\cdot \langle e_1, \hat{z} - y\rangle  =|\langle \hat{z} -y, e_1\rangle|\cdot \langle e_1, z - y\rangle$ and by Cauchy-Schwarz inequality, we have 
\begin{align}
\sum_{i=2}^{d}\langle z -\theta, e_i\rangle \cdot\langle e_i, \hat{z} - y\rangle \leq \|z-\theta\|_{\hat{e}_1} \|\hat{z} - y\|_{\hat{e}_1}.
\end{align}  
Plugging these into the right side of~\eqref{eq:SplitInProd} and applying $\|\hat{z}-y\|_{\hat{e}_1}\leq \tau |\langle\hat{z} -y, e_1\rangle|$ and simplifying yields $0 \ge \langle z-  \theta, e_1\rangle\geq -\tau \|z-\theta\|_{\hat{e}_1}$ which shows that $z\in V_{\tau}$. Hence, $\partial\psi^{*}(U_{\tau})\subset V_{\tau}$.

Since $U_{\tau}$ contains $y$ and $y+e_1$, following a geometric argument of \cite[Theorem~5.1, Figure~1]{FKM09} (see also \cite{CF19}), we see that the Lebesgue measure of $U_{\tau}\cap B_2(y)$ is bounded below by $C_1\tau^{d-1}$ for some positive constant $C_1>0$ when $\tau$ is small. Using a similar argument as in \cite[Theorem~5.1]{FKM09} (see also \cite{Figalli18}), it also follows that the Lebesgue measure of $V_{\tau}$ is bounded above by $C_2 \tau^{d+1}$, where $C_2 >0$. Combining all these, we get, for some $c_1, c_2 >0$, 
\begin{align}
c_2 \tau^{d+1}\geq \mu(V_{\tau})&\geq \mu(\partial \psi^{*}(U_{\tau}\cap B_2(y)))\\&= \mu ((\partial \psi)^{-1}(U_{\tau}\cap B_2(y)))= \nu(U_{\tau}\cap B_2(y))\geq c_1 \tau^{d-1}.\label{eq:string}
\end{align}
The first inequality follows since the Lebesgue measure of $V_{\tau}$ is bounded above by $C_2\tau^{d+1}$ and $\mu$ has bounded density on $\mathcal{S}$. The second inequality follows since $\partial\psi^{*}(U_{\tau})\subset V_{\tau}$. The first equality is a consequence of Lemma~\ref{lem:SubD} and the second equality holds due to $\partial \psi\# \mu=\nu $. The last inequality follows from combining \eqref{eq:sandwitch} with the fact that the Lebesgue measure of $U_{\tau}\cap B_2(y)$ is lower bounded by $C_1\tau^{d-1}$.

Letting $\tau\to 0$, we see that \eqref{eq:string} cannot hold. Thus, $\Sigma$ cannot have a half-line in $\mathcal{Y}$ starting at $y$. This shows that $\Sigma$ is a singleton set for all $y \in \mathrm{Int}(\mathcal{Y})$. Hence, $\psi$ must be a strictly convex function. This completes the proof of $(a)$.

\noindent {\textbf{Proof of}} (b): Since $\nabla\psi^{*}$ is an one-to-one and onto map from $\mathrm{Int}(\Y)$ to $\mathrm{Int}(\s)$, due to Lemma~\ref{lem:SubD}, $\nabla \psi$ is also an one-to-one and onto map from $\mathrm{Int}(\Y)$ to $\mathrm{Int}(\s)$. Moreover, Lemma~\ref{lem:SubD} also implies that $\nabla \psi = (\nabla \psi^{*})^{-1}$ in $\mathrm{Int}(\s)$. From part-$(a)$, we know $\nabla \psi^{*}$ is a homeomorphism from $\mathrm{Int}(\Y)$ to $\mathrm{Int}(\s)$. Hence, $\nabla \psi$ is a homeomorphism from $\mathrm{Int}(\s)$ to $\mathrm{Int}(\Y)$. This completes the proof. \qed

\subsection{Proof of Lemma~\ref{cor:RankProp}}\label{sec:RankProp}
By the definition of the rank map of $\nu$ (w.r.t. $\mu$), there exists a convex function $\psi:\s \to \RR \cup \{+\infty\}$ such that $\partial\psi \# \mu = \nu$ and $R(x)= \nabla \psi^{*}(x)$ where $\psi^{*}$ is the Legendre-Fenchel dual of $\psi$. Since $R$ is a homeomorphism from $\mathrm{Int}(\Y)$ to $\mathrm{Int}(\s)$, so is $\nabla \psi^{*}$. Note that $\psi^{*}$ satisfies the same conditions as the function $\phi$ in Lemma~\ref{lem:UnifConv}. Therefore,  Lemma~\ref{cor:RankProp} follows directly from $(a)$ and $(c)$ of Lemma~\ref{lem:UnifConv}.  
\qed

\subsection{Proof of Lemma~\ref{lem:RAltDef}}\label{sec:RAltDef}
Fix $y\in \RR^d$. Suppose $\hat{R}_n(y)= u_0$. Then, by Lemma~\ref{lem:SubD}, $y\in \partial \hat{\psi}_n(u_0)$. From the definition of the subdifferential of a convex function, we may write 
\begin{align}
\langle y, u_0\rangle - \hat{\psi}_n(u_0)\geq \langle y, u\rangle - \hat{\psi}_n(u)
\end{align}
for any $u \in \s$. Owing to the last inequality, 
\begin{align}
u_0 \in \argmax_{u \in \s}\{\langle y, u\rangle - \hat{\psi}_n(u)\}.
\end{align}
This completes the proof of the lemma. \qed

\subsection{Proof of Lemma~\ref{lem:Char3}}\label{sec:Char3}
As $y= \mathrm{Cl}(W_{i_1}(\hat{h}))\cap \ldots \cap \mathrm{Cl}(W_{i_{d+1}}(\hat{h}))$ we have
\begin{align}\label{eq:haty}
\hat{\psi}_n(y) = \langle y, X_{i_1}\rangle + \hat{h}_{i_1} =\ldots = \langle y, X_{i_{d+1}}\rangle + \hat{h}_{i_{d+1}}  
\end{align}
where $\nabla \hat{\psi}_n\# \mu = \nu_n$.
As a consequence, we get
\begin{align}
\partial \hat{\psi}_n(y) = \mathrm{Conv}(X_{i_1}, \ldots , X_{i_{d+1}}).
\end{align}
Since $x\in \mathrm{Int}\big(\mathrm{Conv}(X_{i_1}, \ldots , X_{i_{d+1}})\big)$, by Lemma~\ref{lem:SubD}, we have $y\in \partial\hat{\psi}^{*}_n(x) $. If $\hat{\psi}^{*}_n$ is affine in a neighborhood of $x$, then, $\nabla \hat{\psi}^{*}_n(x)= y = \hat{R}_n(x)$. Throughout the rest of the proof, we show that $\hat{\psi}^{*}_n$ is affine in $\mathrm{Int}\big(\mathrm{Conv}(X_{i_1}, \ldots , X_{i_{d+1}})\big)$. For this it suffices to show that 
\begin{align}
\hat{\psi}^{*}_n(x) = \sum_{k=1}^{d+1}\theta_k\hat{\psi}^{*}_n(X_{i_k}), \quad \text{such that}\;\; x=\sum_{k=1}^{d+1}\theta_kX_{i_k}  \label{eq:ToProve}
\end{align}
where $\theta_k\geq 0$ for all $k\in \{1,\ldots,d+1\}$ and $\sum_{k=1}^{d+1} \theta_k=1$.
By convexity, we know 
\begin{align}
\hat{\psi}^{*}_n(x) \leq \sum_{k=1}^{d+1}\theta_k\hat{\psi}^{*}_n(X_{i_k}).\label{eq:ForIneq}
\end{align}
We now have to show that the reverse inequality also holds.
 From the definition of Legendre-Fenchel dual,
\begin{align}
\hat{\psi}^{*}_n(x) \geq \langle y, x\rangle - \hat{\psi}_n(y)= \sum_{k=1}^{d+1}\theta_k(\langle y, X_{i_k}\rangle -\hat{\psi}_n(y))= -\sum_{k=1}^{d+1}\theta_k \hat{h}_{i_k}\label{eq:RevIneq}
\end{align}
where the last equality is obtained by using \eqref{eq:haty}. By Lemma~\ref{lem:ExtremePt}, $\hat{\psi}^{*}_n(X_i)= -\hat{h}_i$ for all $i=1,\ldots, n$. Hence, from \eqref{eq:RevIneq}, we have $\hat{\psi}^{*}_n(x) \ge \sum_{k=1}^{d+1}\theta_k\hat{\psi}^{*}_n(X_{i_k})$ a.s. Combining this with \eqref{eq:ForIneq} proves \eqref{eq:ToProve}. This completes the proof. \qed

\subsection{Proof of Lemma~\ref{lem:Uniform}}\label{pf:Uniform}
It suffices to show that for any Borel set $B\subseteq \s$, $\mathbb{P}(\hat{R}_n(X_i)\in B)= \mu(B)$. Let $\mathfrak{S}_{n}$ be the set of all permutations of $\{1,\ldots ,n\}$. Let $\sigma\in \mathfrak{S}_n$ be a random permutation uniformly sampled from $\mathfrak{S}_{n}$ and independent of $X_1, \ldots , X_n$. Note that 
\begin{align*}
(X_1, \ldots , X_n) \stackrel{d}{=}(X_{\sigma(1)}, \ldots , X_{\sigma(n)}).
\end{align*}
Furthermore, $\hat{R}_n$ is a random map which does not depend on the permutation of $X_i$'s. Owing to this, we have $\hat{R}_n(X_i)\stackrel{d}{=} \hat{R}_n(X_{\sigma(i)})$ which yields  
\begin{align}
\mathbb{P}(\hat{R}_n(X_i)\in B) &= \mathbb{P}(\hat{R}_n(X_{\sigma(i)})\in B) = \mathbb{E}[\mathbbm{1}(\hat{R}_n(X_{\sigma(i)})\in B)]\nonumber\\
& = \mathbb{E}_{X}\Big[\mathbb{E}_{\sigma}[\mathbbm{1}(\hat{R}_n(X_{\sigma(i)})\in B)]\Big] \nonumber
\end{align}
where $\mathbb{E}_{X}$ denotes the expectation w.r.t.~$\{X_1, \ldots , X_n\}$ and $\mathbb{E}_\sigma$ denotes the expectation w.r.t.~$\sigma$ and the randomization in $\hat{R}_n(X_{\sigma(i)})$. Taking the expectation w.r.t.~the randomness in the definition of the rank map (see~\eqref{eq:Rank-X_i}), we see that 
\begin{align}
 \mathbb{E}_{X}\Big[\mathbb{E}_{\sigma}[\mathbbm{1}(\hat{R}_n(X_{\sigma(i)})\in B)]\Big]
& = \mathbb{E}_{X}\Big[\frac{1}{n}\sum_{i=1}^{n} n \; \mu(W_i(\hat h) \cap B) \Big] \nonumber \\
& = \mathbb{E}_{X}\Big[\mu\big(B\big)  \Big] = \mu(B) \nonumber 
\end{align}
where $W_i(\hat h)$ is the cell decomposition of $\s$ which corresponds to $X_i$ under the rank map $\hat{R}_n$ and we have used the fact that $\mu(\cup_{i=1}^n W_i(\hat h) \cup B) = \mu(B)$. 
This completes the proof. \qed

\section{Proof of results in Section~4}\label{sec:UnifConv-App}
Theorem~\ref{thm:GCProp} aims towards a unification of the asymptotic results of the Monge-Kantorovich ranks and quantiles of \cite[Theorem 3.1]{Cher17} and the center-outward ranks and quantiles of \cite[Propositions~1.5.1 and 1.5.2]{dCHM}. In~\cite{Cher17}, the authors show the convergence (in probability) of the Monge-Kantorovich ranks and quantiles by relating them to the solutions of the dual problem of the Kantorovich relaxation (to Monge's problem). This correspondence works under the assumption that both the reference and target measures have finite second moments. Hallin et al.~\cite{dCHM} avoided this dependence on the finite second moment assumption by defining the center-outward ranks and quantiles on the basis of McCann's construction~\cite{McCann95} of transport maps as limits of cyclically monotone maps. 
For proving our uniform consistency result, we marry the weak convergence theory of the Monge-Amp\`{e}re measures (see~\cite{G16}) and the recently introduced theory of graphical convergence of transport maps (see~\cite{del2019}); the connections between these two notions of convergence are made precise in Lemmas~\ref{KeyLemma},~\ref{KeyLemma2} and~\ref{lem:UnifConv} below.
By doing so we neither have to impose any finite moment assumptions on the target measure nor have to interpret the quantiles or ranks as limits of cyclically monotone maps.  

\subsection{Some important lemmas needed in the proof of Theorem~\ref{thm:GCProp}}\label{sec:Supp-Lemmas}

One of the key ingredients in the proof of Theorem~\ref{thm:GCProp} is the following lemma which investigates some limiting properties of (sub)-gradients of a sequence of convex functions. Before proceeding to the main statement, let us provide some motivation for such a result. Suppose that $\phi_n:\R^d \to \R \cup \{+\infty\}$, for ${n\ge 1}$, is a sequence of convex functions such that the sets $\{\phi^{*}_n<\infty\}$ are uniformly bounded; here $\phi^{*}_n$ denotes the Legendre-Fenchel dual of $\phi_n$. We may think of $\phi_n$ as $\hat{\psi}^*_n$ (see~\eqref{eq:Sup}) in which case $\phi_n^* = \hat{\psi}_n$. Then, owing to Lemma~\ref{lem:SubD}, $\{\partial \phi_n(\RR^d): n \ge 1\}$,  are uniformly bounded sets. Further, on any bounded set $\mathfrak{K} \subset \R^d$, $\{\phi_n\}_{n\ge 1}$ will be a sequence of uniformly bounded equicontinous\footnote{We say a sequence of functions $\{f_n:\RR^d\to \RR\}_{n\ge 1}$ is uniformly equicontinuous if for any $\epsilon>0$ there exists $\delta>0$ such that $|f_n(x)-f_n(y)|\leq \epsilon$ whenever $\|x-y\|\leq \delta$.} functions on $\mathfrak{K}$; see~Lemma~\ref{lem:Bd}. Therefore, by the Arzela-Ascoli theorem, there exists a convex function $\phi$ such that $\phi_n$ converges to $\phi$ uniformly on $\mathfrak{K}$. However, this does not guarantee the uniform convergence of the corresponding subdifferential sets. The following lemma, proved in Section~\ref{pf:KeyLemma}, addresses the mode of convergence of the subdifferential sets when the underlying convex functions converge uniformly.  The result may be of independent interest.

\bl\label{KeyLemma} 
For $n \ge 1$, let $\phi_n:\RR^d \to \RR$ be a sequence of convex functions with $\phi_n(0)=0$ and $\partial \phi_n(\RR^d)  \subset S$ for some compact set $S\subset \RR^d$. Let $K\subset \RR^d$ be a compact set. Assume that there exists a convex function $\phi:\R^d \to \R \cup \{+\infty\}$ such that $$\sup_{u\in K^{\prime}} |\phi_n(u) - \phi(u) |\to 0 \quad \mbox{ as } \quad n \to \infty,$$ where $K^{\prime}\subset \RR^d$ is a compact set such that $\mathrm{Int}(K') \supset K$. Then, the following hold:
\begin{enumerate}
 
\item[(a)] Fix $\delta>0$. For a set $A\subset \RR^d$, define $B(\delta, A) := \cup_{x\in A} B_{\delta}(x)$. There exists $n_0=n_0(\delta) \in \N$ such that for all $n\geq n_0$, we have $\partial \phi_n(K)\subset B(\delta, \partial \phi(K))$.

\item[(b)] $\partial\phi(K) \subset S$.

%

\item[(c)] Furthermore, if $\phi$ is differentiable everywhere on $K$, then,
\begin{equation}
\sup_{u\in K} \sup_{\xi\in \partial \phi_n(u)} \|\xi- \nabla \phi(u)\|   \to 0, \quad \text{as }n\to \infty. \nonumber
\end{equation}

\item[(d)] If $\phi$ is strictly convex on $K^{\prime}$, then, for any open set $U\subset K$, there exists $n_0=n_0(U)\in \NN$ such that $\partial \phi(U)\subset \partial \phi_n(K^{\prime})$ for all $n\geq n_0$. 

\item[(e)] Let $K^{\prime\prime} \subset S$ be a compact convex set. Suppose that $\phi^{*}_n$, the Legendre-Fenchel dual of $\phi_n$, converges uniformly to another  differentiable convex function $\psi$ on $K^{\prime\prime}$, i.e., $\sup_{u\in K^{\prime\prime}} |\phi_n^*(u) - \psi(u) |\to 0$ as $n \to \infty$. Also suppose that $\nabla \psi(K^{\prime\prime}) \subset K$. For any $x\in K^{\prime\prime}$, if $\nabla\psi(x) =y\in K$, then $x\in \partial \phi(y)$.  
\end{enumerate}

\el

Lemma~\ref{KeyLemma}-(e) shows that if a sequence of convex functions and their Legendre-Fenchel duals converge uniformly to the convex functions $\phi$ and $\psi$ respectively, and if one of these two functions is everywhere differentiable, then,~\eqref{eq:Charac-Sub} of Lemma~\ref{lem:SubD} holds partially between $\phi$ and $\psi$. This raises the following question (which is also relevant to the proof of Theorem~\ref{thm:GCProp}): If \eqref{eq:Charac-Sub} holds partially between any two convex functions $\phi$ and $\psi$, then, is it true that $\phi$ is the Legendre-Fenchel dual of $\psi$? As one can guess, this is not true in general and the main reason is the invariance of \eqref{eq:Charac-Sub} under the addition of constants to $\phi$ or $\psi$. However, one may still hope that the subgradients of $\phi$ will be the same as the subgradients of the Legendre-Fenchel dual of $\psi$.   
Our next result, proved in Section~\ref{pf:KeyLemma2}, provides a sufficient condition under which we are able to validate this claim.

 \bl\label{KeyLemma2} 
Let $\mathcal{S}$ and $\mathcal{Y}$ be two sets in $\RR^d$. Suppose that $\mathcal{S}$ is bounded and there is a differentiable convex function $\psi: \mathcal{S}\to \RR$ such that $\nabla \psi:\mathrm{Int}(\mathcal{S}) \to \mathrm{Int}(\mathcal{Y})$ is a homeomorphism. Then $\psi^{*}$, the convex conjugate of $\psi$, is differentiable everywhere on $\mathrm{Int}(\mathcal{Y})$. Now, let $\phi:\RR^d\to \RR$ be another convex function such that: (i) $\partial \phi(\RR^d)\subset \mathcal{S}$, and (ii) $y= \nabla \psi(x)$ for some $x\in \mathcal{S}\; \Rightarrow \; x\in\partial \phi(y)$. Then, $\phi$ is differentiable everywhere on $\mathrm{Int}(\mathcal{Y})$, and $\nabla\phi(y) = \nabla \psi^{*}(y)$ for all $y\in \mathrm{Int}(\mathcal{Y})$.   
\el

Suppose that a sequence of convex functions $\{\phi_n\}_{n\geq 1}$ converges to another convex function $\phi$ uniformly on every compact subset of $\RR^d$. Our last result of this section Lemma~\ref{lem:UnifConv} (proved in Section~\ref{pf:UnifConv}) provides a sufficient condition under which one can show the uniform convergence of the subgradients of $\phi_n$ to those of $\phi$, over the whole of $\RR^d$.

\bl\label{lem:UnifConv}
Let $\mathcal{S} \subset \RR^d$ be a strictly convex compact set. Let $\Y \subset \RR^d$ have nonempty interior. Let $\phi_n:\RR^d \to \RR$ be a sequence of convex functions such that $\partial \phi_n(\RR^d)\subset \mathcal{S}$ for all $n\geq 1$. Suppose that $\{\phi_n\}_{n\geq 1}$ converges uniformly to a convex function $\phi:\RR^d \to \RR$ on every compact set of $\RR^d$ and $\nabla\phi$ is a homeomorphism from $\mathrm{Int}(\Y)$ to $\mathrm{Int}(\mathcal{S})$. Then, 
\begin{enumerate}
\item[(a)] $\phi$ is everywhere differentiable in $\RR^d$,
\item[(b)] $\sup_{x\in \RR^d} \sup_{y\in \partial\phi_n(x)}\|y - \nabla \phi(x)\|\to 0, \quad \text{ as }n\to \infty,$
\item[(c)] for any $x\in \RR^d$, $\lim_{\lambda \to  +\infty} \nabla \phi(\lambda x) = \argmax_{v\in \mathcal{S}} \langle x, v\rangle$.
\end{enumerate}   
\el
We observe that \cite[Section~3.2.3]{dCHM} has a similar result like Lemma~\ref{lem:UnifConv} when $\mathcal{S}$ is the closed unit ball in $\RR^d$. Lemma~\ref{lem:UnifConv} generalizes \cite{dCHM} by giving a sufficient condition on $\mathcal{S}$ under which a similar conclusion holds.

\subsection{Proof of Theorem~\ref{thm:GCProp}}\label{proof:thm-GC}
We prove part-$(a)$ by contradiction.  Let us assume that there exists $K\subset \mathrm{Int}(\s)$ compact such that 
 \begin{align}\label{eq:Contra}
 \mathbb{P}\big(\A_K\big)> 0 
 \end{align}
where $$\mathcal{A}_K := \big\{\omega \in \Omega: \sup_{u\in K} \|\hat{Q}_n(u)- Q(u)\|(\omega) \not\to 0\big\}.$$
Our goal is to show that \eqref{eq:Contra} does not hold. Let $\mathfrak{F}$ be the set of all 1-Lipschitz continuous functions on $\mathrm{Int}(\Y)$ that are bounded by absolute value 1. Denote by $$\mathfrak{D} := \left\{\omega \in \Omega: \sup_{f\in \mathfrak{F}} \left|\int fd\hat{\nu}_n- \int fd\nu \right|\to 0 \right\}.$$ As $\hat{\nu}_n$ converges weakly to $\nu$ (for a.e.~$\omega$), we know that $\mathbb{P}(\mathfrak{D})=1$. Therefore, $\mathbb{P}(\mathcal{A}_K\cap \mathfrak{D}) = \mathbb{P}(\mathcal{A}_K)$. Thus, it suffices to show that $\mathbb{P}(\mathcal{A}_K\cap \mathfrak{D})>0$ does not hold, which would yield  a contradiction and prove part-$(a)$.

Fix some $\omega \in \mathcal{A}_K\cap \mathfrak{D}$. There exists a subsequence $\{n_k\}_{k \ge 1}$ and $\delta >0$ such that $\sup_{u\in K}\|\hat{Q}_{n_k}(u)-Q(u)\|(\omega)\geq \delta$, for all $k \ge 1$. In the following four steps, we show that there exists a further subsequence $\{n_{k_\ell}\}_{\ell\ge 1}$ such that $\sup_{u\in K}\|\hat{Q}_{n_{k_\ell}}(u) - Q(u)\|(\omega) \to 0$ as $\ell \to \infty$, which will give rise to a contradiction.
 
In the following, for notational simplicity, we will drop using the sample space $\Omega$ and sample element $\omega$ from all probability statements.

 Let $\hat{\psi}_n:\R^d \to \R \cup \{+\infty\}$ be a potential function whose gradient transports $\mu$ to $\hat{\nu}_n$. Let $\hat{\psi}^{*}_n:\R^d \to \R $ be the convex conjugate (Legendre-Fenchel dual) of $\hat{\psi}_n$. Let $x_0\in \mathrm{Int}(\mathcal{S})$. For convenience, without loss of generality (as $\hat{\psi}_n$ is unique up to a constant, by Lemma~\ref{lem:Affinenv} the same is true for $\hat{\psi}^{*}_n$), we define $\hat{\psi}^{*}_n(Q(x_0) ) = 0$ for all $n\in \NN$ (see e.g.,~\cite[Lemma 2.1]{del2019}).  
 
Fix a sequence of convex compact sets $\{\mathcal{F}_{m}\}_{m \ge 1}\subset \R^d$ such that $\mathcal{F}_m\uparrow \R^d$ (see Definition~\ref{bd:SetConv}), $\mathcal{F}_m$ is a continuity set\footnote{Recall that a continuity set of a measure $\rho$ is any Borel set $B$ such that $\rho (\mathrm{Bd}(B))=0$
where $\mathrm{Bd}(B)$ is the boundary set of $B$.} of $\nu$,  and $Q(x_0)\in \mathcal{F}_m$, for all $m \ge 1$. \newline

\noindent \textbf{Step I}: We construct a convex function $\psi^{*}:\R^d\to \RR$ and a subsequence of convex functions $\{\hat{\xi}_{p}\}_{p\geq 1}\subset \{\hat{\psi}^{*}_{n_k}\}_{k\ge 1}$ such that $\hat{\xi}_{p}$ converges  to $\psi^{*}$ uniformly in $\mathcal{F}_m$ as $p\to \infty$, for all $m\geq 1$. 

\noindent {\bf Proof of Step I}: Since $\sup_{x\in \RR^d}\|\partial \hat{\psi}^{*}_n(x)\|\leq \sup_{u\in \mathcal{S}}\|u\|$ (as $\partial \hat{\psi}^{*}_n \# \mu = \hat{\nu}_n$) and $\hat \psi^{*}_n(Q(x_0)) = 0$, $\{\hat{\psi}^{*}_n\}_{n\ge 1}$ is a sequence of uniformly bounded continuous functions when restricted on any compact set $\mathcal{F}_m$, for $m\in \NN$.  Using Arzela-Ascoli's theorem we get a subsequence $\{\xi^{(1)}_{\ell}\}_{\ell}\subset \{\hat{\psi}^{*}_{n_k}\}_{k}$ such that $\sup_{x\in \mathcal{F}_1}|\xi^{(1)}_{\ell}(x)-\psi^{*}_{1}(x)| \to 0$ as $\ell \to \infty$, where $\psi^{*}_1$ is a continuous function on $\mathcal{F}_1$. In fact, $\psi^{*}_1$ is a convex function on $\mathcal{F}_1$, as for $x_1,x_2 \in \mathcal{F}_1$ and $\theta \in [0,1]$, 
\begin{eqnarray}\label{eq:VerConv}
  \psi^{*}_1(\theta x_1 + (1-\theta) x_2) & = & \lim_{\ell \to \infty} \xi^{1}_\ell(\theta x_1 + (1-\theta) x_2) \\ & \le &  \lim_{\ell \to \infty} [\theta  \xi^{1}_\ell(x_1) + (1-\theta) \xi^{1}_\ell(x_2) ] = \theta  \psi^*_{1}(x_1) + (1-\theta) \psi_{1}^*(x_2), 
 \end{eqnarray}
 where we have used the convexity of $\xi^{1}_\ell$.

Using Arzela-Ascoli repeatedly, we construct an ordered collection of sequences $\{\xi^{(1)}_{\ell}\}_{\ell}\supset \{\xi^{(2)}_{\ell}\}_{\ell}\supset \ldots $ such that, for any $m \ge 1$, 
\begin{align}\label{eq:UConv}
\sup_{x\in \mathcal{F}_m} |\xi^{(m)}_{\ell}(x)- \psi^{*}_{m}(x)|\to 0, \quad \text{ as }\ell\to \infty, 
\end{align}
where $\psi^{*}_{m}: \mathcal{F}_m \to \R$ is a convex function. From the construction, it is clear that $\{\psi^{*}_m\}_{m \ge 1}$ satisfies the \emph{tower property}, i.e., $$\psi^{*}_m\big|_{\mathcal{F}_{m-1}}= \psi^{*}_{m-1} \qquad \mbox{ for any }m \geq 2.$$ Using $\{\psi^{*}_{m}\}_{m \ge 1}$, we construct a function $\psi^{*}:\R^d \to \R$ such that $\psi^{*}\big|_{\mathcal{F}_{m}} = \psi^{*}_m$, i.e., for $x \in \mathcal{F}_m$, let $$\psi^*(x) := \psi^*_m(x).$$ Because of the tower property and the continuity, the function $\psi^{*}$ is well-defined everywhere in $\R^d$. Note that $\psi^{*}$ is a convex function. To see this, fix $x_1,x_2\in \RR^d$ and $\theta \in [0,1]$. There exists $m$ such that $x_1, x_2\in \mathcal{F}_m$. To this end, combining \eqref{eq:UConv} and $\psi^{*}\big|_{\mathcal{F}_{m}}= \psi^{*}_m$ yields $\psi^{*}(\theta x_1+(1-\theta)x_2)\leq \theta \psi^{*}(x_1)+ (1-\theta) \psi^{*}(x_2)$ which shows the convexity of $\psi^{*}$.   

For any $p\in \NN$, define $\hat{\xi}_{p} := \xi^{(p)}_p$. Hence, $\{\hat{\xi}_{p}\}_{p\geq m}\subset \{\xi^{(m)}_{p}\}_{p\geq m}$ for all $m\in \NN$. As a consequence, $\hat{\xi}_{p}$ converges uniformly to $\psi^{*}$ in $\mathcal{F}_m$ as $p\to \infty$, for all $m\geq 1$. This completes the construction of {\bf Step I}.     
\newline


\noindent \textbf{Step II:} We will construct a convex function $\widetilde{\psi}:\s \to \R \cup \{+\infty\}$ from the subsequential limit of $\{\hat{\xi}^{*}_p\}_{p\ge 1}$ where $\hat{\xi}^{*}_p$ is the Legendre-Fenchel dual of $\hat{\xi}_p$.

Before proceeding to the details of \textbf{Step II}, let us sketch the overall outline of the main ideas in the proof. For any compact set $\mathfrak{K} \subset \mathrm{Int}(\s)$ and $y_0\in \mathfrak{K}$, using the triangle inequality and Lemma~\ref{lem:Bd}, we get 
\begin{align}
\sup_{y\in \mathfrak{K}}|\hat{\xi}^{*}_p(y)|&\leq |\hat{\xi}^{*}_p(y_{0})|+\sup_{y\in \mathfrak{K}}|\hat{\xi}^{*}_p(y)-\hat{\xi}^{*}_{p}(y_{0})| \\
& \le |\hat{\xi}^{*}_p(y_{0})|+ \sup_{y\in \mathfrak{K}} \|\partial \hat{\xi}^{*}_p(y)\| \; \sup_{y\in \mathfrak{K}}\|y-y_{0}\|. \label{eq:TriangleIneq}
\end{align}
If one can show that there exists $y_0\in \mathfrak{K}$ such that $|\hat{\xi}^{*}_p(y_{0})|$, for $p\geq 1$, is uniformly bounded and $\partial \hat{\xi}^{*}_p(\mathfrak{K})$, for all large $p$, is contained in one compact set, then, the right hand side of \eqref{eq:TriangleIneq} can be uniformly bounded (as $\s$ is compact). If this is possible, then, we can apply Arzela-Ascoli's theorem (in the same way as in \textbf{Step I}) for constructing a subsequential limit of $\{\hat{\xi}^{*}_p\}_{p\geq 1}$.  
 One may notice that $\partial\hat{\xi}^{*}_p$ maps $\mathcal{S}$ to $\mathcal{Y}$ which can potentially be an unbounded set. If $\mathcal{Y}$ is bounded, we can easily bound the second term on the right hand side of \eqref{eq:TriangleIneq}. So, without loss of generality, we assume throughout the rest of \textbf{Step II} that $\Y$ is unbounded. Hence, $\{\hat{\xi}^{*}_p\}_{p\geq 1}$ may not be uniformly bounded in any compact subset of $\mathrm{Int}(\mathcal{S})$. Thus, it is not a priori clear if one can apply Arzela-Ascoli's theorem to the sequence $\{\hat{\xi}^{*}_p\}_{p\geq 1}$. The main challenge of \textbf{Step II} is to overcome this difficulty.

 The pivotal point of \textbf{Step II} is Claim~\ref{cl:MT4}, stated below, where we show that for any compact set $\mathfrak{K}$, there exists $p_0=p_0(\mathfrak{K})$ such that $\partial\hat{\xi}^{*}_{p}(\mathfrak{K})$'s are uniformly bounded, for all $p\geq p_0$. To prove Claim~\ref{cl:MT4}, we mainly need two inputs: (i) $\nabla\psi^{*}$ (see \textbf{Step I} for its definition) is a homeomorphism,  and (ii) for any compact set $\mathfrak{K}\subset \mathrm{Int}(\mathcal{S})$, there exists $m\in \NN$ such that $\mathfrak{K}\subset\mathrm{Int}(\nabla \psi^{*}(\mathcal{F}_m))$. The first of these two facts follows as a consequence of Claim~\ref{cl:MT2} (stated below) and the second is stated as Claim~\ref{cl:MT3} below. Since Claim~\ref{cl:MT4} will guarantee that $\{\partial\hat{\xi}^{*}_p(\mathfrak{K})\}_{p\geq 1}$ is uniformly bounded, for applying Arzela-Ascoli's theorem to the sequence $\{\hat{\xi}^{*}_p\}_{p\geq 1}$, it only remains to show (via \eqref{eq:TriangleIneq}) that $\{\hat{\xi}^{*}_{p}\}_{p\geq 1}$ is uniformly bounded at least at one point in $\mathfrak{K}$. This will be proved in Claim~\ref{cl:MT5} below. In what follows, we formalize the above steps with the precise statements of the claims (whose proofs are deferred to Section~\ref{Appendix-B}).

 We first state Claim~\ref{cl:MT2} (proved in Section~\ref{sec:cl:MT2}) which will help us show that $\nabla\psi^*$ is a homeomorphism from $\mathrm{Int}(\mathcal{Y})$ to $\mathrm{Int}(\mathcal{S})$. 

\begin{claim}\label{cl:MT2}
 Let $\psi$ be the Legendre-Fenchel dual of $\psi^{*}$. Then, $\partial  \psi\# \mu  = \nu$ for all $\omega\in \mathcal{A}_K \cap \mathfrak{D}$. 
 \end{claim}

Combining Claim~\ref{cl:MT2} with Theorem~\ref{thm:Brenier} and Lemma~\ref{lem:CvxEq} shows that $\nabla \psi:\mathrm{Int}(\mathcal{S})\to \mathrm{Int}(\mathcal{Y})$ and $\nabla \psi^{*}:\mathrm{Int}(\mathcal{Y})\to \mathrm{Int}(\mathcal{S})$ remain the same as we vary $\omega \in \mathcal{A}_k\cap \mathfrak{D}$. This is a consequence of the two facts: $(a)$ $\nabla \psi$ and $\nabla\psi^{*}$ are a.e.~equal to the quantile map $Q$ and the rank map $R$ (by Claim~\ref{cl:MT2} and Theorem~\ref{thm:Brenier}), and $(b)$ $Q:\mathrm{Int}(\mathcal{S})\to \mathrm{Int}(\mathcal{Y})$ and $R:\mathrm{Int}(\mathcal{Y})\to \mathrm{Int}(\mathcal{S})$ are homeomorphisms (from the assumption of Theorem~\ref{thm:GCProp}). Therefore, $\nabla\psi:\mathrm{Int}(\mathcal{S})\to \mathrm{Int}(\mathcal{Y})$ and $\nabla\psi^{*}:\mathrm{Int}(\mathcal{Y})\to \mathrm{Int}(\mathcal{S})$ are homeomorphisms and $\nabla \psi^{*}(y)= (\nabla \psi)^{-1}(y)$ for all $y\in \mathrm{Int}(\mathcal{Y})$.

Let us now fix a compact set $\mathfrak{K}\subset \mathrm{Int}(\mathcal{S})$ such that $B_{\delta_0}(x_0)\in \mathrm{Int}(\mathfrak{K})$, for some $\delta_0 >0$. 
We now state Claim~\ref{cl:MT3} which will be proved in Section~\ref{sec:cl:MT3}. 

\begin{claim}\label{cl:MT3}
Let $\mathfrak{K} \subset \mathrm{Int}(\s)$ be a compact set. There exists a compact set  $\mathfrak{J}\subset \mathrm{Int}(\Y)$ such that $\mathfrak{K} \subset \mathrm{Int}(\nabla \psi^{*}(\mathfrak{J}))$. 
\end{claim}

 Now, we are ready to state two pivotal claims (proved in Sections~\ref{sec:cl:MT4} and~\ref{sec:cl:MT5}, respectively) of \textbf{Step II}. Each of these two claims will be followed by a brief outline of the main ideas used in their proofs.   
\noindent \begin{claim}\label{cl:MT4}
Recall that $\hat{\xi}^{*}_{p}$ is the convex conjugate of $\hat{\xi}_{p}$, for all $p\geq 1$, and let $\mathfrak{K} \subset \mathrm{Int}(\s)$ be a compact set. There exist a compact set $\mathfrak{W}\subset \mathrm{Int}(\Y)$ and  $p_0=p_0(\mathfrak{K}) \in \N$  such that $\partial \hat{\xi}^{*}_{p}(\mathfrak{K})\subset \mathfrak{W}$ for all $p\geq p_0$.  
 \end{claim}
 
The proof of Claim~\ref{cl:MT4} has mainly two components, namely, $(a)$ there exists a compact set $\mathfrak{W}\subset \mathrm{Int}(\Y)$ such that $\mathfrak{K}\subset\partial \hat{\xi}_p(\mathfrak{W})$, and $(b)$ $\mathfrak{K}\cap \partial\hat{\xi}_p(\Y\backslash \mathfrak{W})=\emptyset$ for all large $p$. For showing $(a)$, we rely on Claim~\ref{cl:MT3} and the fact that $\partial\hat{\xi}_p$ converges uniformly to $\nabla\psi^{*}$ on compacts (for this we use Lemma~\ref{KeyLemma}). In order to show $(b)$, we again use Lemma~\ref{KeyLemma} along with the fact that $\nabla \psi^{*}$ is injective in $\mathrm{Int}(\mathcal{Y})$.

\begin{claim}\label{cl:MT5}
 Recall that we have fixed $\hat{\xi}_{p}
  (Q(x_0)) = 0$ before \textbf{Step I}. For any $\delta>0$, there exists $p_0=p_0(\delta) \in \N$ and $y=y(x_0,\delta)\in B_{\delta}(x_0)$ such that $$|\hat{\xi}^{*}_{p}(y)|\leq (\|x_0\|+\delta)\|Q(x_0)\| \qquad \mbox{for all }\; p\geq p_0.$$
 \end{claim}
 
To prove Claim~\ref{cl:MT5}, we appeal to the fact that there exists $y\in \partial \hat{\xi}_p(Q(x_0))$ such that $\hat{\xi}^{*}_p(y)+ \hat{\xi}_p(Q(x_0)) = \langle y, Q(x_0)\rangle$ (since $\hat{\xi}^{*}_p$ is the Legendre-Fenchel dual of $\hat{\xi}_p$). Note that $y$ comes closer to $x_0$ as $p\to\infty$ because $\partial \hat{\xi}_p(Q(x_0))$ converges to $R(Q(x_0))=x_0$ via Lemma~\ref{KeyLemma}.

Now, we are ready to complete the construction of \textbf{Step II}. Combining Claims~\ref{cl:MT4} and~\ref{cl:MT5} with Lemma~\ref{lem:Bd} and~\eqref{eq:TriangleIneq} shows that, there exist $p_0=p_0(\mathfrak{K},\delta_0)\in \NN$, a compact set  $\mathfrak{W}\subset \mathrm{Int}(\Y)$ and $y_{x_0}\in B(x_0,\delta)$ such that for all $p\geq p_0$,
  \begin{align}
  \sup_{y\in \mathfrak{K}}|\hat{\xi}^{*}_p(y)|&\leq |\hat{\xi}^{*}_p(y_{x_0})|+\sup_{y\in \mathfrak{K}}|\hat{\xi}^{*}_p(y)-\hat{\xi}^{*}_{p}(y_{x_0})|\\  
  &\leq (\|x_0\|+\delta)\|Q(x_0)\| + \sup_{x\in \mathfrak{W}}\|x\|\times \mathrm{diam}(\mathcal{S}),\label{eq:AAstep}
\end{align}     
  where the second line of \eqref{eq:AAstep} follows since $|\hat{\xi}^{*}_{p}(y)|\leq (\|x_0\|+\delta)\|Q(x_0)\|$ and $|\hat{\xi}^{*}_p(y)-\hat{\xi}^{*}_{p}(y_{x_0})|\leq \sup_{x\in \mathfrak{W}}\|x\|\cdot\|y-y_{x_0}\|$ by Lemma~\ref{lem:Bd}. Note that \eqref{eq:AAstep} implies that $\{\hat{\xi}^{*}_{p}\}_{p \ge 1}$ is uniformly bounded in $\mathfrak{K}$ for all $p\geq p_0$.
  Hence, by applying Arzela-Ascoli's theorem, there exists a further subsequence $\{ \hat{\xi}^{*}_{p_{r}}\}_{r}\subset \{\hat{\xi}^{*}_{p}\}_{p} $ such that $\hat{\xi}^{*}_{p_r}$ converges to some convex function\footnote{Convexity of the $\widetilde{\psi}_{\mathfrak{K}}$ in $\mathfrak{K}$ follows from a similar argument as in \eqref{eq:VerConv}.} $\widetilde{\psi}_{\mathfrak{K}}:\mathfrak{K}\to \RR$ uniformly over $\mathfrak{K}$ as $r \to \infty$. Now, we fix a sequence of compact sets $\{\mathfrak{K}_{m}\}_{m\geq 1}$ such that $\mathfrak{K}_{m}\uparrow \mathrm{Int}(\mathcal{S})$ and there exists $\delta_0>0$ such that $B_{\delta_0}(x_0)\subset \mathfrak{K}_m$ for all $m\geq 1$. Owing to the construction given above, there exists a set of towering subsequences
  $$\{\hat{\xi}^{*}_{p}\}_{p\geq 1} \supset \{\hat{\xi}^{*,(1)}_{p}\}_{p\geq 1}\supset \{\hat{\xi}^{*,(2)}_{p}\}_{p\geq 1}\supset \ldots $$ 
  such that $\hat{\xi}^{*,(m)}_{p}$ converges uniformly to a convex function $\widetilde{\psi}_{m}$ in $\mathfrak{K}_m$ (similar to \textbf{Step I}). Moreover,
 $\{\widetilde{\psi}_{m}\}_{m}$ has the {towering property}, i.e., for any $m_2<m_1$ then, one has $\widetilde{\psi}_{m_1}|_{\mathfrak{K}_{m_2}}= \widetilde{\psi}_{m_2}$. Thus, one can define a convex function $\widetilde{\psi}:\mathrm{Int}(\s) \to \R$ by taking increasing limits of $\{\widetilde{\psi}_{\mathfrak{K}_{m}}\}$ such that $$\widetilde{\psi}\big|_{\mathfrak{K}_m} =\widetilde{\psi}_{m}, \quad \text{ for all }m\geq 1.$$ Further, we extend the definition of $\widetilde{\psi}$ to $\s$ by enforcing l.s.c. This completes the construction of \textbf{Step II}.


\smallskip

\noindent\textbf{Step III:} In the last two steps, we have constructed two convex functions $\psi^{*}:\RR^d\to \RR$ and $\widetilde{\psi}:\mathcal{S}\to \RR$. Here, we claim and prove that 
 $\widetilde{\psi}$ is differentiable everywhere in $\mathrm{Int}(\mathcal{S})$ and $\nabla \widetilde{\psi} = (\nabla \psi^{*})^{-1}$ in $\mathrm{Int}(\mathcal{S})$, 
 for all $\omega \in \mathcal{A}_K\cap \mathfrak{D}$.
 \smallskip 

\noindent \textbf{Proof of Step III:}
 We begin by proving that if $y = \nabla \psi^{*}(x)$ for some $y\in \mathrm{Int}(\mathcal{S})$ and $x\in \mathrm{Int}(\mathcal{Y})$, then, $x\in \partial \widetilde{\psi}(y)$. Since $\nabla\psi:\mathrm{Int}(\mathcal{Y})\to \mathrm{Int}(\mathcal{S})$ is a homeomorphism, there exists two compact convex sets $\mathfrak{K}^{\prime}\subset\mathrm{Int}(\mathcal{S})$  and $\mathfrak{K}^{\prime\prime}\subset \mathrm{Int}(\mathcal{Y})$ such that $x \in \nabla\psi^{*}(\mathfrak{K}^{\prime\prime})\subset \mathfrak{K}^{\prime}$. By the construction of \textbf{Step II}, there exists $m_0\in \NN$ such that $\mathfrak{K}^{\prime}\subset \mathrm{Int}(\mathfrak{K}_{m_0})$ and $\hat{\xi}^{*,(m_0)}_{p}$ converges uniformly to $\widetilde{\psi}$ in $\mathfrak{K}_{m_0}$. From \textbf{Step I}, we know that $\hat{\xi}_p$ converges uniformly to $\psi^{*}$ uniformly in $\mathfrak{K}^{\prime\prime}$ as $p\to \infty$. Combining this uniform convergence with Lemma~\ref{KeyLemma}-(e) yields that $x\in \partial \widetilde{\psi}(y)$. This holds for all $x\in \mathrm{Int}(\mathcal{Y})$ and $y\in \mathrm{Int}(\mathcal{S})$ satisfying $y = \nabla \psi^{*}(x)$. Hence, the claim  of \textbf{Step III} follows by Lemma~\ref{KeyLemma2}.   \smallskip

\noindent\textbf{Step IV:} Here, we combine \textbf{Steps I-III} to complete the proof of Theorem~\ref{thm:GCProp}-(a). 

%

The following claim (proved in Section~\ref{sec:cl:MT6}) summarizes the main properties of the functions $\psi^{*}$ and $\tilde{\psi}$ constructed in \textbf{Step I} and \textbf{Step II}, respectively.

\begin{claim}\label{cl:MT6}
Assume that $Q:\mathrm{Int}(\s)\to \mathrm{Int}(\Y)$ is a homeomorphism. Then, for all $\omega\in \mathcal{A}_{K}\cap \mathfrak{D}$, $\nabla\tilde{\psi}:\mathrm{Int}(\mathcal{S})\to\mathrm{Int}(\mathcal{Y}) $ and $\nabla\psi^{*}:\mathrm{Int}(\mathcal{Y})\to\mathrm{Int}(\mathcal{S})$ are  homeomorphisms. Furthermore, $\nabla\tilde{\psi} =Q$ and $\nabla\psi^{*}=R$ everywhere in $\mathrm{Int}(\mathcal{S})$ and $\RR^d$, respectively. 
\end{claim}

\noindent{\bf Proof of $(a)$:} On the event $\mathcal{A}_K\cap \mathfrak{D}$, for $\{\hat{\psi}_{n_k}\}_{k\geq 1}$ satisfying $\sup_{u\in K}\|\partial\hat{\psi}_{n_k}(u)-Q(u)\|\geq \delta$ and $\hat{\psi}_{n_k}\# \mu= \hat{\nu}_{n_k}$, we have constructed (in \textbf{Steps I-III}) two convex functions $\widetilde{\psi}$ and $\psi^{*}$ such that $\nabla\widetilde{\psi}:\mathrm{Int}(\mathcal{S})\to \mathrm{Int}(\mathcal{Y})$ and $\nabla\psi^{*}:\mathrm{Int}(\mathcal{Y})\to \mathrm{Int}(\mathcal{S})$ are homeomorphisms and for any two compact sets $K \subset \mathrm{Int}(\mathcal{S})$ and $K^{\prime} \subset \mathrm{Int}(\Y)$ there always exists a further subsequence $\{n_{k_\ell}\}_{\ell\geq 1}$ satisfying 
\begin{align}\label{eq:CVXunifConv}
\sup_{u\in K}\|\hat{\psi}_{n_{k_\ell}}(u) - \widetilde{\psi}(u)\| \to 0, \quad \sup_{v\in K^{\prime}}\|\hat{\psi}^{*}_{n_{k_\ell}}(v) - \psi^{*}(v)\|\to 0,
\end{align}
as $\ell\to \infty$, with $\nabla \widetilde{\psi} =Q$ and $\nabla \psi^{*} = R$ everywhere in $\mathcal{S}$ and $\mathcal{Y}$ respectively. Combining \eqref{eq:CVXunifConv}  with Lemma~\ref{KeyLemma}-(c), we get, as $\ell\to \infty$,
\begin{align}
\sup_{u\in K}\|\partial\hat{\psi}_{n_{k_\ell}}(u) - \nabla\widetilde{\psi}(u)\| \to 0, \qquad \sup_{v\in K^{\prime}}\|\partial\hat{\psi}^{*}_{n_{k_\ell}}(v) - \nabla\psi^{*}(v)\|\to 0.
\end{align}   
This implies $\mathcal{A}_K\cap \mathfrak{D}$ is empty. Hence, $\mathbb{P}(\mathcal{A}_K)= \mathbb{P}(\mathcal{A}_K\cap \mathfrak{D})=0$. 
This completes the proof.

\noindent{\bf Proof of $(b)$:}
On the event $\mathfrak{D}$, for any subsequence  $\{n_k\}_{k\geq 1}$ satisfying $\sup_{v\in K^{\prime}} \|\partial \hat{\psi}^{*}_{n_k}(v)- R(v)\|\geq \delta$ (or, satisfying $\sup_{v\in \RR^d} \|\partial \hat{\psi}^{*}_{n_k}(v)- R(v)\|\geq \delta$), one can construct a  further subsequence $\{n_{k_\ell}\}_{\ell\geq 1}$ (using \textbf{Step I} and \textbf{Step II}) such that   $\hat{\psi}^{*}_{n_{k_\ell}}$ converges uniformly to a convex function $\psi^{*}$ on compacts such that $\nabla\psi^{*}:\mathrm{Int}(\s)\to \mathrm{Int}(\Y)$ is a homeomorphism. 
 
So, $\nabla\psi^{*}$ is equal to $ R$ a.e. in $\mathrm{Int}(\Y)$. Due to the continuity of $R$, $\nabla \psi^{*}$ is equal to $R$ everywhere in $\mathrm{Int}(\Y)$ by Lemma~\ref{lem:CvxEq}. Hence, by Lemma~\ref{KeyLemma}-(c), 
 \begin{align}
 \sup_{v\in \mathfrak{K}}\|\partial \hat{\psi}^{*}_{n_{k_\ell}}(v)-\nabla \psi^{*}(v)\|\to 0 \quad \text{ as }\ell\to \infty, \label{eq:CmptUConv}
 \end{align}
for all compact sets $\mathfrak{K}\subset \mathrm{Int}(\Y)$. By taking $\mathfrak{K}= K^{\prime}$, one gets 
\begin{align}
0&=\mathbb{P}\Big(\mathfrak{D}\cap \big\{\exists \{n_k\}_{k\geq 1} \text{ s.t. }\sup_{v\in K^{\prime}} \|\partial \hat{\psi}^{*}_{n_k}(v)- R(v)\|\geq \delta\big\}\Big)\\&=\mathbb{P}\Big(\big\{\exists \{n_k\}_{k\geq 1} \text{ s.t. }\sup_{v\in K^{\prime}} \|\partial \hat{\psi}^{*}_{n_k}(v)- R(v)\|\geq \delta\big\}\Big)\label{eq:IntStep}
\end{align}
where the second equality follows since $\mathbb{P}(\mathfrak{D})=1$. Taking $\{\delta_m\}_{m\geq 1}\subset (0,\infty)$ such that $\delta_m\downarrow 0$ as $m\to \infty$ and applying these in \eqref{eq:IntStep}, we see that 
\begin{align}
\mathbb{P}&\Big(\big\{\sup_{v\in K^{\prime}} \|\partial \hat{\psi}^{*}_{n}(v)- R(v)\|\nrightarrow 0\big\}\Big)\\&= \mathbb{P}\Big(\bigcup_{m\geq 1}\big\{\exists \{n_k\}_{k\geq 1}\text{ s.t. }\sup_{v\in K^{\prime}} \|\partial \hat{\psi}^{*}_{n_k}(v)- R(v)\|\geq \delta_m\big\}\Big) =0.\label{eq:GCLast}
\end{align}
This shows \eqref{eq:GC-Rank}. 

\noindent{\bf Proof of $(c)$:}
By Claim~\ref{cl:MT6}, $\nabla\psi^{*} =R$ everywhere in $\RR^d$ and by Lemma~\ref{lem:UnifConv}-$(a)$, $R$ is continuous on entire $\RR^d$. Then, \eqref{eq:GC-Rank} holds for any compact set $K\subset \RR^d$ in a similar way as in \eqref{eq:GCLast}.  We may now apply Lemma~\ref{lem:UnifConv}-$(b)$ which yields
\begin{align}\label{eq:UnifConvInside}
\sup_{v\in \RR^d}\|\partial \hat{\psi}^{*}_{n_{k_\ell}}(v)-\nabla \psi^{*}(v)\|\to 0\quad \text{ as }\ell\to \infty,
\end{align}
for all $\omega\in \mathfrak{D}$. In a similar way as in \eqref{eq:IntStep}, owing to \eqref{eq:UnifConvInside}, we get
\begin{align}\label{eq:Prob0}
\mathbb{P}\Big(\big\{\exists \{n_k\}_{k\geq 1} \text{ s.t }\sup_{v\in \RR^d} \|\partial \hat{\psi}^{*}_{n_k}(v)- R(v)\|\geq \delta\big\}\Big) =0.
\end{align}
Now, applying \eqref{eq:Prob0} in the same way as in \eqref{eq:GCLast}, yields 
\begin{align}
\mathbb{P}&\Big(\big\{\sup_{v\in \RR^d} \|\partial \hat{\psi}^{*}_{n}(v)- R(v)\|\nrightarrow 0\big\}\Big)\\&= \mathbb{P}\Big(\bigcup_{m\geq 1}\big\{\exists \{n_k\}_{k\geq 1}\text{ s.t }\sup_{v\in \RR^d} \|\partial \hat{\psi}^{*}_{n_k}(v)- R(v)\|\geq \delta_m\big\}\Big) =0.\label{eq:GCLastFiner}
\end{align}  
 Note that \eqref{eq:GCLastFiner} shows \eqref{eq:GC-RankFiner}. For any sequence $\{\lambda_n\}_{n\in \NN}$ such that $\lambda_n\to \infty$ as $n\to \infty$, if $\lim_{\lambda_n\uparrow \infty}R(\lambda_n x)$ exists, then, by \eqref{eq:GC-RankFiner}, 
 \begin{align}
 \lim_{\lambda_n\uparrow \infty} \hat{R}_n(\lambda_n x) \stackrel{a.s.}{=} \lim_{\lambda_n\uparrow \infty}R(\lambda_n x), \quad \forall x\in \RR^d. \label{eq:LimLim}
 \end{align}
 Owing to Lemma~\ref{lem:UnifConv}-$(c)$, $\lim_{\lambda_n\uparrow \infty}R(\lambda_n x)$ is equal to $\argmax_{v\in \s}\langle x,v\rangle$ for all $x\in \s$. Combining this with~\eqref{eq:LimLim} proves \eqref{eq:Asymptot}.    \qed


%
%
%
%
%

\subsection{Proof of Lemma~\ref{KeyLemma}}\label{pf:KeyLemma}
(a) We prove this by contradiction. Assume that for some $\delta>0$, there exists  a sequence of points $\{u_n\}_{n\ge 1} \subset K$ and $y_n\in \partial \phi_n(u_n)$ such that $\inf_{x\in \partial \phi(K)}\|x-y_n\|\geq \delta$ for all $n\in \NN$. Since $\{u_n\}_{n\ge 1}\subset K$ and $\{y_n\}_{n\ge 1} \subset S$, and both $K$ and $S$ are compact, there exists a subsequence $\{n_k\}_{k\ge 1}$ such that $u_{n_k}\to u\in K$ and $y_{n_k}\to y\in S$ as $k\to \infty$. 

\noindent \textbf{Claim:} $y\in \partial \phi(u)$.

\noindent \textsc{Proof of Claim:} Recall that $K\subset \mathrm{Int}(K^{\prime})$ is compact. So, there exists $r_0>0$ such that $B_{r_0}(x)\subset K^{\prime}$ for all $x\in K$ where $B_{r_0}(x)$ is a ball of radius $r_0$ around $x$. Owing to the uniform convergence of $\phi_n$ to $\phi$ on $K^{\prime}$ and $u_{n_k}\to u$, $y_{n_k}\to y$ and the subgradient inequality, one has 
\begin{align}\label{eq:Hypo}
\phi(z)\geq \phi(u) +\langle y, z-u\rangle, \quad \forall z\in B_{r_0}(u). 
\end{align}
Combining~\eqref{eq:Hypo} with Lemma~\ref{lem:SubGIneq} shows that the inequality in \eqref{eq:Hypo} holds for all $z \in \R^d$. This proves the claim.

Since, $y_{n_k}\to y$ as $k \to \infty$, so there exists $k_0$ such that $\|y_{n_{k}}-y\|\leq \delta$ for all $k\geq k_0$. This contradicts $\inf_{x\in \partial \phi(K)}\|x-y_{n_k}\|\geq \delta$ for all $k\in \NN$ as $y \in  \partial \phi(u) \subset  \partial \phi(K)$. Hence, the result follows.

(b) Fix $y\in K$ and $x\in \partial \phi(y)$. We will show that $x\in S$. Let us first assume that $\partial \phi(y)=\{x\}$. Let $\{x_n\}\subset S$ be a sequence such that $x_n\in \partial \phi_n(y)$ for all $n\in \NN$. Since $S$ is a compact set, so, there exists a subsequence $\{n_k\}_k$ such that $x_{n_k}$ converges to some point $z$ in $S$. In what follows, we show that $z$ is in fact equal to $x$. To prove this, recall that 
\begin{align}\label{eq:SubGradCond3}
\phi_{n_k}(w)\geq \phi_{n_k}(y) + \langle x_{n_k}, w-y\rangle
\end{align} 
for all $w\in \RR^d$. For any $w\in K^{\prime}$, letting $k \to \infty$ on both sides of \eqref{eq:SubGradCond3}, we see 
\begin{align}
\phi(w)\geq \phi(y) + \langle z, w- y\rangle. 
\end{align}
Repeating the argument in the proof of the claim of $(a)$, we get $z\in \partial \phi(y)$. But, $\partial \phi(y)=\{x\}$. This implies $z=x$.

Now, we prove the result when $\partial \phi(y)$ has more than one element. Fix $x\in \partial \phi(y)$. Assume $x\notin S$. Appealing to the convexity of $S$, by the separating hyperplane theorem, one can find $\theta\in \RR^d$ and $c_1 <c_2$ such that 
\begin{align}\label{eq:Sep}
\langle x, \theta - y \rangle > c_2, \quad\qquad  \langle w,\theta -y \rangle < c_1 \;\; \mbox{for all }\; w\in S.   
\end{align} 
Now note that the set of point where $\phi$ is differentiable is a dense set. Hence, one can find a point $z$ in a small neighborhood of $\theta$ such that $\langle x, z - y\rangle > c_2$ and $\partial \phi(z) = \{\nabla \phi(z)\}$. Using the convexity of $\phi$ one has 
\begin{align}
\phi(z)\geq \phi(y)+ \langle x, z-y\rangle  \quad \mbox{and} \quad \phi(y)\geq \phi(z) + \langle \nabla \phi(z), y-z\rangle  
\end{align}
which after combining shows
\begin{align}\label{eq:Combine}
\langle \nabla \phi(z), z-y\rangle \geq \phi(z) -\phi(y)\geq \langle x,z-y\rangle > c_2.
\end{align}
This implies $\langle \nabla \phi(z), z-y\rangle>c_2$ whereas we have proved that $\nabla \phi(z)\in S$ (from the first part of (b) as $\phi(\cdot)$ is differentiable at $z$) which indicates $\langle \nabla \phi(z), z-y\rangle < c_1$ via \eqref{eq:Sep}. This gives rise to a contradiction and thus proves the result.

(c) This result follows directly from \cite[Lemma~3.10]{SS11}.

(d) We denote the boundary of $K^{\prime}$ by $\mathrm{Bd}(K^{\prime})$. Let us define 
\begin{align}
\delta:= \inf_{u \in U, v\in \partial\phi(u)}\inf_{x\in \mathrm{Bd}(K^{\prime})}\Big\{\phi(x)-\phi(u) - \langle v, x-u\rangle\Big\}.
\end{align}

\noindent \textbf{Claim:} $\delta>0$. 

\noindent \textsc{Proof of Claim:} If $\delta=0$, then, there exists a sequence $\{(x_k, u_k, v_k)\}_{k\geq 0}$ such that $x_k\in \mathrm{Bd}(K^{\prime})$, $u_k\in U$, $v_k\in \partial \phi(u_k)$ and 
\begin{align}\label{eq:phiineq}
\phi(x_k)-\phi(u_k) - \langle v, x_k-u\rangle \to 0, \quad \text{ as }k\to \infty.
\end{align}
Due to the compactness of $K^{\prime}$, $K$ and $\mathcal{S}$, there exists a subsequence $(x_{k_\ell}, u_{k_{\ell}}, v_{k_{\ell}})$ such that $$x_{k_\ell}\to x_0\in \mathrm{Bd}(K^{\prime}), \quad u_{k_{\ell}} \in u_0\in K,\quad v_{k_{\ell}}\to v_0 \in \mathcal{S}.$$ 
Furthermore, since $v_{k_{\ell}}\in \partial \phi(u_{k_{\ell}})$ for all $\ell\geq 1$, we have $v_0\in \partial \phi(u_0)$. Thanks to \eqref{eq:phiineq} and the continuity of $\phi$, we have 
\begin{align}
\phi(x_0)-\phi(u_0) - \langle v_0, x_0-u_0\rangle =0\label{eq:Strict}
\end{align}
Note that $x_0\neq u_0$ because $x_0\in \mathrm{Bd}(K^{\prime})$, $u_0\in K$ and $K\subset \mathrm{Int}(K^{\prime})$. Hence, \eqref{eq:Strict} contradicts the strict convexity of $\phi$. Therefore, $\delta > 0$.     

Returning to complete the proof, due to uniform convergence of $\phi_n$ to $\phi$ on $K^{\prime}$, there exists $n_0=n_0(\delta)$ such that $|\phi_{n}(x)-\phi(x)|\leq \delta/3$ for all $x\in K^{\prime}$ and $n\geq n_0$. Hence, for all $u\in U$ and $v\in \partial \phi(u)$, 
\begin{align}
\inf_{x\in \mathrm{Bd}(K^{\prime})}\Big\{\phi_n(x) - \phi_{n}(u) -\langle v, x-u\rangle \Big\}\geq \frac{\delta}{3}, \qquad \forall n\geq n_0.\label{eq:deltageq}
\end{align}
Fix $u\in U$ and $v\in \partial \phi(u)$. Define 
\begin{align}
\theta^{(n)}_{u,v}:= \mathrm{arginf}_{x\in K^{\prime}}\Big\{\phi_n(x) - \phi_{n}(u) -\langle v, x-u\rangle \Big\}.
\end{align}
Note that the minimum value of $\phi_n(x)-\phi_{n}(u) -\langle v, x-u\rangle$ over all $x\in K^{\prime}$ is less than or equal to $0$ because $u\in K^{\prime}$. Combining this with \eqref{eq:deltageq} yields $\theta^{(n)}_{u,v}\in \mathrm{Int}(K^{\prime})$. So there exists a open ball $\Xi^{(n)}_{u,v}\subset \mathrm{Int}(K^{\prime})$ around $\theta^{(n)}_{u,v}$ such that for all $n\geq n_0$,
\begin{align}
\phi_n(x)\geq \phi_n(\theta^{(n)}_{u,v}) + \langle v, x- \theta^{(n)}_{u,v}\rangle, \quad \forall x\in \Xi^{(n)}_{u,v}. \label{eq:ConIneq}
\end{align} 
Now, by Lemma~\ref{lem:SubGIneq}, \eqref{eq:ConIneq} holds for all $x\in \RR^d$. Therefore, $v\in \partial \phi_n(\theta^{(n)}_{u,v})\subset \partial \phi_n(K^{\prime})$ for all $n\geq n_0$. This completes the proof.

(e) Suppose that $\nabla\psi(x) =y\in K$ for some $x\in K^{\prime\prime}$. By an application of $(c)$, $\sup_{z\in K^{\prime\prime}}d_H(\partial \phi_n^{*}(z), \nabla \psi(z))\to 0$ as $n\to \infty$. Hence, if $\{y_n\}_{n\ge 1}$ is a sequence such that $y_n \in \partial \phi^{*}_n(x) $, then, $\|y_n-y\|\to 0$ as $n\to \infty$, and as $K \subset \mathrm{Int}(K')$, there exists $n_0$ such that $y_n \in K'$ for all $n\geq n_0$. Furthermore, $x\in \partial \phi_n(y_n)$ (thanks to Lemma~\ref{lem:SubD}) which implies 
\begin{align}
\phi_n(w) \geq \phi_n(y_n) + \langle x, w-y_n\rangle, \qquad \mbox{for all $w\in \RR^d$ and $n\ge 1$.}
\end{align}
 Letting $n \to \infty$, using the uniform convergence of $\phi_n$ to $\phi$ on $K^{\prime}$, we see 
\begin{align}\label{eq:SubGradCond2}
\phi(w)\geq \phi(y) + \langle x, w-y\rangle
\end{align}
for all $w\in K^{\prime}$. Since $K$ is embedded inside the interior of $K^{\prime}$, so, \eqref{eq:SubGradCond2} holds for all $w$ in an open ball $B_{\delta}(y)$ around $y$ for some $\delta>0$. Now, by Lemma~\ref{lem:SubGIneq}, \eqref{eq:SubGradCond2} for all $w\in \RR^d$. Hence, $x\in \partial \phi(y)$. This completes the proof.  \qed

\subsection{Proof of Lemma~\ref{KeyLemma2}}\label{pf:KeyLemma2}
We first show that $\psi^{*}(\cdot)$ is differentiable for all $y \in \mathcal{Y}$. To see this, let us assume that $\{z_1, z_2\} \in \partial \psi^{*}(y)$. Then, by Lemma~\ref{lem:SubD} and the fact that $\psi(\cdot)$ is differentiable for all $z \in \mathcal{S}$, for $i=1,2$, we have $z_i \in\partial  \psi^{*}(y)  \Longleftrightarrow y = \nabla\psi(z_i).$ However, as $\nabla\psi$ is a homeomorphism from $\s$ to $\mathcal{Y}$, and thus one-to-one, $z_1=z_2$. Therefore,  for $y \in \mathcal{Y}$, $z = \nabla \psi^{*}(y)  \Longleftrightarrow y = \nabla \psi(z),$ i.e., $(\nabla\psi)^{-1} = \nabla \psi^{*}$ on $\mathcal{Y}$.

Let us now prove that $\nabla\phi(y)= \nabla \psi^{*}(y)$ for all $y\in \mathrm{Int}(\mathcal{Y})$ where $\nabla \phi(y)$ exists (i.e., $\phi$ is differentiable).  As $\nabla \psi$ is a homeomorphism from $\mathcal{S}$ and $\mathcal{Y}$, for any $y\in \mathcal{Y}$, there exists a unique $x\in \mathcal{S}$ such that $y= \nabla \psi(x)$. By the assumption in the lemma this implies that $x\in \partial \phi(y)$. If $\phi$ is differentiable at $y$, then, $\partial \phi(y)= \{\nabla \phi(y)\}$ which implies that $\nabla \phi(y)= x =  (\nabla\psi)^{-1}(y) = \nabla \psi^{*}(y)$.

Now, we turn to prove the result when $y\in \mathrm{Int}(\mathcal{Y})$ is a non-differentiable point of $\phi$. Suppose that $y\in \mathrm{Int}(\mathcal{Y})$ be such that $\partial \phi(y)$ contains more than one element. Note that $\partial \phi(y)$ is a closed bounded convex set (as $\partial \phi(\RR^d)\subset \mathcal{S}$). Let $z_1$ and $z_2$ be two points  in $\partial \phi(y)$ such that $\|z_1-z_2\|= \mathrm{diam}(\partial \phi(y)) >0$. Note that two such points exist as $\partial \phi(y)$ is a compact convex set.
\smallskip 

\begin{claim}\label{cl:EqCvx}
 For any given $ \delta>0$  there exist $\epsilon=\epsilon(\delta) >0$ and $w_i \in B_{\epsilon}(y) \subset \mathcal{Y}$ such that $\phi$ is differentiable at $w_i$ and $\|\nabla \phi(w_i)-z_i\|< \delta$, for $i=1,2$.
\end{claim}
\smallskip 

\noindent \textsc{Proof of Claim:} We prove the claim for $z_1$. The proof for the case of $z_2$ is similar. If $z_1 =\nabla \psi^{*}(y)$, then, the claim follows from the continuity of $\nabla \psi^{*}(\cdot)$ and the fact that the set of points where $\phi$ is differentiable is a dense set of $\RR^d$. Thus, without loss of generality, we may assume that $z_1\neq \nabla \psi^{*}(y)$. Owing to this, there exists $\delta^{\prime} >0$ such that $B_{\delta^{\prime}}(\nabla \psi^{*}(y))\cap B_{\delta^{\prime}}(z_1)= \emptyset$. Thanks to the continuity of $\nabla \psi^{*}$, there exists $\epsilon>0$ such that $ \nabla \psi^{*}(B_{\epsilon}(y))\subset B_{\delta^{\prime}}(\nabla \psi^{*}(y))$. Now, by the separating hyperplane theorem \cite[Section~11]{Rockf}, there exist $\theta \in B_{\epsilon}(y)$ and $c_1>c_2\in \RR$ such that 
\begin{align}\label{eq:CVXSep}
\qquad \langle u, \theta -y\rangle >c_1, \; \forall u\in B_{\delta^{\prime}}(z_1), \;  \langle v, \theta -y\rangle \leq c_2, \; \forall v\in B_{\delta^{\prime}}(\nabla \psi^{*}(y)).  
\end{align} 
Recalling that the set of points where $\phi$ is differentiable is a dense set, we obtain $\theta_1\in B_{\epsilon}(y)$ arbitrarily close to $\theta$ such that $\phi$ is differentiable at $\theta_1$ and $\langle u, \theta_1 -y\rangle >c_1$ for all $u\in B_{\delta^{\prime}}(z_1)$. Due to the convexity of $\phi$ and the fact that $z_1\in \partial \phi(y)$, we know that $\langle \nabla\phi(\theta_1) - z_1, \theta_1-y\rangle\geq 0$
which implies 
\begin{align}\label{eq:PhiTheta}
\langle \nabla \phi(\theta_1), \theta_1 - y\rangle>c_1.
\end{align}
However, 
$$\langle \nabla \phi(\theta_1), \theta_1 - y\rangle= \langle \nabla\psi^{*}(\theta_1), \theta_1-y\rangle \leq c_2$$ where the equality follows since $\nabla \psi^{*}(y)= \nabla \phi(y)$ for all $y$ where $\phi$ is differentiable (and $\phi$ is differentiable at $\theta_1$) and the inequality follows from the second inequality of \eqref{eq:CVXSep} combined with the facts: (i) $\theta_1\in B_{\epsilon}(y)$, and (ii) $\nabla \psi^{*}(B_{\epsilon}(y))\subset B_{\delta^{\prime}}(\nabla \psi^{*}(y))$. This contradicts \eqref{eq:PhiTheta} and hence, completes the proof of the claim. 
\smallskip 

Now, we return to complete the proof of this lemma. Owing to the last claim, we observe that $\min_{u\in B_{\epsilon}(y)}\|z_1- \nabla \psi^{*}(u)\|$ and $\min_{u\in B_{\epsilon}(y)}\|z_2- \nabla \psi^{*}(u)\|$ converge to $0$ as $\epsilon \to 0$. However, appealing to the continuity of $\nabla \psi^{*}$ yields $\sup_{u\in B_{\epsilon}(y)}\|\nabla\psi^{*}(u)- \nabla\psi^{*}(y)\|\to 0$ as $\epsilon\to 0$. Hence, the distance between $z_1$ (respectively $z_2$) and $\nabla \psi^{*}(y)$ decreases as $\epsilon$ goes to $0$. However, this contradicts $\|z_1-z_2\|=\mathrm{diam}(\partial \phi(y))>0$. Hence, $\partial \phi(y)$ is a singleton set and consequently, $\nabla\phi(y)= \nabla \psi^{*}(y)$. This completes the proof.   \qed

\subsection{Proof of Lemma~\ref{lem:UnifConv}}\label{pf:UnifConv}
\noindent{Proof of (a):} Since $\phi_n$ converges uniformly to $\phi$ on every compact set of $\RR^d$ and $\partial \phi_n(\RR^d)\subset \s$ for all $n\geq 1$, by Lemma~\ref{KeyLemma}-(b), $\partial \phi(\RR^d)\subset \s$. Hence, to show that $\phi$ is differentiable everywhere, it suffices to show $\partial \phi(x)$ is a singleton set for all $x\in \RR^d$. We prove this by contradiction. Suppose that there exists $x\in \RR^d, z_1\neq z_2\in \s$ such that $z_1,z_2\in \partial \phi(x)$. Let $\phi^{*}$ be the Legendre-Fenchel dual of $\phi$. Since $\nabla\phi$ is a homeomorphism from $\mathrm{Int}(\Y)$ to $\mathrm{Int}(\s)$, by Lemma~\ref{lem:SubD}, $\nabla\phi^{*}$ is a homeomorphism from $\mathrm{Int}(\s)$ to $\mathrm{Int}(\Y)$. As a consequence, we get if $z_1,z_2\in \mathrm{Int}(\s)$, then, $\nabla\phi^{*}(z_1)\neq \nabla\phi^{*}(z_2)$. However, this contradicts $z_1,z_2\in \partial \phi(x)$ because $\nabla\phi^{*} = (\nabla \phi)^{-1}$ in $\mathrm{Int}(\s)$ implying $x\in \nabla \phi^{*}(z_1)\cap \nabla \phi^{*}(z_2)$. So, both of $z_1$ and $z_2$ cannot belong to $\mathrm{Int}(\s)$. Now, suppose that $z_1\in \mathrm{Bd}(\s)$ and $z_2\in \mathrm{Int}(\s)$. As $z_1,z_2\in \partial \phi(x)$, so, $x\in \partial \phi^{*}(z_1)\cap \partial \phi^{*}(z_2)$ by Lemma~\ref{lem:SubD}. This implies $x$ belongs to the subgradient set of $\phi^{*}$ at all $v$ in the line segment $\mathbf{s}_{z_1,z_2}$ joining $z_1$ and $z_2$ (this follows from the definition of a convex function and the subgradient inequality). Note that $(\mathbf{s}_{z_1,z_2}\cap \mathrm{Int}(\s)\big)\backslash \{z_2\}$ is not empty because $z_2\in \mathrm{Int}(\s)$. But this contradicts injectivity of $\nabla \phi^{*}$ in $\mathrm{Int}(\s)$. Hence, none of $z_1,z_2$ belongs to $\mathrm{Int}(\s)$. To complete the proof, it suffices now to show that $z_1$ and $z_2$ cannot be on $\mathrm{Bd}(\s)$. Suppose $z_1,z_2\in \mathrm{Bd}(\s)$. Consider the line segment $\mathbf{s}_{z_1,z_2}$ joining $z_1$ and $z_2$. If $\mathbf{s}_{z_1,z_2}\cap\mathrm{Int}(\s)= \emptyset$, then, $\mathbf{s}_{z_1,z_2}$ belongs to a supporting hyperplane which touch $\mathrm{Bd}(\s)$ at more than one points. This contradicts our assumption that $\s$ is strictly convex, i.e., every supporting hyperplane of $\s$ touches $\mathrm{Bd}(\s)$ just once. Therefore, $\mathbf{s}_{z_1,z_2}\cap \mathrm{Int}(\s)$ is not empty. Fix $v\in \mathbf{s}_{z_1,z_2}\cap \mathrm{Int}(\s)$. Notice that $x\in \partial \phi^{*}(v)$ because $x\in \partial \phi^{*}(z_1)\cap \partial \phi^{*}(z_2)$. However, we have proved before that this cannot be true. Hence, the result follows.

\noindent{Proof of (b):} Since $\phi$ is differentiable everywhere in $\RR^d$, by Lemma~\ref{KeyLemma}-(c), for any compact set $K\subset \RR^d$, 
\begin{align}\label{eq:InCompacts}
\sup_{x\in K}\sup_{y\in \partial \phi_n(x)} \|y - \nabla \phi(x)\|\to 0, \quad \text{ as }n\to \infty.
\end{align}
Now, we extend the uniform convergence of $\partial \phi_n(\cdot)$  to the whole of $\RR^d$, by contradiction. Suppose that there exists $\epsilon>0$, two sequences $\{x_n\}_{n\geq 1}\subset \RR^d$ and $\{y_n\}_{n\geq 1}\subset \mathcal{S}$  such that $y_n\in \partial \phi_n(x_n)$ and
\begin{align}
\|y_n -\nabla \phi(x_n)\|\geq \epsilon, \qquad \mbox{for all }\; n\geq 1.\label{eq:EpDiff}
\end{align} 
Suppose that $x_n=\lambda_n u_n$, for some $\lambda_n\ge0$ and $\|u_n\|=1$. Then, $\lambda_n\to \infty$ as $n\to \infty$, otherwise, we will reach a contradiction to \eqref{eq:InCompacts}. As $\mathcal{S}$ and the unit sphere in $\RR^d$ are compact sets, there exists a converging subsequence of $(u_n,y_n,\nabla \phi(x_n))$. Here, we abuse notation by taking $(u_n,y_n,\nabla \phi(x_n))$ as such a converging sequence. Let $(u,y,z)$ be the limit of $(u_n,y_n, \nabla\phi(x_n))$ as $n\to \infty$.  Now, we claim and prove that 
\begin{align}
z \in \argmax_{v\in \mathcal{S}} \langle u, v\rangle,\qquad \mbox{and} \qquad   y \in \argmax_{v\in \mathcal{S}} \langle u, v\rangle.\label{eq:ZYrep} 
\end{align}

Before proceeding to the proof of \eqref{eq:ZYrep}, let us explain how \eqref{eq:ZYrep} will contradict \eqref{eq:EpDiff}. We first show that if \eqref{eq:ZYrep} holds, then $z=y$. Let us suppose that $z\neq y$. Due to the convexity of $\mathcal{S}$ and the convexity of the functional $\langle u, \cdot\rangle $ in~\eqref{eq:ZYrep}, $z,y$ belong to the boundary of $\mathcal{S}$ and $\langle u, z-y \rangle =0$. Owing to this last equality, we observe that $\langle u, z\rangle = \langle u, v\rangle$ for all $v$ in the line segment $\vec{s}_{z,y}$ joining $z$ and $y$. As all the supporting hyperplanes of $\mathcal{S}$ touch the boundary of $\mathcal{S}$ at most at one point, therefore, $\vec{s}_{z,y}\cap \mathrm{Int}(\mathcal{S})$ is  nonempty. Fix $v\in \vec{s}_{z,y}\cap \mathrm{Int}(\mathcal{S})$. There exists $\delta>0$ such that $B_{\delta}(v)\subset \mathrm{Int}(\mathcal{S})$. But then $$\sup_{w\in B_{\delta}(v)}\langle u, w\rangle> \langle u,v\rangle.$$
Combining this with \eqref{eq:ZYrep} implies $\sup_{w\in B_{\delta}(v)}\langle u, w\rangle$ is greater than $\sup_{w\in \mathcal{S}}\langle u, w\rangle$ which is a contradiction. Hence, $z=y$. However, according to \eqref{eq:EpDiff}, we have $\|z-y\|\geq \epsilon$. This leads to a contradicts our assumption. Therefore, to complete the proof of the lemma, it suffices to verify \eqref{eq:ZYrep} which is given below. 

 By the convexity of $\phi$, for any fixed $x\in \RR^d$, 
\begin{align}
\langle x_n- x, \nabla \phi(x_n) - \nabla \phi(x)\rangle \geq 0.
\end{align}   
Diving both sides of the above display by $\lambda_n$, and letting $n\to \infty$, we see that
\begin{align}
\qquad \lim_{n\to \infty}\frac{1}{\lambda_n}\langle x_n- x, \nabla \phi(x_n) - \nabla \phi(x)\rangle =\langle u, z- \nabla\phi(x)\rangle \geq 0, \; \forall x\in \RR^d.\label{eq:Monotonocity}
\end{align}
Since $z\in \mathcal{S}$ and $\nabla \phi$ is a homeomorphism from $\mathrm{Int}(\Y)$ to $\mathrm{Int}(\mathcal{S})$,~\eqref{eq:Monotonocity} implies that
\begin{align}
\langle u, z\rangle\geq \sup_{v\in \mathcal{S}} \langle u, v\rangle \quad \Rightarrow \quad z\in \argmax_{v\in \mathcal{S}}\langle u,v\rangle. 
\end{align}
Now, it remains to show the second part of \eqref{eq:ZYrep}. Fix a compact set $\mathfrak{K} \subset \mathrm{Int}(\mathcal{S})$. As $\nabla \phi$ is homeomorphism from $\RR^d$ to $\mathrm{Int}(\mathcal{S})$, there exist a compact set $K\subset \mathrm{Int}(\Y)$ such that $\nabla \phi(K)=\mathfrak{K}$. For any $x\in K$ and $v_n\in \partial \phi_n(x)$, due to the convexity of $\phi_n$, 
\begin{align}
\langle x_n -x, y_n - v_n\rangle \geq 0. \label{eq:yineq}
\end{align}
Owing to \eqref{eq:InCompacts}, $v_n$ converges to $\nabla \phi(x)$ as $n\to \infty$. Hence, dividing both sides of \eqref{eq:yineq} by $\lambda_n$ and letting $n\to \infty$ yields 
\begin{align}
\lim_{n\to \infty}\frac{1}{\lambda_n}\langle x_n- x,y_n-v_n\rangle = \langle u, y-\nabla \phi(x)\rangle\geq 0,\quad  \forall x\in K
\end{align}
which implies 
\begin{align}
\langle u, y\rangle\geq \sup_{v\in \mathfrak{K}} \langle u, v\rangle. \label{eq:EachComp}
\end{align}
Taking a sequence of compact sets $\{\mathfrak{K}_n\}_{n\geq 1}$ with $\mathfrak{K}_n\subset \mathrm{Int}(\mathcal{S})$ for all $n\geq 1$ such that $\mathfrak{K}_n \uparrow \mathrm{Int}(\mathcal{S})$ and using those in \eqref{eq:EachComp}, we observe that $\langle u, y\rangle$ is greater than or equal to the maximum value of $\langle u, v\rangle$ over all $v\in \mathcal{S}$. Note that $y\in \mathcal{S}$. Therefore, $$y\in \argmax_{v\in \mathcal{S}}\langle u, v\rangle. $$
This completes the proof. 

\noindent{Proof of (c):} We parallel the proof of \eqref{eq:ZYrep}. Thanks to the convexity of $\phi$, 
\begin{align}
\langle \lambda x - y, \nabla \phi(\lambda x) - \nabla \phi(y) \rangle \geq 0
\end{align}
for any $y \in \RR^d$. Dividing both sides by $\lambda$ and letting $\lambda\to \infty$ yields 
\begin{align}\label{eq:convCons}
\liminf_{\lambda\to \infty}\langle x,  \nabla \phi(\lambda x)\rangle\geq \langle x,\nabla \phi(y) \rangle.
\end{align}
Since $\nabla\phi(\RR^d)\subset \mathcal{S}$ by Lemma~\ref{KeyLemma}-$(b)$, we have $\max_{v\in \mathcal{S}}\langle x,v \rangle \geq \langle x,  \nabla \phi(\lambda x)\rangle$ for all $\lambda \in \RR$. Combining this with \eqref{eq:convCons} and the fact $\nabla \phi:\mathrm{Int}(\Y)\to \mathrm{Int}(\s)$ is a homeomorphism, we arrive at
\begin{align}
\qquad \;\;\;\max_{v\in \mathcal{S}}\langle x,v \rangle\geq \limsup_{\lambda \to \infty}\langle x,  \nabla \phi(\lambda x)\rangle\geq \liminf_{\lambda\to \infty}\langle x,  \nabla \phi(\lambda x)\rangle\geq \max_{v\in \mathcal{S}}\langle x,v \rangle. \label{eq:UpLowBd}
\end{align}
Note that \eqref{eq:UpLowBd} implies any limit points of the sequence $\{\nabla \phi(\lambda_{n} x)\}_{n\in \NN}$ belongs to the set $\argmax_{v\in \s}\langle x, v\rangle$ where $\lambda_n\uparrow \infty$ as $n\to\infty$. Since all the supporting hyperplanes of $\s$ touch the boundary of $\s$ at most once, due to a reason explained in part-$(b)$, $\argmax_{v\in \s}\langle x, v\rangle$ is a singleton set. This shows $\{\nabla \phi(\lambda_{n} x)\}_{n\in \NN}$ can have only one limit point which is equal to $\argmax_{v\in \s}\langle x, v\rangle$. Hence, we get $\lim_{\lambda\to \infty}\langle x,\nabla \phi(\lambda x)\rangle = \argmax_{v\in \s}\langle x, v\rangle$. 
\qed

\section{Proofs of the results in Sections~5.1,~5.2}\label{pf:Global-Rate}

\subsection{Proof of Lemma~\ref{ppn:RateProp2}}\label{pf:RateProp2}
To simplify notation let $T \equiv \nabla \psi$ and $\tilde T \equiv \nabla {\tilde \psi}$. First note that  
\begin{align}
\int \psi^* \,d ( {\tilde \nu} - \nu) =  \int [\psi^*(\tilde T(u)) - \psi^*(T(u))] d \mu(u). \label{eq:step-1}
\end{align}
Now, by~\eqref{eq:step-1}, the strong convexity of $\psi^*$, and the fact that $\nabla \psi^*(T(u)) = u$ for $\mu$-a.e.~$u$, we have 
\begin{align}
\int \psi^* \,d ( {\tilde \nu} - \nu) & \ge \int \left\{\nabla \psi^*(T(u))^\top (\tilde T(u) - T(u)) + \frac{\lambda}{2} \|\tilde T(u)-T(u)\|^2 \right\} d \mu(u)  \nonumber \\
        & = \int u^\top (\tilde T(u) - T(u)) d \mu(u) + \frac{\lambda}{2} \int  \|\tilde T(u)-T(u)\|^2  d \mu(u). \label{eq:v1}
\end{align}
Further, using the facts $$\int \|\tilde T(u)\|^2 d \mu(u) = \int \|x\|^2 d \tilde \nu(x) \quad \mbox{and} \quad \int \|T(u)\|^2 d \mu(u) = \int \|x\|^2 d \nu(x),$$ note that $2 \int u^\top (\tilde T(u) - T(u)) d \mu(u)$ equals
\begin{align}\label{eq:v2}
\quad \int  \|u - T(u)\|^2  d \mu(u) -\int  \|u - \tilde T(u)\|^2  d \mu(u)  - \int \|x\|^2 d (\nu - \tilde \nu)(x).
\end{align}
Further, as $\int \|x\|^2 d\nu(x) <+\infty$ and $\int \|x\|^2 d {\tilde \nu}(x) <+\infty$ the corresponding Wasserstein distances (w.r.t.~$\mu$) exist and we have 
\begin{equation}\label{eq:Wass-nu-t-nu}
\int  \|u - T(u)\|^2  d \mu(u) = W_2^2(\mu,\nu) \;\, \mbox{ and } \, \int  \|u - \tilde T(u)\|^2  d \mu(u) = W_2^2(\mu,\tilde \nu).
\end{equation}
Thus, combining~\eqref{eq:v1} and~\eqref{eq:v2}, and noting $g(x) =  \frac{\|x\|^2}{2} - \psi^*(x)$, we get 
\begin{equation}\label{eq:Global-Diff}
{\lambda}\int  \|\tilde T(u)-T(u)\|^2  d \mu(u) \le \left[W_2^2(\mu,\tilde \nu) - W_2^2(\mu, \nu)\right] + 2 \int g \, d (\nu - \tilde \nu)
\end{equation}
which yields the desired result. \qed

\subsection{Proof of Theorem~\ref{thm:Q-Rate}}\label{pf:Q-Rate}
Note that, from convex analysis, the assumption of strong convexity of $\psi^*$ (with parameter $\lambda>0$) implies that $\nabla \psi$ is $\lambda^{-1}$-Lipschitz. As $\nabla \psi \# \mu = \nu$ where $\mu$ has compact support $\s$, this necessarily implies that $\nabla \psi$ (restricted to $\s$) takes values in a bounded set, i.e., $\nu$ is compactly supported. 

We first use Lemma~\ref{ppn:RateProp2} to note that the left hand side of~\eqref{eq:Q-Rate} can be controlling by controlling the expectation of the right hand side of~\eqref{eq:Global-Diff} (as $\nabla {\hat \psi}_n \equiv \hat{Q}_n$, and $\nabla {\psi} \equiv Q$). From the proof of~\cite[Theorem 2]{Chizat2020} we can bound 
\begin{equation}\label{eq:Q_0stBound}
\E\left[ \big|W_2^2(\mu,\hat \nu_n) - W_2^2(\mu, \nu)\big| \right]  \le C \, r_{d,n},
\end{equation}
where $C >0$ is a constant and $\hat \nu_n$ denotes the empirical measure of $X_1,\ldots, X_n$. Moreover, as $\nu$ has a compact support and $g$ is bounded inside the support, $$\left|\E \left[ \int g \,d(\hat \nu_n - \nu) \right]  \right|<C n^{-1/2},$$ for some constant $C >0$. Adding the last two displays complete the proof of the  result.

Now, we proceed to prove \eqref{eq:Prob_bound}. By Lemma~\ref{ppn:RateProp2}, it suffices to show that there exists $c>0$ such that 
\begin{align}
\mathbb{P}\Big(\big|W^2_2(\mu, \hat \nu_n)- W^2_2(\mu, \nu)\big|>Cr_{d,n}+ n^{-1/2}s\Big) &\leq \exp(-cs^2),\label{eq:Q_1stBound}\\ \mathbb{P}\Big(\big|\int g \, d (\hat \nu_n - \nu)\big|> n^{-1/2}s\Big) &\leq \exp(-cs^2).\label{eq:Q_2ndBound}
\end{align} 
From \cite[Theorem~2]{Chizat2020} it also follows that, for some $c >0$,
\begin{align}
\mathbb{P}\Big(\big|W^2_2(\mu, \hat \nu_n)- \mathbb{E}[W^2_2(\mu, \hat \nu_n)]\big|\geq n^{-1/2}s\Big)\leq \exp(-cs^2).
\end{align}
Combining the last inequality with \eqref{eq:Q_0stBound} shows \eqref{eq:Q_1stBound}, whereas \eqref{eq:Q_2ndBound} follows from Hoeffding's inequality. This completes the proof of this result.
\qed

\subsection{Proof of Theorem~\ref{thm:RateProp1}}\label{pf:RateProp1}
Let $\hat \nu_n$ be the empirical distribution of the $X_i$'s. By the strong convexity of $\psi$ on $\s$ (with parameter $\lambda >0$), we have
\begin{align}
D_1 := & \int \left[\psi(\hat R_n(x)) - \psi(R(x)) \right] d \hat  \nu_n(x)  \nonumber \\
    \ge & \int \left\{\nabla \psi(R(x))^\top \left(\hat R_n(x) - R(x) \right) + \frac{\lambda}{2} \|\hat R_n(x)-R(x)\|^2 \right\} d \hat \nu_n(x)  \nonumber \\
      = & \int x^\top (\hat R_n(x) - R(x)) d \hat \nu_n(x) + \frac{\lambda}{2} \int  \|\hat R_n(x)-R(x)\|^2  d \hat \nu_n(x)  \label{eq:L2-v1}
\end{align}
where we have used and the fact that $\nabla \psi(R(x)) = x$ for $\hat \nu_n$-a.e.~$x$. Recall the definition of $\hat \psi_n$; see~\eqref{eq:SampQ_n}. Then, 
\begin{equation}
 \;\;\; D_2 := \int \left[\hat \psi_n(R(x)) - \hat \psi_n(\hat R_n(x)) \right] d \hat \nu_n(x) \ge  \int x^\top (R(x) - \hat R_n(x)) d \hat \nu_n(x)  \label{eq:L2-v2}
\end{equation}
where we have used the facts that: (i) $\hat \psi_n$ is convex and thus $\hat \psi_n(u) - \hat \psi_n(v) \ge \xi_v^\top (u-v)$ for any $\xi_v \in \partial \hat \psi_n(v)$,  and (ii) $x \in \partial \hat \psi_n(\hat R_n(x))$ for $\hat \nu_n$-a.e.~$x$.

Now, using~\eqref{eq:L2-v1} and~\eqref{eq:L2-v2}, we get $$D_1 + D_2 \ge \frac{\lambda}{2} \int  \|\hat R_n(x)-R(x)\|^2  d \hat \nu_n(x).$$ As $\hat \nu_n$ is the empirical distribution of the $X_i$'s, then using the randomized choice of $\hat R_n(X_i)$ as in~\eqref{eq:Rank-X_i}, we have 
$$\int \psi(\hat R_n(x)) d \hat \nu_n(x) = \frac{1}{n} \sum_{i=1}^n \psi(U_i) \quad \mbox{where} \quad \hat R_n(X_i)  \equiv U_i|X_1,\ldots, X_n \sim \hat \mu_i,$$ and thus, $$\E_U \left[\int \psi(\hat R_n(x)) \,d \hat \nu_n(x) \right] = \int \psi(u) d\mu(u),$$
where the expectation is taken only with respect to the $U_i$'s (independent conditional on the data). Similarly, $$\E_U \left[\int \hat \psi_n(\hat R_n(x)) \, d \hat \nu_n(x) \right] = \int \hat \psi_n(u) d\mu(u).$$
Thus, 
\begin{align}
\E_U [D_1+ D_2] & =  \int \left(\psi(u) - \hat \psi_n(u) \right) d \mu(u)  +  \int \left[\hat \psi_n(R(x)) - \psi(R(x))  \right] d \hat \nu_n(x) \nonumber \\
 &=  \int \left(\psi(u) - \hat \psi_n(u) \right) d (\mu - \hat \mu_n)(u).
\end{align}
where by $\hat \mu_n$ we denote the empirical distribution of $V_i :=R(X_i)$ (for $i=1,\ldots, n$) which are i.i.d.~$\mu$ (i.e., $\hat \mu_n := R \#\hat \nu_n$). As $\hat \psi_n$ is a convex function defined on $\s$ such that $\hat \psi_n(u_0) = 0$ (for some $u_0 \in \s$; w.l.o.g.) and $\nabla \hat \psi_n$ takes values inside the support of $\nu$, the absolute value of the right hand side of the above display can be bounded by $$\sup_{\varphi \in \mathcal{F}} \left| \int \varphi \, d (\hat \mu_n -\mu) \right| + \left| \int \psi \, d (\hat \mu_n - \mu)  \right|$$ where the supremum is over the class of all convex functions $\mathcal{F}$ defined on $\s$ (a compact convex set) and $\beta$-Lipschitz (where we assume that the support of $\nu$ is contained in a ball of radius $\beta$ around 0). Thus, by a similar analysis as in~\cite[Lemma 4]{Chizat2020} allows us to control the expected value of the first term in the above display. Note that, the expectation of the second term can be controlled easily as $\mu$ is compactly supported, in fact, $\E \left| \int \psi \,d(\hat \mu_n - \mu) \right|<C n^{-1/2},$ for some constant $C >0$. This completes the proof of the result.  \qed

\section{Proofs of Results in Section~5.3}\label{pf:Loc-RateSec}

As a first step to finding the local uniform rate of convergence of the empirical quantile and rank functions, we give the following result (proved in Section~\ref{sec:RateProp1}) that provides a deterministic upper bound on the local uniform rate of convergence of the subdifferentials (of a sequence of convex functions and their Legendre-Fenchel duals) by a local $L_2$-loss of the subdifferentials. 
\begin{proposition}\label{ppn:RateProp1}
Let $\mu$ be an absolutely continuous probability measure supported on a compact convex set $\mathcal{S} \subset \R^d$. Let $\nu$ be a probability measure supported on $\mathcal{Y}\subset \RR^d$ and let $\psi$ be a convex function such that $\nabla \psi \# \mu = \nu$. Suppose that $\{\hat{\nu}_n\}_{n \ge 1}$ is a sequence of probability distributions on $\RR^d$. Let $\{\hat{\psi}_n\}_{n \ge 1}$ be a sequence of convex functions such that $\nabla \hat{\psi}_n\# \mu = \hat{\nu}_n$, for all $n \ge 1$. Fix $u_0\in \mathrm{Int}(\s)$ and $\delta_0\equiv\delta_0(u_0)>0$ such that $B_{\delta_0}(u_0) \subset \s$. Suppose that $\mu$ has a bounded (from below and above) nonvanishing density on $B_{\delta_0}(u_0)$. Suppose that $\psi$ is differentiable everywhere in $\mathrm{Int}(\s)$ and let $\nabla\psi$ be locally uniformly Lipschitz  in $B_{\delta_0}(u_0)$ with Lipschitz constant $K$. Define 
 \begin{align}\label{eq:DeltaDef}
\delta_n:= \Big(\int_{B_{\delta_0}(u_0)} \|\nabla \hat{\psi}_n(u)- \nabla\psi(u)\|^2 d\mu(u)\Big)^{\frac{1}{d+2}}. 
\end{align}
  Then, there exists $C=C(\mu, d, K) >1/2$ such that 
{\small \begin{align}\label{eq:SupBd}
\quad \;\; \sup_{u\in B_{\delta_0/3}(u_0)} \sup_{y \in \partial \hat{\psi}_n(u)}\|y- \nabla\psi(u) \|  \leq \begin{cases}
C\delta_n & \text{if }\delta_n \leq \delta_0/3, \\
C \delta_n^{d+2}\delta^{-(d+1)}_{0}  + \frac{1}{2}\delta_0& \text{if } \delta_n> \delta_0/3.
\end{cases}
\end{align}}
Now, assume that $\mu$ has bounded nonvanishing density everywhere on $\s$ and that $\nabla \psi$ is uniformly Lipschitz on $\mathcal{S}$ with Lipschitz constant $K$. Suppose that there exists $\tilde{\delta}_0>0$\footnote{The existence of $\tilde \delta_0$ is guaranteed if $\nabla \psi$ is a homeomorphism from $\mathrm{Int}(\s)$ to $\mathrm{Int}(\Y)$.} such that $\nabla\psi^{*}(\mathrm{Cl}(B_{\tilde{\delta}_0}(\nabla \psi(u_0))))\subset B_{\delta_{0}/6}(u_0)$. Define 
\begin{align*}
\delta^{\prime}_n:= \Big(\int_{\mathcal{S}}\|\nabla \hat{\psi}_n(u) - \nabla \psi(u)\|^2 d\mu(u)\Big)^{\frac{1}{d+2}},
\end{align*}
and
{\small $$ \tilde{\delta}_n: = \sup_{v \in \mathrm{Cl}(B_{\tilde{\delta}_0}(\nabla \psi(u_0))), \; w\in \nabla\psi(\mathrm{Cl}(B_{\delta_0/3}(u_0)))}\Big\{\|\nabla\psi^{*}(v)-\nabla\psi^{*}(w)\|: \|v-w\|\leq C\delta^{\prime}_n\Big\} $$}
where $C$ is the same constant as in \eqref{eq:SupBd}.
If $\max\{C \delta^{\prime}_n, \tilde{\delta}_n\}<\delta_0/6$, then, 
\begin{align}\label{eq:SupBd1}
\sup_{x\in B_{\tilde{\delta}_0}(\nabla \psi(u_0))} \sup_{w \in \partial \hat{\psi}^{*}_n(x)}\|w -\nabla\psi^{*}(x)\| \leq  
\tilde{\delta}_n.
\end{align}
\end{proposition}

In~\eqref{eq:SupBd} of Proposition~\ref{ppn:RateProp1} we provide a deterministic upper bound on the pointwise difference between $\partial \hat{\psi}_n(\cdot)$ and $\nabla\psi(\cdot)$, uniformly on the local ball $B_{\delta_0/3}(u_0)$, in terms of the $L_2$-loss~\eqref{eq:DeltaDef}. Similarly,~\eqref{eq:SupBd1} bounds the pointwise difference between $\partial \hat{\psi}_n^*(\cdot)$ and $\nabla\psi^*(\cdot)$, uniformly over a local ball around $\nabla\psi(u_0)$, in terms of $\tilde{\delta}_n$.

\begin{remark}[$\tilde \delta_n$ adapts to the local smoothness of $\nabla \psi$ and $\nabla \psi^{*}$] 
As $\tilde{\delta}_n$ is defined  implicitly in terms of $\delta_0,\tilde{\delta}_0$ and $\delta^{\prime}_n$ it is natural to ask how $\tilde{\delta}_n$ varies with $\delta^{\prime}_n$. If $\nabla \psi^{*}$ is uniformly $\beta$-H\"older continuous, for $\beta \in (0,1]$,  in the neighborhood $\mathrm{Cl}(B_{\tilde{\delta}_0}(\nabla \psi(u_0)))\cup \nabla\psi(\mathrm{Cl}(B_{\delta_0/3}(u_0)))$, then, for all   $v,w\in \mathrm{Cl}(B_{\tilde{\delta}_0}(\nabla \psi(u_0)))\cup \nabla\psi(\mathrm{Cl}(B_{\delta_0/3}(u_0)))$
\begin{align}
\|\nabla \psi^{*}(v)- \nabla \psi^{*}(w)\| \leq \kappa (\delta^{\prime}_n)^{\beta}, \quad \text{such that } \|v-w\|\leq C\delta^{\prime}_n, \nonumber
\end{align} 
where $\kappa = \kappa(C,\beta, \mu, \delta_0)>0$ is a constant. As a consequence, $\tilde{\delta}_n$ is also less than $\kappa (\delta^{\prime}_n)^{\beta}$. This shows that the local rate of convergence of the rank map adapts to the local smoothness of the transport maps $\nabla \psi$ and $\nabla \psi^{*}$.
\end{remark}

\subsection{Proof of Theorem~\ref{thm:RateTheo}}\label{pf:RateTheo}
Note that $R \equiv \nabla \psi^*$ and $Q \equiv \nabla \psi$. Let  $\hat{Q}_n \equiv \nabla \hat{\psi}_n$ where $\hat{\psi}_n$ is the convex potential function. As $\hat{Q}_n(\cdot)$ and $\hat{R}_n(\cdot)$ can be any point in the corresponding subdifferential sets when the underlying potential functions ($\hat{\psi}_n$ and $\hat{\psi}_n^*$) are not differentiable, we will upper bound the following two quantities: $$\mathbb{E}\Big[\sup_{u\in B_{\delta_0/3}(u_0)}\sup_{y \in \partial \hat{\psi}_n(u)}\|y- \nabla \psi(u)\|\Big]$$ and  $$\mathbb{E}\Big[\sup_{x\in B_{\delta_0/6}(\nabla u_0)}\sup_{u \in \partial \hat{\psi}^{*}_n(x)}\|u-\nabla \psi^{*}(x)\|\Big].$$
 Define 
\begin{align}\label{eq:DefT}
\mathcal{T}(\hat \nu_n, \nu) := \int_{\RR^d}\|\nabla\hat{\psi}_n(u) -\nabla \psi(u)\|^2 d\mu(u). 
\end{align}
Under the assumption that $\psi^{*}$ is strongly convex function with parameter $L>0$, Theorem~\ref{thm:Q-Rate} shows that, for some $C \equiv C(\mu,\nu) >0$,
\begin{align}\label{eq:MKduality}
 \mathbb{E}[\mathcal{T}(\hat \nu_n, \nu)]\leq Cr_{d,n}. 
\end{align}
As $\psi^{*}$ is strongly convex in $\mathcal{Y}$ and $\nabla\psi$ is a homeomorphism from $\mathrm{Int}(\mathcal{S})$ to $\mathrm{Int}(\mathcal{Y})$, $\nabla\psi$ is uniquely defined and is locally uniformly Lipschitz in $B_{\delta_0}(u_0)$, say with Lipschitz constant $K >0$. Recall the definition of $\delta_n$ from \eqref{eq:DeltaDef}. Thus,  $\delta_n^{d+2} \le \mathcal{T}(\hat \nu_n, \nu)$. Now, by Proposition~\ref{ppn:RateProp1} (in particular, by~\eqref{eq:SupBd}), we get 
\begin{align}\label{eq:DB}
\mathbb{E}\Big[\sup_{u\in B_{\delta_0/3}(u_0)}\sup_{y \in \partial \hat{\psi}_n(u)}\|y- \nabla \psi(u)\|\Big] \hspace{3in} &\nonumber \\
\leq C\mathbb{E}\Big[\delta_n\mathbbm{1}(\delta_n\leq 3^{-1}\delta_0)\Big] + C\delta^{-(d+1)}_0\mathbb{E}\big[\delta^{d+2}_n\mathbbm{1}(\delta_n> 3^{-1}\delta_0)\big] + \frac{\delta_0}{2}\mathbb{P}(\delta_n> 3^{-1}\delta_0) &\nonumber \\ 
\leq C\mathbb{E}\big[\mathcal{T}(\hat{\nu}_n,\nu)^{\frac{1}{d+2}}\big]+C\delta^{-(d+1)}_0\mathbb{E}\big[\mathcal{T}(\hat{\nu}_n,\nu)\big]+ \frac{\delta_0}{2}\mathbb{P}\Big(\mathcal{T}(\hat{\nu}_n, \nu)\geq (3^{-1}\delta_0)^{d+2}\Big) \hspace{.15in}&\nonumber \\
\leq \tilde K \Big(\mathbb{E}\big[\mathcal{T}(\hat{\nu}_n,\nu)^{\frac{1}{d+2}}\big]+\delta^{-(d+1)}_0\mathbb{E}\big[\mathcal{T}(\hat{\nu}_n,\nu)\big]\Big) \hspace{1in}&
\end{align}
where $C$ is the same constant as in~\eqref{eq:SupBd} and $\tilde K >0$ depends on $u_0,\mu$ and $d$. The last inequality follows by upper bounding $\mathbb{P}\big(\mathcal{T}(\hat{\nu}_n, \nu)\geq (3^{-1}\delta_0)^{d+2}\big)$ by $ 3^{d+2}\delta^{-(d+2)}_0\mathbb{E}[\mathcal{T}(\hat{\nu}_n, \nu)]$.
 Using Jensen's inequality, 
\begin{align}
\mathbb{E}\big[\mathcal{T}(\hat{\nu}_n, \nu)^{\frac{1}{d+2}}\big]\leq \Big(\mathbb{E}\big[\mathcal{T}(\hat{\nu}_n, \nu)\big]\Big)^{\frac{1}{d+2}}. 
\end{align}
Applying this in the first term of the right hand side of \eqref{eq:DB}, we arrive at 
\begin{align}
\mathbb{E}\Big[\sup_{u\in B_{\delta_0/3}(u_0)}\sup_{y \in \partial \hat{\psi}_n(u)}\|y- \nabla \psi(u)\|\Big]\leq K^{\prime} \big(\mathbb{E}\big[\mathcal{T}(\hat{\nu}_n, \nu)\big]\big)^{\frac{1}{d+2}}\label{eq:Coeff}
\end{align}
for some constant $K^{\prime}=K^{\prime}(u_0,\delta_0,\mu,\nu)>0$. Now, the first part of Theorem~\ref{thm:RateTheo} follows by combining \eqref{eq:Coeff} with \eqref{eq:MKduality} and observing that $Q \equiv \nabla \psi$.

Now, we turn to prove the second part of Theorem~\ref{thm:RateTheo}. For this, we intend to apply \eqref{eq:SupBd1} of Proposition~\ref{ppn:RateProp1}. Suppose that $\psi^{*}$ and $\psi$ are strongly convex with parameters $L>0$ and $L'>0$ respectively. There exists $\tilde{\delta}_0 >0$ such that $\nabla\psi^{*}(\mathrm{Cl}(B_{\tilde{\delta}_0}(\nabla \psi(u_0))))\subset B_{\delta_{0}/6}(u_0)$ since $\nabla \psi^{*}$ is a homeomorphism and, in fact, the inverse of $\nabla \psi$. By our assumption, $\nabla\psi$ and $\nabla \psi^{*}$ are uniformly Lipschitz continuous with parameters $1/L$ and $1/L'$ respectively. So, we have $\tilde \delta_n \leq C \delta'_n/L$. Also, note that $\delta_n' \equiv \mathcal{T}(\hat{\nu}_n,\nu)^{\frac{1}{d+2}}$. Let us denote $C':= C \max\{1,1/L\}$. Using \eqref{eq:SupBd1} in the same way as in \eqref{eq:DB}, we may write
\begin{eqnarray*}
\mathbb{E}\Big[\sup_{x\in B_{\delta_0/6}(\nabla u_0)}\sup_{u \in \partial \hat{\psi}^{*}_n(x)}\|u-\nabla \psi^{*}(x)\|\Big]  \hspace{3in} & \\ \leq \E \left[\tilde \delta_n \mathbbm{1}(\max\{C \delta^{\prime}_n, \tilde{\delta}_n\}<\delta_0/6)\right] + \Big(\|u_0\|+\frac{1}{6}\delta_0+\sup_{y\in \mathrm{supp}(\mu)}\|y\|\Big)\mathbb{P} \left(\max\{C \delta^{\prime}_n, \tilde{\delta}_n\}\ge \delta_0/6 \right)&\\
\leq K\mathbb{E}\Big[ \delta'_n\mathbbm{1}(C'\delta'_n< 6^{-1}\delta_0)\Big] + \Big(\|u_0\|+\frac{1}{6}\delta_0+\sup_{y\in \mathrm{supp}(\mu)}\|y\|\Big)\mathbb{P}(C' \delta'_n\geq 6^{-1}\delta_0) \hspace{1in} &\\
\leq K\Big(\mathbb{E}\big[\mathcal{T}(\hat{\nu}_n,\nu)^{\frac{1}{d+2}}\big]+\delta^{-(d+2)}_0\mathbb{E}\big[\mathcal{T}(\hat{\nu}_n,\nu)\big]\Big) \hspace{3in} & 
\label{eq:DB2}
\end{eqnarray*}
for some $K \equiv K(u_0,\mu,\nu)>0$. Now, the result follows from the above display by applying~\eqref{eq:MKduality} in a similar way as in~\eqref{eq:Coeff}. \qed

\subsection{Proof of Proposition~\ref{ppn:RateProp1}}\label{sec:RateProp1}

We first prove \eqref{eq:SupBd}. For any $u,v\in \mathcal{S}$, suppose $z_1\in \partial \hat{\psi}_{n}(u)$ and $z_2 \in \partial \hat{\psi}_{n}(v)$. Then, due to the convexity of $\hat{\psi}_n$, $\langle z_1-z_2, u-v\rangle \geq 0$, which implies 
\begin{align}\label{eq:Bound}
\qquad \langle z_1-\nabla \psi(u)-z_2+\nabla \psi(v), u-v\rangle \geq -\langle \nabla \psi(u) - \nabla \psi(v), u-v\rangle.
\end{align}
Let us define, for every $x \in \s$, 
 \begin{align}
 u(x) := \{z-\nabla \psi(x): z\in \partial \hat{\psi}_{n}(x)\}, \;\; \|u(x)\|_{\infty}:= \sup\{\|z-\nabla\psi(x)\|: z\in \partial \hat{\psi}_{n}(x)\}. \nonumber
 \end{align}
Fix $x = (x^{(1)},\ldots, x^{(d)}) \in B_{\delta_0/3}(u_0)$ and $y \in B_{\delta_0}(u_0)$. Then, using the $K$-Lipschitzness of $\psi$ on $B_{\delta_0}(u_0)$ and the Cauchy-Schwartz inequality, we have, $$\langle \nabla \psi(x) - \nabla \psi(y), x-y\rangle \le K\| x-y\|^2.$$ Using the above,  for some $z_x \in u(x)$ and $z_y \in u(y)$, we can rewrite~\eqref{eq:Bound} as 
 \begin{align}
 \langle z_x, y-x\rangle\leq \langle z_y, y-x\rangle + K\|y-x\|^2.\label{eq:ReIneq}
 \end{align}
Let us fix some $z_x = \big(z^{(1)}_x,\ldots ,z^{(d)}_x\big)\in u(x)$. We will bound $\|z_x\|$ from above. From~\eqref{eq:ReIneq}, we get 
 \begin{equation}
\quad\;  \langle z_x, y-x\rangle\leq \sup_{z_y \in u(y)}|\langle z_y, y-x\rangle| + K\|y-x\|^2 \le \frac{1}{2}\|u(y)\|_\infty^2 +(K+\frac{1}{2})\|y-x\|^2 \label{eq:MonotoneCons}
 \end{equation}
where the last inequality follows from the bound $|\langle z_y, y-x\rangle| \le \|z_y\|^2/2 + \|y-x\|^2/2$. Let $\theta_{x,\pm} := x\pm\delta e_1/3$ (here $e_1 =(1,0,\ldots, 0) \in \R^d$). Since $B_{\delta_0}(u_0)$ is inside the support of $\mu$ and $\mu$ has a nonvanishing bounded density in $B_{\delta_0}(u_0)$, for each $i=1,\ldots ,d$ and $\delta<\delta_0$, there exists an interval $\mathcal{B}^{(i)}_{x,\delta,\pm}$ around $\theta^{(i)}_{x,\pm}$ and $c=c(u_0,\delta_0,\mu)>0$ such that the following conditions hold\footnote{To see how we find such intervals, we fix the length of the interval $\mathcal{B}^{(i)}_{x,\delta,+}$ on the left hand side of $\theta^{(i)}_x$ to be $c\delta$ for some number $c\leq (3\sqrt{d+1})^{-1}\delta$. Then, there exists $\omega^{(i)}=\omega(c,i,u_0,\mu)>0$ such that \eqref{eq:ExpVanish} will be satisfied if we define $\mathcal{B}^{(i)}_{x,\delta,+}= (\theta^{(i)}_x-c\delta,\theta^{(i)}_x+\omega^{(i)}\delta)$. One can give an upper bound on $\omega^{(i)}$ which will depend on $c$ and the upper and lower bounds on the density of $\mu$ in $B_{\delta_0}(u_0)$. Denote this upper bound by  $\xi(c)$. Then, we will find the optimal choice of $c$ from the constraint $\max\{c,\xi(c)\}\leq (3\sqrt{d+1})^{-1} \delta$.}:
\begin{align}
\int_{y^{(i)}\in \mathcal{B}^{(i)}_{x,\delta,\pm}} (y^{(i)}-\theta^{(i)}_{x,\pm}) d\mu(y)&=0, \qquad \mbox{for} \; i = 1,\ldots, d,\label{eq:ExpVanish}\\ c\delta= \sup_{w\in \mathcal{B}^{(i)}_{x,\delta,\pm}}|w-\theta^{(i)}_{x,\pm}|& \leq \frac{1}{3\sqrt{d+1}}\delta. \label{eq:LBd} 
\end{align}
Define $\mathcal{B}_{x,\delta,\pm} := \mathcal{B}^{(1)}_{x,\delta,\pm}\times \ldots \times \mathcal{B}^{(d)}_{x,\delta,\pm}$, i.e., $\mathcal{B}_{x,\delta,\pm}$ is the Cartesian product of $\mathcal{B}^{(i)}_{x,\delta,\pm}$ for $i=1,\ldots, d$. Due to \eqref{eq:LBd}, $\mathcal{B}_{x,\delta,\pm} \subset B_{\delta/3}(\theta_{x,\pm})$. Observe that, for $y = (y^{(1)}, \ldots, y^{(d)}) \in \mathcal{B}_{x,\delta,+}$,
\begin{align}
  \langle z_x, y-x\rangle = & \sum_{i=1}^d z_x^{(i)} (y^{(i)} - x^{(i)}) = \sum_{i=1}^d z_x^{(i)} (y^{(i)} - \theta^{(i)}_{x,+}) + \sum_{i=1}^d z_x^{(i)} (\theta^{(i)}_{x,+} - x^{(i)}) \nonumber \\
= & \sum_{i=1}^d z_x^{(i)} (y^{(i)} - \theta^{(i)}_{x,+})  + z_x^{(1)} \frac{\delta}{3}.
  \end{align}
Note that, for $y \in \mathcal{B}_{x,\delta,+}$, we have   \begin{align*}
\|y - x\|^2 = (y^{(1)} - x^{(1)})^2 + \sum_{j=2}^d (y^{(j)} - x^{(j)})^2 & \le 2 (y^{(1)} - \theta^{(1)}_{x,+})^2 + 2 (\theta^{(j)}_{x,+} - x^{(1)})^2  \\
\hspace{1.5in} + \frac{\delta^2(d-1)}{9(d+1)} \;\;\;\; & \le \frac{\delta^2}{9} + 2 \frac{\delta^2}{9} = \frac{\delta^2}{3}.
  \end{align*} 
  as $x_j = \theta^{(j)}_{x,+}$ for $j = 2,\ldots, d$, and we have used~\eqref{eq:LBd}.
Therefore, integrating both sides of~\eqref{eq:MonotoneCons} w.r.t.~$\mu$ as $y$ varies over $\mathcal{B}_{x,\delta,+}$ yields 
  \begin{align}
\qquad \; \frac{\delta}{3} z^{(1)}_x \mu(\mathcal{B}_{x,\delta,+})& \leq \frac{1}{2}  \int_{B_{\delta_0}(u_0)}\|u(y)\|^2_\infty d\mu(y) + \frac{(2K+1)}{6}\delta^2\mu(\mathcal{B}_{x,\delta,+})\label{eq:UpIneq}
  \end{align}
where the left side is a consequence of~\eqref{eq:ExpVanish}, and the inequality follows from combining~\eqref{eq:MonotoneCons} with the fact $\mathcal{B}_{x,\delta,+}\subset B_{\delta/3}(\theta_{x,+})\subset B_{\delta_0}(u_0)$. 

Similarly integrating  \eqref{eq:MonotoneCons} w.r.t.~$\mu$ over $\mathcal{B}_{x,\delta,-}$, we get 
 \begin{align}
\qquad \;\; -\frac{\delta}{3} z^{(1)}_x \mu(\mathcal{B}_{x,\delta,-})\leq \frac{1}{2} \int_{B_{\delta_0}(u_0)}\|u(y)\|^2_\infty d\mu(y) + \frac{(2K+1)}{6}\delta^2 \mu(\mathcal{B}_{x,\delta,-}).\label{eq:DownIneq}
 \end{align}
 Since $\mu$ has a nonvanishing bounded density (w.r.t.~Lebesgue measure) in $B_{\delta_0}(u_0)$, $\mathcal{B}_{x,\delta,\pm}$ belongs to the support of $\mu$ for all $\delta<\delta_0$ and by \eqref{eq:LBd}, we know
\begin{align}
(c \delta)^d\leq \lambda_d(\mathcal{B}_{x,\delta,\pm}) \le  \frac{1}{3^d(d+1)^{d/2}} \delta^d.
\end{align} 
Therefore, there exists $0<c_2 \equiv c_2(u_0,\delta_0,\mu)<c_1 \equiv c_1(u_0, \delta_0,\mu)$ such that  
\begin{align}\label{eq:MeasureUpDown}
c_2\delta^{d}\leq \mu\big(\mathcal{B}_{x,\delta,\pm}\big)\leq c_1\delta^{d}, \quad \forall \; x\in B_{\delta_0/3}(u_0).
\end{align}
 Combining \eqref{eq:UpIneq} with \eqref{eq:DownIneq} and \eqref{eq:MeasureUpDown} yields
 \begin{align}
 |z^{(1)}_x|\leq \frac{3}{2c_2\delta^{d+1}}\int_{B_{\delta_0}(u_0)}\|u(y)\|^2_\infty d\mu(y) + \frac{(2K+1)}{2}\delta.
 \end{align}
Note that $\int_{B_{\delta_0}(u_0)}\|u(y)\|^2_\infty d\mu(y) = \int_{B_{\delta_0}(u_0)} \|\partial \hat{\psi}_n(y) - \nabla \psi(y)\|^2  d\mu(y)$. Optimizing the right side of the above display w.r.t.~$\delta$, we see
 \begin{align}\label{eq:1stBd}
\;\; \qquad |z^{(1)}_x| \leq \begin{cases} 
C_0\Big(\int_{B_{\delta_0}(u_0)} \|\partial \hat{\psi}_n(y) - \nabla \psi(y)\|^2  d\mu(y)\Big)^{\frac{1}{d+2}}=C\delta_n & \text{if } \delta_n\leq \delta_0/3\\
 \frac{C_0}{\delta^{d+1}_0}  \delta_n^{d+2} + \frac{(2K+1)\delta_0}{6} & \text{if }\delta_n > \delta_0/3 
\end{cases} 
\end{align}
for some constant $C_0= C_0(\mu,d, K)>0$; we further assume that $C_0 d > 1/2$. In a similar way, $|z^{(i)}_x|$ can be bounded above by the right side of \eqref{eq:1stBd}, for $i=2,\ldots,d$. As a consequence
\begin{align}\label{eq:0thBd}
\;\;\qquad \|z_x\| \leq d\times \begin{cases} 
C_0\delta_n & \text{if } \delta_n\leq \delta_0/3\\
 \frac{C_0}{\delta^{d+1}_0} \delta_n^{d+2} + \frac{(2K+1)\delta_0}{6} & \text{if }\delta_n > \delta_0/3. 
\end{cases} 
\end{align}
 Note that the right side of \eqref{eq:0thBd} does not depend on $z_x$ or $x$, for all $x\in B_{\delta_{0}/3}(u_0)$. Therefore, $\sup_{x\in B_{\delta_0/3}(u_0)} \sup_{z\in u(x)} \|z\|$ is bounded above by the right side of the above display.
Letting $C= C_0d$ and substituting into \eqref{eq:0thBd} proves \eqref{eq:SupBd}.

Now, we turn to verify the bound in \eqref{eq:SupBd1} when $\max\{C\delta^{\prime}_n,\tilde{\delta}_n\}<\delta_0/6$. We prove this by contradiction. Suppose that there exists $x\in B_{\tilde{\delta}_0}(\nabla\psi(u_0))$ and $u \in \partial \hat{\psi}^{*}_n(x)$ such that $\|u-\nabla\psi^{*}(x)\|>\tilde{\delta}_n$. We first show the contradiction under the assumption $\|u-\nabla\psi^{*}(x)\|\leq \delta_0/6$. Since we know that $\nabla\psi^{*}(\mathrm{Cl}(B_{\tilde{\delta}_0}(\nabla\psi(u_0))))\subset B_{\delta_0/6}(u_0)$ and $x\in B_{\tilde{\delta}_0}(\nabla\psi(u_0))$, therefore, $\|\nabla\psi^{*}(x) - u_0\|\leq \delta_0/6$. Combining this fact with the assumption $\|u-\nabla\psi^{*}(x)\|\leq \delta_0/6$  and using the subadditivity of the Euclidean norm yields $\|u-u_0\|\leq \delta_0/3$. However, we know that $x\in \partial \hat{\psi}_n(u)$ by Lemma~\ref{lem:SubD}, and by \eqref{eq:SupBd}, $\|x-\nabla\psi(u)\|\leq C\delta^{\prime}_n$. Further, since $x\in  B_{\tilde{\delta}_0}(\nabla\psi(u_0))$ and $\nabla\psi(u) \in \nabla\psi(\mathrm{Cl}(B_{\delta_0/3}(u_0)))$, we get $\|\nabla\psi^{*}(x)-u\|\leq \tilde{\delta}_n$ by using the definition of $\tilde{\delta}_n$. This contradicts $\|u-\nabla\psi^{*}(x)\|>\tilde{\delta}_n$. 

Thus, we may now assume $\|u-\nabla\psi^{*}(x)\|> \delta_0/6$. We first suppose that $\nabla\psi(y)=y$ for all $y \in \s$ (which holds when $\mu=\nu$). As a consequence, we get $\nabla \psi^{*}(x)=x$. Fix $v$ on the line joining $u$ and $x$ such that $\|v-x\|=\delta_0/6$. Note that $v\in B_{\delta_0/3}(u_0)$ because $\nabla \psi^{*}(x)=x \in B_{\delta_0/6}(u_0)$. Fix $z\in \partial \hat{\psi}_n(v)$. Since $x\in \partial \hat{\psi}_n(u)$, we have $\langle x-z, u-v\rangle \geq 0$.  Now, we note 
\begin{align}
\langle x-z, u-v\rangle &= \langle x-v, u-v\rangle+ \langle v-z, u-v\rangle \\& \leq -\frac{\delta_0}{6}\|u-v\| + \|v-z\|\|u-v\|\\
 &\leq \|u-v\|\big(-\frac{\delta_0}{6}+\sup_{y\in B_{\delta_0/3}(u_0)} \sup_{w\in \partial \hat{\psi}_n(y)} \|w-y\|\big)\\
 &\leq \|u-v\|(-\frac{\delta_0}{6}+C\delta^{\prime}_n) < 0\label{eq:SeqIneq}
\end{align}
where the first inequality follows since $\langle x-v, u-v\rangle = -\|u-v\|\delta_0/6 $ as $x,v,u$ are collinear. The second inequality is obtained by recalling that $z\in \partial \hat{\psi}_n(v)$ for $v\in B_{\delta_0/3}(u_0)$ and the last inequality holds because $C\delta^{\prime}_n$ is assumed to be less than $\delta_0/6$. Note that \eqref{eq:SeqIneq} contradicts $\langle x-z, u-v\rangle \geq 0$. Hence, the result follows when $\nabla \psi(y)=y$ for all $y\in \s$.

Thus, we may assume $\|u-\nabla\psi^{*}(x)\|> \delta_0/6$ and $\nabla \psi^{*}(x)\neq x$. Since $\nabla \psi$ is uniformly Lipschitz continuous in $\s$ with the Lipschitz constant $K$ and $\mu$ has uniformly bounded nonvanishing density everywhere $\s$, by repeating the same argument as in the proof of \eqref{eq:SupBd} yields 
\begin{align}
\sup_{y\in \partial\hat{\psi}_n(u)}\|y-\nabla \psi(u)\|\leq \begin{cases}
C\delta^{\prime}_n & \text{if }\delta^{\prime}_n\leq \delta_0/3,\\
\frac{3C}{\delta^{(d+1)}_0}(\delta^{\prime}_n)^{d+2} + \frac{d(2K+1)}{2}\delta_0, & \text{if }\delta^{\prime}_n>\delta_0/3.
\end{cases} 
\end{align}
Since $C>\frac{1}{2}$, under the assumption $C\delta^{\prime}_n \leq \frac{\delta_0}{6}$, we have 
\begin{align}\label{eq:Conseq}
\sup_{y \in \partial\hat{\psi}_n(u)}\|y - \nabla \psi(u)\|\leq C\delta^{\prime}_n
\end{align} 
Due to the fact that $u \in \partial \hat{\psi}^{*}_n(x)$ and Lemma~\ref{lem:SubD}, we get $x\in \partial \hat{\psi}_n(u)$. Combining this with \eqref{eq:Conseq} shows $\|x-\nabla\psi(u)\|\leq C\delta^{\prime}_n$. From the definition of $\tilde{\delta}_n$, we get $\|\nabla \psi^{*}(x)- u\|\leq \tilde{\delta}_n$. According to the condition of Theorem~\ref{ppn:RateProp1}, we have $\tilde{\delta}_n\leq {\delta_0}/{6}$. This contradicts our assumption $\|\nabla \psi^{*}(x)-u\|>{\delta_0}/{6}$ and hence, completes the proof.   \qed

\section{Proofs of the results in Section~6}\label{pf:Goodness-Fit-Test}
\subsection{Proof of Lemma~\ref{lem:Distfree}}\label{pf:Distfree}

In this proof $U$ will denote an independent random vector with distribution $\mu$ with support $\s$ (where $\s \subset \R^d$ is a convex compact set) and independent of the observed data. Let $\hat Q_{X,Y}$ be the sample quantile map obtained from the pooled samples $X_1,\ldots, X_m, Y_1,\ldots, Y_n$. Recall that $\hat Q_{X,Y}$ induces a cell decomposition of $\s$ (see e.g.,~\eqref{eq:W_i}) into $m+n$ sets $\{\mathcal{C}^{X,Y}_{X,i}\}_{i=1}^{n}\bigcup \{\mathcal{C}^{X,Y}_{Y,j}\}_{j=1}^{m}$  such that $\hat Q_{X,Y} (u)= X_i$ for any $u\in \mathcal{C}^{X,Y}_{X,i}$ and $\hat Q_{X,Y}(u) = Y_j$ for any $u\in \mathcal{C}^{X,Y}_{Y,j}$. Further, for $\mu$-a.e.~$u$, $\hat Q_{X,Y} (u) \in \{X_1,\ldots, X_m, Y_1,\ldots, Y_n\}$. Denote $\mathbf{X}:= \{X_1, \ldots , X_m\}$ and $\mathbf{Y}:= \{Y_1, \ldots , Y_n\}$. Now, for any Borel set $B \subset \s$,
 \begin{align}
  \mathbb{P}&\Big(\hat R_{X,Y}(\hat Q_{X}(U))\in B \Big|\mathbf{X}\cup \mathbf{Y} \Big) \nonumber \\
  & =  \sum_{i=1}^{m}\mathbb{P}(\hat Q_{X}(U)=X_i|\mathbf{X})\times \mathbb{P} \left(\hat R_{X,Y}(X_i)\in B|\mathbf{X}\cup \mathbf{Y}\right) \nonumber \\
  & = \sum_{i=1}^{m} \frac{1}{m} \left[ (n+m) \, \mu\big(\mathcal{C}^{X,Y}_{X,i}\cap B\big) \right] \nonumber\\
 &   = \frac{m+n}{m} \, \mu\left(\Big\{\bigcup_{i=1}^{m}\mathcal{C}^{X,Y}_{X,i}\Big\}\cap B \right) \nonumber 
 \end{align}
 where the second equality follows by noting that $\mathbb{P}(\hat Q_{X}(U)= X_i|\{X_1, \ldots , X_m\}) =\frac{1}{m}$ and 
 \begin{align*}
 \mathbb{P}\Big(\hat R_{X,Y}(X_i)\in B | \mathbf{X}\cup \mathbf{Y}\Big) = (n+m) \, \mu\big(\mathcal{C}^{X,Y}_{X,i}\cap B\big) 
\end{align*} 
as $\hat R_{X,Y}(X_i)$ is a randomly chosen point supported on the cell $\mathcal{C}^{X,Y}_{X,i}$ and $\mu(\mathcal{C}^{X,Y}_{X,i}) = (m+n)^{-1}$; see~\eqref{eq:Rank-X_i}. Thus, we see that 
{\small \begin{align}\label{eq:ProbB}
\qquad \mathbb{P}\Big(\hat R_{X,Y}(\hat Q_{X}(U))\in B\Big) = \frac{n+m}{m} \;\mathbb{E}_{X,Y}\left[ \mu \Big(B \cap \Big\{\bigcup_{i=1}^{m} \mathcal{C}^{X,Y}_{X,i}\Big\}\Big)\right],
\end{align}}
where $\mathbb{E}_{X,Y}$ denotes the expectation w.r.t.~$\{X_1, \ldots , X_m\}\cup \{Y_1, \ldots , Y_n\}$.

When the law of $X_i$'s are the same as that of $Y_j$'s, then, $\{X_1, \ldots , X_m\}$ could be any random permutation of the joint data set $\{X_1, \ldots , X_m, Y_1, \ldots , Y_n\}$. Let us denote $Z_i:= X_i$ for $i=1, \ldots ,m$, and $Z_{m+j}:= Y_{j}$ for $j=1, \ldots , n$. Let $\mathfrak{S}_{m+n}$ be the set of all permutations of $\{1, \ldots ,m+n\}$. If $\sigma$ is a random permutation uniformly chosen from  $\mathfrak{S}_{m+n}$, then, $(Z_1, \ldots , Z_m) \stackrel{d}{=} (Z_{\sigma(1)}, \ldots , Z_{\sigma(m)})$. Owing to this, 
  \begin{align}
\qquad  \quad\mu \Big(B\cap \Big\{\bigcup_{i=1}^{m}\mathcal{C}^{X,Y}_{X,i}\Big\}\Big) &= \mu \Big(B\cap \Big\{\bigcup_{i=1}^{m}\mathcal{C}^{Z}_{i}\Big\}\Big) \stackrel{d}{=} \mu \Big(B\cap \Big\{\bigcup_{i=1}^{m}\mathcal{C}^{Z}_{\sigma(i)}\Big\}\Big).  \label{eq:Disteq}
\end{align}      
Employing~\eqref{eq:Disteq}, we notice that 
\begin{align}
 \mathbb{E}_{X,Y}\Big[\mu \Big(B\cap \Big\{\bigcup_{i=1}^{m}\mathcal{C}^{X,Y}_{X,i}\Big\}\Big)\Big] &= \mathbb{E}_{Z,\sigma}\Big[\mu \Big(B\cap \Big\{\bigcup_{i=1}^{m}\mathcal{C}^{Z}_{\sigma(i)}\Big\}\Big)\Big] = \mathbb{E}_{Z}\Big[\mathbb{E}_{\sigma}\Big[\mu \Big(B\cap \Big\{\bigcup_{i=1}^{m}\mathcal{C}^{Z}_{\sigma(i)}\Big\}\Big)\Big]\Big] \nonumber \\
 &=\mathbb{E}_Z\Big[\frac{\binom{n+m-1}{m-1}}{\binom{n+m}{m}} \sum_{\ell=1}^{n+m} \mu \Big(B\cap  \mathcal{C}^{Z}_{\ell}\Big)\Big] = \frac{m}{m+n} \mu(B).\label{eq:ExpSame} 
\end{align} 
To see the third equality in \eqref{eq:ExpSame}, owing to the independence of $Z$ and $\sigma$, we note
\begin{align}\label{eq:SumExpan}
\mathbb{E}_{\sigma}&\Big[\mu \Big(B\cap \Big\{\bigcup_{i=1}^{m}\mathcal{C}^{Z}_{\sigma(i)}\Big\}\Big)\Big] = \frac{1}{(n+m)!} \sum_{\tau\in \mathfrak{S}_{n+m}} \Big[\sum_{i=1}^{m} \mu(B\cap \mathcal{C}^{Z}_{\tau(i)}) \Big] \nonumber\\
& \;\;\;\;\;\; = \frac{1}{(n+m)!} \sum_{j=1}^{n+m}\#\Big\{\tau\in \mathfrak{S}_{n+m}:j\in \big\{\tau(1), \ldots , \tau(m)\big\}\Big\} \times \mu \big(B\cap \mathcal{C}^{Z}_{j}\big). 
\end{align}   
 For any $j\in \{1, \ldots , m+n\}$, the total number of permutations $\tau \in \mathfrak{S}_{m+n}$ in which $j\in \{\tau(1), \ldots , \tau(m)\}$ is $\binom{m+n-1}{m-1}m!n!$.  Plugging this into the right hand side of \eqref{eq:SumExpan} and taking the sum inside the integral, we get~\eqref{eq:ExpSame}.

Plugging \eqref{eq:ExpSame} into \eqref{eq:ProbB}, we notice that $\mathbb{P}(\hat R_{X,Y}(\hat Q_{X}(U))\in B)$ does not depend on the distribution of $\{X_1, \ldots, X_m\}$ and $\{Y_1, \ldots , Y_n\}$. This shows the distribution-freeness of $\mathbb{P}(\hat R_{X,Y}(\hat Q_{X}(U))\in B)$. In fact, we have shown that $\hat R_{X,Y}(\hat Q_{X}(U))$ has distribution $\mu$. Using a similar argument, one can show that $\hat R_{X,Y}(\hat Q_{Y}(U)) \sim \mu$ and thus is also distribution-free.   \qed

\subsection{Proof of Proposition~\ref{lem:Power1}}\label{pf:Power1}
Assume that $\nu_X = \nu_Y = \nu$ which is absolutely continuous. Let us define $\hat{\mathcal{F}}:\s\to \R$ as $\hat{\mathcal{F}}(u) := \|\hat R_{X,Y}(\hat Q_{X}(u))- \hat R_{X,Y}(\hat Q_{Y}(u))\|^2$. Since $\hat{R}_{X,Y}$ maps $\RR^d$ to $\s$, therefore, for all $u\in \s$, 
\begin{align}\label{eq:OutBd}
\hat{\mathcal{F}}(u)\leq \mathrm{diam}(\s). 
\end{align}
Fix $\epsilon >0$. Choose a compact set $K\subset \mathrm{Int}(\s)$ such that $\int \mathbbm{1}(u\in \s \backslash K) d\mu(u)\leq \epsilon$. Since $Q$ is a homeomorphism from $\mathrm{Int}(\s)$ to $\mathrm{Int}(\Y)$, so, $J:=Q(K)$ is compact subset of $\mathrm{Int}(\Y)$. Let $\delta_0>0$ be such that $\mathrm{Cl}(J+B_{\delta_0}(0))\subset\mathrm{Int}(\Y)$.

Owing to Theorem~\ref{thm:GCProp}, for any two compact sets $\mathfrak{K}\subset \mathcal{S}$ and $\mathfrak{K}^{\prime}\subset \mathcal{Y}$, we have  
\begin{align}
\sup_{u\in \mathfrak{K}}\max\{\|\hat Q_{X}(u)- Q(u)\|,\|\hat Q_{Y}(u)- Q(u)\|\} & \stackrel{a.s.}{\longrightarrow} 0, \nonumber \\
\sup_{u\in \mathfrak{K}^{\prime}}\|\hat R_{X,Y}(u)- R(u)\| & \stackrel{a.s.}{\longrightarrow} 0 \label{eq:GlivCantelli}
\end{align}  
as $m, n\to \infty$, where $Q$ and $R$ are the quantile and rank maps for $\nu$ w.r.t.~$\mu$ (i.e., they are gradients of convex functions such that $Q\# \mu = \nu$ and $R\#\nu= \mu$). Due to \eqref{eq:GlivCantelli}, w.p.~$1$, $\hat Q_{X}(K), \hat Q_{Y}(K)$ will be contained in $\mathrm{Cl}(J+B_{\delta}(0))$ as $m,n\to \infty$ for any $\delta<\delta_0$. Combining this with the continuity of the map $R(\cdot)$ in $\mathrm{Int}(\Y)$ yields 
{\begin{align}\label{eq:CombineCon}
\qquad \sup_{u\in K}\max\{\|\hat R_{X,Y}(\hat Q_{X}(u))- R(Q(u))\|,\|\hat R_{X,Y}(\hat Q_{Y}(u))-R(Q(u))\|\}\stackrel{a.s.}{\longrightarrow} 0.
\end{align}
This implies, by the dominated convergence theorem,  
\begin{align}\label{eq:insIntCon}
\int_{K}\hat{\mathcal{F}}(u) d\mu(u) \stackrel{a.s.}{\longrightarrow} 0, \quad \text{as }m,n\to \infty.
\end{align}
 Then, 
\begin{eqnarray*}
\limsup_{m,n\to \infty}T_{X, Y} & \stackrel{a.s.}{\le} & \limsup_{m,n\to \infty} \int_{\mathcal{S}\backslash K} \hat{\mathcal{F}}(u) d \mu(u) \le \epsilon \,\mathrm{diam}(\s)
\end{eqnarray*}
where the first inequality follows from \eqref{eq:insIntCon} and the second inequality follows from \eqref{eq:OutBd}. Letting $\epsilon\to 0$ completes the proof of the first part of Proposition~\ref{lem:Power1}.

Now, we turn to the next part of Proposition~\ref{lem:Power1}. Let $Q_{X}, Q_{Y}$ and $R_{X,Y}$ be the gradients of convex functions such that $Q_{X}\# \mu = \nu_X$, $Q_{Y}\#\mu =\nu_Y$ and $R_{X,Y}\#\big(\theta \nu_X +(1-\theta)\nu_Y\big)=\mu$. In a similar way as in \eqref{eq:CombineCon}, for any compact set $\mathcal{K}\subset \mathrm{Int}(\s)$,
\begin{align}\label{eq:SepASConv}
\sup_{u\in \mathcal{K}}\max\{\|\hat R_{X,Y}(\hat Q_{X}(u))- R_{X,Y}(Q_{X}(u))\|,\|\hat R_{X,Y}(\hat Q_{Y}(u))- R_{X,Y}(Q_{Y}(u))\|\}\stackrel{a.s.}{\to} 0, \qquad
\end{align}
as $m,n \to \infty$. Fix $\epsilon>0$. Recall that $K\subset \mathrm{Int}(\s)$ is a compact set satisfying $\int \mathbbm{1}(u\in \s \backslash K) d\mu(u)\leq \epsilon$. Also recall that $\hat{\mathcal{F}}(u) = \|\hat R_{X,Y}(\hat Q_{X}(u))- \hat R_{X,Y}(\hat Q_{Y}(u))\|^2$ and define $\mathcal{F}:\s \to \RR$ as $\mathcal{F}(u) := \|R_{X,Y}(Q_{X}(u))- R_{X,Y}(Q_{Y}(u))\|^2$. Then, using \eqref{eq:SepASConv}, we get 
\begin{align}
\int_{K} \hat{\mathcal{F}}(u) d \mu(u) \stackrel{a.s.}{\to} \int_{K}{\mathcal{F}}(u) d \mu(u), \qquad \mbox{as } n,m\to \infty.\label{eq:limit}
\end{align}
   Owing to \eqref{eq:OutBd} and \eqref{eq:limit}, as $n,m\to\infty$, 
\begin{align}
 \limsup_{n\to \infty}T_{X, Y} &\stackrel{a.s.}{\leq} \epsilon \mathrm{diam}(\s) + \int_{K} \mathcal{F}(u)d \mu(u) ,\label{eq:conv1}\\
 \liminf_{n\to \infty} T_{X, Y} &\stackrel{a.s.}{\geq}  \int_{K} \mathcal{F}(u)d \mu(u).\label{eq:conv2}
\end{align} 
 Letting $\epsilon\to 0$ in \eqref{eq:conv1} and combining it with \eqref{eq:conv2} yields
{\small \begin{align}
\qquad  T_{X, Y} \stackrel{a.s.}{\longrightarrow} \int_{\s} \mathcal{F}(u) d \mu(u) =\int_{\s}\|R_{X,Y}(Q_{X}(u))-R_{X,Y}(Q_{Y}(u))\|^2 d\mu(u). \label{eq:limit2}
 \end{align}}
  Since $Q_{Y}\neq Q_{X}$ and both are continuous functions, there exists an open set $U\subset \mathrm{Int}(\s)$ such that $Q_{Y}(U)\cap  Q_{X}(U)=\emptyset$. We have assumed that $R_{X,Y}:\mathrm{Int}(\Y_X \cup \Y_Y)\to \mathrm{Int}(\s)$ is a homeomorphism. Hence, the right side of~\eqref{eq:limit2} is lower bounded by 
  \begin{align}
 \int_{U}\|R_{X,Y}(Q_{X}(u))-R_{X,Y}(Q_{Y}(u))\|^2 du >0,
\end{align}   
which implies the desired result.   \qed

\subsection{Proof of Proposition~\ref{thm:Two-S-Rate}}\label{pf:Two-S-Rate}
Let $R_{X,Y}$ be the rank map of $\nu$. Let
\begin{align}
A_{m,n} & := \int_{\mathcal{S}}  \|\hat{R}_{X,Y}(\hat{Q}_{X}(u))- R_{X,Y}(\hat{Q}_{X}(u))\|^2 \, d\mu(u), \\
B_{m,n} & := \int_{\mathcal{S}}  \|R_{X,Y}(\hat{Q}_{X}(u))- R_{X,Y}(\hat{Q}_{Y}(u))\|^2 \, d \mu(u), \quad \mbox{and}\\
C_{m,n} & := \int_{\mathcal{S}}  \|R_{X,Y}(\hat{Q}_{Y}(u))- \hat{R}_{X,Y}(\hat{Q}_{Y}(u))\|^2 \, d \mu(u).
\end{align}
By using the fact that $(a+ b + c)^2 \le 4 (a^2 + b^2 + c^2)$, we may write 
\begin{align}\label{eq:BreakItDown}
\mathbb{E}[T_{X,Y}] &\leq 4 \, \E \big[A_{m,n} + B_{m,n} + C_{m,n} \big].
\end{align}
In what follows, we will show that
\begin{equation} \label{eq:Bound1}
\E[A_{m,n} + C_{m,n}]  \leq C \, r_{d,m+n}, 
\end{equation}
and 
\begin{equation}\label{eq:Bound2}
\E[B_{m,n}]  \leq C(r_{d,m}+ r_{d,n}) 
\end{equation}
for some constant $C\equiv C(\mu,\nu, \theta)>0$. Owing to \eqref{eq:BreakItDown}, combining \eqref{eq:Bound1} and \eqref{eq:Bound2} would then complete the proof. 

We first show \eqref{eq:Bound1}. Note that $\hat{Q}_{X}(u)= X_i$ if $u\in \mathrm{Int}\big(W_i^X(\hat{h}_{X})\big)$ and $\hat{Q}_{Y}(u)= Y_j$ if $u\in \mathrm{Int}\big(W_j^Y(\hat{h}_{Y})\big)$, where $\{W_i^X\}_{i=1}^m$ and $\{W_j^Y\}_{j=1}^n$ are the cell-decompositions of $\s$ induced by the empirical distributions of $X_i$'s and $Y_j$'s respectively. Further, recall that $\int_{W_i^X(\hat{h}_{X})} d \mu  = m^{-1}$ and $\int_{W_j^Y(\hat{h}_{Y})} d \mu  =n^{-1}$, for all $i=1,\ldots, m$ and $j=1,\ldots, n$. Thus, 
\begin{align}
A_{m,n} & = \sum_{i=1}^{m}\int_{W_i^X(\hat{h}_X)} \|\hat{R}_{X,Y}(\hat{Q}_{X}(u))- R_{X,Y}(\hat{Q}_{X}(u))\|^2 d \mu(u) \\ & = \frac{1}{m} \sum_{i=1}^{m}\|\hat{R}_{X,Y}(X_i)- R_{X,Y}(X_i)\|^2.\label{eq:AmnBd}
\end{align}
Using a similar reduction for $C_{m,n}$ shows that $A_{m,n} + C_{m,n}$ can be upper bounded by
\begin{align*}
& \frac{1}{m} \sum_{i=1}^{m}\|\hat{R}_{X,Y}(X_i)- R_{X,Y}(X_i)\|^2 + \frac{1}{n}\sum_{j=1}^{n}\|\hat{R}_{X,Y}(Y_i)- R_{X,Y}(Y_i)\|^2 \\
\le & \;\;  \frac{\theta^{-1}}{m + n} \left[\sum_{i=1}^{m} \|\hat{R}_{X,Y}(X_i)- R_{X,Y}(X_i)\|^2 +\sum_{j=1}^{n}\|\hat{R}_{X,Y}(Y_i)- R_{X,Y}(Y_i)\|^2\right].
\end{align*}
Note that the right side of the above display is just (a constant times) the $L^2$-deviation of the empirical rank map $\hat{R}_{X,Y}$ (based on the pooled sample $\{X_1,\ldots, X_m, Y_1,\ldots, Y_n\}$) from its population version. Thus, taking expectation  and applying Theorem~\ref{thm:RateProp1} yields~\eqref{eq:Bound1}.

Now we turn to show \eqref{eq:Bound2}. For $u \in \s\vspace{-0.15in}$,
\begin{align}
\|R_{X,Y}(\hat{Q}_{X}(u))- R_{X,Y}(\hat{Q}_{Y}(u))\|^2 &\leq 2\big(\|R_{X,Y}(\hat{Q}_{X}(u)) - R_{X,Y}(Q_{X,Y}(u))\|^2\\
& + \|R_{X,Y}(\hat{Q}_{Y}(u)) - R_{X,Y}(Q_{X,Y}(u))\|^2\big).  \label{eq:Rbound}
\end{align} 
By our assumption, the convex potential $\psi_{X,Y}$ is strongly convex with parameter $L >0$. Thus, from convex analysis, $R_{X,Y}$ is Lipschitz continuous with parameter $\frac{1}{L}$. Owing to this, we may write 
\begin{align}
\|R_{X,Y}(\hat{Q}_{X}(u)) - R_{X,Y}(Q_{X,Y}(u))\|& \leq \frac{1}{L}\|\hat{Q}_{X}(u)- Q_{X,Y}(u)\|, \\  \|R_{X,Y}(\hat{Q}_{Y}(u))- R_{X,Y}(Q_{X,Y}(u))\| &\leq \frac{1}{L}\|\hat{Q}_{Y}(u)- Q_{X,Y}(u)\|. 
\end{align}
Applying these inequalities into the right hand side of \eqref{eq:Rbound}, we notice that 
$$B_n \leq \frac{2}{L} \left[\int_{\mathcal{S}} \|\hat{Q}_{X}(u)- Q_{X,Y}(u)\|^2 \, d \mu(u)+ \int_{\mathcal{S}}  \|\hat{Q}_{Y}(u)- Q_{X,Y}(u)\|^2 d \mu(u)\right]. $$ Taking expectation on both sides of the above display and applying Theorem~\ref{thm:Q-Rate} shows that $\mathbb{E}[B_n]$ is bounded above by $C (r_{d,n}+r_{d,m})$, for some $C >0$, depending on $\mu$ and $\nu$. This completes showing \eqref{eq:Bound2}.
\qed

\subsection{Proof of Proposition~\ref{thm:Two-S-Alt}}\label{pf:Two-S-Alt}
Let us define $\mathcal{R}_{m,n}:\mathcal{S}\to \mathcal{S}$ as 
\begin{equation}
\mathcal{R}_{m,n}(u) := \hat{R}_{X,Y}(\hat{Q}_{X}(u))- \hat{R}_{X,Y}(\hat{Q}_{Y}(u)) -R_{X,Y}(Q_X(u))+ R_{X,Y}(Q_Y(u)).
\end{equation}
We will first show that   
\begin{equation}\label{eq:TxyBound} 
\mathbb{E}[|T_{X,Y}-c|] \leq \mathbb{E}\Big[\int \|\mathcal{R}_{m,n} \|^2 d\mu \Big]+ \mathrm{diam}(\mathcal{S})  \Big(\mathbb{E}\Big[\int  \|\mathcal{R}_{m,n} \|^2 d\mu \Big]\Big)^{1/2}. 
\end{equation}
Next we will show that there exists constant $C \equiv C(\mu,\nu,\theta)>0$, depending only on $\mu, \nu$, and $\theta$ such that
\begin{equation}\label{eq:TxyBoundFinal}
\mathbb{E}\Big[\int_\s \|\mathcal{R}_{m,n}\|^2 d\mu \Big] \leq C r_{d,N}.
\end{equation}
Note that we get \eqref{eq:TBound} by combining \eqref{eq:TxyBound} and \eqref{eq:TxyBoundFinal}. 

We may write 
\begin{align}
T_{X,Y}&= \int \|\mathcal{R}_{m,n}(u)\|^2  d\mu(u) + c\\&+ 2\int \langle \mathcal{R}_{m,n}(u), R_{X,Y}(Q_X(u))- R_{X,Y}(Q_Y(u))\rangle d\mu(u).
\end{align}
Subtracting $c$ from both sides of the above display and using Cauchy-Schwarz inequality to obtain $\langle \mathcal{R}_{m,n}(u), R_{X,Y}(Q_X(u))- R_{X,Y}(Q_Y(u))\rangle\leq \|\mathcal{R}_{m,n}(u)\|\|R_{X,Y}(Q_X(u))- R_{X,Y}(Q_Y(u))\|$, we get 
\begin{align}
|T_{X,Y}-c|&\leq \int \|\mathcal{R}_{m,n}(u)\|^2 d\mu(u) \\&\quad+ 2\max_{u\in \mathcal{S}} \Big\|R_{X,Y}(Q_X(u))- R_{X,Y}(Q_Y(u)) \Big\| \cdot\int \|\mathcal{R}_{m,n}(u) \| d\mu(u)\\&\leq \int \|\mathcal{R}_{m,n} \|^2 d\mu  + 2\mathrm{diam}(\mathcal{S}) \int \|\mathcal{R}_{m,n} \| d\mu \label{eq:TXYBound}
\end{align} 
where the last inequality follows since $\max_{u\in \mathcal{S}}\|R_{X,Y}(Q_X(u))- R_{X,Y}(Q_Y(u))\|$ is bounded above by $\mathrm{diam}(\mathcal{S}) $. Due to Cauchy-Schwarz inequality, we have the following:  \begin{align}
\int \|\mathcal{R}_{m,n} \| d\mu  &\leq \Big(\int \|\mathcal{R}_{m,n} \|^2 d\mu \Big)^{1/2},\\ \mathbb{E}\Big[\Big(\int \|\mathcal{R}_{m,n} \|^2 d\mu \Big)^{1/2}\Big]&\leq \Big(\mathbb{E}\Big[\int \|\mathcal{R}_{m,n}\|^2 d\mu\Big]\Big)^{1/2}. 
\end{align}
Taking expectations on both sides of \eqref{eq:TXYBound} and applying the above inequalities into the right hand side of \eqref{eq:TXYBound} yields \eqref{eq:TxyBound}.  

Now, we turn our attention to~\eqref{eq:TxyBoundFinal}. Let us define 
\begin{align}
A_{m,n} & := \int_{\mathcal{S}}  \|\hat{R}_{X,Y}(\hat{Q}_{X}(u))- R_{X,Y}(\hat{Q}_{X}(u))\|^2 \, d \mu(u), \\
B_{m,n} & := \int_{\mathcal{S}}  \|R_{X,Y}(\hat{Q}_{X}(u))- R_{X,Y}(Q_{X}(u))\|^2 \, d \mu(u),\\
C_{m,n} & := \int_{\mathcal{S}}  \|\hat{R}_{X,Y}(\hat{Q}_{Y}(u))- R_{X,Y}(\hat{Q}_{Y}(u))\|^2 \, d \mu(u), \quad \mbox{and} \\
D_{m,n} & := \int_{\mathcal{S}}  \|R_{X,Y}(\hat{Q}_{Y}(u))- R_{X,Y}(Q_{Y}(u))\|^2 \, d \mu(u).
\end{align}

By using the fact that $(a+ b + c+d)^2 \le 4 (a^2 + b^2 + c^2+d^2)$, we may write 
\begin{align}\label{eq:BreakItDownNext}
\mathbb{E}\Big[\int \|\mathcal{R}_{m,n} \|^2 d\mu \Big] &\leq 4 \, \E \big[A_{m,n} + B_{m,n} + C_{m,n} +D_{m,n} \big].
\end{align}
In what follows, we will show that
\begin{equation} \label{eq:Bound1Next}
\E[A_{m,n} + C_{m,n}]  \leq C \, r_{d,N}, 
\end{equation}
and 
\begin{equation}\label{eq:Bound2Next}
\E[B_{m,n}+D_{m,n}]  \leq C \,r_{d,N}
\end{equation}
for some constant $C\equiv C(\mu,\nu_X,\nu_Y,\theta)>0$. Owing to \eqref{eq:BreakItDownNext}, combining \eqref{eq:Bound1Next} and \eqref{eq:Bound2Next} would then complete the proof. 

  
The proof of \eqref{eq:Bound1Next} is similar to the proof of \eqref{eq:Bound1} in Section~\ref{pf:Two-S-Rate}. 
By using a similar argument as in \eqref{eq:AmnBd}, we get 
$$A_{m,n} = \frac{1}{m}\sum_{i=1}^{m}\|\hat{R}_{X,Y}(X_i)- R_{X,Y}(X_i) \|^2, \quad C_{m,n} = \frac{1}{n}\sum_{j=1}^{n}\|\hat{R}_{X,Y}(Y_j)- R_{X,Y}(Y_j) \|^2.$$
Since $\|\hat{R}_{X,Y}(X_i)- R_{X,Y}(X_i) \|$ and $\|\hat{R}_{X,Y}(Y_i)- R_{X,Y}(Y_i) \|$ are bounded by $\mathrm{diam}(\mathcal{S})$, we upper bound $A_{m,n} + C_{m,n}$ by
\begin{align} \label{eq:ACUpBd}
\mathrm{diam}^2(\mathcal{S})&\Big[\big|\frac{N}{m}-\frac{1}{\theta}\big|+\big|\frac{N}{n}-\frac{1}{1-\theta}\big|\Big] +\max\left\{\frac{1}{\theta}, \frac{1}{1-\theta} \right\} \frac{1}{N}\Big[m A_{m,n}+ n C_{m,n}\Big].
\end{align} 
Recall that $m|N \sim \mathrm{Bin}(N,\theta)$. Thus, we have $\mathbb{E}[|\frac{N}{m}-\theta^{-1}|]\leq K N ^{-1/2}$, and $\mathbb{E}[|\frac{N}{n}-(1-\theta)^{-1}|]\leq K N^{-1/2}$ for some $K\equiv K(\theta)>0$ and $\{X_1,\ldots , X_m,Y_1,\ldots, Y_n\}$ is a set of i.i.d. random variables from the distribution $\theta\nu_{X}+ (1-\theta)\nu_Y$. Due to the last fact and Theorem~\ref{thm:RateProp1}, there exists $C >0$ such that
$$\frac{1}{N}\mathbb{E}\Big[m A_{m,n}+ n C_{m,n}\Big]\leq Cr_{d,N}.$$
Combining the above inequality with the bounds in the expectations of the other terms of \eqref{eq:ACUpBd} shows \eqref{eq:Bound1Next}.

The proof of \eqref{eq:Bound2Next} is similar to the proof of \eqref{eq:Bound2}. By our assumption, the convex potential $\psi_{X,Y}$ is strongly convex with parameter $L >0$ (say). This shows that $R_{X,Y}$ is Lipschitz continuous with parameter $\frac{1}{L}$. Owing to this, we may write 
\begin{align}
\|R_{X,Y}(\hat{Q}_{X}(u)) - R_{X,Y}(Q_{X}(u))\|& \leq \frac{1}{L}\|\hat{Q}_{X}(u)- Q_{X}(u)\|, \\  \|R_{X,Y}(\hat{Q}_{Y}(u))- R_{X,Y}(Q_{Y}(u))\| &\leq \frac{1}{L}\|\hat{Q}_{Y}(u)- Q_{Y}(u)\|. 
\end{align}
Applying these inequalities to control the integrand of $B_{m,n}$ and $D_{m,n}$, we notice that 
$$B_{m,n}+D_{m,n} \leq \frac{2}{L^2} \left[\int \|\hat{Q}_{X} - Q_{X} \|^2 \, d \mu + \int   \|\hat{Q}_{Y} - Q_{Y} \|^2 d \mu \right]. $$ Taking expectations on both sides of the above display and applying Theorem~\ref{thm:Q-Rate} shows that $\mathbb{E}[B_{m,n}+D_{m,n}|N]$ is bounded above by $C (r_{d,n}+r_{d,m})$ (which is further bounded by $C r_{d,N}$), for some $C=C(\mu,\nu_X, \nu_Y) >0$. This completes showing \eqref{eq:Bound2Next}. \qed

\subsection{Proof of Proposition~\ref{lem:Power2}}\label{sec:Pf:lem:Power2}
Note that it suffices to show the distributions of $\hat R(\hat Q(U))$ and $\tilde R(\hat Q(U))$ do not depend on the distribution of $(X_1,Y_1)$. 
Here, $U$ is an independent random vector with distribution $\mathrm{Uniform}(\s)$ where $\s = [0,1]^d$.

Let us denote by $Z_i = (X_i,Y_i)$, for $i=1,\ldots, n$. Recall that $\hat Q$ induces a cell decomposition of $\s = [0,1]^{d}$ into $n$ polyhedral sets $\{\mathcal{C}_{i} \}_{i=1}^{n}$ such that $\hat Q(u)= Z_i$ for any $u\in \mathcal{C}_{i}$.  Now, under $H_0$, for Borel $B \subset \s$,
 \begin{align*}
  \mathbb{P}&\Big(\hat R(\hat Q(U))\in B \Big| Z_1, \ldots , Z_n \Big)\nonumber\\& =  \sum_{i=1}^{n}\mathbb{P}(\hat Q(U)=Z_i |Z_1, \ldots, Z_n)\; \mathbb{P} (\hat R(Z_i)\in B|Z_1, \ldots , Z_n) \nonumber \\ 
  & = \frac{1}{n} \sum_{i=1}^{n} \mathbb{P} (\hat R(Z_i)\in B|Z_1, \ldots , Z_n)\nonumber
 \end{align*}
 where the second equality follows as $\mathbb{P}(\hat Q(U)= Z_i|Z_1, \ldots ,Z_n) = n^{-1}$, for all $i=1,\ldots, n$. Thus, $$ \mathbb{P}\Big(\hat R(\hat Q(U))\in B \Big) = \mathbb{P} (\hat R(Z_i)\in B) = \mu(B)$$
by Lemma~\ref{lem:Uniform} (as $\hat R(Z_i) \sim \mu$ for every $i=1,\ldots, n$).

Let us next show that $\tilde  R(\hat Q(U)) \sim \mu$. 
Let us denote by $\hat Q_X$ ($\hat Q_Y$) the sample quantile function obtained only from the $X_i$'s ($Y_i$'s) and let the corresponding a cell decomposition of $[0,1]^{d_X}$ ($[0,1]^{d_Y}$) into $n$ polyhedral sets be denotes by $\{\mathcal{C}_{X,i} \}_{i=1}^{n}$ ($\{\mathcal{C}_{Y,i} \}_{i=1}^{n}$), i.e., $\hat Q_X(u)= X_i$ for any $u\in \mathcal{C}_{X,i}$ (and $\hat Q_Y(u)= Y_i$ for any $u\in \mathcal{C}_{Y,i}$).

Also, let $\mu_X := $ Uniform$([0,1]^{d_X})$ and $\mu_Y := $ Uniform$([0,1]^{d_Y})$. Now, for $B_X \subset [0,1]^{d_X}$ and $B_Y \subset [0,1]^{d_Y}$, and under $H_0$,
 \begin{align}
  \mathbb{P}&\Big(\tilde R(\hat Q(U))\in B_X \times B_Y \Big| Z_1, \ldots , Z_n \Big) \nonumber \\
  & =  \sum_{i=1}^{n}\mathbb{P}(\hat Q(U)=(X_i,Y_i)|Z_1, \ldots , Z_n)\times \mathbb{P} (\tilde R(X_i,Y_i)\in B_X \times B_Y|Z_1, \ldots , Z_n) \nonumber \\  
  & = \sum_{i=1}^{n} \frac{1}{n} \mathbb{P} (\hat R_X(X_i) \in B_X, \hat R_Y(Y_i) \in B_Y|Z_1, \ldots , Z_n)  \nonumber\\
  & = \frac{1}{n} \sum_{i=1}^{n} \left[ n \mu_X(\mathcal{C}_{X,i}\cap B_X) \right] \left[ n \mu_Y(\mathcal{C}_{Y,i}\cap B_Y) \right]  \nonumber
 \end{align}
 where the second equality follows from: (i) $\mathbb{P}(\hat Q(U)= (X_i,Y_i)|Z_1, \ldots ,Z_n) = n^{-1}$ and 
\begin{eqnarray*}
&& \mathbb{P} (\hat R_X(X_i) \in B_X, \hat R_Y(Y_i) \in B_Y|Z_1, \ldots , Z_n) \\
& = & \mathbb{P} (\hat R_X(X_i) \in B_X| X_1,\ldots, X_n)  \times \mathbb{P} (\hat R_Y(Y_i) \in B_Y|Y_1, \ldots , Y_n).
\end{eqnarray*} 
The third equality follows from the fact that $\hat R_{X}(X_i)$ is a randomly chosen point uniformly drawn from the polytope $\mathcal{C}_{X, i}$ where $\mu_X(\mathcal{C}_{X,i}) = \lambda_{d_X}(\mathcal{C}_{X, i})= n^{-1}$. Thus, we see that 
{\small \begin{align}\label{eq:ProbB-Ind}
\qquad  \mathbb{P}\Big(\tilde R(\hat Q(U))\in B_X \times B_Y\Big) = n \;\mathbb{E}_{X,Y}\Big[ \sum_{i=1}^{n}  \mu_X(\mathcal{C}_{X,i}\cap B_X) \mu_Y(\mathcal{C}_{Y,i}\cap B_Y)\Big],
\end{align}} 
where $\mathbb{E}_{X,Y}$ denotes the expectation with respect to $(X_1,Y_1), \ldots , (X_n, Y_n)$.

As the $X_i$'s are independent of the $Y_i$'s, then, letting $\mathfrak{S}_{n}$ be the set of all permutations of $\{1, \ldots ,n\}$, $(Z_1, \ldots , Z_n) \stackrel{d}{=} ((X_1,Y_{\sigma(1)}), \ldots ,(X_n, Y_{\sigma(n)}))$ for any random permutation $\sigma$ chosen uniformly  from $\mathfrak{S}_{n}$. Owing to this, for any random $\sigma$ chosen uniformly from $\mathfrak{S}_{n}$, $$\sum_{i=1}^{n}  \mu_X(\mathcal{C}_{X,i}\cap B_X) \mu_Y(\mathcal{C}_{Y,i}\cap B_Y) \stackrel{d}{=} \sum_{i=1}^{n}  \mu_X(\mathcal{C}_{X,i}\cap B_X) \mu_Y(\mathcal{C}_{Y,\sigma(i)}\cap B_Y).$$
Using the above display, we notice that 
  \begin{align}
 \mathbb{E}_{X,Y} &\Big[\sum_{i=1}^{n}  \mu_X(\mathcal{C}_{X,i}\cap B_X) \mu_Y(\mathcal{C}_{Y,i}\cap B_Y) \Big] \\&= \mathbb{E}_{Z,\sigma}\Big[\sum_{i=1}^{n}  \mu_X(\mathcal{C}_{X,i}\cap B_X) \mu_Y(\mathcal{C}_{Y,\sigma(i)}\cap B_Y)\Big]\nonumber  \\ 
 &= \mathbb{E}_{Z}\Big[\mu_X(\mathcal{C}_{X,i}\cap B_X)  \sum_{i=1}^{n}  \mu_Y(\mathcal{C}_{Y,\sigma(i)}\cap B_Y)\Big] \nonumber \\
  &= \mathbb{E}_{Z}\Big[\mu_X(\mathcal{C}_{X,i}\cap B_X) \mu_Y(\cup_{i=1}^{n}  \mathcal{C}_{Y,\sigma(i)}\cap B_Y)\Big] \nonumber \\
    &=  \mu_Y(B_Y) \; \mathbb{E}_{X}\Big[\mu_X(\mathcal{C}_{X,i}\cap B_X)\Big] \nonumber \\
 &= \mu_Y(B_Y) \; n^{-1} \P(\hat R_X(X_i) \in B_X) \nonumber \\
  &=n^{-1} \mu_Y(B_Y) \; \mu_X(B_X),\label{eq:ExpSame-Ind} 
\end{align}    
where the fifth equality above holds as 
\begin{align}\label{eq:SumExpan-Ind}
\P(\hat R_X(X_i) \in B_X) = \E [\mathbbm{1}(\hat{R}_n(X_i)\in B)|X_1,\ldots, X_n]= n\mathbb{E}[\mu_X(\mathcal{C}_{X,i}\cap B)].
\end{align}   
Therefore, using~\eqref{eq:ProbB-Ind}, we have 
 $$ \mathbb{P}\Big(\tilde R(\hat Q(U))\in B_X \times B_Y\Big) = \mu_Y(B_Y) \; \mu_X(B_X) = \mu(B_X \times B_Y),$$ thereby implying that $\tilde R(\hat Q(U)) \sim \mu$.
 
The rest of the proof follows from a similar argument as in the proof of Proposition~\ref{lem:Power1}. 
\qed

\subsection{Proof of Proposition~\ref{lem:Rate2}}\label{pf:Rate2} 
Note that 
\begin{equation}\label{eq:T_n-Ind}
T_{n} \le \frac{2 }{n} \left[\sum_{i=1}^n \|\hat R(Z_i) - R(Z_i)\|^2 + \sum_{i=1}^n \| R(Z_i) - \tilde R(Z_i)\|^2 \right].
\end{equation}
Observe that the expectation of the first term on the right side of the above display can be bounded above as $$\E \left[ \frac{1}{n} \sum_{i=1}^n \|\hat R(Z_i) - R(Z_i)\|^2  \right]
\le C \,r_{d,n},$$ by an application of Theorem~\ref{thm:RateProp1}, for some constant $C \equiv C(\mu,\nu) >0$. As, 
$$ \sum_{i=1}^n \| R(Z_i) - \tilde R(Z_i)\|^2= \sum_{i=1}^n \| \hat R_X(X_i) - R_X(X_i)\|^2  + \sum_{i=1}^n \| \hat R_Y(Y_i) - R_Y(Y_i)\|^2,$$ the expectation of the second term in~\eqref{eq:T_n-Ind} can be upper bounded as $$\E \left[  \sum_{i=1}^n \| R(Z_i) - \tilde R(Z_i)\|^2 \right] \le C  (r_{d_X,n} + r_{d_Y,n})$$ for some constant $C \equiv C(\mu,\nu) >0$ (by using Theorem~\ref{thm:RateProp1}). This completes the proof. \qed

\subsection{Proof of Proposition~\ref{lem:Rate2-Alt}}\label{pf:Rate2-Alt}
Recall that $T_n =\frac{1}{n}\sum_{i=1}^n \|\hat R(Z_i) - \tilde R(Z_i)\|^2$. Let $\mathcal{R}_n :=\frac{1}{n}\sum_{i=1}^n \|\hat R(Z_i) - R(Z_i)\|^2$. The idea of the proof is to bound (the absolute values of) the following terms separately:
\begin{align}
T_{n}-\frac{1}{n}\sum_{i=1}^{n}  \| R(Z_i) - \tilde R(Z_i)\|^2, \hspace{0.4in} \label{eq:3-Terms-1}\\
 \frac{1}{n}\sum_{i=1}^{n}  \| R(Z_i) - \tilde R(Z_i)\|^2 - \frac{1}{n}\sum_{i=1}^n \| R(Z_i) - \bar R(Z_i)\|^2, \quad \mbox{and} \label{eq:3-Terms-2}\\
\frac{1}{n} \sum_{i=1}^n \| R(Z_i) - \bar R(Z_i)\|^2  - \int_{[0,1]^d} \| u - \bar R(Q(u))\|^2 du. \hspace{0.4in} \label{eq:3-Terms-3}
\end{align}
Combining these would yield the desired result. 

Let us start with the first term in the above display, i.e.,~\eqref{eq:3-Terms-1}. By using a similar argument as in the proof of \eqref{eq:TxyBound}, we have
\begin{align}\label{eq:T_n-Ind2}
\quad\mathbb{E}\left[\Big|T_{n}-\frac{1}{n}\sum_{i=1}^{n}  \| R(Z_i) - \tilde R(Z_i)\|^2 \Big|\right] &\le  \mathbb{E}\big[\mathcal{R}_n\big] + 2\mathrm{diam}(\mathcal{S})\Big(\mathbb{E}\big[\mathcal{R}_n\big]\Big)^{\frac{1}{2}}.
\end{align}
 By an application of Theorem~\ref{thm:RateProp1}, we observe that $\E \left[ \mathcal{R}_n  \right]
\le C \,r_{d,n},$ for some constant $C \equiv C(\mu,\nu) >0$. To bound~\eqref{eq:3-Terms-2}, we write 
\begin{eqnarray*} 
 \sum_{i=1}^n \| R(Z_i) - \tilde R(Z_i)\|^2 &= &\sum_{i=1}^n \| R(Z_i) - \bar R(Z_i)\|^2 + \sum_{i=1}^n \| \bar R(Z_i) - \tilde R(Z_i)\|^2
 \\ && \qquad + 2 \sum_{i=1}^n \langle R(Z_i) - \bar R(Z_i), \bar R(Z_i) - \tilde R(Z_i)\rangle.  
\end{eqnarray*} 
Now, letting $\mathcal{S}_n :=  \frac{1}{n}\sum_{i=1}^n \| \bar R(Z_i) - \tilde R(Z_i)\|^2 $, and using the above relation and a similar argument to derive \eqref{eq:TxyBound},
\begin{align}
\mathbb{E}\Big[\Big|\frac{1}{n}&\sum_{i=1}^{n}  \| R(Z_i) - \tilde R(Z_i)\|^2 - \frac{1}{n} \sum_{i=1}^n \| R(Z_i) - \bar R(Z_i)\|^2\Big|\Big]\\
&\leq \mathbb{E}\big[\mathcal{S}_n \big]+ 2\mathrm{diam}(\mathcal{S})\Big(\mathbb{E}\big[\mathcal{S}_n \big]\Big)^{\frac{1}{2}}.
\end{align}
As, 
$$\mathcal{S}_n = \sum_{i=1}^n \| \hat R_X(X_i) - R_X(X_i)\|^2+ \sum_{i=1}^n \| \hat R_Y(Y_i) - R_Y(Y_i)\|^2$$
the expectation of the right hand side of the previous display can be upper bounded yielding $$\E \left[\mathcal{S}_n  \right] \le C  (r_{d_X,n} + r_{d_Y,n}) \le C \, r_{d,n}$$ for some constant $C \equiv C(\nu_X,\nu_Y) >0$ (by using Theorem~\ref{thm:RateProp1}). We can easily bound~\eqref{eq:3-Terms-3} by an application of Hoeffding's inequality as we are dealing with an average of bounded random variables with expectation $\int_{[0,1]^d} \| u - \bar R(Q(u))\|^2 du$. Combining all this yields the desired conclusion. \qed

\section{Simulation study}\label{sec:Simul}
\subsection{Two-sample goodness-of-fit testing}\label{sec:Simul-2S}
Consider the testing problem~\eqref{eq:2-Sample-Test} using the test statistic $T_{X,Y}$ defined in~\eqref{eq:2-S-TS}. Let us denote by $\{W_i^X\}_{i=1}^m$ the cell decomposition of $\s$ induced by the sample $X_1,\ldots, X_m$ (see~\eqref{eq:W_i}); similarly, let $\{W_j^Y\}_{j=1}^n$ be the cell decomposition of $\s$ induced by $Y_1,\ldots, Y_n$. Then, the test statistic $T_{X,Y}$ can be expressed as: $$T_{X,Y} = \sum_{i=1}^m \sum_{j=1}^n \|\hat R_{X,Y}(X_i) - \hat R_{X,Y}(Y_j)\|^2 p_{ij}$$ where $p_{ij}$ is the volume of $W_i^X \cap W_j^Y$. However, it is computationally intensive to calculate $T_{X,Y}$ using the above exact expression as finding $p_{ij}$, for $i,j \in \{0,1,\ldots, n\}$, in non-trivial. A natural alternative is to use Monte Carlo simulations to approximate $T_{X,Y}$, as defined in~\eqref{eq:2-S-TS}, by drawing random samples from the distribution $\mu$. This is the strategy we adopt in this subsection where we use $10^4$ Monte Carlo samples to approximate $T_{X,Y}$.

\begin{figure}
\includegraphics{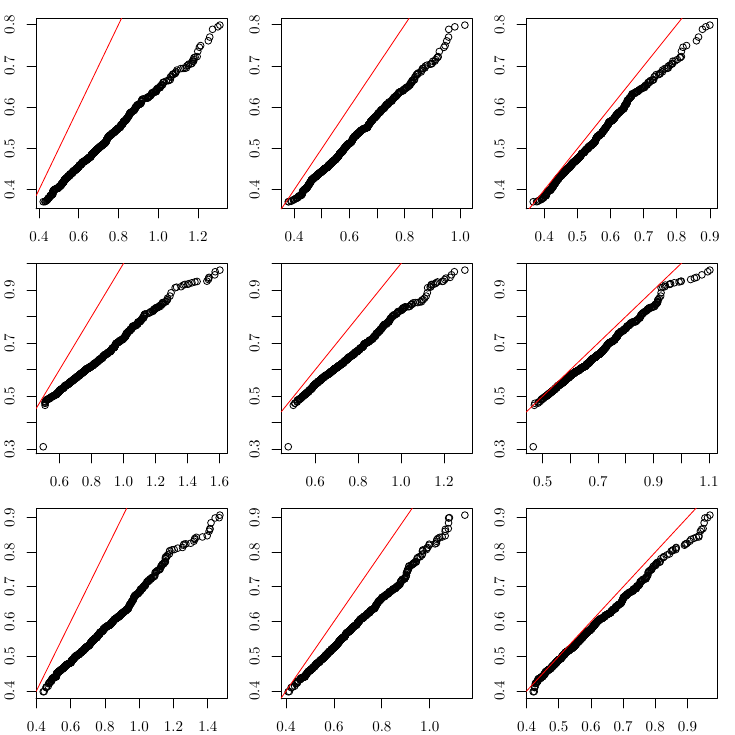}
\caption{The three rows correspond to three different simulation settings --- (i), (ii), and (iii). As we move from left to right in each row, the QQ-plots are obtained by comparing the sampling distribution of $(n/\log n) T_{X,Y}$ for $m \equiv n=100$ (left-most plots), $m \equiv n=1000$ (center plots), $m \equiv n=10^4$ (right-most plots) with that of $m \equiv n=5 \times 10^4$.}  
\label{fig:T-X-Y-Conv}
\end{figure}

In this subsection we focus on the case $d=2$, $m \equiv n$, and $\nu_X = \nu_Y$ (i.e., the null hypothesis holds), and study the following three questions:
\begin{enumerate}
	\item[(a)] Does $T_{X,Y}$, properly normalized, converge to a nondegenerate asymptotic weak limit as $m \equiv n$ grows?
	
	\item[(b)] Is the above asymptotic limit distribution-free, i.e., free of the model parameters, and thus, universal?
	
	\item[(c)] Is $T_{X,Y}$ distribution-free for every finite $n, m$?
\end{enumerate}
To try to answer the above questions we conduct a simulation study. We take $\nu_X = \nu_Y$ and we consider the following settings:
\begin{enumerate}
\item[(i)] The multivariate Gaussian distribution \(N_2({0}, {I})\).

\item[(ii)] The Gaussian mixture distribution \[\frac{1}{2}\cdot N_2\left[\begin{pmatrix}
5\\0
\end{pmatrix}, \begin{pmatrix}
2 & 0\\
0 & 2
\end{pmatrix}\right]+\frac{1}{2}\cdot N_2\left[\begin{pmatrix}
0\\5
\end{pmatrix}, \begin{pmatrix}
5 & 2\\
2 & 5
\end{pmatrix}\right].\]
\item[(iii)] The banana-shaped distribution as in~\cite{Cher17} generated as \[\begin{pmatrix}
X+R\cos\Phi\\
X^2+R\sin\Phi
\end{pmatrix},\] where \(X\sim U[-1,1]\), \(\Phi\sim U[0, 2\pi]\), and \(R=0.2 Z[1+(1-|X|)/2]\) for \(Z\in U[0,1]\). Here $X$, $\Phi$ and $Z$ are drawn independently. 
\end{enumerate}

\begin{figure}
\includegraphics{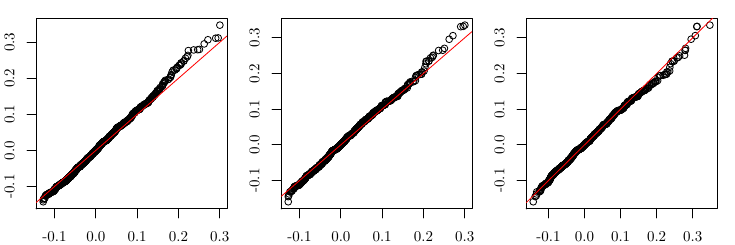}
\caption{The three QQ-plots are obtained when comparing the asymptotic weak limits of $(n/\log n) \{T_{X,Y} - \E (T_{X,Y})\}$ (approximated by $m \equiv n=5 \times 10^4$) for settings (i) versus (ii), (i) versus (iii), and (ii) versus (iii) (from left to right).}  
\label{fig:T-X-Y-Asym}
\end{figure}

\begin{figure}
\includegraphics{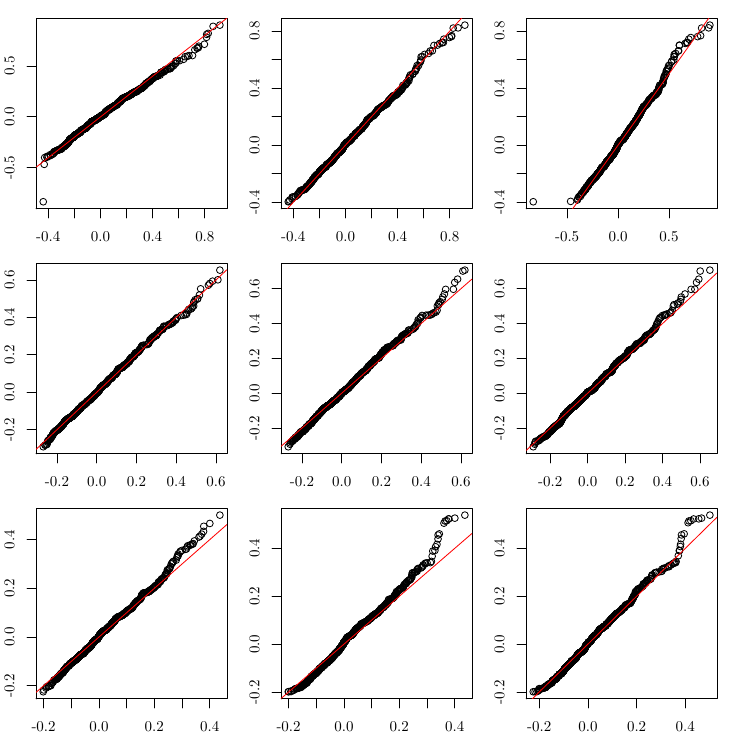}
\caption{The three rows correspond to three different sample size settings: $m \equiv n = 25$, $m \equiv n = 100$ and $m \equiv n = 1000$. As we move from left to right in each row, the three QQ-plots are obtained when comparing the asymptotic limits of $(n/\log n) \{T_{X,Y} - \E (T_{X,Y})\}$ (for the three settings of $m= n$) for settings (i) versus (ii), (i) versus (iii), and (ii) versus (iii).}  
\label{fig:T-X-Y-FS}
\end{figure}

We draw i.i.d.~samples from the above distributions with sample sizes $m \equiv  n= 100,  1000, 10^4$ and $5 \times 10^4$ and compute the test statistic $T_{X,Y}$. To approximate the distribution of $T_{X,Y}$, for every sample size and simulation setting, we use 1000 independent replications. From~\cite{AKT} it seems natural to consider the distribution of the normalized statistic $(n/\log n) T_{X,Y}$.

Figure~\ref{fig:T-X-Y-Conv} shows the QQ-plots obtained from comparing the distribution of $(n/\log n) T_{X,Y}$, for $m \equiv n=100, 1000, 10^4$, with that for $m \equiv n=5 \times 10^4$, for simulation settings (i)-(iii). The plots illustrate that as the sample sizes $m \equiv n$ increase, the sampling distribution of $(n/\log n) T_{X,Y}$ converges to a weak limit.

Figure~\ref{fig:T-X-Y-Asym} gives the QQ-plots for pairwise comparisons of the asymptotic limits of $(n/\log n) \{T_{X,Y} - \E (T_{X,Y})\}$ for the simulation settings (i)-(iii). The QQ-plots (and the corresponding two-sample Kolmogorov-Smirnov tests; not provided here) illustrate that the limiting distributions (as approximated by taking $m \equiv n=10^5$) are the same for the three settings considered. Note that the centering by $\E (T_{X,Y})$ is necessary, without which there seems to be a difference in the means of the asymptotic weak limits. This is, by itself, an interesting phenomenon that needs further study.

Figure~\ref{fig:T-X-Y-FS} shows the QQ-plots of the distribution of $(n/\log n) \{T_{X,Y} - \E (T_{X,Y})\}$, for $m \equiv n = 25$, $m\equiv n = 100$ and $m\equiv n = 1000$. The QQ-plots (and the corresponding two-sample Kolmogorov-Smirnov tests; not given here) show remarkable resemblance with the $y=x$ line (in red) which suggests that the finite sample distributions are probably the same for the three settings considered. We plan to thoroughly investigate this phenomenon in a future paper. As above, the centering is necessary, without which the QQ-plots differ by a constant shift.

\subsection{Testing for mutual independence}\label{sec:Simul-Ind}
\begin{figure}
\includegraphics{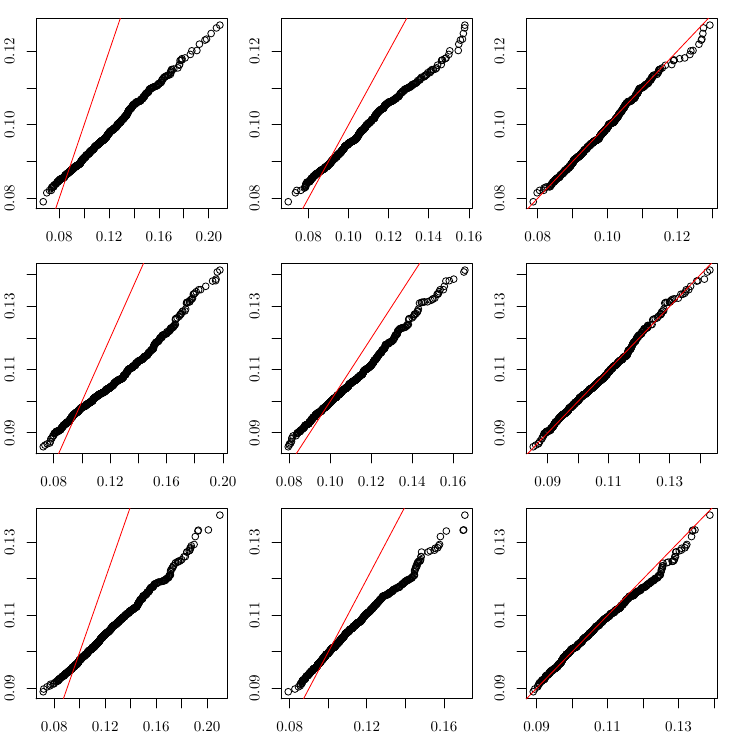}
\caption{The three rows correspond to the three different simulation settings (i)-(iii). As we move from left to right in each row, the QQ-plots are obtained by comparing the sampling distribution of $(n/\log n) T_{n}$ for $n=100$ (left-most plots), $n=1000$ (center plots), $n=5 \times 10^4$ (right-most plots) with that of $n=10^5$.}  
\label{fig:T-n-Conv}
\end{figure}

\begin{figure}
\includegraphics{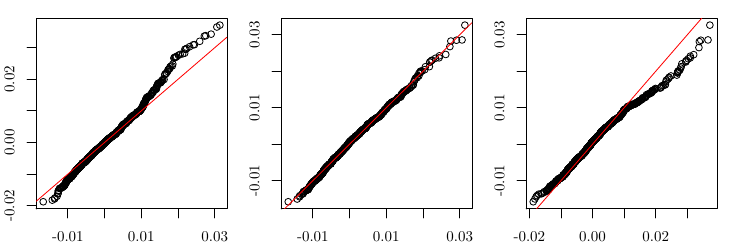}
\caption{The three QQ-plots are obtained when comparing the asymptotic weak limits of $(n/\log n) \{T_{n} - \E (T_{n})\}$ (approximated by $n=10^5$) for settings (i) versus (ii), (i) versus (iii), and (ii) versus (iii) (from left to right).}  
\label{fig:T-n-Asym}
\end{figure}

Let us consider the testing problem~\eqref{eq:Test-2} using the test statistic $T_n$ defined in~\eqref{eq:Ind-TS}. We focus on the case $d=2$ and illustrate the behavior of $T_n$ under different simulation settings. In particular, we use $T_n$ to answer the three questions (a)-(c) outlined in Section~\ref{sec:Simul-2S}. In our study we consider the following simulation settings:
\begin{enumerate}
	\item[(i)] Two independent standard Gaussian: \(X\sim N(0, 1)\) and \(Y\sim N(0,1)\).
	\item[(ii)] Two independent random variables having log-normal and Gamma distribution: \(X\sim\ln N(0, 0.5)\) and \(Y\sim\Gamma(3, 2)\).
	\item[(iii)] Two independent Gaussian mixtures: \[X\sim\frac{1}{4}\cdot N(-1, 2)+\frac{3}{4}\cdot N(5,3),\] and \[Y\sim\frac{3}{10}\cdot N\left(0,\frac{1}{2}\right)+\frac{3}{10}\cdot N(5,2)+\frac{2}{5}\cdot N(-5,1).\]
\end{enumerate}
We draw i.i.d.~samples from the above distribution with sample sizes $n=25, 100, 1000, 5\times 10^4$ and $10^5$ and compute the test statistic $T_n$. To approximate the distribution of $T_n$, for every sample size and simulation setting we use 1000 independent replications. From~\cite{AKT} it seems natural to consider the distribution of the normalized statistic $(n/\log n) T_n$.

\begin{figure}
\includegraphics{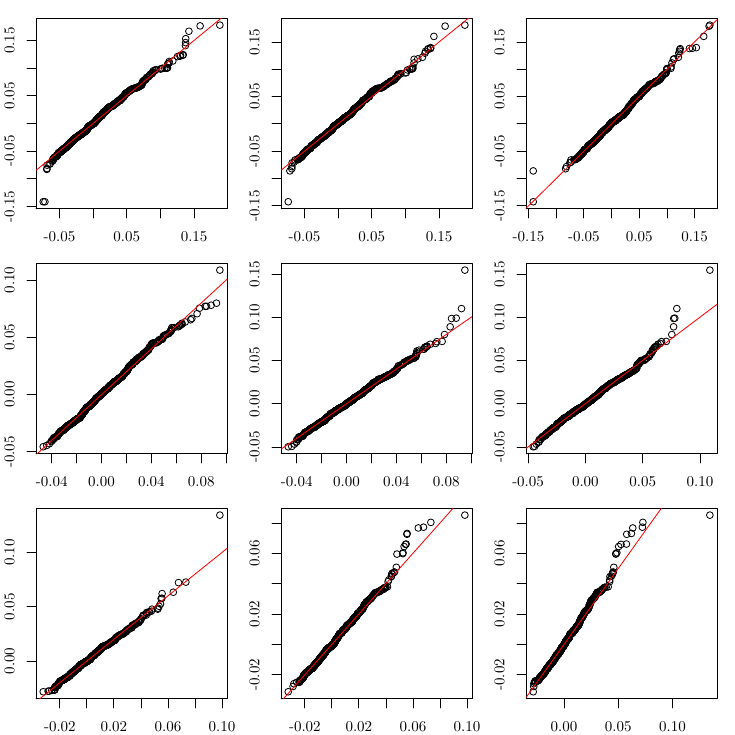}
\caption{The three rows correspond to three different sample size settings: $n = 25$, $n = 100$ and $n = 1000$. As we move from left to right in each row, the three QQ-plots are obtained when comparing the asymptotic limits of $(n/\log n) \{T_{n} - \E (T_{n})\}$ for settings (i) versus (ii), (i) versus (iii), and (ii) versus (iii).}  
\label{fig:T-n-FS}
\end{figure}

Figure~\ref{fig:T-n-Conv} gives the QQ-plots of the distribution of $(n/\log n) T_n$, for $n=100, 1000, 5 \times 10^4$, compared to that when $n=10^5$, for settings (i)-(iii). It illustrates that indeed, as the sample size increases, the sampling distribution of $(n/\log n) T_n$ converges to a weak limit, for each of the simulation settings. 

As illustrated in the QQ-plots in Figure~\ref{fig:T-X-Y-Asym}, we observe a similar phenomenon in Figure~\ref{fig:T-n-Asym}, which shows that the asymptotic distribution of $(n/\log n) \{T_{n} - \E (T_{n})\}$ (approximated by its sampling distribution when $n=10^5$) in the three different settings are probably the same. The corresponding two-sample Kolmogorov-Smirnov tests (not provided here) are all accepted under 5\% nominal level. Finally, Figure~\ref{fig:T-n-FS} illustrates that the finite sample distributions of $(n/\log n) \{T_{n} - \E (T_{n})\}$ are also very similar. 

Our simulation studies indicate that in both the testing problems (see~\eqref{eq:2-Sample-Test} and~\eqref{eq:Test-2}) considered in this paper, the (asymptotic) distribution of $(n/\log n) \{T_{n} - \E (T_{n})\}$ is universal. This justifies the form of the test statistics considered in this paper. However, a formal treatment of these issues is beyond the scope of the present paper and would be an interesting future research direction.

\section{Proofs of the Claims in the Appendix~D.2}\label{Appendix-B}

\subsection{Proof of Claim~\ref{cl:MT2}}\label{sec:cl:MT2}
Since $\psi$ is the Legendre-Fenchel dual of $\psi^{*}$, hence, $\partial\psi^{*}(B) = (\partial\psi)^{-1}(B)$, for any Borel $B \subset \R^d$, by Lemma~\ref{lem:SubD-1}. Since $\psi^{*}$ is a continuous function on $\RR^d$, by~\cite[Theorem~1.1.13]{G16}, $\partial \psi^*(B)$ is   Lebesgue measurable for any Borel set $B\subset \RR^d$. Therefore, we may write $\mu((\partial \psi)^{-1}(B))= \mu(\partial \psi^{*}(B))$. In order to show that $\partial  \psi\# \mu  = \nu$, it is enough to show that $\mu(\partial \psi^{*}(B))= \nu(B)$ for all Borel sets $B \subset \RR^d$. This will be proved hereafter.  

Now, suppose $G\subset \RR^d$ is an open set. We claim and prove that 
\begin{align}\label{eq:LiminfClaim}
\liminf_{p\to \infty}\mu(\partial \hat{\xi}_{p}(G))\geq \mu(\partial \psi^{*}(G)) 
\end{align}
which, by the Portmanteau theorem, will show that $\mu(\partial \hat{\xi}_{p}(\cdot))\stackrel{d}{\to} \mu(\partial \psi^{*}(\cdot))$ weakly. We first show \eqref{eq:LiminfClaim} when $G$ is a bounded set and hence, the closure of $G$ is a compact set. Define $\mathfrak{B}:= \{u\in \mathcal{S}: \partial\psi(u) \text{ is not a singleton set}\}$. Since $\psi$ is a convex function and $\mu$ is absolutely continuous, by Alexandrov's differentiability theorem (\cite{Ambro08}), we know that $\mu(\mathfrak{B})=0$. 
It suffices to show \eqref{eq:LiminfClaim} with $\mu(\partial \psi^{*}(G)\backslash \mathfrak{B})$ on the right hand side. Fix $u_0\in \partial \psi^{*}(G)\backslash \mathfrak{B}$. We show that there exists $p_0= p_0(u_0)\in \NN$ such that $u_0\in \partial \hat{\xi}_{p}(y_p)$ for some $y_{p}\in G$, for all $p\geq p_{0}$. Due to Lemma~\ref{lem:SubD},
\begin{align}
\mathfrak{B} = \{u\in \mathcal{S}: u\in \partial \psi^{*}(x_1)\cap \partial \psi^{*}(x_2), \text{ for some }x_1\neq x_2\in \RR^d\}.
\end{align} 
As $u_0\notin \mathfrak{B}$, there exists a unique $x_0\in G$ such that $u_0\in \partial \psi^{*}(x_0)$ and $u_0\notin \partial \psi^{*}(x)$ for any $x\in \RR^d$. Since $G$ is open, there exists $\eta=\eta(u_0)>0$ such that Cl(${B_{\eta}(x_0)}$) (closure of the open ball $B_{\eta}(x_0)$) belongs to $G$. Let 
\begin{align}
\Theta_{\eta}:= \inf_{x\in \mathrm{Bd}(B_{\eta}(x_0))}\Big\{\psi^{*}(x) - \psi^{*}(x_0)- \langle u_0, x-x_0\rangle\Big\}.
\end{align}
Note that $\Theta_{\eta}>0$ because $\psi^{*}$ is convex and $u_0 \in \partial \psi^{*}(G)\backslash\mathfrak{B}$. Since $G\subset \mathcal{F}_{m}$ for some $m\in \NN$ (as $G$ is bounded), $\hat{\xi}_{p}$ converges uniformly to $\psi^{*}$ in Cl($B_{\eta}(x_0)$). Hence, there exists $p_0=p_0(u_0)$ such that $| \hat{\xi}_{p}(x)-\psi^{*}(x)|\leq \Theta_{\eta}/3$ for all $x\in \mathrm{Cl}(B_{\eta}(x_0))$ and $p\geq p_0$. Therefore, for all $p\geq p_0$, 
\begin{align}
\inf_{x\in \mathrm{Bd}(B_{\eta}(x_0))}\Big\{\hat{\xi}_{p}(x) - \hat{\xi}_{p}(x_0) -\langle u_0, x-x_0\rangle\Big\}\geq \frac{\Theta_{\eta}}{3}>0.\label{eq:Infimum}
\end{align}  
Let 
\begin{align}
y_{p}:= \arg \min_{x\in \mathrm{Cl}(B_{\eta}(x_0)) }\Big\{\hat{\xi}_{p}(x) - \hat{\xi}_{p}(x_0) -\langle u_0, x-x_0\rangle\Big\}.\label{eq:Defyl}
\end{align}
Note that the minimum value of $g(x) := \hat{\xi}_{p}(x) - \hat{\xi}_{p}(x_0) -\langle u_0, x-x_0\rangle$ in $\mathrm{Cl}(B_{\eta}(x_0))$ is bounded above by $0$ because $\hat{\xi}_{p}(x_0) - \hat{\xi}_{p}(x_0) -\langle u_0, x_0-x_0\rangle =0$. Owing to \eqref{eq:Infimum}, $y_{p}$ cannot be on the boundary of the ball $B_{\eta}(x_0)$. Hence, using the definition of $y_p$ in \eqref{eq:Defyl}, for all $p\geq p_0$,
$$g(x) \ge g(y_p) \quad \forall x\in B_{\eta}(x_0) \quad \Leftrightarrow \quad \hat{\xi}_{p}(x) \geq \hat{\xi}_{p}(y_{p}) +\langle u_0,x-y_{p}\rangle, \quad \forall x\in B_{\eta}(x_0). $$
Since $B_{\eta}(x_0)$ is an open set containing $y_p$ and $\hat{\xi}_{p}$ is a convex function, the above display along with Lemma~\ref{lem:SubGIneq} implies $u_0\in \partial\hat{\xi}_{p}(y_{p})$, for all $p\geq p_0$. Hence, for the functions $f(u):=\mathbbm{1}(u\in \partial\psi^{*}(G)\backslash\mathfrak{B} )$ and $f_{p}(u):= \mathbbm{1}(u\in \partial \hat{\xi}_{p}(G))$, we get $\liminf_{p\to \infty} f_{p}(u)\geq f(u)$ for all $u\in \mathcal{S}$. Applying Fatou's lemma, 
\begin{align}
\liminf_{p\to \infty}\mu(\partial \hat{\xi}_{p}(G)) &= \liminf_{p\to \infty}\int f_{p}(u)\mu(du) \nonumber\\&\geq \int f(u)\mu(du)= \mu(\partial\psi^{*}(G)\backslash \mathfrak{B}) = \mu(\partial \psi^{*}(G)).
\end{align} 
This proves \eqref{eq:LiminfClaim} when $G$ is an open bounded set. Now, we suppose that $G$ is open, but unbounded. There exists a sequence of open bounded sets $\{G_{j}\}_{j\geq 1}$ such that $G_j\uparrow G$ as $j\to \infty$.  Note that $\hat{\xi}_p$ converges uniformly to $\psi^{*}$ in Cl(${G}_j$), for all $j\geq 1$. Therefore,   
\begin{align}
\liminf_{p\to \infty} \mu(\partial\hat{\xi}_p(G))\geq \liminf_{p\to \infty} \mu(\partial\hat{\xi}_p(G_j)) \geq \mu(\partial\psi^{*}(G_j)) \label{eq:partialIneq}
\end{align} 
for all $j\geq 1$ where the first inequality follows since $G_j\subset G$ and the second inequality holds since we have proved \eqref{eq:LiminfClaim} for all open bounded sets. As $j\to \infty$, the right hand side of \eqref{eq:partialIneq} converges to $\mu(\partial\psi^{*}(G))$ by an application of the monotone convergence theorem. Hence, \eqref{eq:LiminfClaim} holds for all open sets $G\subseteq \RR^d$.

Returning to the proof of the claim, owing to Lemma~\ref{lem:SubD-1}, for any Borel set $B \subset \R^d$, $$\mu\big(\partial \hat{\psi}^{*}_n (B)\big) = \mu\big( (\partial\hat{\psi}_n)^{-1} (B)\big)= \hat{\nu}_n(B), \qquad \mbox{for all}\;\; n \ge 1.$$ Combining this observation with \eqref{eq:LiminfClaim} shows that $\hat{\nu}_n(\cdot)$ weakly converges to the measure $\mu(\partial \psi^{*}(\cdot))$, along the subsequence corresponding to $\{\partial \hat{\xi}_{p}\}_{p \ge 1}$. However, we know that $\hat{\nu}_n$ converges weakly to $\nu$, for all $\omega \in \mathcal{A}_K\cap \mathcal{D}$. By the uniqueness of the weak limit, $\nu(B)= \mu(\partial\psi^{*}(B))$ for all Borel sets $B\subset \RR^d$. This completes the proof of the claim. \qed 

\subsection{Proof of Claim~\ref{cl:MT3}}\label{sec:cl:MT3}
We first claim and prove that there exists an open set $U$ such that \begin{align}\label{eq:Construction}
\mathfrak{K}\subset U\subset\mathrm{Cl}(U)\subset \mathrm{Int}(\s).
\end{align}
Owing to the compactness of $\mathfrak{K}$ and $\mathrm{Bd}(\s)$, there exists $x_0\in \mathfrak{K}$ and $y_0\in \mathrm{Bd}(\s)$ such that 
\begin{align}
\|x_0-y_0\| = \argmin_{x\in \mathfrak{K}, y\in \mathrm{Bd}(\s)}\|x-y\|>0
\end{align}
Denote $\delta:=\|x_0-y_0\|$. Define $U :=\mathfrak{K}+B_{\delta/2}(0)$. We will now show that $U \subset \mathrm{Int}(\s)$, i.e., for any $y\in \mathfrak{K}$, $B_{\delta/2}(y) \in \mathrm{Int}(\s)$. We prove this by contradiction. Suppose that there exists $y\in \mathfrak{K}\subset \mathrm{Int}(\s)$ and $z\in B_{\delta/2}(y)$ such that $z\notin \mathrm{Int}(\s)$. Then, there must exist $z^{\prime}$ in the line joining $y$ and $z$ such that $z^{\prime}\in \mathrm{Bd}(\s)$. This implies $\|y-z^{\prime}\|<\delta/2$ and hence, contradicts the definition of $\delta$. Therefore, $B_{\delta/2}(y)\subset \mathrm{Int}(\s)$ for all $y\in \mathfrak{K}$ and hence, $U \subset \mathrm{Int}(\s)$. Since the distance between $\mathrm{Cl}(U)$ and $\mathrm{Bd}(\s)$ is equal to $\delta/2$, therefore, $\mathrm{Cl}(U)\subset \mathrm{Int}(\s)$. This shows \eqref{eq:Construction}. 
 
 Now, we return to prove Claim~\ref{cl:MT3}. As $\nabla\psi^{*}$ is a homeomorphism from $\mathrm{Int}(\Y)$ to $\mathrm{Int}(\s)$, hence, $(\nabla\psi^{*})^{-1}$ is a homeomorphism from $\mathrm{Int}(\s)$ to $\mathrm{Int}(\Y)$. This implies $(\nabla\psi^{*})^{-1}$ maps compact subsets of $\mathrm{Int}(\s)$ to compact subsets of $\mathrm{Int}(\Y)$. Therefore, $\mathfrak{J}:=(\nabla\psi^{*})^{-1}(\mathrm{Cl}(U))$ is a compact subset of $\mathrm{Int}(\Y)$ and $$\mathfrak{K}\subset U=\mathrm{Int}(\nabla\psi^{*}(\mathfrak{J})).$$ 
 This completes the proof of the claim. \qed

 \subsection{Proof of Claim~\ref{cl:MT4}}\label{sec:cl:MT4}
  Since $\mathfrak{K}\subset \mathrm{Int}(\s)$ is a compact set, there exists $\delta_0>0$ such that $\mathrm{Cl}\big(\mathfrak{K}+B_{\delta_0}(0)\big)\subset \mathrm{Int}(\mathcal{S})$. Owing to Claim~\ref{cl:MT3}, there exists a compact set $\mathfrak{J}\subset \mathrm{Int}(\Y)$ such that $\mathrm{Cl}\big(\mathfrak{K}+B_{\delta_0}(0)\big)\subset \mathrm{Int}(\nabla \psi^{*}(\mathfrak{J}))$. Let $U$ be an open set and $\mathfrak{W}$ be a compact set in $\mathrm{Int}(\Y)$ such that $$\mathfrak{J}\subset U \subset \mathrm{Cl}(U)\subset \mathrm{Int}(\mathfrak{W}).$$  By the construction in \textbf{Step I}, $\hat{\xi}_{p}$ converges to $\psi^{*}$ uniformly on any compact subset of $\RR^d$, as $p\to \infty$. Moreover, we have $\mathrm{Cl}\big(\mathfrak{K}+B_{\delta_0}(0)\big)\subset \nabla \psi^{*}(U)$.
 Due to the uniform convergence of $\{\hat{\xi}_{p}\}_{p\ge 1}$ to $\psi^*$ in $\mathfrak{W}$, the strict convexity of $\psi^{*}$ in $U$ (since $\nabla \psi^{*}$ is a homeomorphism in $\mathrm{Int}(\Y)$) and Lemma~\ref{KeyLemma}-(d), there exists $p_0=p_0(\mathfrak{K}) \in \N$ such that for all $p\geq p_0$ we have $\mathfrak{K}\subset \partial \hat{\xi}_{p}(\mathfrak{W})$. But, this does not ensure that $\partial \hat{\xi}^{*}_{p}(\mathfrak{K})\subset \mathfrak{W}$. Note that $\partial\hat{\xi}^{*}_p(\mathrm{Int}(\mathcal{S}))\subset \mathrm{Int}(\mathcal{Y})$ w.p.~$1$. Hence, to complete the proof, due to Lemma~\ref{lem:SubD} and the fact that $\mathfrak{K}\subset \mathrm{Int}(\s)$, one needs to further show that there exists $p_0$ such that $\mathfrak{K}\cap \partial \hat{\xi}_{p}(\Y\backslash \mathfrak{W}) = \emptyset$ for all $p\geq p_0$. This will be showed in the rest of the proof. To this end, we will first find $p_0$ such that $\mathfrak{K}$ will be embedded inside the $\delta$-core of $\partial \hat{\xi}_{p}(\mathfrak{W})$, for some $\delta>0$, for all $p\geq p_0$, where for any set $B\subset \RR^d$, the $\delta$-core of $B$ is defined as  
   \begin{align}
   B^{\delta} := \Big\{y\in B : y+ B_{\delta}(0) \subset B \Big\}. 
\end{align}  
 Once again owing to Claim~\ref{cl:MT3}, the uniform convergence of $\hat{\xi}_{p}$ to $\psi^*$ in $\mathfrak{W}$, the strict convexity of $\psi^{*}$ in $U$ and Lemma~\ref{KeyLemma}-(d), there exists $p_0$ such that $\nabla\psi^{*}(U)\subset \partial \hat{\xi}_{p}(\mathfrak{W})$ for all $p\geq p_0$. Since, $\mathrm{Cl}(\mathfrak{K}+B_{\delta_0}(0))\subset \nabla\psi^{*}(U)$, therefore, $\mathfrak{K}\subset (\nabla\psi^{*}(U))^{\delta_0}\subset (\partial \hat{\xi}_{p}(\mathfrak{W}))^{\delta_0}$ for all $p\geq p_0$.   
   
   Fix $R>0$ such that $\mathfrak{W}+B_R(0)\subset \mathrm{Int}(\Y)$. Define $\mathcal{G}:=\mathrm{Cl}\Big((\mathfrak{W}+B_R(0))\backslash \mathfrak{W}\Big).$   Owing to Lemma~\ref{KeyLemma}-(a), there exists $p_1=p_1(\delta_0)$ such that,     for all $p\geq p_1$, 
   \begin{align}\label{eq:Incl1}
   \partial \hat{\xi}_{p}(\mathcal{G})\subset \nabla \psi^{*}(\mathcal{G})+ B_{\delta_0/2}(0).
   \end{align}
Since $\nabla \psi^{*}$ is an injective map in $\mathrm{Int}(\Y)$ (as it is a homeomorphism), we have 
    \begin{align}\label{eq:Incl2}
     \big(\nabla \psi^{*}(\mathcal{G})+ B_{\delta_0/2}(0)\big)\cap (\nabla\psi^{*}(\mathfrak{W}))^{\delta_0} =\emptyset.
     \end{align}
      On the other hand, we know  
   \begin{align}\label{eq:Incl3}
   \mathfrak{K}\subset (\nabla \psi^{*}(U))^{\delta_0}\subset (\nabla \psi^{*}(\mathfrak{W}))^{\delta_0} .
    \end{align}   
     Combining \eqref{eq:Incl1}, \eqref{eq:Incl2} and \eqref{eq:Incl3}, we have $\mathfrak{K}\cap \partial \hat{\xi}_{p}(\mathcal{G}) = \emptyset$ for all $p\geq \max\{p_0 , p_1\}$. Now, we claim and prove that $\mathfrak{K}\cap \partial \hat{\xi}_{p}(\Y \backslash (\mathfrak{W}+B_R(0)))= \emptyset$ for all $p\geq \max\{p_0 , p_1\}$. Suppose that there exist $x\in \mathfrak{W}$, $y\in \big(\Y\backslash (\mathfrak{W}+B_R(0))\big)$ and $z\in \mathfrak{K}$ such that $z\in \partial \hat{\xi}_{p}(x)\cap \partial \hat{\xi}_{p}(y)$. Then, $z$ belongs to the subdifferential set of the functions $\hat{\xi}_{p}$ at any point in the line segment joining $x$ and $y$. However, the line joining $x$ and $y$ passes through $\mathcal{G}$ and we know that $\mathfrak{K}\cap \partial\hat{\xi}_{p}(\mathcal{G}) = \emptyset$. This contradicts the existence of $y$. Hence, the claim follows. Combining this claim with  $\mathfrak{K}\cap \partial\hat{\xi}_{p}(\mathcal{G}) = \emptyset$ yields
     \begin{align}
     \Big(\mathfrak{K}\cap \partial \hat{\xi}_{p}(\Y\backslash \mathfrak{W})\Big)&\subset \mathfrak{K}\cap \Big(\partial \hat{\xi}_{p}\big(\Y\backslash (\mathfrak{W}+B_{R}(0))\big)\bigcup \partial\hat{\xi}_{p}(\mathcal{G})\Big) \\ 
     &= \Big(\mathfrak{K}\cap \partial \hat{\xi}_p\big(\Y\backslash (\mathfrak{W}+B_{R}(0))\big)\Big)\bigcup \Big(\mathfrak{K}\cap \partial\hat{\xi}_{p}(\mathcal{G})\Big)\\
     &= \emptyset, \qquad \forall p\geq \max\{p_0, p_1\}.
     \end{align} 
     Owing to this, $\partial \hat{\xi}^{*}_{p}(\mathfrak{K})\subset \mathfrak{W}$ for all $p\geq \max\{p_0, p_1\}$. This completes the proof. \qed

  \subsection{Proof of Claim~\ref{cl:MT5}}\label{sec:cl:MT5}
    Recall that $
  \nabla \psi^{*}(Q(x_0))=x_0 $ (because $\nabla \psi^{*}=Q^{-1}$ in $\mathrm{Int}(\Y)$ via Claim~\ref{cl:MT2}) and $B_{\delta_0}(x_0)\subset \mathfrak{K}$ where $\mathfrak{K}$ is a compact set in $\mathrm{Int}(\s)$. Owing to the 
  uniform convergence of $\{\partial\hat{\xi}_{p}\}_{p\ge 1}$ to $\nabla\psi^{*}$ (by \textbf{Step I} and Lemma~\ref{KeyLemma}-(c)), for any $
  \delta>0$ there exists $p_0=p_0(\delta)$ such that 
  \begin{equation}\label{eq:UConCons}
  \sup_{y\in \partial \hat{\xi}_{p}(Q(x_0))}\|y- x_0\|\leq \delta, \qquad\forall p\geq p_0.
  \end{equation}
   Furthermore, for any $p \ge 1$, from the definition of subdifferentials, there exists $y\in\partial \hat{\xi}_{p}(Q(x_0))$ such that 
  \begin{align}\label{eq:Duality}
  \hat{\xi}_{p}(Q(x_0)) + \hat{\xi}^{*}_{p}(y) = \langle y, Q(x_0)\rangle. 
\end{align}    
Note that~\eqref{eq:UConCons} shows that for any $\delta>0$ there exists $p_0=p_0(\delta)$ such that $y\in B_{\delta}(x_0)$. Combining this with~\eqref{eq:Duality} and the fact that $  \hat{\xi}_{p}(Q(x_0)) =0$, we have $|\hat{\xi}^{*}_{p}(y)|\leq (\|x_0\|+\delta)\|Q(x_0)\|$. \qed

\subsection{Proof of Claim~\ref{cl:MT6}}\label{sec:cl:MT6}
  By Claim~\ref{cl:MT2} and Lemma~\ref{lem:CvxEq}, for all $\omega\in \mathcal{A}_{K}\cap \mathfrak{D}$, $\nabla\psi^{*}:\mathrm{Int}(\mathcal{Y})\to\mathrm{Int}(\mathcal{S})$ is a homeomorphism and $\nabla\psi^{*}=R$ everywhere in $\mathrm{Int}(\mathcal{Y})$. Since $\nabla \widetilde{\psi} = (\nabla \psi^{*})^{-1}$ in $\mathrm{Int}(\mathcal{S})$ 
 for all $\omega \in \mathcal{A}_K\cap \mathfrak{D}$ via \textbf{Step III} and $Q=R^{-1}$ in $\mathrm{Int}(\s)$, thus, $\nabla\tilde{\psi}:\mathrm{Int}(\mathcal{S})\to\mathrm{Int}(\mathcal{Y}) $ is a homeomorphism and $\nabla\tilde{\psi}=Q$ is same as $Q$ in $\mathrm{Int}(\mathcal{S})$ 
 for all $\omega \in \mathcal{A}_K\cap \mathfrak{D}$. This completes the proof of the first part of the claim.

 Now, we turn to show that $\nabla\psi^{*}=R$ everywhere in $\RR^d$. Recall from \textbf{Step I} that $\hat{\xi}_p$ converges uniformly to $\psi^{*}$ on compacts of $\RR^d$, as $p\to \infty$. Since $\partial\hat{\xi}_p(\RR^d)\subset \mathcal{S}$, by Lemma~\ref{KeyLemma}-(b), $\partial \psi^{*}(\RR^d)\subset \mathcal{S}$. Note that $\psi$ is the Legendre-Fenchel dual of $\psi^{*}$. From the definition of the Legendre-Fenchel dual of a convex function, $$\{u:\psi(u)<\infty\}= \{y\in \RR^d:y\in \partial \psi^{*}(x) \text{ for some }x\in \RR^d\}.$$     
Hence, $\psi(u)=+\infty$ for all $u\in \RR^d \backslash \mathcal{S}.$ 
 Since $\nabla \psi^{*}$ is a homeomorphism from $\mathrm{Int}(\Y)$ to $\mathrm{Int}(\s)$, by Lemma~\ref{lem:SubD}, $\nabla\psi (u)$ is equal to $(\nabla \psi^{*})^{-1}(u)$  and hence, equal to $Q(u)$ for all $u\in \mathrm{Int}(\s)$. Let $\phi$ be the convex function such that $\nabla\phi (u)=Q(u)$ for all $u\in \mathrm{Int}(\s)$. Then $\psi-\phi$ is a differentiable function in $\mathrm{Int}(\s)$ with $\nabla(\psi-\phi)(u)=0$ for all $u\in \mathrm{Int}(\s)$. This implies $\psi-\phi =c$ for some constant $c$ in $\mathrm{Int}(\s)$. As both $\psi$ and $\phi$ are l.s.c., we have $\phi-\psi=c$ in $\s$. Furthermore, from the definition of the quantile map, we know $\phi=+\infty$ in $\RR^d \backslash \s$. Let $\phi^{*}$ be the Legendre-Fenchel dual of $\phi$. Now, by applying Lemma~\ref{lem:Affinenv}, we get $\psi^{*}(y)-\phi^{*}(y)=-c$ for all $y\in \RR^d$. As a consequence, $\partial \phi^{*}=\partial \psi^{*}$ in $\RR^d$. Since $\nabla\phi =Q$, we know that $R= \nabla \phi^{*}$. This implies $R=\nabla\psi^{*}$ everywhere in $\RR^d$. \qed

\bibliographystyle{imsart-number}
\bibliography{OT}

\begin{thebibliography}{107}

\bibitem{AKT}
\begin{barticle}[author]
\bauthor{\bsnm{Ajtai},~\bfnm{M.}\binits{M.}},
  \bauthor{\bsnm{Koml\'{o}s},~\bfnm{J.}\binits{J.}} \AND
  \bauthor{\bsnm{Tusn\'{a}dy},~\bfnm{G.}\binits{G.}}
(\byear{1984}).
\btitle{On optimal matchings}.
\bjournal{Combinatorica}
\bvolume{4}
\bpages{259--264}.
\end{barticle}
\endbibitem

\bibitem{Ambro08}
\begin{bbook}[author]
\bauthor{\bsnm{Ambrosio},~\bfnm{Luigi}\binits{L.}},
  \bauthor{\bsnm{Gigli},~\bfnm{Nicola}\binits{N.}} \AND
  \bauthor{\bsnm{Savar\'{e}},~\bfnm{Giuseppe}\binits{G.}}
(\byear{2008}).
\btitle{Gradient flows in metric spaces and in the space of probability
  measures},
\bedition{second} ed.
\bseries{Lectures in Mathematics ETH Z\"{u}rich}.
\bpublisher{Birkh\"{a}user Verlag, Basel}.
\bmrnumber{2401600}
\end{bbook}
\endbibitem

\bibitem{Auren87}
\begin{barticle}[author]
\bauthor{\bsnm{Aurenhammer},~\bfnm{F.}\binits{F.}}
(\byear{1987}).
\btitle{Power diagrams: properties, algorithms and applications}.
\bjournal{SIAM J. Comput.}
\bvolume{16}
\bpages{78--96}.
\bdoi{10.1137/0216006}
\bmrnumber{873251}
\end{barticle}
\endbibitem

\bibitem{Auren98}
\begin{barticle}[author]
\bauthor{\bsnm{Aurenhammer},~\bfnm{Franz}\binits{F.}},
  \bauthor{\bsnm{Hoffmann},~\bfnm{Friedrich}\binits{F.}} \AND
  \bauthor{\bsnm{Aronov},~\bfnm{Boris}\binits{B.}}
(\byear{1998}).
\btitle{Minkowski-type theorems and least-squares clustering}.
\bjournal{Algorithmica}
\bvolume{20}
\bpages{61--76}.
\end{barticle}
\endbibitem

\bibitem{Auren13}
\begin{bbook}[author]
\bauthor{\bsnm{Aurenhammer},~\bfnm{Franz}\binits{F.}},
  \bauthor{\bsnm{Klein},~\bfnm{Rolf}\binits{R.}} \AND
  \bauthor{\bsnm{Lee},~\bfnm{Der-Tsai}\binits{D.-T.}}
(\byear{2013}).
\btitle{Voronoi diagrams and {D}elaunay triangulations}.
\bpublisher{World Scientific Publishing Co. Pte. Ltd., Hackensack, NJ}.
\bdoi{10.1142/8685}
\bmrnumber{3186045}
\end{bbook}
\endbibitem

\bibitem{BF04}
\begin{barticle}[author]
\bauthor{\bsnm{Baringhaus},~\bfnm{L.}\binits{L.}} \AND
  \bauthor{\bsnm{Franz},~\bfnm{C.}\binits{C.}}
(\byear{2004}).
\btitle{On a new multivariate two-sample test}.
\bjournal{J. Multivariate Anal.}
\bvolume{88}
\bpages{190--206}.
\bdoi{10.1016/S0047-259X(03)00079-4}
\bmrnumber{2021870}
\end{barticle}
\endbibitem

\bibitem{Berman2018}
\begin{barticle}[author]
\bauthor{\bsnm{Berman},~\bfnm{Robert~J}\binits{R.~J.}}
(\byear{2018}).
\btitle{Convergence rates for discretized {M}onge-{A}mp\`{e}re equations and
  quantitative stability of optimal transport}.
\bjournal{arXiv preprint arXiv:1803.00785}.
\end{barticle}
\endbibitem

\bibitem{Berrett19}
\begin{barticle}[author]
\bauthor{\bsnm{Berrett},~\bfnm{T.~B.}\binits{T.~B.}} \AND
  \bauthor{\bsnm{Samworth},~\bfnm{R.~J.}\binits{R.~J.}}
(\byear{2019}).
\btitle{Nonparametric independence testing via mutual information}.
\bjournal{Biometrika}
\bvolume{106}
\bpages{547--566}.
\bdoi{10.1093/biomet/asz024}
\bmrnumber{3992389}
\end{barticle}
\endbibitem

\bibitem{B15}
\begin{barticle}[author]
\bauthor{\bsnm{Bhattacharya},~\bfnm{Bhaswar~B}\binits{B.~B.}}
(\byear{2015}).
\btitle{Two-Sample Tests Based on Geometric Graphs: Asymptotic Distribution and
  Detection Thresholds}.
\bjournal{arXiv preprint arXiv:1512.00384}.
\end{barticle}
\endbibitem

\bibitem{Bickel68}
\begin{barticle}[author]
\bauthor{\bsnm{Bickel},~\bfnm{P.~J.}\binits{P.~J.}}
(\byear{1968}).
\btitle{A distribution free version of the {S}mirnov two sample test in the
  {$p$}-variate case}.
\bjournal{Ann. Math. Statist.}
\bvolume{40}
\bpages{1--23}.
\bdoi{10.1214/aoms/1177697800}
\bmrnumber{0256519}
\end{barticle}
\endbibitem

\bibitem{Blomqvist50}
\begin{barticle}[author]
\bauthor{\bsnm{Blomqvist},~\bfnm{Nils}\binits{N.}}
(\byear{1950}).
\btitle{On a measure of dependence between two random variables}.
\bjournal{Ann. Math. Statistics}
\bvolume{21}
\bpages{593--600}.
\bdoi{10.1214/aoms/1177729754}
\bmrnumber{0039190}
\end{barticle}
\endbibitem

\bibitem{BlumEtAl61}
\begin{barticle}[author]
\bauthor{\bsnm{Blum},~\bfnm{J.~R.}\binits{J.~R.}},
  \bauthor{\bsnm{Kiefer},~\bfnm{J.}\binits{J.}} \AND
  \bauthor{\bsnm{Rosenblatt},~\bfnm{M.}\binits{M.}}
(\byear{1961}).
\btitle{Distribution free tests of independence based on the sample
  distribution function}.
\bjournal{Ann. Math. Statist.}
\bvolume{32}
\bpages{485--498}.
\bdoi{10.1214/aoms/1177705055}
\bmrnumber{0125690}
\end{barticle}
\endbibitem

\bibitem{BL19}
\begin{barticle}[author]
\bauthor{\bsnm{Bobkov},~\bfnm{Sergey}\binits{S.}} \AND
  \bauthor{\bsnm{Ledoux},~\bfnm{Michel}\binits{M.}}
(\byear{2019}).
\btitle{One-dimensional empirical measures, order statistics, and {K}antorovich
  transport distances}.
\bjournal{Mem. Amer. Math. Soc.}
\bvolume{261}
\bpages{v+126}.
\bdoi{10.1090/memo/1259}
\bmrnumber{4028181}
\end{barticle}
\endbibitem

\bibitem{BSS18}
\begin{barticle}[author]
\bauthor{\bsnm{Boeckel},~\bfnm{Melf}\binits{M.}},
  \bauthor{\bsnm{Spokoiny},~\bfnm{Vladimir}\binits{V.}} \AND
  \bauthor{\bsnm{Suvorikova},~\bfnm{Alexandra}\binits{A.}}
(\byear{2018}).
\btitle{Multivariate Brenier cumulative distribution functions and their
  application to non-parametric testing}.
\bjournal{arXiv preprint arXiv:1809.04090}.
\end{barticle}
\endbibitem

\bibitem{B91}
\begin{barticle}[author]
\bauthor{\bsnm{Brenier},~\bfnm{Yann}\binits{Y.}}
(\byear{1991}).
\btitle{Polar factorization and monotone rearrangement of vector-valued
  functions}.
\bjournal{Comm. Pure Appl. Math.}
\bvolume{44}
\bpages{375--417}.
\bdoi{10.1002/cpa.3160440402}
\bmrnumber{1100809}
\end{barticle}
\endbibitem

\bibitem{Cafa90}
\begin{barticle}[author]
\bauthor{\bsnm{Caffarelli},~\bfnm{Luis~A.}\binits{L.~A.}}
(\byear{1990}).
\btitle{Interior {$W^{2,p}$} estimates for solutions of the {M}onge-{A}mp\`ere
  equation}.
\bjournal{Ann. of Math. (2)}
\bvolume{131}
\bpages{135--150}.
\bdoi{10.2307/1971510}
\bmrnumber{1038360}
\end{barticle}
\endbibitem

\bibitem{Ca1}
\begin{barticle}[author]
\bauthor{\bsnm{Caffarelli},~\bfnm{Luis~A.}\binits{L.~A.}}
(\byear{1992}).
\btitle{The regularity of mappings with a convex potential}.
\bjournal{J. Amer. Math. Soc.}
\bvolume{5}
\bpages{99--104}.
\end{barticle}
\endbibitem

\bibitem{Ca2}
\begin{barticle}[author]
\bauthor{\bsnm{Caffarelli},~\bfnm{Luis~A.}\binits{L.~A.}}
(\byear{1992}).
\btitle{Boundary regularity of maps with convex potentials}.
\bjournal{Comm. Pure Appl. Math.}
\bvolume{45}
\bpages{1141--1151}.
\end{barticle}
\endbibitem

\bibitem{Ca3}
\begin{barticle}[author]
\bauthor{\bsnm{Caffarelli},~\bfnm{Luis~A.}\binits{L.~A.}}
(\byear{1996}).
\btitle{Boundary regularity of maps with convex potentials. {II}}.
\bjournal{Ann. of Math. (2)}
\bvolume{144}
\bpages{453--496}.
\end{barticle}
\endbibitem

\bibitem{Caffarelli99}
\begin{bincollection}[author]
\bauthor{\bsnm{Caffarelli},~\bfnm{Luis~A.}\binits{L.~A.}},
  \bauthor{\bsnm{Kochengin},~\bfnm{Sergey~A.}\binits{S.~A.}} \AND
  \bauthor{\bsnm{Oliker},~\bfnm{Vladimir~I.}\binits{V.~I.}}
(\byear{1999}).
\btitle{On the numerical solution of the problem of reflector design with given
  far-field scattering data}.
In \bbooktitle{Monge {A}mp\`ere equation: applications to geometry and
  optimization ({D}eerfield {B}each, {FL}, 1997)}.
\bseries{Contemp. Math.}
\bvolume{226}
\bpages{13--32}.
\bpublisher{Amer. Math. Soc., Providence, RI}.
\bdoi{10.1090/conm/226/03233}
\bmrnumber{1660740}
\end{bincollection}
\endbibitem

\bibitem{Chaud96}
\begin{barticle}[author]
\bauthor{\bsnm{Chaudhuri},~\bfnm{Probal}\binits{P.}}
(\byear{1996}).
\btitle{On a geometric notion of quantiles for multivariate data}.
\bjournal{J. Amer. Statist. Assoc.}
\bvolume{91}
\bpages{862--872}.
\bdoi{10.2307/2291681}
\bmrnumber{1395753}
\end{barticle}
\endbibitem

\bibitem{Cher17}
\begin{barticle}[author]
\bauthor{\bsnm{Chernozhukov},~\bfnm{Victor}\binits{V.}},
  \bauthor{\bsnm{Galichon},~\bfnm{Alfred}\binits{A.}},
  \bauthor{\bsnm{Hallin},~\bfnm{Marc}\binits{M.}} \AND
  \bauthor{\bsnm{Henry},~\bfnm{Marc}\binits{M.}}
(\byear{2017}).
\btitle{Monge-{K}antorovich depth, quantiles, ranks and signs}.
\bjournal{Ann. Statist.}
\bvolume{45}
\bpages{223--256}.
\bdoi{10.1214/16-AOS1450}
\bmrnumber{3611491}
\end{barticle}
\endbibitem

\bibitem{Chizat2020}
\begin{barticle}[author]
\bauthor{\bsnm{Chizat},~\bfnm{Lenaic}\binits{L.}},
  \bauthor{\bsnm{Roussillon},~\bfnm{Pierre}\binits{P.}},
  \bauthor{\bsnm{L{\'e}ger},~\bfnm{Flavien}\binits{F.}},
  \bauthor{\bsnm{Vialard},~\bfnm{Fran{\c{c}}ois-Xavier}\binits{F.-X.}} \AND
  \bauthor{\bsnm{Peyr{\'e}},~\bfnm{Gabriel}\binits{G.}}
(\byear{2020}).
\btitle{Faster {W}asserstein Distance Estimation with the Sinkhorn Divergence}.
\bjournal{Advances in Neural Information Processing Systems}
\bvolume{33}.
\end{barticle}
\endbibitem

\bibitem{CF19}
\begin{barticle}[author]
\bauthor{\bsnm{Cordero-Erausquin},~\bfnm{Dario}\binits{D.}} \AND
  \bauthor{\bsnm{Figalli},~\bfnm{Alessio}\binits{A.}}
(\byear{2019}).
\btitle{Regularity of monotone transport maps between unbounded domains}.
\bjournal{Discrete Contin. Dyn. Syst.}
\bvolume{39}
\bpages{7101--7112}.
\bdoi{10.3934/dcds.2019297}
\bmrnumber{4026183}
\end{barticle}
\endbibitem

\bibitem{Cuesta-2013}
\begin{barticle}[author]
\bauthor{\bsnm{Cuesta-Albertos},~\bfnm{J.~A.}\binits{J.~A.}},
  \bauthor{\bsnm{R\"{u}schendorf},~\bfnm{L.}\binits{L.}} \AND
  \bauthor{\bsnm{Tuero-D\'{\i}az},~\bfnm{A.}\binits{A.}}
(\byear{1993}).
\btitle{Optimal coupling of multivariate distributions and stochastic
  processes}.
\bjournal{J. Multivariate Anal.}
\bvolume{46}
\bpages{335--361}.
\bdoi{10.1006/jmva.1993.1064}
\bmrnumber{1240428}
\end{barticle}
\endbibitem

\bibitem{PF13}
\begin{barticle}[author]
\bauthor{\bsnm{De~Philippis},~\bfnm{Guido}\binits{G.}} \AND
  \bauthor{\bsnm{Figalli},~\bfnm{Alessio}\binits{A.}}
(\byear{2013}).
\btitle{{$W^{2,1}$} regularity for solutions of the {M}onge-{A}mp\`ere
  equation}.
\bjournal{Invent. Math.}
\bvolume{192}
\bpages{55--69}.
\bdoi{10.1007/s00222-012-0405-4}
\bmrnumber{3032325}
\end{barticle}
\endbibitem

\bibitem{PF14}
\begin{barticle}[author]
\bauthor{\bsnm{De~Philippis},~\bfnm{Guido}\binits{G.}} \AND
  \bauthor{\bsnm{Figalli},~\bfnm{Alessio}\binits{A.}}
(\byear{2014}).
\btitle{The {M}onge-{A}mp\`ere equation and its link to optimal
  transportation}.
\bjournal{Bull. Amer. Math. Soc. (N.S.)}
\bvolume{51}
\bpages{527--580}.
\bdoi{10.1090/S0273-0979-2014-01459-4}
\bmrnumber{3237759}
\end{barticle}
\endbibitem

\bibitem{PF15}
\begin{barticle}[author]
\bauthor{\bsnm{De~Philippis},~\bfnm{Guido}\binits{G.}} \AND
  \bauthor{\bsnm{Figalli},~\bfnm{Alessio}\binits{A.}}
(\byear{2015}).
\btitle{Partial regularity for optimal transport maps}.
\bjournal{Publ. Math. Inst. Hautes \'{E}tudes Sci.}
\bvolume{121}
\bpages{81--112}.
\bdoi{10.1007/s10240-014-0064-7}
\bmrnumber{3349831}
\end{barticle}
\endbibitem

\bibitem{VS18}
\begin{barticle}[author]
\bauthor{\bparticle{de} \bsnm{Valk},~\bfnm{Cees}\binits{C.}} \AND
  \bauthor{\bsnm{Segers},~\bfnm{Johan}\binits{J.}}
(\byear{2018}).
\btitle{Tails of optimal transport plans for regularly varying probability
  measures}.
\bjournal{arXiv preprint arXiv:1811.12061}.
\end{barticle}
\endbibitem

\bibitem{deb2021efficiency}
\begin{barticle}[author]
\bauthor{\bsnm{Deb},~\bfnm{Nabarun}\binits{N.}},
  \bauthor{\bsnm{Bhattacharya},~\bfnm{Bhaswar~B}\binits{B.~B.}} \AND
  \bauthor{\bsnm{Sen},~\bfnm{Bodhisattva}\binits{B.}}
(\byear{2021}).
\btitle{Efficiency Lower Bounds for Distribution-Free Hotelling-Type Two-Sample
  Tests Based on Optimal Transport}.
\bjournal{arXiv preprint arXiv:2104.01986}.
\end{barticle}
\endbibitem

\bibitem{DebSen2019}
\begin{barticle}[author]
\bauthor{\bsnm{Deb},~\bfnm{Nabarun}\binits{N.}} \AND
  \bauthor{\bsnm{Sen},~\bfnm{Bodhisattva}\binits{B.}}
(\byear{2019}).
\btitle{Multivariate rank-based distribution-free nonparametric testing using
  measure transportation}.
\bjournal{arXiv preprint arXiv:1909.08733}.
\end{barticle}
\endbibitem

\bibitem{D14}
\begin{barticle}[author]
\bauthor{\bsnm{Decurninge},~\bfnm{Alexis}\binits{A.}}
(\byear{2014}).
\btitle{Multivariate quantiles and multivariate L-moments}.
\bjournal{arXiv preprint arXiv:1409.6013}.
\end{barticle}
\endbibitem

\bibitem{BGU05}
\begin{barticle}[author]
\bauthor{\bparticle{del} \bsnm{Barrio},~\bfnm{Eustasio}\binits{E.}},
  \bauthor{\bsnm{Gin\'{e}},~\bfnm{Evarist}\binits{E.}} \AND
  \bauthor{\bsnm{Utzet},~\bfnm{Frederic}\binits{F.}}
(\byear{2005}).
\btitle{Asymptotics for {$L_2$} functionals of the empirical quantile process,
  with applications to tests of fit based on weighted {W}asserstein distances}.
\bjournal{Bernoulli}
\bvolume{11}
\bpages{131--189}.
\bdoi{10.3150/bj/1110228245}
\bmrnumber{2121458}
\end{barticle}
\endbibitem

\bibitem{delBarrio2020}
\begin{barticle}[author]
\bauthor{\bparticle{del} \bsnm{Barrio},~\bfnm{Eustasio}\binits{E.}},
  \bauthor{\bsnm{Gonz\'{a}lez-Sanz},~\bfnm{Alberto}\binits{A.}} \AND
  \bauthor{\bsnm{Hallin},~\bfnm{Marc}\binits{M.}}
(\byear{2020}).
\btitle{A note on the regularity of optimal-transport-based center-outward
  distribution and quantile functions}.
\bjournal{J. Multivariate Anal.}
\bvolume{180}
\bpages{104671, 13}.
\bdoi{10.1016/j.jmva.2020.104671}
\bmrnumber{4147635}
\end{barticle}
\endbibitem

\bibitem{delB-JMA-19}
\begin{barticle}[author]
\bauthor{\bparticle{del} \bsnm{Barrio},~\bfnm{Eustasio}\binits{E.}},
  \bauthor{\bsnm{Gordaliza},~\bfnm{Paula}\binits{P.}},
  \bauthor{\bsnm{Lescornel},~\bfnm{H\'{e}l\`ene}\binits{H.}} \AND
  \bauthor{\bsnm{Loubes},~\bfnm{Jean-Michel}\binits{J.-M.}}
(\byear{2019}).
\btitle{Central limit theorem and bootstrap procedure for {W}asserstein's
  variations with an application to structural relationships between
  distributions}.
\bjournal{J. Multivariate Anal.}
\bvolume{169}
\bpages{341--362}.
\bdoi{10.1016/j.jmva.2018.09.014}
\bmrnumber{3875604}
\end{barticle}
\endbibitem

\bibitem{delB-AoP-19}
\begin{barticle}[author]
\bauthor{\bparticle{del} \bsnm{Barrio},~\bfnm{Eustasio}\binits{E.}} \AND
  \bauthor{\bsnm{Loubes},~\bfnm{Jean-Michel}\binits{J.-M.}}
(\byear{2019}).
\btitle{Central limit theorems for empirical transportation cost in general
  dimension}.
\bjournal{Ann. Probab.}
\bvolume{47}
\bpages{926--951}.
\bdoi{10.1214/18-AOP1275}
\bmrnumber{3916938}
\end{barticle}
\endbibitem

\bibitem{del2019}
\begin{barticle}[author]
\bauthor{\bsnm{Del~Barrio},~\bfnm{Eustasio}\binits{E.}} \AND
  \bauthor{\bsnm{Loubes},~\bfnm{Jean-Michel}\binits{J.-M.}}
(\byear{2019}).
\btitle{Central limit theorems for empirical transportation cost in general
  dimension}.
\bjournal{The Annals of Probability}
\bvolume{47}
\bpages{926--951}.
\end{barticle}
\endbibitem

\bibitem{Dudley}
\begin{bbook}[author]
\bauthor{\bsnm{Dudley},~\bfnm{R.~M.}\binits{R.~M.}}
(\byear{2002}).
\btitle{Real analysis and probability}.
\bpublisher{Cambridge University Press, New York}.
\end{bbook}
\endbibitem

\bibitem{F10}
\begin{barticle}[author]
\bauthor{\bsnm{Figalli},~\bfnm{Alessio}\binits{A.}}
(\byear{2010}).
\btitle{Regularity properties of optimal maps between nonconvex domains in the
  plane}.
\bjournal{Comm. Partial Differential Equations}
\bvolume{35}
\bpages{465--479}.
\bdoi{10.1080/03605300903307673}
\bmrnumber{2748633}
\end{barticle}
\endbibitem

\bibitem{Figalli}
\begin{bbook}[author]
\bauthor{\bsnm{Figalli},~\bfnm{A.}\binits{A.}}
(\byear{2017}).
\btitle{The Monge-Amp\'ere equation and its application}.
\bseries{Zurich lectures in advanced mathematics}
\bvolume{58}.
\bpublisher{European Mathematical Society}.
\end{bbook}
\endbibitem

\bibitem{Figalli18}
\begin{barticle}[author]
\bauthor{\bsnm{Figalli},~\bfnm{Alessio}\binits{A.}}
(\byear{2018}).
\btitle{On the continuity of center-outward distribution and quantile
  functions}.
\bjournal{Nonlinear Anal.}
\bvolume{177}
\bpages{413--421}.
\bdoi{10.1016/j.na.2018.05.008}
\bmrnumber{3886582}
\end{barticle}
\endbibitem

\bibitem{FKM09}
\begin{barticle}[author]
\bauthor{\bsnm{Figalli},~\bfnm{Alessio}\binits{A.}},
  \bauthor{\bsnm{Kim},~\bfnm{Young-Heon}\binits{Y.-H.}} \AND
  \bauthor{\bsnm{McCann},~\bfnm{Robert~J.}\binits{R.~J.}}
(\byear{2013}).
\btitle{H\"{o}lder continuity and injectivity of optimal maps}.
\bjournal{Arch. Ration. Mech. Anal.}
\bvolume{209}
\bpages{747--795}.
\bdoi{10.1007/s00205-013-0629-5}
\bmrnumber{3067826}
\end{barticle}
\endbibitem

\bibitem{FRV11}
\begin{barticle}[author]
\bauthor{\bsnm{Figalli},~\bfnm{Alessio}\binits{A.}},
  \bauthor{\bsnm{Rifford},~\bfnm{Ludovic}\binits{L.}} \AND
  \bauthor{\bsnm{Villani},~\bfnm{C\'{e}dric}\binits{C.}}
(\byear{2011}).
\btitle{Necessary and sufficient conditions for continuity of optimal transport
  maps on {R}iemannian manifolds}.
\bjournal{Tohoku Math. J. (2)}
\bvolume{63}
\bpages{855--876}.
\bdoi{10.2748/tmj/1325886291}
\bmrnumber{2872966}
\end{barticle}
\endbibitem

\bibitem{FG15}
\begin{barticle}[author]
\bauthor{\bsnm{Fournier},~\bfnm{N.}\binits{N.}} \AND
  \bauthor{\bsnm{Guillin},~\bfnm{A.}\binits{A.}}
(\byear{2015}).
\btitle{On the rate of convergence in {W}asserstein distance of the empirical
  measure}.
\bjournal{Probab. Theory Related Fields}
\bvolume{162}
\bpages{707--738}.
\bdoi{10.1007/s00440-014-0583-7}
\bmrnumber{3383341}
\end{barticle}
\endbibitem

\bibitem{FR79}
\begin{barticle}[author]
\bauthor{\bsnm{Friedman},~\bfnm{Jerome~H.}\binits{J.~H.}} \AND
  \bauthor{\bsnm{Rafsky},~\bfnm{Lawrence~C.}\binits{L.~C.}}
(\byear{1979}).
\btitle{Multivariate generalizations of the {W}ald-{W}olfowitz and {S}mirnov
  two-sample tests}.
\bjournal{Ann. Statist.}
\bvolume{7}
\bpages{697--717}.
\bmrnumber{532236}
\end{barticle}
\endbibitem

\bibitem{Galichon16}
\begin{bbook}[author]
\bauthor{\bsnm{Galichon},~\bfnm{Alfred}\binits{A.}}
(\byear{2016}).
\btitle{Optimal transport methods in economics}.
\bpublisher{Princeton University Press, Princeton, NJ}.
\bdoi{10.1515/9781400883592}
\bmrnumber{3586373}
\end{bbook}
\endbibitem

\bibitem{Gangbo99}
\begin{bincollection}[author]
\bauthor{\bsnm{Gangbo},~\bfnm{Wilfrid}\binits{W.}}
(\byear{1999}).
\btitle{The {M}onge mass transfer problem and its applications}.
In \bbooktitle{Monge {A}mp\`ere equation: applications to geometry and
  optimization ({D}eerfield {B}each, {FL}, 1997)}.
\bseries{Contemp. Math.}
\bvolume{226}
\bpages{79--104}.
\bpublisher{Amer. Math. Soc., Providence, RI}.
\bdoi{10.1090/conm/226/03236}
\bmrnumber{1660743}
\end{bincollection}
\endbibitem

\bibitem{GM96}
\begin{barticle}[author]
\bauthor{\bsnm{Gangbo},~\bfnm{Wilfrid}\binits{W.}} \AND
  \bauthor{\bsnm{McCann},~\bfnm{Robert~J.}\binits{R.~J.}}
(\byear{1996}).
\btitle{The geometry of optimal transportation}.
\bjournal{Acta Math.}
\bvolume{177}
\bpages{113--161}.
\bdoi{10.1007/BF02392620}
\bmrnumber{1440931}
\end{barticle}
\endbibitem

\bibitem{Gigli2011}
\begin{barticle}[author]
\bauthor{\bsnm{Gigli},~\bfnm{Nicola}\binits{N.}}
(\byear{2011}).
\btitle{On {H}\"{o}lder continuity-in-time of the optimal transport map towards
  measures along a curve}.
\bjournal{Proc. Edinb. Math. Soc. (2)}
\bvolume{54}
\bpages{401--409}.
\bdoi{10.1017/S001309150800117X}
\bmrnumber{2794662}
\end{barticle}
\endbibitem

\bibitem{GO17}
\begin{barticle}[author]
\bauthor{\bsnm{Goldman},~\bfnm{Michael}\binits{M.}} \AND
  \bauthor{\bsnm{Otto},~\bfnm{Felix}\binits{F.}}
(\byear{2017}).
\btitle{A variational proof of partial regularity for optimal transportation
  maps}.
\bjournal{arXiv preprint arXiv:1704.05339}.
\end{barticle}
\endbibitem

\bibitem{GrettonKernelMeasInd05}
\begin{barticle}[author]
\bauthor{\bsnm{Gretton},~\bfnm{A.}\binits{A.}},
  \bauthor{\bsnm{Herbrich},~\bfnm{R.}\binits{R.}},
  \bauthor{\bsnm{Smola},~\bfnm{A.}\binits{A.}},
  \bauthor{\bsnm{Bousquet},~\bfnm{O.}\binits{O.}} \AND
  \bauthor{\bsnm{Sch{\"o}lkopf},~\bfnm{B.}\binits{B.}}
(\byear{2005}).
\btitle{Kernel methods for measuring independence}.
\bjournal{J. Mach. Learn. Res.}
\bvolume{6}
\bpages{2075--2129 (electronic)}.
\bmrnumber{2249882}
\end{barticle}
\endbibitem

\bibitem{Gu12}
\begin{barticle}[author]
\bauthor{\bsnm{Gu},~\bfnm{Xianfeng}\binits{X.}},
  \bauthor{\bsnm{Luo},~\bfnm{Feng}\binits{F.}},
  \bauthor{\bsnm{Sun},~\bfnm{Jian}\binits{J.}} \AND
  \bauthor{\bsnm{Yau},~\bfnm{Shing-Tung}\binits{S.-T.}}
(\byear{2016}).
\btitle{Variational principles for {M}inkowski type problems, discrete optimal
  transport, and discrete {M}onge-{A}mpere equations}.
\bjournal{Asian J. Math.}
\bvolume{20}
\bpages{383--398}.
\bdoi{10.4310/AJM.2016.v20.n2.a7}
\bmrnumber{3480024}
\end{barticle}
\endbibitem

\bibitem{G16}
\begin{bbook}[author]
\bauthor{\bsnm{Guti\'{e}rrez},~\bfnm{C.~E.}\binits{C.~E.}}
(\byear{2016}).
\btitle{The {M}onge-{A}mp\`ere equation}.
\bseries{Progress in Nonlinear Differential Equations and their Applications}
\bvolume{89}.
\bpublisher{Birkh\"{a}user/Springer, [Cham]}
\bnote{Second edition [of MR1829162]}.
\end{bbook}
\endbibitem

\bibitem{dCHM}
\begin{barticle}[author]
\bauthor{\bsnm{Hallin},~\bfnm{Marc}\binits{M.}}, \bauthor{\bparticle{del}
  \bsnm{Barrio},~\bfnm{Eustasio}\binits{E.}},
  \bauthor{\bsnm{Cuesta-Albertos},~\bfnm{Juan}\binits{J.}} \AND
  \bauthor{\bsnm{Matr{\'a}n},~\bfnm{Carlos}\binits{C.}}
(\byear{2021}).
\btitle{Distribution and quantile functions, ranks and signs in dimension $d$:
  A measure transportation approach}.
\bjournal{The Annals of Statistics}
\bvolume{49}
\bpages{1139--1165}.
\end{barticle}
\endbibitem

\bibitem{Hallin2020}
\begin{barticle}[author]
\bauthor{\bsnm{Hallin},~\bfnm{Marc}\binits{M.}},
  \bauthor{\bsnm{Mordant},~\bfnm{Gilles}\binits{G.}} \AND
  \bauthor{\bsnm{Segers},~\bfnm{Johan}\binits{J.}}
(\byear{2020}).
\btitle{Multivariate goodness-of-Fit tests based on Wasserstein distance}.
\bjournal{arXiv preprint arXiv:2003.06684}.
\end{barticle}
\endbibitem

\bibitem{HPS10}
\begin{barticle}[author]
\bauthor{\bsnm{Hallin},~\bfnm{Marc}\binits{M.}},
  \bauthor{\bsnm{Paindaveine},~\bfnm{Davy}\binits{D.}} \AND
  \bauthor{\bsnm{\v{S}iman},~\bfnm{Miroslav}\binits{M.}}
(\byear{2010}).
\btitle{Multivariate quantiles and multiple-output regression quantiles: from
  {$L_1$} optimization to halfspace depth}.
\bjournal{Ann. Statist.}
\bvolume{38}
\bpages{635--669}.
\bdoi{10.1214/09-AOS723}
\bmrnumber{2604670}
\end{barticle}
\endbibitem

\bibitem{HW03}
\begin{barticle}[author]
\bauthor{\bsnm{Hallin},~\bfnm{Marc}\binits{M.}} \AND
  \bauthor{\bsnm{Werker},~\bfnm{Bas J.~M.}\binits{B.~J.~M.}}
(\byear{2003}).
\btitle{Semi-parametric efficiency, distribution-freeness and invariance}.
\bjournal{Bernoulli}
\bvolume{9}
\bpages{137--165}.
\bdoi{10.3150/bj/1068129013}
\bmrnumber{1963675}
\end{barticle}
\endbibitem

\bibitem{Hoeffding1952}
\begin{barticle}[author]
\bauthor{\bsnm{Hoeffding},~\bfnm{Wassily}\binits{W.}}
(\byear{1952}).
\btitle{The large-sample power of tests based on permutations of observations}.
\bjournal{Ann. Math. Statistics}
\bvolume{23}
\bpages{169--192}.
\bdoi{10.1214/aoms/1177729436}
\bmrnumber{57521}
\end{barticle}
\endbibitem

\bibitem{HW99}
\begin{bbook}[author]
\bauthor{\bsnm{Hollander},~\bfnm{Myles}\binits{M.}} \AND
  \bauthor{\bsnm{Wolfe},~\bfnm{Douglas~A.}\binits{D.~A.}}
(\byear{1999}).
\btitle{Nonparametric statistical methods},
\bedition{second} ed.
\bseries{Wiley Series in Probability and Statistics: Texts and References
  Section}.
\bpublisher{John Wiley \& Sons, Inc., New York}
\bnote{A Wiley-Interscience Publication}.
\bmrnumber{1666064}
\end{bbook}
\endbibitem

\bibitem{H94}
\begin{bbook}[author]
\bauthor{\bsnm{H\"{o}rmander},~\bfnm{Lars}\binits{L.}}
(\byear{1994}).
\btitle{Notions of convexity}.
\bseries{Progress in Mathematics}
\bvolume{127}.
\bpublisher{Birkh\"{a}user Boston, Inc., Boston, MA}.
\bmrnumber{1301332}
\end{bbook}
\endbibitem

\bibitem{HR19}
\begin{barticle}[author]
\bauthor{\bsnm{H{\"u}tter},~\bfnm{Jan-Christian}\binits{J.-C.}} \AND
  \bauthor{\bsnm{Rigollet},~\bfnm{Philippe}\binits{P.}}
(\byear{2019}).
\btitle{Minimax rates of estimation for smooth optimal transport maps}.
\bjournal{arXiv preprint arXiv:1905.05828}.
\end{barticle}
\endbibitem

\bibitem{Kim2020}
\begin{barticle}[author]
\bauthor{\bsnm{Kim},~\bfnm{Ilmun}\binits{I.}},
  \bauthor{\bsnm{Balakrishnan},~\bfnm{Sivaraman}\binits{S.}} \AND
  \bauthor{\bsnm{Wasserman},~\bfnm{Larry}\binits{L.}}
(\byear{2020}).
\btitle{Minimax optimality of permutation tests}.
\bjournal{arXiv preprint arXiv:2003.13208}.
\end{barticle}
\endbibitem

\bibitem{Kitagawa14}
\begin{barticle}[author]
\bauthor{\bsnm{Kitagawa},~\bfnm{Jun}\binits{J.}}
(\byear{2014}).
\btitle{An iterative scheme for solving the optimal transportation problem}.
\bjournal{Calc. Var. Partial Differential Equations}
\bvolume{51}
\bpages{243--263}.
\bdoi{10.1007/s00526-013-0673-x}
\bmrnumber{3247388}
\end{barticle}
\endbibitem

\bibitem{KM17}
\begin{barticle}[author]
\bauthor{\bsnm{Kitagawa},~\bfnm{Jun}\binits{J.}} \AND
  \bauthor{\bsnm{McCann},~\bfnm{Robert}\binits{R.}}
(\byear{2019}).
\btitle{Free discontinuities in optimal transport}.
\bjournal{Arch. Ration. Mech. Anal.}
\bvolume{232}
\bpages{1505--1541}.
\bdoi{10.1007/s00205-018-01348-3}
\bmrnumber{3928755}
\end{barticle}
\endbibitem

\bibitem{Kita-Meri-19}
\begin{barticle}[author]
\bauthor{\bsnm{Kitagawa},~\bfnm{Jun}\binits{J.}},
  \bauthor{\bsnm{M\'{e}rigot},~\bfnm{Quentin}\binits{Q.}} \AND
  \bauthor{\bsnm{Thibert},~\bfnm{Boris}\binits{B.}}
(\byear{2019}).
\btitle{Convergence of a {N}ewton algorithm for semi-discrete optimal
  transport}.
\bjournal{J. Eur. Math. Soc. (JEMS)}
\bvolume{21}
\bpages{2603--2651}.
\bdoi{10.4171/JEMS/889}
\bmrnumber{3985609}
\end{barticle}
\endbibitem

\bibitem{Klatt20}
\begin{barticle}[author]
\bauthor{\bsnm{Klatt},~\bfnm{Marcel}\binits{M.}},
  \bauthor{\bsnm{Tameling},~\bfnm{Carla}\binits{C.}} \AND
  \bauthor{\bsnm{Munk},~\bfnm{Axel}\binits{A.}}
(\byear{2020}).
\btitle{Empirical regularized optimal transport: statistical theory and
  applications}.
\bjournal{SIAM J. Math. Data Sci.}
\bvolume{2}
\bpages{419--443}.
\bdoi{10.1137/19M1278788}
\bmrnumber{4105566}
\end{barticle}
\endbibitem

\bibitem{Kol97}
\begin{barticle}[author]
\bauthor{\bsnm{Koltchinskii},~\bfnm{V.~I.}\binits{V.~I.}}
(\byear{1997}).
\btitle{{$M$}-estimation, convexity and quantiles}.
\bjournal{Ann. Statist.}
\bvolume{25}
\bpages{435--477}.
\bdoi{10.1214/aos/1031833659}
\bmrnumber{1439309}
\end{barticle}
\endbibitem

\bibitem{Lehmann75}
\begin{bbook}[author]
\bauthor{\bsnm{Lehmann},~\bfnm{E.~L.}\binits{E.~L.}}
(\byear{1975}).
\btitle{Nonparametrics: statistical methods based on ranks}.
\bpublisher{Holden-Day, Inc., San Francisco, Calif.; McGraw-Hill International
  Book Co., New York-D\"{u}sseldorf}
\bnote{With the special assistance of H. J. M. d'Abrera, Holden-Day Series in
  Probability and Statistics}.
\bmrnumber{0395032}
\end{bbook}
\endbibitem

\bibitem{TSH-2005}
\begin{bbook}[author]
\bauthor{\bsnm{Lehmann},~\bfnm{E.~L.}\binits{E.~L.}} \AND
  \bauthor{\bsnm{Romano},~\bfnm{Joseph~P.}\binits{J.~P.}}
(\byear{2005}).
\btitle{Testing statistical hypotheses},
\bedition{third} ed.
\bseries{Springer Texts in Statistics}.
\bpublisher{Springer, New York}.
\bmrnumber{2135927}
\end{bbook}
\endbibitem

\bibitem{Levy15}
\begin{barticle}[author]
\bauthor{\bsnm{L\'{e}vy},~\bfnm{Bruno}\binits{B.}}
(\byear{2015}).
\btitle{A numerical algorithm for {$L_2$} semi-discrete optimal transport in
  3{D}}.
\bjournal{ESAIM Math. Model. Numer. Anal.}
\bvolume{49}
\bpages{1693--1715}.
\bdoi{10.1051/m2an/2015055}
\bmrnumber{3423272}
\end{barticle}
\endbibitem

\bibitem{LN2020}
\begin{barticle}[author]
\bauthor{\bsnm{Li},~\bfnm{Wenbo}\binits{W.}} \AND
  \bauthor{\bsnm{Nochetto},~\bfnm{Ricardo~H}\binits{R.~H.}}
(\byear{2020}).
\btitle{Quantitative Stability and Error Estimates for Optimal Transport
  Plans}.
\bjournal{arXiv preprint arXiv:2004.05299}.
\end{barticle}
\endbibitem

\bibitem{Liu92}
\begin{barticle}[author]
\bauthor{\bsnm{Liu},~\bfnm{Regina~Y.}\binits{R.~Y.}}
(\byear{1990}).
\btitle{On a notion of data depth based on random simplices}.
\bjournal{Ann. Statist.}
\bvolume{18}
\bpages{405--414}.
\bdoi{10.1214/aos/1176347507}
\bmrnumber{1041400}
\end{barticle}
\endbibitem

\bibitem{Lyons13}
\begin{barticle}[author]
\bauthor{\bsnm{Lyons},~\bfnm{Russell}\binits{R.}}
(\byear{2013}).
\btitle{Distance covariance in metric spaces}.
\bjournal{Ann. Probab.}
\bvolume{41}
\bpages{3284--3305}.
\bdoi{10.1214/12-AOP803}
\bmrnumber{3127883}
\end{barticle}
\endbibitem

\bibitem{McCann95}
\begin{barticle}[author]
\bauthor{\bsnm{McCann},~\bfnm{Robert~J.}\binits{R.~J.}}
(\byear{1995}).
\btitle{Existence and uniqueness of monotone measure-preserving maps}.
\bjournal{Duke Math. J.}
\bvolume{80}
\bpages{309--323}.
\bdoi{10.1215/S0012-7094-95-08013-2}
\bmrnumber{1369395}
\end{barticle}
\endbibitem

\bibitem{Merigot11}
\begin{binproceedings}[author]
\bauthor{\bsnm{M{\'e}rigot},~\bfnm{Quentin}\binits{Q.}}
(\byear{2011}).
\btitle{A multiscale approach to optimal transport}.
In \bbooktitle{Computer Graphics Forum}
\bvolume{30}
\bpages{1583--1592}.
\bpublisher{Wiley Online Library}.
\end{binproceedings}
\endbibitem

\bibitem{Merigot2020}
\begin{barticle}[author]
\bauthor{\bsnm{Merigot},~\bfnm{Quentin}\binits{Q.}} \AND
  \bauthor{\bsnm{Thibert},~\bfnm{Boris}\binits{B.}}
(\byear{2020}).
\btitle{Optimal transport: discretization and algorithms}.
\bjournal{Handbook of Numerical Analysis 22 -- Geometric PDES (arXiv preprint
  arXiv:2003.00855)}.
\end{barticle}
\endbibitem

\bibitem{Monge1781}
\begin{barticle}[author]
\bauthor{\bsnm{Monge},~\bfnm{Gaspard}\binits{G.}}
(\byear{1781}).
\btitle{M\'{e}moire sur la th\'{e}orie des d\'{e}blais et des remblais}.
\bjournal{Histoire de l'Acad\'{e}mie Royale des Sciences de Paris}
\bpages{666--704}.
\end{barticle}
\endbibitem

\bibitem{4-Axioms-MS19}
\begin{barticle}[author]
\bauthor{\bsnm{M\'{o}ri},~\bfnm{Tam\'{a}s~F.}\binits{T.~F.}} \AND
  \bauthor{\bsnm{Sz\'{e}kely},~\bfnm{G\'{a}bor~J.}\binits{G.~J.}}
(\byear{2019}).
\btitle{Four simple axioms of dependence measures}.
\bjournal{Metrika}
\bvolume{82}
\bpages{1--16}.
\bdoi{10.1007/s00184-018-0670-3}
\bmrnumber{3897521}
\end{barticle}
\endbibitem

\bibitem{Oja83}
\begin{barticle}[author]
\bauthor{\bsnm{Oja},~\bfnm{Hannu}\binits{H.}}
(\byear{1983}).
\btitle{Descriptive statistics for multivariate distributions}.
\bjournal{Statist. Probab. Lett.}
\bvolume{1}
\bpages{327--332}.
\bdoi{10.1016/0167-7152(83)90054-8}
\bmrnumber{721446}
\end{barticle}
\endbibitem

\bibitem{OP88}
\begin{barticle}[author]
\bauthor{\bsnm{Oliker},~\bfnm{V.~I.}\binits{V.~I.}} \AND
  \bauthor{\bsnm{Prussner},~\bfnm{L.~D.}\binits{L.~D.}}
(\byear{1988}).
\btitle{On the numerical solution of the equation {$(\partial^2z/\partial
  x^2)(\partial^2z/\partial y^2)-((\partial^2z/\partial x\partial y))^2=f$} and
  its discretizations. {I}}.
\bjournal{Numer. Math.}
\bvolume{54}
\bpages{271--293}.
\bdoi{10.1007/BF01396762}
\bmrnumber{971703}
\end{barticle}
\endbibitem

\bibitem{PZ19}
\begin{barticle}[author]
\bauthor{\bsnm{Panaretos},~\bfnm{Victor~M.}\binits{V.~M.}} \AND
  \bauthor{\bsnm{Zemel},~\bfnm{Yoav}\binits{Y.}}
(\byear{2019}).
\btitle{Statistical aspects of {W}asserstein distances}.
\bjournal{Annu. Rev. Stat. Appl.}
\bvolume{6}
\bpages{405--431}.
\bdoi{10.1146/annurev-statistics-030718-104938}
\bmrnumber{3939527}
\end{barticle}
\endbibitem

\bibitem{Peyre2019}
\begin{barticle}[author]
\bauthor{\bsnm{Peyr{\'e}},~\bfnm{Gabriel}\binits{G.}},
  \bauthor{\bsnm{Cuturi},~\bfnm{Marco}\binits{M.}} \betal{et~al.}
(\byear{2019}).
\btitle{Computational Optimal Transport: With Applications to Data Science}.
\bjournal{Foundations and Trends{\textregistered} in Machine Learning}
\bvolume{11}
\bpages{355--607}.
\end{barticle}
\endbibitem

\bibitem{Pfister18}
\begin{barticle}[author]
\bauthor{\bsnm{Pfister},~\bfnm{Niklas}\binits{N.}},
  \bauthor{\bsnm{B\"{u}hlmann},~\bfnm{Peter}\binits{P.}},
  \bauthor{\bsnm{Sch\"{o}lkopf},~\bfnm{Bernhard}\binits{B.}} \AND
  \bauthor{\bsnm{Peters},~\bfnm{Jonas}\binits{J.}}
(\byear{2018}).
\btitle{Kernel-based tests for joint independence}.
\bjournal{J. R. Stat. Soc. Ser. B. Stat. Methodol.}
\bvolume{80}
\bpages{5--31}.
\bdoi{10.1111/rssb.12235}
\bmrnumber{3744710}
\end{barticle}
\endbibitem

\bibitem{Ramdas17}
\begin{barticle}[author]
\bauthor{\bsnm{Ramdas},~\bfnm{Aaditya}\binits{A.}},
  \bauthor{\bsnm{Garc\'{\i}a~Trillos},~\bfnm{Nicol\'{a}s}\binits{N.}} \AND
  \bauthor{\bsnm{Cuturi},~\bfnm{Marco}\binits{M.}}
(\byear{2017}).
\btitle{On {W}asserstein two-sample testing and related families of
  nonparametric tests}.
\bjournal{Entropy}
\bvolume{19}
\bpages{Paper No. 47, 15}.
\bdoi{10.3390/e19020047}
\bmrnumber{3608466}
\end{barticle}
\endbibitem

\bibitem{RW18-MLE}
\begin{barticle}[author]
\bauthor{\bsnm{Rigollet},~\bfnm{Philippe}\binits{P.}} \AND
  \bauthor{\bsnm{Weed},~\bfnm{Jonathan}\binits{J.}}
(\byear{2018}).
\btitle{Entropic optimal transport is maximum-likelihood deconvolution}.
\bjournal{C. R. Math. Acad. Sci. Paris}
\bvolume{356}
\bpages{1228--1235}.
\bdoi{10.1016/j.crma.2018.10.010}
\bmrnumber{3907589}
\end{barticle}
\endbibitem

\bibitem{RW19}
\begin{barticle}[author]
\bauthor{\bsnm{Rigollet},~\bfnm{Philippe}\binits{P.}} \AND
  \bauthor{\bsnm{Weed},~\bfnm{Jonathan}\binits{J.}}
(\byear{2019}).
\btitle{Uncoupled isotonic regression via minimum {W}asserstein deconvolution}.
\bjournal{Inf. Inference}
\bvolume{8}
\bpages{691--717}.
\bdoi{10.1093/imaiai/iaz006}
\bmrnumber{4045481}
\end{barticle}
\endbibitem

\bibitem{Rockf}
\begin{bbook}[author]
\bauthor{\bsnm{Rockafellar},~\bfnm{R.~Tyrrell}\binits{R.~T.}}
(\byear{1970}).
\btitle{Convex analysis}.
\bseries{Princeton Mathematical Series, No. 28}.
\bpublisher{Princeton University Press, Princeton, N.J.}
\bmrnumber{0274683}
\end{bbook}
\endbibitem

\bibitem{Rosenbaum05}
\begin{barticle}[author]
\bauthor{\bsnm{Rosenbaum},~\bfnm{Paul~R.}\binits{P.~R.}}
(\byear{2005}).
\btitle{An exact distribution-free test comparing two multivariate
  distributions based on adjacency}.
\bjournal{J. R. Stat. Soc. Ser. B Stat. Methodol.}
\bvolume{67}
\bpages{515--530}.
\bdoi{10.1111/j.1467-9868.2005.00513.x}
\bmrnumber{2168202}
\end{barticle}
\endbibitem

\bibitem{RR1999}
\begin{barticle}[author]
\bauthor{\bsnm{Rousseeuw},~\bfnm{Peter~J.}\binits{P.~J.}} \AND
  \bauthor{\bsnm{Ruts},~\bfnm{Ida}\binits{I.}}
(\byear{1999}).
\btitle{The depth function of a population distribution}.
\bjournal{Metrika}
\bvolume{49}
\bpages{213--244}.
\bdoi{10.1007/PL00020903}
\end{barticle}
\endbibitem

\bibitem{Sc86}
\begin{barticle}[author]
\bauthor{\bsnm{Schilling},~\bfnm{Mark~F.}\binits{M.~F.}}
(\byear{1986}).
\btitle{Multivariate two-sample tests based on nearest neighbors}.
\bjournal{J. Amer. Statist. Assoc.}
\bvolume{81}
\bpages{799--806}.
\bmrnumber{860514}
\end{barticle}
\endbibitem

\bibitem{SS11}
\begin{barticle}[author]
\bauthor{\bsnm{Seijo},~\bfnm{Emilio}\binits{E.}} \AND
  \bauthor{\bsnm{Sen},~\bfnm{Bodhisattva}\binits{B.}}
(\byear{2011}).
\btitle{Nonparametric least squares estimation of a multivariate convex
  regression function}.
\bjournal{Ann. Statist.}
\bvolume{39}
\bpages{1633--1657}.
\end{barticle}
\endbibitem

\bibitem{EquivRKHS13}
\begin{barticle}[author]
\bauthor{\bsnm{Sejdinovic},~\bfnm{Dino}\binits{D.}},
  \bauthor{\bsnm{Sriperumbudur},~\bfnm{Bharath}\binits{B.}},
  \bauthor{\bsnm{Gretton},~\bfnm{Arthur}\binits{A.}} \AND
  \bauthor{\bsnm{Fukumizu},~\bfnm{Kenji}\binits{K.}}
(\byear{2013}).
\btitle{Equivalence of distance-based and {RKHS}-based statistics in hypothesis
  testing}.
\bjournal{Ann. Statist.}
\bvolume{41}
\bpages{2263--2291}.
\bdoi{10.1214/13-AOS1140}
\bmrnumber{3127866}
\end{barticle}
\endbibitem

\bibitem{Serfling10}
\begin{barticle}[author]
\bauthor{\bsnm{Serfling},~\bfnm{Robert}\binits{R.}}
(\byear{2010}).
\btitle{Equivariance and invariance properties of multivariate quantile and
  related functions, and the role of standardisation}.
\bjournal{J. Nonparametr. Stat.}
\bvolume{22}
\bpages{915--936}.
\bdoi{10.1080/10485250903431710}
\bmrnumber{2738875}
\end{barticle}
\endbibitem

\bibitem{Shi2019}
\begin{barticle}[author]
\bauthor{\bsnm{Shi},~\bfnm{Hongjian}\binits{H.}},
  \bauthor{\bsnm{Drton},~\bfnm{Mathias}\binits{M.}} \AND
  \bauthor{\bsnm{Han},~\bfnm{Fang}\binits{F.}}
(\byear{2019}).
\btitle{Distribution-free consistent independence tests via {H}allin's
  multivariate rank}.
\bjournal{arXiv preprint arXiv:1909.10024}.
\end{barticle}
\endbibitem

\bibitem{Shi2020}
\begin{barticle}[author]
\bauthor{\bsnm{Shi},~\bfnm{Hongjian}\binits{H.}},
  \bauthor{\bsnm{Hallin},~\bfnm{Marc}\binits{M.}},
  \bauthor{\bsnm{Drton},~\bfnm{Mathias}\binits{M.}} \AND
  \bauthor{\bsnm{Han},~\bfnm{Fang}\binits{F.}}
(\byear{2020}).
\btitle{Rate-optimality of consistent distribution-free tests of independence
  based on center-outward ranks and signs}.
\bjournal{arXiv preprint arXiv:2007.02186}.
\end{barticle}
\endbibitem

\bibitem{SzekelyBDCov09}
\begin{barticle}[author]
\bauthor{\bsnm{Sz{\'e}kely},~\bfnm{G.~J.}\binits{G.~J.}} \AND
  \bauthor{\bsnm{Rizzo},~\bfnm{M.~L.}\binits{M.~L.}}
(\byear{2009}).
\btitle{Brownian distance covariance}.
\bjournal{Ann. Appl. Stat.}
\bvolume{3}
\bpages{1236--1265}.
\bdoi{10.1214/09-AOAS312}
\bmrnumber{2752127}
\end{barticle}
\endbibitem

\bibitem{SR13}
\begin{barticle}[author]
\bauthor{\bsnm{Sz\'{e}kely},~\bfnm{G\'{a}bor~J.}\binits{G.~J.}} \AND
  \bauthor{\bsnm{Rizzo},~\bfnm{Maria~L.}\binits{M.~L.}}
(\byear{2013}).
\btitle{Energy statistics: a class of statistics based on distances}.
\bjournal{J. Statist. Plann. Inference}
\bvolume{143}
\bpages{1249--1272}.
\bdoi{10.1016/j.jspi.2013.03.018}
\bmrnumber{3055745}
\end{barticle}
\endbibitem

\bibitem{SzekelyCorrDist07}
\begin{barticle}[author]
\bauthor{\bsnm{Sz{\'e}kely},~\bfnm{G.~J.}\binits{G.~J.}},
  \bauthor{\bsnm{Rizzo},~\bfnm{M.~L.}\binits{M.~L.}} \AND
  \bauthor{\bsnm{Bakirov},~\bfnm{N.~K.}\binits{N.~K.}}
(\byear{2007}).
\btitle{Measuring and testing dependence by correlation of distances}.
\bjournal{Ann. Statist.}
\bvolume{35}
\bpages{2769--2794}.
\bdoi{10.1214/009053607000000505}
\bmrnumber{2382665 (2009a:62206)}
\end{barticle}
\endbibitem

\bibitem{vdV98}
\begin{bbook}[author]
\bauthor{\bparticle{van~der} \bsnm{Vaart},~\bfnm{A.~W.}\binits{A.~W.}}
(\byear{1998}).
\btitle{Asymptotic statistics}.
\bseries{Cambridge Series in Statistical and Probabilistic Mathematics}
\bvolume{3}.
\bpublisher{Cambridge University Press, Cambridge}.
\bdoi{10.1017/CBO9780511802256}
\bmrnumber{1652247}
\end{bbook}
\endbibitem

\bibitem{V03}
\begin{bbook}[author]
\bauthor{\bsnm{Villani},~\bfnm{C\'{e}dric}\binits{C.}}
(\byear{2003}).
\btitle{Topics in optimal transportation}.
\bseries{Graduate Studies in Mathematics}
\bvolume{58}.
\bpublisher{American Mathematical Society, Providence, RI}.
\bdoi{10.1007/b12016}
\bmrnumber{1964483}
\end{bbook}
\endbibitem

\bibitem{V09}
\begin{bbook}[author]
\bauthor{\bsnm{Villani},~\bfnm{C\'{e}dric}\binits{C.}}
(\byear{2009}).
\btitle{Optimal transport}.
\bseries{Grundlehren der Mathematischen Wissenschaften [Fundamental Principles
  of Mathematical Sciences]}
\bvolume{338}.
\bpublisher{Springer-Verlag, Berlin}
\bnote{Old and new}.
\bdoi{10.1007/978-3-540-71050-9}
\bmrnumber{2459454}
\end{bbook}
\endbibitem

\bibitem{WB19}
\begin{barticle}[author]
\bauthor{\bsnm{Weed},~\bfnm{Jonathan}\binits{J.}} \AND
  \bauthor{\bsnm{Bach},~\bfnm{Francis}\binits{F.}}
(\byear{2019}).
\btitle{Sharp asymptotic and finite-sample rates of convergence of empirical
  measures in {W}asserstein distance}.
\bjournal{Bernoulli}
\bvolume{25}
\bpages{2620--2648}.
\bdoi{10.3150/18-BEJ1065}
\bmrnumber{4003560}
\end{barticle}
\endbibitem

\bibitem{Weihs18}
\begin{barticle}[author]
\bauthor{\bsnm{Weihs},~\bfnm{L.}\binits{L.}},
  \bauthor{\bsnm{Drton},~\bfnm{M.}\binits{M.}} \AND
  \bauthor{\bsnm{Meinshausen},~\bfnm{N.}\binits{N.}}
(\byear{2018}).
\btitle{Symmetric rank covariances: a generalized framework for nonparametric
  measures of dependence}.
\bjournal{Biometrika}
\bvolume{105}
\bpages{547--562}.
\bdoi{10.1093/biomet/asy021}
\bmrnumber{3842884}
\end{barticle}
\endbibitem

\bibitem{Weiss60}
\begin{barticle}[author]
\bauthor{\bsnm{Weiss},~\bfnm{Lionel}\binits{L.}}
(\byear{1960}).
\btitle{Two-sample tests for multivariate distributions}.
\bjournal{Ann. Math. Statist.}
\bvolume{31}
\bpages{159--164}.
\bdoi{10.1214/aoms/1177705995}
\bmrnumber{0119305}
\end{barticle}
\endbibitem

\bibitem{OTM}
\begin{bmanual}[author]
\bauthor{\bsnm{Xu},~\bfnm{Peng}\binits{P.}}
(\byear{2019}).
\btitle{testOTM: Multivariate Ranks and Quantiles using Optimal Transportation}
\bnote{R package version 1.00.0}.
\end{bmanual}
\endbibitem

\bibitem{Zemel2019}
\begin{barticle}[author]
\bauthor{\bsnm{Zemel},~\bfnm{Yoav}\binits{Y.}} \AND
  \bauthor{\bsnm{Panaretos},~\bfnm{Victor~M.}\binits{V.~M.}}
(\byear{2019}).
\btitle{Fr\'{e}chet means and {P}rocrustes analysis in {W}asserstein space}.
\bjournal{Bernoulli}
\bvolume{25}
\bpages{932--976}.
\bdoi{10.3150/17-bej1009}
\bmrnumber{3920362}
\end{barticle}
\endbibitem

\bibitem{Zou03}
\begin{barticle}[author]
\bauthor{\bsnm{Zuo},~\bfnm{Yijun}\binits{Y.}}
(\byear{2003}).
\btitle{Projection-based depth functions and associated medians}.
\bjournal{Ann. Statist.}
\bvolume{31}
\bpages{1460--1490}.
\bdoi{10.1214/aos/1065705115}
\bmrnumber{2012822}
\end{barticle}
\endbibitem

\end{thebibliography}
\end{document}